\documentclass[11pt]{article}
\usepackage[utf8x]{inputenc}
\usepackage[toc,page]{appendix}
\usepackage{microtype}
\usepackage[T1]{fontenc}
\usepackage{enumitem}
\usepackage{tikz-cd}
\usepackage{slashed}
\usepackage{mathrsfs} 
\usepackage{ae,aecompl}
\usepackage{times}
\usepackage{float}
\usepackage{verbatim,amsmath,amsthm,amsfonts,amssymb,latexsym,graphicx,mathtools,extpfeil,color,mathabx,scalerel}
\usepackage{stmaryrd}
\usepackage{epstopdf,pinlabel} 
\usepackage[all]{xy}
\usepackage{graphicx}
\usepackage{caption}
\usepackage{subcaption}
\DeclareMathAlphabet{\mathpzc}{OT1}{pzc}{m}{it}
\usepackage[colorlinks,pagebackref,hypertexnames=false]{hyperref} 
\usepackage{hyperref}
\hypersetup{
  colorlinks = true,
  linkcolor  = black
} 
\usepackage[alphabetic,backrefs,msc-links]{amsrefs}
\usepackage{accents}

\usepackage{aliascnt}
\numberwithin{equation}{section}
\newcommand{\wt}{\widetilde}

\newcommand{\bA}{{\bf A}}

\newcommand{\bE}{{\bf E}}

\newcommand{\bR}{{\bf R}}

\newcommand{\fA}{{\mathfrak A}}

\newcommand{\fC}{{\mathfrak C}}

\newcommand{\fZ}{{\mathfrak Z}}

\newcommand{\R}{\bR}

\DeclareMathOperator{\ad}{ad}
\DeclareMathOperator{\codim}{codim}

\DeclareMathOperator{\tr}{tr}

\renewcommand{\det}{\operatorname{det}}

\renewcommand{\epsilon}{\varepsilon}

\def\({\mathopen{}\left(}
\def\){\right)\mathclose{}}
\def\<{\mathopen{}\left<}
\def\>{\right>\mathclose{}}

\usepackage{multicol, color}

\definecolor{gold}{rgb}{0.85,.66,0}
\definecolor{cherry}{rgb}{0.9,.1,.2}
\definecolor{burgundy}{rgb}{0.8,.2,.2}
\definecolor{orangered}{rgb}{0.85,.3,0}
\definecolor{orange}{rgb}{0.85,.4,0}
\definecolor{olive}{rgb}{.45,.4,0}
\definecolor{lime}{rgb}{.6,.9,0}
\definecolor{green}{rgb}{.2,.7,0}
\definecolor{grey}{rgb}{.4,.4,.2}
\definecolor{brown}{rgb}{.4,.3,.1}

\def\makeautorefname#1#2{\AtBeginDocument{\expandafter\def\csname#1autorefname\endcsname{#2}}}

\newcommand{\mynewtheorem}[2]{
  \newaliascnt{#1}{equation}          
  \newtheorem{#1}[#1]{#2}
  \aliascntresetthe{#1}
  \makeautorefname{#1}{#2}
}
\mynewtheorem{theorem}{Theorem}
\mynewtheorem{prop}{Proposition}
\mynewtheorem{cor}{Corollary}
\mynewtheorem{lemma}{Lemma}
\mynewtheorem{conjecture}{Conjecture}
\mynewtheorem{question}{Question}
\mynewtheorem{problem}{Problem}

\numberwithin{substep}{step}
\makeautorefname{step}{Step}
\makeautorefname{substep}{Step}

\numberwithin{subcase}{case}
\makeautorefname{case}{Case}
\makeautorefname{subcase}{case}

\theoremstyle{remark}
\mynewtheorem{remarkx}{Remark}
\newenvironment{remark}
  {\pushQED{\qed}\remarkx}
  {\popQED\endremarkx}

\theoremstyle{definition}
\mynewtheorem{definitionx}{Definition}
\newenvironment{definition}
  {\pushQED{\qed}\definitionx}
  {\popQED\enddefinitionx}
\mynewtheorem{constructionx}{Construction}
\newenvironment{construction}
  {\pushQED{\qed}\constructionx}
  {\popQED\endconstructionx}
\mynewtheorem{examplex}{Example}
\newenvironment{example}
  {\pushQED{\qed}\examplex}
  {\popQED\endexamplex}
\mynewtheorem{exercise}{Exercise}
\mynewtheorem{conventionx}{Convention}
\newenvironment{convention}
  {\pushQED{\qed}\conventionx}
  {\popQED\endconventionx}
\newtheorem*{remark*}{Remark}
\newtheorem*{convention*}{Convention}
\newtheorem*{conventions*}{Conventions}
\newtheorem*{assumption*}{Assumption}
        
\makeautorefname{chapter}{Chapter}
\makeautorefname{section}{Section}
\makeautorefname{subsection}{Section}
\makeautorefname{subsubsection}{Section}

\makeatletter
\renewcommand*\env@matrix[1][*\c@MaxMatrixCols c]{%
  \hskip -\arraycolsep
  \let\@ifnextchar\new@ifnextchar
  \array{#1}}
\makeatother

\theoremstyle{plain}
\newtheorem{theorem-intro}{Theorem}
\newtheorem{cor-intro}{Corollary}
\newtheorem{prop-intro}{Proposition}

\title{\large \bf Instantons and Rational Homology Spheres}
\author{Aliakbar Daemi\thanks{AD was partially supported by NSF Grants DMS-1812033, DMS-2208181 and NSF FRG Grant DMS-1952762. This work was completed while AD was in residence at the Simons Laufer Mathematical Sciences Institute during Fall 2023 and supported by NSF Grant DMS-1928930.} \hspace{1cm} Mike Miller Eismeier\thanks{MME was partially supported by NSF FRG Grant DMS-1952762.}}

\begin{document}
\date{}
\maketitle

\begin{abstract}
In previous work, the second author defined `equivariant instanton homology groups' $I^\bullet(Y,\pi;R)$ for a rational homology 3-sphere $Y$, a set of auxiliary data $\pi$, and a PID $R$. These objects are modules over the cohomology ring $H^{-*}(BSO_3;R)$.

We prove that the equivariant instanton homology groups $I^\bullet(Y;R)$ are independent of the auxiliary data $\pi$, and thus define topological invariants of rational homology spheres. Further, we prove that these invariants are functorial under cobordisms of 3-manifolds with a path between the boundary components.

For any rational homology sphere $Y$, we may also define an analogue of Floer's irreducible instanton homology group of integer homology spheres $I_*(Y, \pi; R)$ which now depends on the auxiliary data $\pi$, unlike the equivariant instanton homology groups. However, our methods allow us to prove a precise ``wall-crossing formula'' for $I_*(Y, \pi; R)$ as the auxiliary data $\pi$ moves between adjacent chambers. We use this to define an instanton invariant $\lambda_I(Y) \in \Bbb Q$ of rational homology spheres, conjecturally equal to the Casson--Walker invariant.

Our approach to invariance uses a novel technique known as a suspended flow category. Given an obstructed cobordism $W: Y \to Y'$, which supports reducible instantons which can neither be cut out transversely nor be removed by a small change to the perturbation, we remove and replace a neighborhood of obstructed solutions in the moduli space of instantons. The resulting moduli spaces have a new type of boundary component, so do not define a chain map between the instanton chain complexes of $Y$ and $Y'$. However, they do define a chain map between the instanton chain complex of $Y$ and a sort of suspension of the instanton chain complex of $Y'$.
\end{abstract}

\tableofcontents

\section{Introduction}
%!TEX root = equivariant-functoriality.tex

Instanton homology originates in Floer's landmark paper \cite{Fl}, where he associates to a 3-manifold $Y$ with the same integral homology groups as the 3-sphere (an integer homology 3-sphere)  a finite-dimensional $\Bbb Z/8$-graded abelian group, $I_*(Y)$. This group has many desiderata: its Euler characteristic returns twice the Casson invariant \cite{Taubes}, while a four-manifold $X$ with boundary $Y$ gives elements $D_X \in I_*(Y)$ generalizing Donaldson's polynomial invariants for closed 4-manifolds; if a closed 4-manifold $X$ can be decomposed as $X = X_1 \cup_Y X_2$ along an integer homology 3-sphere $Y$, then the Donaldson invariants of $X$ may be computed by pairing $D_{X_1} \in I_*(Y)$ against $D_{X_2} \in I_*(-Y) =: I^*(Y)$.

Three-manifolds with trivial first Betti number (rational homology 3-spheres) provide an important extension of integer homology 3-spheres. The naive generalization of Floer's instanton homology to rational homology spheres does not produce a well-defined invariant (but see Theorem \ref{introthm:wallcross} below). Instead, one needs to introduce a variation called the \emph{equivariant} instanton homology groups. Such a construction was first attempted in \cite{AB1}, though it was not established that the result is an invariant of the homeomorphism type of 3-manifolds. An analogous and more transparently well-defined construction for integer homology spheres was carried out in \cite[Section 7.3.4]{Don}.

Our first theorem is that a variant of the Austin and Braam construction does indeed produce well-defined invariants of rational homology spheres.

\begin{theorem-intro}\label{introthm:main}
	For any PID $R$ and any rational homology 3-sphere $Y$, there are three well-defined $\Bbb Z/8$-graded groups 
	\[
	I^+(Y;R), \hspace{1cm} I^-(Y;R), \hspace{1cm} I^\infty(Y;R),
	\]
	which are invariants of the homeomorphism type of $Y$. These three groups are graded modules over the ring $H^{-*}(BSO(3); R)$ and fit into an exact triangle 
	\[
	  \cdots \xrightarrow{[-4]} I^+(Y;R) \xrightarrow{[3]} I^-(Y;R) \to I^\infty(Y;R) \xrightarrow{[-4]} \cdots
	\]   
	in which the morphism $I^- \to I^\infty$ is an $H^{-*}(BSO(3); R)$-module homomorphism, as are the other two when $2 = 0 \in R$ or $2 \in R^\times$, and the grading shifts are indicated above each arrow.
\end{theorem-intro}

More generally, we define invariants for 3-manifolds equipped with $SO(3)$-bundles. The more general statement appears as Theorem \ref{thm:main-ex}.

Previously, instantons have been used to define {\it framed Floer homology} $I^\sharp(Y)$ for any 3-manifold $Y$ \cite{YAFT}. In the case that $Y$ is a rational homology sphere, the invariants of the present paper are expected to strictly contain more information than framed Floer homology. To be a bit more specific and analogous to the case of integer homology spheres in \cite{Scaduto}, a connected-sum theorem would show that our package of invariants recovers Kronheimer and Mrowka's $I^\sharp(Y)$. A more precise statement appears as Remark \ref{rmk:compare}. For readers familiar with Heegaard Floer homology \cite{OzSz:HF}, we remark that $I^\bullet(Y)$ is the formal counterpart of ${\rm HF}^\bullet(Y)$ for $\bullet \in \{+,-,\infty\}$ whereas $I^\sharp(Y)$ is formally similar to $\widehat{\rm HF}(Y)$.

The equivariant instanton Floer homology groups of Theorem \ref{introthm:main} are constructed in \cite{M}\footnote{The definition of $I^\infty$ and of the exact triangle above should be corrected following \cite[Section 5.4]{Helle}}, by a two-step process. First, a chain complex $\widetilde C(Y, \pi)$ depending on a rational homology sphere $Y$ and some auxiliary data $\pi$ is constructed, so that this chain complex carries an action of the Lie group $SO(3)$. Next, an algebraic construction is applied to produce equivariant homology groups $I^\bullet(Y, \pi)$ for $\bullet \in \{+,-,\infty\}$. The auxiliary data $\pi$ consists of a Riemannian metric on $Y$ and a {\it regular} perturbation of the Chern-Simons functional. The novelty of this article is that the groups $I^\bullet(Y, \pi)$ are 3-manifold invariants and do not depend on $\pi$.

Each choice of $\pi$ gives rise to a function $\sigma_\pi:\fA(Y)\to 2\Bbb Z$, called the \emph{signature data} of $\pi$, where $\fA(Y)$ is the free orbits of the action of $\Bbb Z/2$ on $H^2(Y;\Bbb Z)$ given by multiplication by $\pm 1$. The values of this function modulo 4 are independent of $\pi$, and are given by the dimension of certain cohomology groups. It is shown in \cite{M} that $I^\bullet(Y, \pi)$ is an invariant of $(Y,\sigma_\pi)$, and in this paper we show that indeed $I^\bullet(Y, \pi)$ does not depend on $\sigma_\pi$. See Section \ref{section:overview} for more details on the strategy of the proof.

There is also an extension of Floer's instanton homology to rational homology spheres which is defined only in terms of irreducible (perturbed) flat connections. Unlike the invariants in Theorem \ref{introthm:main}, this invariant depends on the choice of signature data; we obtain a $\Bbb Z/8$-graded abelian group $I_*(Y,\sigma)$ for any signature data function $\sigma:\fA(Y)\to 2\Bbb Z$. The method of the proof of Theorem \ref{introthm:main} gives the following relation among $I_*(Y,\sigma)$ for different choices of $\sigma$. 

\begin{theorem-intro}\label{introthm:wallcross}
	Suppose $\sigma_0,\,\sigma_1: \fA(Y)\to 2\Bbb Z$ are adjacent signature data functions in the following sense: $\sigma_1(\rho)=\sigma_0(\rho)$ for any $\rho\in \fA(Y)$ except one element $\rho_0$ where we have $\sigma_1(\rho_0)=\sigma_0(\rho_0)+4$. Then there exists an exact triangle

	\begin{equation*}\label{top-equiv-triangle-intro}
	\begin{tikzcd}[column sep=1ex, row sep=5ex, fill=none, /tikz/baseline=-10pt]
	I_*(Y,\sigma_0) \arrow{rr} & & I_*(Y, \sigma_1) \arrow{dl}{[-1]} \\
	& \Bbb Z[\wt i(\rho_0)] \arrow{ul} & \\
	\end{tikzcd}
	\vspace{-1cm}
	\end{equation*}
	where $\wt i(\rho_0)$ is the integer given in Definition \ref{def:mod-2-grading}.
\end{theorem-intro}
\noindent 
In Section \ref{sec:CW}, we state a conjecture which proposes a relation between the Euler characteristic of this invariant (plus a spectral correcion term) and the Casson--Walker invariant \cite{Walker}. We verify this conjecture for lens spaces. For integer homology spheres, our claim is the main result of \cite{Taubes}. 

\subsection*{Functoriality of the invariants}
One of the most important features of Floer's instanton homology for integer homology spheres is its connection to 4-manifolds. Given a negative-definite cobordism $W: Y \to Y'$ between integer homology spheres with $b_1(W) = 0$, there is a well-defined induced map $I(W): I_*(Y) \to I_*(Y')$, which gives $I_*$ the structure of a functor over an appropriate cobordism category of 3- and 4-manifolds \cite[Chapter 7]{Don} and \cite[Section 3.2]{Fr1}.

The second major result of the present paper is related to the functoriality of the invariants of Theorem \ref{introthm:main} with respect to cobordisms. Suppose $\mathsf{Cob}_{3+1}^{\Bbb Q}$ is the category whose objects are rational homology 3-spheres with basepoints and whose morphisms between rational homology 3-spheres $Y_1$ and $Y_2$ are given by 4-dimensional smooth, oriented, and connected cobordisms $W:Y_1 \to Y_2$ together with a path $\gamma$ running between basepoints of $Y_i$ and an orientation of $H^1(W) \oplus H^{2,+}(W)$ (a homology orientation of $W$). The composition of morphisms is defined using composition of cobordisms. Let us also define $\mathsf{Cob}_{3+1}^{\Bbb Q,-}$ to be the subcategory with the same objects but whose morphisms are the negative-definite cobordisms. In this case, the vector space $H^{2,+}$ is trivial, and a homology orientation reduces to a choice of orientation on $H^1(W)$.

\begin{theorem-intro}\label{introthm:functorial}
	Let $R$ be a field with $\textup{char}(R) \ne 2$. Then there is a functor $I^\bullet$ from the category $\mathsf{Cob}_{3+1}^{\Bbb Q}$ to the category of $\Bbb Z/8$-graded modules over the ring $H^{-*}(BSO(3); R)$
	that specializes to the invariants of Theorem \ref{introthm:main} at the level of objects. 
	For an arbitrary PID $R$, there is a functor $I^\bullet$ from the category $\mathsf{Cob}_{3+1}^{\Bbb Q,-}$ to the category of $\Bbb Z/8$-graded modules over the ring $H^{-*}(BSO(3); R)$ satisfying a similar property.
\end{theorem-intro}

As with Theorem \ref{introthm:main}, there is an extension of this theorem to a more general class of 3-manifolds and $SO(3)$-bundles, which is stated in Section \ref{section:overview}. A special case of Theorem \ref{introthm:functorial} is established in \cite{M} for {\it unobstructed cobordisms}; in the next section, we review the scope of the already-known functoriality results and describe our strategy to extend this to the level of generality stated above. In fact, this strategy is closely related to the method used in the proof of Theorem \ref{introthm:main}. We also expect that the assumption on invertibility of $2$ in Theorem \ref{introthm:functorial} can be dropped. In general, one expects that instanton theory provides a functor $I^\bullet$ from $\mathsf{Cob}_{3+1}^{\Bbb Q}$ for any PID $R$.

\begin{remark*}
The methods in this paper are mostly unrelated to the specifics of the instanton theory. We expect that a similar set of arguments and constructions could be applied in the Seiberg-Witten setting, which would give a new way to define the $S^1$-equivariant Seiberg-Witten Floer homology groups, at least in the case of rational homology spheres. The counterparts of the equivariant homology groups from Theorem \ref{introthm:main} most likely coincide with those defined in \cite{KM}. For a rational homology sphere $Y$, there is an analogue of irreducible instanton homology, $\textrm{HM}^{\text{irr}}_*(Y, \mathfrak s, \pi)$, defined as the homology of a complex generated by irreducible solutions to the 3-dimensional Seiberg-Witten equations. We expect that an analogue of Theorem \ref{introthm:wallcross} holds for these irreducible monopole Floer groups; the authors are not aware of such a result in the literature. 
\end{remark*}

\subsection*{Future directions}
These results leave open a handful of foundational questions which are beyond the scope of this article. First off, the functoriality result in Theorem \ref{introthm:functorial} asserts that $I^\bullet(Y;R)$ is functorial under negative-definite cobordisms for any PID $R$, and for arbitrary cobordisms when $1/2 \in R$.

\begin{problem}\label{problem:any-R}
For any PID $R$, show that $I^\bullet(Y;R)$ is functorial for arbitrary cobordisms.
\end{problem}

We expect that this problem can be attacked using the ideas of this article (finding an algebraic home for obstructed gluing theory). On the other hand, extending functoriality on the object level seems significantly more difficult.

\begin{problem}
Extend the definition of $I^\bullet(Y; R)$ to allow for 3-manifolds with $b_1(Y) > 0$. 
\end{problem}

If $b_1(Y) \le 2$ or $H_1(Y;\Bbb Z)$ has no torsion except 2-torsion, then as in \cite[Proposition 3.9]{AB1} one can achieve transversality for moduli spaces on the cylinder $\Bbb R \times Y$, and it seems plausible that our proofs of invariance and functoriality extend to this setting. For $b_1(Y) > 2$, there will be obstructed reducibles on the cylinder; it already becomes difficult to define the differential on $\widetilde C(Y)$. It is plausible that the ideas of this paper, used here to construct cobordism maps for obstructed cobordisms, could be helpful to make progress toward this goal. This will require a better understanding of the Uhlenbeck compactification. Understanding the Uhlenbeck compactification is also the main difficulty in making sense of Donaldson's polynomial invariants \cite{Don:polynomial, DK} for negative definite cobordisms. In the case of non-negative definite cobordism between integer homology spheres, polynomial cobordism maps are constructed in \cite{Don}.

\begin{problem}
Extend $I^\bullet(Y; R)$ to a \emph{polynomial functor}.
\end{problem}

There still remains the problem of developing good tools for calculation.  Complete computations of the invariants of Theorem \ref{introthm:main} only exist for some spherical manifolds. Two fundamental tools in Floer theory would help to remedy this situation.

%There still remains the problem of developing good tools for calculation. Even with the recent calculations for binary polyhedral 3-manifolds \cite[Section 6]{Helle}, complete computations only exist for spherical manifolds (and only some of them at that). Two fundamental tools in Floer theory would help to remedy this situation.

\begin{problem}
Construct a connected-sum quasi-isomorphism $\widetilde C(Y) \otimes \widetilde C(Y') \simeq \widetilde C(Y \# Y')$, and use this to derive connected-sum formulas for $I^\bullet$ and other related invariants.
\end{problem}
\begin{problem}
	Construct surgery exact triangles relating the equivariant instanton homology of Dehn surgeries on a knot. 
\end{problem}

Fr\o yshov used instanton homology for integer homology spheres to define a numerical invariant $h$ of integer homology spheres that is invariant with respect to homology cobordism relation \cite{Fr1}. Similar to its counterparts in monopole Floer homology \cite{Fr:SW-4-man-bdry,KM} and Heegaard Floer homology \cite{OzSz:d-inv}, Fr\o yshov's $h$-invariant has many interesting topological applications. In compare to $h$, it is easier to work with the invariant of \cite{OzSz:d-inv}, called the correction term, mainly because the correction term can be generalized to rational homology spheres and its behavior with respect to Dehn surgery is better understood. The constructions of this paper provide the essential ingredients to address the following question.

\begin{problem}\label{problem:h-invariant}
	Find analogues of Fr\o yshov's $h$-invariant for rational homology spheres.
\end{problem}

The Chern--Simons functionals on the space of $SU(2)$-connections can be used to define an additional algebraic structure on equivariant instanton homology. In the absence of a perturbation of the Chern--Simons, this gives rise to a filtration of equivariant instanton complexes\footnote{In general we need to pick a non-trivial perturbation of the Chern--Simons functional to define instanton homology, and this additional structure might not be quite a filtration. A more precise formulation of a closely related structure is given in \cite[Subsection 7.3]{DS:equiv-asp-sing}.}. In this paper, we do not attempt to flesh out the details of this construction. However, it is straightforward to see that the method of this paper implies that the additional structure on $\widetilde C(Y)$ given by the Chern--Simons functional is also an invariant of $Y$. In particular, we expect that this structure can be employed to generalize $\Gamma_Y$ from \cite{D}, which is a refinement of Fr\o yshov's invariant, to any rational homology sphere. Similarly, one should be able to generalize homology cobordism invariants of \cite{NST:fil-HCG} using a similar strategy.

In the spirit of Problem \ref{problem:any-R}, we should mention that the structure of the invariants are especially intricate over the ground ring $R = \Bbb F_2$; to compare, when $1/2 \in R$, the invariants described above are modules over $H^{-*}(BSO(3); R) \cong R[U]$, where $|U| = 4$. On the other hand, we have $H^{-*}(BSO(3); \Bbb F_2) \cong \Bbb F_2[w_2, w_3]$, where $|w_i| = i$. 

\begin{problem}
Find appropriate modifications of the previous invariants to the setting of $\Bbb F_2$ coefficients, and develop their properties.
\end{problem}

We hope to address some of these problems in future work.

{\bf Acknowledgements.} The authors would like to thank Christopher Scaduto and Piotr Suwara for interesting conversations about the present work. The authors are grateful to the Simons Center for Geometry and Physics and the organizers of the 2019 PCMI Research Program “Quantum Field Theory and Manifold Invariants” where parts of the work were carried out. 
\newpage

%!TEX root = equivariant-functoriality.tex

\section{An overview of the article}\label{section:overview}

In this section, we review the previous work on (equivariant) instanton homology of rational homology spheres \cite{M,AB2}. Then we review the novel techniques of the present paper and how they can be combined with the previous works to establish the theorem in the introductions. In fact, we can prove these results in the slightly more general context of {\it weakly admissible bundles}. 

\subsection*{Weakly admissible bundles}
A $U(2)$-bundle on a 3-manifold $Y$ is determined by its first Chern class. If $w$ is an {\it oriented 1-cycle} (an embedded oriented 1-manifold) in $Y$, we may form in a canonical way a $U(2)$-bundle $E = E_w$ over $Y$ such that $c_1(E)= {\rm PD}(w)$ and so that $E$ carries a canonical trivialization over $Y \setminus w$. A pair $(Y,w)$ is admissible if there is an embedded surface $\Sigma$ in $Y$ such that $w\cdot \Sigma$ is odd. This pair is weakly admissible if either $b_1(Y)=0$ or $(Y,w)$ is admissible.

We fix a connection $B_0$ on $\det(E)$ and define $\mathcal A(Y,w)$ as the space of all $U(2)$ connections on $E$ that induce the connection $B_0$ on $\det(E)$. This is an affine space modeled on $\Omega^1(\text{ad} (E))$. Here $\text{ad} (E)$ is the vector bundle associated with $E$ corresponding to the adjoint action of $U(2)$ on the Lie algebra $\frak{su}(2)$, and $\Omega^1(\text{ad} (E))$ denotes the space of 1-forms on $Y$ with values in $\text{ad} (E)$. We fix a basepoint $y$ on $Y$ and a trivialization of $\det(E)$ over $y$. A {\it framed connection} $\widetilde A = (A, \Phi)$ on $E$ is a pair of a connection $A$ and a $U(2)$-equivariant trivialization $\Phi$ of the fiber of $\text{ad} (E)$ above $y$. We write $\widetilde {\mathcal A}(Y,w)$ for the space of all framed connections on $E$. Let $\mathcal G(Y,w)$ be the space of all automorphisms of $E$ whose fiberwise determinant is $1$. This space acts on $\mathcal A(Y,w)$ and $\widetilde {\mathcal A}(Y,w)$ by taking pullback of connections (and their framing), and the quotient spaces are respectively denoted by $\mathcal B(Y,w)$ and $\widetilde {\mathcal B}(Y,w)$; the gauge group action is free on $\widetilde{\mathcal A}(Y, w)$, so that $\widetilde{\mathcal B}(Y, w)$ is an infinite-dimensional manifold. These quotients are called the configuration spaces of connections and framed connections associated with $(Y,w)$. Changing the framing of framed connections induces an action of $SO(3)=SU(2)/\{\pm I\}$ on $\widetilde {\mathcal B}(Y,w)$ with the quotient space being $\mathcal B(Y,w)$. 

There are three different possibilities for the stabilizer $\Gamma_A$ of a connection $A \in \mathcal A(Y,w)$ with respect to the action of $\mathcal G(Y,w)$:
\begin{itemize}
	\item[(i)] $\Gamma_A = \{\pm I\}$ where $I$ is the trivial automorphism of $E$.  Any such connection is called {\it irreducible}. The framed connections with the underlying connection $A$ determines an orbit of the $SO(3)$ action of 
	$\widetilde {\mathcal B}(Y,w)$ which is diffeomorphic to $SO(3)$. We denote such orbits by the initial letters of the Greek alphabet $\alpha$, $\beta$ and $\gamma$. 
	\item[(ii)] $\Gamma_A \cong S^1$.  Any such connection is called {\it abelian}. The framed connections with the underlying connection $A$ determines an $SO(3)$-orbit diffeomorphic to $S^2=SU(2)/S^1$ in 
	$\widetilde {\mathcal B}(Y,w)$. We usually denote such orbits by $\rho$. 
	\item[(iii)] $\Gamma_A \cong SU(2)$.  Any such connection is called {\it central}. The framed connections with the underlying connection $A$ determines an $SO(3)$-orbit consisting of one point in 
	$\widetilde {\mathcal B}(Y,w)$. We usually denote such orbits by $\theta$. 
\end{itemize}
In the case that $A$ is not irreducible, we say $A$ is reducible.

We may define a Chern-Simons functional on the configuration spaces of connections associated with $(Y,w)$. Fix a connection $A_0\in \mathcal A(Y,w)$ and let $A$ be any other connection in $\mathcal A(Y,w)$. Let $\bA$ be a connection on the pullback of $E$ to $\R\times Y$ such that the restrictions of $\bA$ to $(-\infty,-1]\times Y$ and $[1,\infty)\times Y$ are respectively equal to the pullbacks of $A$ and $A_0$. Then we defined the value of the Chern--Simons functional at $A$ to be the Chern--Weil integral 
\[
  \frac{1}{8\pi^2}\int_{\R\times Y} \tr(F_{\bA}\wedge F_{\bA}),
\]
where $F_\bA$ denotes the curvature of $\bA$. Varying the connection $\bA$ by a homotopy does not change the value of the above integral, and changing $A$ by an element of $\mathcal G(Y,w)$ changes the above expression by an integer. In particular, we have a well-defined functional $CS:{\mathcal B}(Y,w)\to\Bbb R/\Bbb Z$ that is independent of $A_0$ up to an overall translation. Precomposing $CS$ with the map $\widetilde {\mathcal B}(Y,w)\to {\mathcal B}(Y,w)$ which forgets the framing defines a functional on $\widetilde {\mathcal B}(Y,w)$, which is still called the Chern--Simons functional.

A connection $A$ represents a critical point of the Chern--Simons functional if the induced adjoint connection $\ad(A)$ on the $SO(3)$-bundle $\ad(E)$ is flat. With a slight abuse of language, we call any such connection a flat connection associated with $(Y,w)$. 
We denote the space of all such flat connections by $\fC(Y,w)$, and write this set as a disjoint union based on the three possibilities for the stabilizers of flat connections:
\[
  \fC(Y,w) = \fC^*(Y,w) \sqcup \fA(Y,w) \sqcup \fZ(Y,w).
\]
We can identify $\fC(Y,w)$ with the character variety
\[
  \chi(Y,w):=\{\phi:\pi_1(Y\setminus w)\to SU(2) \mid \text{$\phi$ maps any meridian of $w$ to $-I$}\}/SU(2).  
\]
Here the quotient is given by the conjugation action. 

The following lemma asserts that reduicble flat connections associated with $(Y,w)$ are determined by the cohomology groups of $Y$.

\begin{lemma}\label{lemma:enum-red-3d}
	If $(Y,w)$ is an admissible pair, then there is no reducible flat connection associated with $(Y,w)$. If $b_1(Y) = 0$, then the set of reducible flat connections form a discrete subset of $\mathcal B(Y,w)$, and are characterized as follows.
	\begin{itemize}
		\item[(i)] For any element $\theta$ of $\mathfrak Z(Y, w)$, there is a $\theta$-parallel splitting of $E$ into a direct sum of complex line bundles $L_\theta^{\oplus 2}$, where the two summands are canonically isomorphic. The map $\theta \mapsto c_1(L_\theta)$ defines a canonical bijection
		\[
		  \mathfrak Z(Y,w) \cong \{\xi \in H^2(Y;\Bbb Z) \mid 2\xi ={\rm PD}(w)\}.
		\]
		Moreover, when $\mathfrak Z(Y,E)$ is nonempty, the group $H^1(Y;\Bbb Z/2)$ acts on it freely and transitively via $x \cdot \xi = \xi + \beta x$ with $\beta$ being Bockstein homomorphism.

		\item[(ii)] For any element $\rho$ of $\mathfrak A(Y, w)$, there is a $\rho$-parallel splitting of $E$ as a direct sum of $S^1$-bundles $L_\rho$ and $L_\rho'$. The map $\rho \mapsto \{c_1(L_\rho), c_1(L_\rho')\}$ defines a canonical bijection
		\[
		  \mathfrak A(Y,w) \cong \big\{\{x, y\}\subset H^2(Y;\Bbb Z) \mid x + y = {\rm PD}(w), \;\; x \neq y \big\}.
		\]
	\end{itemize}
\end{lemma}

\begin{comment}
		\item[(i)] For any element $\theta$ of $\mathfrak Z(Y, w)$, there is a canonical $S^1$-bundle $L_\theta$ with an $S^1$-connection $B_\theta$ such that $L_\theta^{\otimes 2}$ is isomorphic to $\det(E)$
		and $B_\theta$ is the unique connections that the induced connection on $L_\theta^{\otimes 2}$ is $B_0$. 
		Then the connection $\theta$ is equal to $B_\theta\oplus B_\theta$ on $L_\theta\oplus L_\theta\cong E$. The map $\theta \mapsto c_1(L_\theta)$ defines a canonical bijection
		\[
		  \mathfrak Z(Y,w) \cong \{\xi \in H^2(Y;\Bbb Z) \mid 2\xi ={\rm PD}(w)\}.
		\]
		Moreover, when $\mathfrak Z(Y,E)$ is nonempty, the group $H^1(Y;\Bbb Z/2)$ acts on it freely and transitively via $x \cdot \xi = \xi + \beta x$ with $\beta$ being Bockstein homomorphism.

		\item[(ii)] For any element $\rho$ of $\mathfrak A(Y, w)$, there are canonical non-isomorphic $S^1$-bundles $L_\rho$ and $L_\rho'$ with $S^1$-connections $B_\rho$ and $B_\rho'$ such that $L_\rho\oplus L_\rho'$ is isomorphic 
		to $E$, and 
		$B_\rho$, $B_\rho'$ are the unique connections (up to automorphisms of $S^1$-bundles) that the induced connection on $L_\rho\otimes L_\rho'$ is $B_0$ and the induced connection on $L_\rho^*\oplus L_\rho'$ is flat.
		In particular, we have the canonical bijection
		\[
		  \mathfrak A(Y,w) \cong \big\{\{z_1, z_2\}\subset H^2(Y;\Bbb Z) \mid z_1 + z_2 = {\rm PD}(w), \;\; z_1 \neq z_2 \big\}.
		\]
\end{comment}

\subsection*{Perturbations and signature data}
To define something like a Morse complex for an $SO(3)$-invariant functional, we need to assume its Hessian at any critical orbit has kernel equal to the tangent space at that orbit, in which case we say it is non-degenerate. This is usually not true for the Chern--Simons functional, so we will need to add a perturbation term.

The standard perturbations of the Chern--Simons functional is given by {\it holonomy perturbations} \cite[Formula (1b.9)]{Fl}. Suppose $CS+\pi:\mathcal B(Y,w) \to \Bbb R/\Bbb Z$ is a perturbation of the Chern--Simons functional. A critical point of $CS+\pi$ is non-degenerate if the Hessian of the functional $CS+\pi$ at this connection, which can be regarded as an elliptic operator, has trivial kernel. (See Subsection \ref{hessian} for more details on the definition of this operator.) 

When $b_1(Y) = 0$, every element of $\fZ(Y,w)$ is non-degenerate. Moreover, the restriction of the Chern--Simons functional to the locus of abelian connections is non-degenerate. Motivated by this, we only consider perturbations $\pi$ of the Chern--Simons functional which vanish on the abelian locus and on a neighborhood of the central locus. For any such perturbation $\pi$, reducible critical points of $CS+\pi$ agree with that of $CS$, and hence we have the following splitting of the critical points $\fC_\pi(Y,w)$ of $CS+\pi$:
\[
  \fC_\pi(Y,w) = \fC_\pi^*(Y,w) \sqcup \fA(Y,w) \sqcup \fZ(Y,w).
\]
We can always find a perturbation term $\pi$ as above such that all elements of $\fC_\pi(Y,w)$ are non-degenerate \cite[Theorem 3.8]{M}. We call any such perturbation {\it nice}.

We associate a {\it signature data function} $\sigma_\pi:\fA(Y,w) \to 2\Bbb Z$ to any nice perturbation $\pi$. Any $\rho$ in $\fA(Y,w)$ is a critical point of $CS+t\pi$ for any $t\in [0,1]$. Thus the Hessian $H_t(\rho)$ of $CS+t\pi$ at the critical point $\rho$ defines a 1-parameter family of Fredholm operators that are perturbations of elliptic operators with compact operators and act on the same space. Then $\sigma_\pi(\rho)$ is defined as 
\[\sigma_\pi(\rho) = \dim \ker H_0(\rho) - 2\text{sf } \{H_{1-t}(\rho)\}_{t\in [0,1]},\]
the latter term being the spectral flow of the family ${H_{1-t}(\rho)}_{t \in [0,1]}$. When $\pi$ is a small perturbation, this quantity can be understood as the signature of the restriction of $H_1(\rho)$ to $\ker H_0(\rho)$, which is a nondegenerate symmetric bilinear form. The quantity $\sigma_\pi(\rho)$ is equal modulo $4\Bbb Z$ to $\dim H_0(\rho)$.

\subsection*{Review of equivariant instanton homology}
The basic object that instanton theory associates to a weakly admissible pair $(Y,w)$ is the \emph{tilde complex}\footnote{This object is called the `framed instanton complex' of $(Y,w)$ in \cite{M}. We use the term `tilde-complex' in this paper to avoid confusion with the Kronheimer and Mrowka's $I^\sharp(Y)$.}\cite{M}. Suppose $\pi$ is a nice perturbation for $(Y,w)$ and $R$ is a ring. We may associate a relatively $\Bbb Z/8$-graded chain complex of $R$-modules 
\begin{equation}\label{SO-com-pair}
	\widetilde C(Y, w, \pi; R),
\end{equation}	
equipped with an action of the dg-algebra 
\begin{equation}\label{SO(3)-alg}
	C_*(SO(3); R)
\end{equation} 
of chains on the Lie group $SO(3)$. When there is no risk of confusion we drop $\pi$ and $R$ from our notation. Any chain complex with this dg-module structure is called an $SO(3)$-complex.

The underlying $SO(3)$-module of $\widetilde C(Y, w, \pi;R)$ is defined in terms of the critical points of $CS+\pi$. The complex $\widetilde C(Y, w, \pi)$ is the direct sum of $C_*(\alpha)$, where $\alpha$ runs over the critical orbits of $CS+\pi$ in the framed configuration space $\widetilde {\mathcal B}(Y,w)$. Thus $\widetilde C(Y, w, \pi)$ contains a summand $C_*(SO(3))$ for any element of $\fC_\pi^*(Y,w)$, a summand of $C_*(S^2)$ for any element of $\fA(Y,w)$ and a summand of $C_*({\rm point})$ for any element of $\fZ(Y,w)$. To define these $R$-modules, we use a certain geometric models for chains where $C_*(\alpha)$ is generated by $\phi$ where $\phi:P\to \alpha$ is a smooth map from a space $P$ into $\alpha$ with $P$ being a structured space analogous to a smooth manifold with corners. The definition of {\it geometric chains} used in our construction will be discussed in Section \ref{sec:gchain-intro} and the foundational material carefully developed in Appendix \ref{A}. The action of $C_*(SO(3))$ on $\widetilde C(Y, w, \pi)$ is induced by the $SO(3)$-action on $\widetilde {\mathcal B}(Y,w)$.

For each two orbits $\alpha, \beta$, (perturbed) {\it anti-self-duality} (ASD) equation determines a compact moduli space $\breve M^+_\zeta(\alpha, \beta)$, equipped with an $SO(3)$-action and $SO(3)$-equivariant {\it restriction} maps to $\alpha$ and $\beta$. For a generic nice perturbation, these moduli spaces are smooth(ish), and any such perturbation is called {\it regular}. The appropriate analogue of the usual point-counting differential in Morse theory is instead a \emph{fiber product differential}; for any $\phi:P\to \alpha$ representing an element of $C_*(\alpha)$, we have
\[
  \widetilde \partial \phi = \partial_{gm} \phi + \sum_\beta \pm P \times_\alpha \breve M^+_\zeta(\alpha, \beta),
  \]
 where the first term is the geometric boundary of $\phi$ (again a chain in $\alpha$), while each term in the sum is a chain in some $C_*(\beta)$. The fact that this operator squares to zero is encoded in boundary relations for the moduli spaces $\breve M^+_\zeta(\alpha, \beta)$; the boundary components of $\breve M^+_\zeta(\alpha, \beta)$ may be described as fiber products of smaller-dimensional moduli spaces. The definition of $\breve M^+_\zeta(\alpha, \beta)$ and the determination of its boundary strata are recalled in Subsection \ref{comp-moduli}, and the precise definition of the differential $\widetilde \partial$ is given in Subsection \ref{inst-cx-unob-cob-map}.  This construction deserves to be called the \emph{Morse-Bott complex} of $CS+\pi$ on $\wt{\mathcal B}(Y, w)$; a de Rham version was initially developed in \cite{AB2}, together with its equivariant homology.

A limited version of functoriality of tilde complexes with respect to cobordisms is established in \cite{M}. First let $(Y, w,\pi)$ and $(Y', w',\pi')$ be a weakly admissible pairs with chosen regular perturbations. Suppose $W:Y \to Y'$ is a cobordism and $c$ is a properly oriented embedded surface in $W$ such that the intersections of $c$ with $Y$ and $Y'$ are respectively $w$ and $w'$. Then $c$ canonically determines a $U(2)$ bundle $\bE_c$ over $W$ such that $c_1(\bE_c)= {\rm PD}(c)$. In particular, the restrictions of $\bE$ to $Y$ and $Y'$ are the $U(2)$-bundles $E_w$ and $E_{w'}$, respectively. The cobordism is also equipped with a path $\gamma$ in $W$ running between the chosen basepoints of $Y$ and $Y'$. As with the basepoints, we usually suppress $\gamma$ from our notation for cobordisms.

The ASD equation on $W$ (after adding cylindrical ends) gives rise to a moduli space $M_z^+(W,c; \alpha, \alpha')$ for any pair of an orbit of $\pi$-flat connection $\alpha$ on $Y$ and an orbit of $\pi'$-flat connection $\alpha'$ on $Y'$. An element of this moduli space is the gauge equivalence of a pair $(A, \Phi)$ of a (possibly broken) connection $A$ on $W$, which satisfies a perturbation of the ASD equation compatible with the perturbation terms on the ends, and a framing $\Phi$ along $\gamma$ which is parallel with respect to $A$. One might hope that after picking an appropriate perturbation of the ASD equation,  $M_z^+(W,c; \alpha, \alpha')$ is a compact smooth(ish) space. However, there is an obstruction which can be expressed in terms of $\sigma_{\pi}$, $\sigma_{\pi'}$, the topology of $W$ and the cohomology class ${\rm PD}(c)$. when this obstruction is trivial, we say $(W,c,\pi,\pi')$ is {\it unobstructed}. Otherwise, it is obstructed. The precise definition of $M_z^+(W,c; \alpha, \alpha')$ and unobstructedness of $(W,c,\pi,\pi')$ will be given in Subsections \ref{cob-basics} and \ref{red-reg} respectively. As an important example, let $\pi$ and $\pi'$ be two regular perturbations for a weakly admissible pair $(Y,w)$. Then $([0,1]\times Y,[0,1]\times w)$ as a cobordism from $(Y,w)$ with perturbation term $\pi$ to $(Y,w)$ with perturbation term $\pi'$ is obstructed if and only if $\sigma_{\pi}\leq \sigma_{\pi'}$.

The second main result of \cite{M} is that when $(W,c,\pi,\pi')$ is an unobstructed cobordism with a path between the baespoints, a choice of an appropriate perturbation of the ASD equation and a homology orientation for $W$ (an orientation of the vector space $H^1(W) \oplus H^{2,+}(W) \oplus H^1(Y)$), we may define a homomorphism 
\[
  \widetilde C_*(W,c): \widetilde C_*(Y,w,\pi) \to \widetilde C_*(Y',w',\pi').
\]
This homomorphism is $SO(3)$-equivariant and is independent of the Riemannian metric on $W$ and the perturbation of the ASD equation up to equivariant chain homotopy. 

Finally, the third main result of \cite{M} is that given an algebraic object resembling $\widetilde C_*(Y, w, \pi)$ (more precisely, a relatively graded dg-module equipped with a \emph{periodic filtration}), one may construct \emph{equivariant homology groups} which fit into an exact triangle. Applying this to $\widetilde C(Y, w, \pi)$ one obtains the equivariant instanton homology groups $I^+, \; I^-, \; I^\infty$ of $(Y, w,\pi)$. 

The definition of $I^\infty$ in \cite[Section A.8]{M} is flawed; the complex $CI^\infty$ is defined in terms of a completion of an algebraic construction. However, the completion process described there does not have all the properties stated. Instead, one should use the more robust \emph{full completion}, as in \cite[Section 5.4]{Helle}; using this modified notion of completion, all the stated results of \cite[Section A.8]{M} hold.

Usually, one proves that Floer homology groups are independent of perturbation, and hence topological invariants, by employing functoriality with respcet to cobordisms. Unfortunately, the unobstructedness assumption is very restrictive. Such an argument can be used to show that equivariant instanton homology groups are invariants of $(Y,w,\sigma)$, where $\sigma$ is the signature data of the perturbation. To show that these Floer homology groups do not depend on the choice of signature data in some subtle way, one also needs cobordism maps associated with the product cobordism $([0,1]\times Y,[0,1]\times w)$ when the signature data does not necessarily increase from the incoming end to the outgoing end. The primary purpose of the present paper is to extend the relevant functoriality properties to remedy this issue.

\subsection*{Results of this article}
Suppose $(W,c):(Y,w)\to (Y',w')$ is a cobordism with $b_1(W)=b^+(W)=0$. We fix a regular perturbation $\pi_i$ for $(Y_i,w_i)$. We also fix a Riemannian metric with cylindrical ends on $W$. The issue of obstructedness for $(W,c,\pi,\pi')$ is related to the {\it reducible instantons}, the elements $\Lambda$ in $M_z^+(W,c; \alpha, \alpha')$ which have non-trivial stabilizer with respect to the $SO(3)$ action. As is explained in Proposition \ref{prop:enum-red}, any such reducible $\Lambda$ is determined by some cohomological data. Any such reducible $\Lambda$ gives rise to a quantity $N(\Lambda,\pi,\pi') \in 2\Bbb Z$ called the \emph{normal index}. The expected dimension of the moduli space $M_z^+(W,c; \alpha, \alpha')$ containing $\Lambda$ is equal to $N(\Lambda,\pi,\pi')+2$. The quantity $N(\Lambda,\pi,\pi')$ can be computed in terms of some cohomological data and the signature data of $\pi$ and $\pi'$ (see Proposition \ref{prop:normal-ind}). In particular, if $\Lambda$ is asymptotic to the abelian reducibles $\rho$ and $\rho'$ on the two ends, then for regular perturbations $\pi_1$ and $\pi_2$ for $(Y,w)$ and regular perturbations $\pi_1'$ and $\pi_2'$ for $(Y',w')$ we have
\begin{equation}\label{nor-ind-change-per}
  N(\Lambda,\pi_1,\pi_1')-N(\Lambda,\pi_2,\pi_2')= \(\sigma_{\pi_1}(\rho)-\sigma_{\pi_2}(\rho)\)-\(\sigma_{\pi_1'}(\rho')-\sigma_{\pi_2'}(\rho')\).
\end{equation}
In the case that $\rho$ or $\rho'$ is central, the corresponding term on the right-hand side of the above identity vanishes. The tuple $(W,c,\pi,\pi')$ is unobstructed exactly when $N(\Lambda,\pi,\pi')$ is non-negative for all possible choices of $\Lambda$.

\begin{definition}\label{def:nearly-unobst}
	We say $(W,c,\pi,\pi')$ is {\it nearly unobstructed} if there is one reducible instanton $\Lambda$ with $N(\Lambda,\pi,\pi')=-2$ and all the remaining reducibles have non-negative normal index.
\end{definition}
\begin{example}\label{adjacent signature data}
	Suppose $(\pi_1,\pi_2)$ are two regular perturbations of $(Y,w)$ with {\it adjacent signature data}, meaning that $\sigma_{\pi_1}=\sigma_{\pi_2}$ except at one abelian flat connection $\rho_0$, where 
	$\sigma_{\pi_1}(\rho_0)=\sigma_{\pi_2}(\rho_0)+4$. Then $(I \times Y, I \times w,\pi_1,\pi_2)$ is nearly unobstructed. On the other hand, $(I \times Y, I \times w,\pi_2,\pi_1)$ is unobstructed.
\end{example}

The following is our main technical theorem which allows us to obtain the results advertised in the introduction. The proof of this theorem will be given in Section \ref{obscob}.
\begin{theorem}\label{intro:induced-map}
	Any weakly admissible pair $(Y,w)$ with a regular perturbation $\pi$ and a choice of abelian flat connection $\rho$ on $Y$ gives rise to an $SO(3)$-complex $S_{\rho}\widetilde C(Y,w, \pi)$ and an $SO(3)$-equivariant quasi-isomorphism 
	\[
	  S_*: \widetilde C(Y,w, \pi) \to S_{\rho}\widetilde C(Y,w, \pi).
	\]  
	Suppose $(W,c): (Y, w) \to (Y',w')$ is a cobordism as above with path between the basepoints and $\pi$, $\pi'$ are regular perturbations for $(Y,w)$, $(Y',w')$ such that $(W,c,\pi,\pi')$ is nearly unobstructed with exactly one obstructed reducible $\Lambda$ running between abelian 
	flat connections $\rho$ and $\rho'$. Then there is an equivariant chain map 
	\[
	  \widetilde{C}(Y,w,\pi) \to S_{\rho'} \widetilde C(Y',w', \pi'),
	\] 
	which is well-defined up to equivariant homotopy.
\end{theorem}

Since $SO(3)$-equivariant quasi-isomorphisms induce isomorphisms on equivariant homology groups, the map $S_*$ may be inverted at the level of equivariant homology. It follows that nearly unobstructed cobordisms give rise to maps on equivariant instanton homology. In Section \ref{sec:CW} we use the above theorem to prove Theorem \ref{introthm:main} from the introduction.

The theorem above, together with an analysis of composites, is nearly enough to prove the invariance of equivariant instanton homology. Suppose $(\pi_1,\pi_2)$ are two regular perturbations of $(Y,w)$ with {\it adjacent signature data}. Then Theorem \ref{intro:induced-map} gives 
\begin{equation}\label{wrong-way-map}
   \widetilde C_*(Y,w, \pi_1) \to S_\rho \widetilde C_*(Y,w, \pi_2).
\end{equation}
Moreover, $(I \times Y, I\times w,\pi_2,\pi_1)$ is unobstructed and the construction of \cite{M} gives a cobordism map $\widetilde C_*(Y,w,\pi_1) \to \widetilde C_*(Y,w,\pi_2)$ and in Subsection \ref{cob-flow} we shall show that this map can be lifted to a homomorphism 
\[
   S_\rho \widetilde C_*(Y,w,\pi_1) \to \widetilde C_*(Y,w,\pi_2);
\]
in Subsection \ref{sec:invt} we will show that it is an equivariant homotopy inverse to the homomorphism in \eqref{wrong-way-map}. Finally, for any two regular perturbations $\pi, \pi'$ of $(Y,w)$, one may always find a sequence 
\[
  \pi = \pi_0 < \cdots < \pi_n = \pi_m' > \cdots > \pi_0' = \pi'
\] 
of regular perturbations such that $(\pi_i,\pi_{i+1})$ and $(\pi_{j+1}',\pi_j')$ have adjacent signature data. This completes the proof of the following generalization of Theorem \ref{introthm:main}.

\begin{theorem}\label{thm:main-ex}
	For any weakly admissible pair $(Y,w)$ and any PID $R$, there are relatively $\Bbb Z/8$-graded groups, 
	\[
	  \widetilde I(Y, w;R), \hspace{1cm} I^+(Y,w;R), \hspace{1cm} I^-(Y,w;R), \hspace{1cm} I^\infty(Y,w;R),
	\] 
	which are invariants of the oriented diffeomorphism type of $Y$ and the element of $H^2(Y;\Bbb Z/2)$ determined by $w$.
	The first is a graded module over the ring $H_*(SO(3);R)$. The next three are graded modules over the ring $H^{-*}(BSO(3); R)$ and fit into an exact triangle 
	\[
	  \cdots \xrightarrow{[-4]} I^+(Y,w;R) \xrightarrow{[3]} I^-(Y,w;R) \to I^\infty(Y,w;R) \xrightarrow{[-4]} \cdots
	\] 
	in which the morphism $I^- \to I^\infty$ is an $H^{-*}(BSO(3); R)$-module homomorphism, as are the other two when $2 = 0 \in R$ or $2 \in R^\times$, and the grading shifts are indicated above each arrow. 
\end{theorem}

We will suppress $R$ from notation when the particular choice of ring is irrelevant to us.

\begin{remark}\label{rmk:compare}
When $(Y, w)$ is admissible, we have $I^\infty(Y, w) = 0$ and $I^+(Y,w) \cong I^-(Y,w)[-3]$ are canonically isomorphic to Floer's instanton homology groups $I_*(Y,w)$ for admissible bundles, while $\widetilde I(Y, w) = 0$.

When $Y$ is a rational homology sphere, we anticipate a relationship to Kronheimer and Mrowka's $I^\sharp(Y,w)$ --- roughly, $I_*$ of $(Y \# T^3, w \sqcup S^1)$ --- which cannot be extracted directly from the homology groups above, but rather from the $SO(3)$-complex $\widetilde C$. Precisely, writing $u$ for the action of $[SO(3)] \in C_3(SO(3);R)$ on $\widetilde C$, a connected sum theorem should give an isomorphism between $\Bbb Z/4$-graded groups $H\left(\widetilde C(Y,w), \widetilde d -4u\right) \cong I^\sharp(Y, w)$. 
\end{remark}

Our next goal is to define cobordism maps for more general cobordisms (arbitrary cobordisms when $1/2 \in R$, and arbitrary negative-definite cobordisms otherwise). Were it not for difficulties with the Uhlenbeck compactification, it seems plausible that the techniques used in the proof of Theorem \ref{intro:induced-map} could be used to construct induced maps for arbitrary cobordisms. Instead, we opt for a different approach. First we must be slightly more specific about the sort of issues one runs into in the definition of cobordism maps. When enumerating the reducible instantons on $W$, they come in three essentially different types:
\begin{itemize}
\item Abelian instantons $\Lambda$ on $W$ which are not central on the boundary. 
\item Abelian instantons $\Lambda$ which restrict to central connections on the boundary. In this case the normal index of $\Lambda$ is independent of the perturbations $\pi$ and $\pi'$. If $\Lambda$ is flat, then the normal index of $\Lambda$ is $-2$.  Otherwise the normal index is non-negative and it is not a cause for concern. Flat abelian instantons which are central on the boundary are called \emph{pseudocentral instantons}.
\item Central connections, which are necessarily flat. These are cut out transversely when $b^+(W) = 0$ and not otherwise. 
\end{itemize}

The first key insight is that because $I^\bullet(Y,w)$ and $I^\bullet(Y',w')$ are independent of the choice of perturbations $\pi$ and $\pi'$, \emph{we are free to choose the perturbation on the ends however we like}. In practice, this implies that the first type of abelian instanton is irrelevant: we may now choose perturbation data $\pi$ and $\pi'$ so that $N(\Lambda,\pi,\pi') \ge 0$ for any reducible $\Lambda$ of the first type (see the identity in \eqref{nor-ind-change-per}). Therefore, the perturbation term over the cobordism can be chosen so that all such reducibles are cut out transversely. To resolve the issue related to pseudocentral reducibles, we remove a regular neighborhood of the elements of the moduli spaces containing the pseudocentral reducibles. An analysis of the boundary of these spaces allows to define an induced map between equivariant instanton homology groups. Resolving the issues related to the first two types of reducible instantons allows us to prove the second part of Theorem \ref{introthm:functorial}.

Central connections complicate the definition of cobordism maps in the case that $b^+(W)>0$. It seems likely that a variation on the construction used in Theorem \ref{intro:induced-map} will allow us to define the induced maps of such cobordism maps and hence complete the proof of Theorem \ref{introthm:functorial}. However, we do not attempt to carry this out here. Instead, we resolve the issue of obstructed reducibles with a blowup trick used to define Donaldson invariants past the realm they are usually defined; for instance, see \cite{MM,KM-Structure,K-HigherRank}.
More precisely, we replace moduli spaces associated with $(W,c)$ with the moduli spaces $(W \# \overline{\Bbb{CP}}^2,c\cup D)$ cutting down by $\mu(D)$. Here $D$ denotes the exceptional divisor in $\overline{\Bbb{CP}}^2$. The new moduli spaces do not have any central flat connections. Moreover, if the cobordism map for $(W,c)$ is already defined, then the above blowup process gives the same cobordism map up to a multiple of $2$. In particular, if $2$ is a unit in the ring $R$, then this multiple of $2$ is harmless and hence it justifies that the blowup process gives a sensible definition of cobordism maps. This is why we require that $2$ is invertible in the first part of Theorem \ref{introthm:functorial}.

In fact, we can slightly extend Theorem \ref{introthm:functorial} to allow cobordisms between weakly admissible pairs. First we extend the definition of categories $\mathsf{Cob}_{3+1}^{\Bbb Q}$ and $\mathsf{Cob}_{3+1}^{\Bbb Q,-}$. 

\begin{definition}\label{def:CobCat-w}
The category $\mathsf{Cob}_{3+1}^{w}$ has as its objects $(Y, w, y)$, where $(Y,w)$ is a weakly admissible pair and $y \in Y$ is a basepoint. The morphisms $(Y, w, y) \to (Y', w', y')$ are $(W, c, \gamma, o)$, where $(W,c)$ is an oriented connected cobordism $(Y, w) \to (Y', w')$, with $\gamma$ a path in $W$ between the basepoints, and $o$ a homology orientation, all considered up to diffeomorphism. Composition is given by gluing along the common boundary of two cobordisms.

\noindent The category $\mathsf{Cob}_{3+1}^{w,-}$ has the same objects, and those morphisms for which either $b^+(W) = 0$ or $(W,c)$ admits no central connections.
\end{definition}

Then the functors $I^\bullet$ extend to these categories, and the exact triangle is functorial.

\begin{theorem}\label{thm:Functor-w}
	Let $R$ be a PID. 
	\begin{enumerate}
	\item[(i)] There is a functor $I^\bullet$ from the category $\mathsf{Cob}_{3+1}^{w,-}$ to the category of relatively $\Bbb Z/8$-graded modules over the ring $H^{-*}(BSO(3);R)$ that specializes to the invariants of Theorem \ref{thm:main-ex} at the level of objects; the induced maps commute with the structure maps in the exact triangle of Theorem \ref{thm:main-ex}.
	\item[(ii)] If $R$ is a field with $\textup{char}(R) \ne 2$, there is a functor $I^\bullet$ from the category $\mathsf{Cob}_{3+1}^w$ to the category of relatively $\Bbb Z/8$-graded modules over $H^{-*}(BSO(3);R)$ which extends the functor from (i).
	\end{enumerate}
\end{theorem}

\subsection*{Organization and strategy}
In Section \ref{GeneralSetup} below, we recall the relevant facts about ASD connections, perturbations, and moduli spaces thereof. To define the tilde-complex in \eqref{SO-com-pair}, we need to have chain-level fiber product maps with respect to whatever kind of smooth structure the moduli spaces have. Such an apparatus is called a theory of geometric chains. In Appendix A, we define and present the foundations of {\it stratified-smooth spaces}. Our moduli spaces carry the structure of stratified-smooth spaces, and transverse intersections of stratified-smooth spaces are again stratified-smooth. We conclude this Appendix by verifying in \ref{CheckAxioms} that we can use stratified-smooth spaces to construct a homology theory of smooth manifolds canonically isomorphic to singular homology; for technical reasons, we truncate it to a bounded complex in \ref{gm-trunc}.

In Section \ref{geo-chain-inst-hom}, we use this apparatus to define the notion of {\it flow categories}, {\it bimodules between flow categories}, and {\it homotopy between bimodules}. These amount to collections of orbits $\alpha$, as well as collections of geometric chains (corresponding to the moduli spaces $\breve M^+_\zeta(\alpha, \beta)$ in the previous discussion) so that each moduli space satisfies an appropriate boundary relation. Using the above geometric chain apparatus, we associate a {\it flow complex} to a flow category, a chain map of flow complexes to a bimodule and a homotopy of chain maps to a homotopy between bimodules. We then explain how the tilde-complex of \eqref{SO-com-pair} and the homomorphisms associated with unobstructed cobordisms arise from this construction.

Section \ref{obscob} is the technical heart of the paper. Given a nearly unobstructed cobordism $(W,c,\pi,\pi')$, the moduli spaces $M_z^+(W,c; \alpha, \alpha')$ will no longer be smooth manifolds with corners. However, if we delete a neighborhood of the `non-smooth points', we are left with a stratified-smooth space which has new boundary strata, that we must investigate. While these new boundary strata are not explicit, we may correct them by adding a term interpolating to some more standard expressions (in terms of zero sets of sections of an \emph{obstruction line bundle}). Doing so, we define \emph{modified moduli spaces} which satisfy an explicit boundary relation --- though not the right boundary relation to define a chain map between tilde-complexes. 

Investigating these boundary relation leads us to the notion of \emph{suspended} flow category, which we explore in Section \ref{cob-flow}; its flow complex is the $SO(3)$-complex $S_\rho \widetilde C$ discussed above. After doing so, we verify in Section \ref{cob-bimods} that the modified moduli spaces do define bimodules with respect to these suspended flow categories, resulting in Theorem \ref{intro:induced-map}. In Section \ref{sec:invt} we prove that these compose in a reasonable way and thus prove Theorem \ref{thm:main-ex}. 

With this theorem in hand, only a little work remains. We conclude by verifying that negative-definite cobordisms (or cobordisms without central connections) with pseudocentral reducibles still give chain maps well-defined up to homotopy in Section \ref{badred}. In Section \ref{sec:composite} we show that these compose functorially. We then discuss the blowup trick in Section \ref{MM-blowup}, which gives us functoriality for arbitrary cobordisms when our coefficient ring contains $1/2$. 

In the final section, we explain how to use our constructions to extract results about the extension of Floer's {\it irreducible instanton homology} to rational homology spheres, and how to use this to construct a numerical invariant which conjecturally agrees with the Casson-Walker invariant.
\newpage

\section{Moduli spaces of instantons}\label{GeneralSetup}
\subsection{Linear analysis on 3-manifolds}\label{hessian}
%!TEX root = equivariant-functoriality.tex

Let $(Y, w)$ be a weakly admissible pair with the associated $U(2)$-bundle $E$. For a non-negative integer $k$, let $\mathcal A_k(Y, w)$ and $\mathcal G_{k}(Y, w)$ be respectively the $L^2_k$-completions of the space of connections $\mathcal A(Y, w)$ and the gauge group $\mathcal G(Y, w)$ introduced in Section \ref{section:overview}. In particular, $\mathcal A_k(Y, w)$ can be identified with the affine space $A_0 + \Omega^1_k (\text{ad} (E))$ where $A_0$ is a smooth reference connection and $ \Omega^1_k (\text{ad} (E))$ denotes the space of $L^2_k$ sections of the bundle $\Lambda^1\otimes \text{ad}(E)$ over $Y$. Fix an integer $k\geq 2$ for the rest of the paper. The space $\mathcal G_{k+1}(Y, w)$ is a Banach Lie group acting smoothly on $\mathcal A_k(Y, w)$. We may thus define $\mathcal B_k(Y, w)$ to be the quotient of this action. We similarly define the space of framed connections $\widetilde {\mathcal A}_k(Y, w)$ and the configuration space of framed connections $\widetilde {\mathcal B}_k(Y, w)$.

The framed configuration space is a Banach manifold with smooth $SO(3)$-action, and hence has a natural $SO(3)$-equivariant tangent bundle $T_{[A,\Phi]} \wt{\mathcal B}_k(Y, w)$. This may be identified with the quotient of $\Omega^1_k(Y, \text{ad}(E)) \oplus \mathfrak{so}(3)$ by the image of 
\[(d_A, \text{ev}_y): \Omega^0_{k+1}(Y, \text{ad}(E)) \to \Omega^1_k(Y, \text{ad}(E)) \oplus \mathfrak{so}(3).\]
Notice that this tangent bundle has natural Sobolev completions $T_j \wt{\mathcal B}_k(Y, w)$ for all $1 \le j \le k$. (The point-evaluation component requires that our sections of $\text{ad}(E)$ are continuous, and hence the restriction that they are of regularity at least $L^2_2$.)

We are now able to give a careful discussion of the Hessian operator. The formal gradient of the perturbed Chern-Simons functional defines an $SO(3)$-equivariant section 
\[\nabla(CS+\pi): \wt{\mathcal B}_k \to T_{k-1} \wt{\mathcal B}_k.\] 
When the perturbation is zero, its zeroes are the orbits of framed flat connections. At a given critical point $[A, \Phi]$ with stabilizer $\Gamma$, there is further a formal $\Gamma$-equivariant Hessian operator $\text{Hess}(CS+\pi): T_k \wt{\mathcal B}_k \to T_{k-1} \wt{\mathcal B}_k$. This gives an acceptable notion of \emph{nondegenerate} critical point: the kernel of this operator should be the tangent space to the orbit $\mathcal O$ through $[A, \Phi]$ (in other words, as small as possible; by $SO(3)$-invariance of $CS+\pi$ the tangent space $T_{[A, \Phi]} \mathcal O$ will always be a subspace of the kernel). 

Alternatively we can define nondegeneracy of critical points in terms of the \emph{extended Hessian} operator. The extended Hessian operator for a connection $A$
\[\widehat{\text{Hess}}_A(CS): L^2_k(\Omega^0 \oplus \Omega^1)(\text{ad}(E)) \to L^2_{k-1}(\Omega^0 \oplus \Omega^1)(\text{ad}(E))\] 
is defined by the matrix $\begin{pmatrix}0 & -d_A^* \\ -d_A & *d_A\end{pmatrix}$. This is a self-adjoint elliptic operator, and has discrete spectrum diverging to $\pm \infty$.
If one instead works with $\widehat{\text{Hess}}_A(CS+\pi)$, then there is an additional term, given by a compact operator which factors as $L^2_k \xrightarrow{\text{Hess}(\pi)} L^2_k \hookrightarrow L^2_{k-1}$; then $\widehat{\text{Hess}}_A(CS+\pi)$ is again a self-adjoint Fredholm operator with discrete spectrum diverging to $\pm \infty$. The extended Hessian is now defined for all connections $A$ and perturbations $\pi$, not just those for which $\nabla_A(CS+\pi) = 0$, and is defined with the same domain and codomain for all $(A, \pi)$. It is thus a suitable home for a notion of spectral flow. It follows from a straightforward analysis that 
\[\ker \widehat{\text{Hess}}_A(CS+\pi) = \mathfrak g_A \iff \ker \text{Hess}_{[A, \Phi]}(CS+\pi) = TO_{[A,\Phi]},\] 
where $\mathfrak g_A$ is the space of $A$-parallel sections of $\text{ad}(E)$, and thus the notion of nondegeneracy outlined above agrees with that of Section 2. This analysis is carried out in detail in \cite[Section 3.2]{M}.

As discussed in Section 2, we choose our perturbations to vanish on the reducible locus and in a neighborhood of the central connections, and therefore the set $\mathfrak A(Y, w)$ of gauge equivalence classes of abelian flat connections coincides with the set of gauge equivalence classes of abelian critical orbits of $CS+\pi$. However, the spectral behavior at these abelian flat connections depends on $\pi$. The following definition is meant to capture this dependence.

\begin{definition}\label{def:spectral-shift}
Let $(Y,w)$ be a rational homology sphere with an orinted 1-cycle, equipped with a nice perturbation $\pi$. Its {\it spectral shift} is a function $S_\pi: \mathfrak A(Y, w) \to 2\Bbb Z,$ defined to be the spectral flow 
\[S_\pi(\rho) = \text{sf}\left(\left\{\widehat{\text{Hess}}_\rho(CS+(1-t)\pi)\right\}_{t\in [0,1]}\right).\]
When defining the spectral flow between possibly degenerate operators $A_0$ and $A_1$ via the path $\{A_t\}_{t\in [0,1]}$, we have 
\[\text{sf}(\{A_t\}) = \lim_{\epsilon \to 0^+} \text{sf}(\{A_t+\epsilon I\}).\] 
That is, we compute the intersection number of the eigenvalues with the line $\lambda = -\epsilon$ for sufficiently small $\epsilon > 0$.
\end{definition}

Each of these Hessian operators splits as a direct sum of an operator on the reducible locus (which is the same for all $t$ because $\pi = 0$ on the reducible locus) and an operator normal to the reducible locus. The latter is $S^1$-equivariant with respect to a weight two $S^1$-action, and hence complex linear with respect to an appropriate complex structure. This spectral flow can be computed as twice the complex spectral flow of a family of complex-linear operators. In particular, it is an even integer.

There is a topological interpretation of the kernel of the unperturbed extended Hessian operator $\widehat{\text{Hess}}_\rho(CS)$. If $\rho$ is a flat abelian connection, the corresponding adjoint representation $\text{ad}_\rho: \pi_1(Y) \to SU(2)/\pm I = SO(3)$ has image inside the subgroup $SO(2)$ induced by diagonal matrices $U(1) \subset SU(2)$, and thus defines a 1-dimensional complex representation $\Bbb C_{\rho}$ of $\pi_1(Y)$. Then we have an isomorphism \[\text{ker} \; \widehat{\text{Hess}}_\rho(CS) \cong \Bbb R \oplus H^1(Y; \Bbb C_{\rho}),\]
where the term $H^1(Y; \Bbb C_\rho)$ is the first cohomology with coefficients in the local system $\Bbb C_\rho$. 

We use this to define slight variation on the spectral shift.

\begin{definition}\label{def:sigdata}
A {\it signature data function} is a function $\sigma: \mathfrak A(Y, w) \to 2\Bbb Z$ so that 
\[\sigma(\rho) \equiv \dim_{\Bbb R} H^1(Y; \Bbb C_{\rho}) \pmod 4.\]
Given a spectral shift $S_\pi: \mathfrak A(Y, w) \to 2 \Bbb Z$, the corresponding signature data function is denoted by $\sigma_\pi$ and is defined as
\[\sigma_\pi(\rho) = \dim_{\Bbb R} H^1(Y; \Bbb C_{\rho}) - 2S_\pi(\rho).\qedhere\]
\end{definition}

\begin{remark}\label{rmk:relationship}
The quantity $\sigma_\pi$ defined here is the same as the quantity $\sigma_\pi$ in \cite[Definition 3.16]{M}. When $\pi$ is a small perturbation, the perturbed Hessian restricts to a nondegenerate $S^1$-invariant bilinear form on $H^1(Y; \Bbb C_{\rho})$, and $\sigma_\pi(\rho)$ is the signature of this bilinear form. In what follows, the perturbations $\pi$ are not necessarily small, but we use the name `signature data' by analogy.
\end{remark}

Later, we will need to know that every signature data function arises as the signature data of some nice perturbation $\pi$. The remainder of this section will be devoted to a proof of this fact, recorded as a proposition below.

\begin{prop}\label{arbitrary-specshift}
Let $(Y, w)$ be a weakly admissible pair. For any given function $f: \mathfrak A(Y, w) \to 2\Bbb Z$, there is a nice perturbation $\pi$ with $S_\pi = f$. Correspondingly, for any given signature data function $\sigma$, there is a nice perturbation $\pi$ with $\sigma_\pi = \sigma$.
\end{prop}

On the way to proving this claim, we will need the following (a priori weaker) fact. 

\begin{lemma}\label{lemma:large-specshift}Let $(Y, w)$ be as above and let $\rho$ be a reducible flat connection. Then for any positive integer $k$ there exists a perturbation $\pi$ such that $S_\pi(\rho) \leq -k$, and $\pi$ vanishes on the reducible locus and in a neighborhood of every other reducible flat connection. A similar claim holds if we replace $S_\pi(\rho) \leq -k$ with $S_\pi(\rho) \geq k$.
\end{lemma}

\begin{proof}
Write $T_\rho^\nu \widetilde{\mathcal B}$ for the normal space to a point in the orbit $\rho$, and $\text{Hess}^\nu_\rho$ for the Hessian of an $SO(3)$-invariant functionl on $\widetilde{\mathcal B}$ restricted to this normal space. 

We will construct this perturbation as a non-negative function which vanishes on the reducible locus (and in a neighborhood of each reducible flat connection but $\rho$). Adding this perturbation can only make eigenvalues larger, and we construct it so that the $k$ smallest negative eigenvalues become positive. As a result, we will have that \[\text{sf}(\text{Hess}^\nu_\rho, \text{Hess}^\nu_\rho + \text{Hess}^\nu(\pi)) \ge k,\] so that $S_\pi(\rho) \le -k$. Let $V_k$ be the $SO(2)$-invariant subspace of $T_\rho^\nu \widetilde{\mathcal B}$ spanned by the $k$ smallest nonpositive eigenvalues of $\text{Hess}^\nu_\rho$, with $-\lambda$ being the least of these eigenvalues; one may exponentiate this to a small $SO(2)$-invariant disc $D_k \subset \widetilde{\mathcal B}$, and consider the $SO(3)$-invariant subset $\tilde D_k$ containing $D_k$.

To construct such a perturbation, first recall the definition of holonomy perturbations. We construct these using a family of thickened loops $S^1 \times D^2 \to Y$ which all coincide on a neighborhood of $\{1\} \times D^2$; taking holonomy along $n$ such thickened loops gives an $SO(3) = \text{Inn}(SU(2))$-equivariant map $D^2 \times \widetilde{\mathcal B}(Y,w) \to SU(2)^N$. We then apply a smooth $SO(3)$-invariant function $f: SU(2)^N \to \Bbb R$ which vanishes on the reducible locus and in a neighborhood of $\{\pm I\}^N$, and then we integrate over the disc.

Given any compact $SO(3)$-invariant submanifold of $\widetilde{\mathcal B}$, such as $\tilde D_k$, we may choose the family of loops so that the holonomy map embeds $\tilde D_k$ into $SU(2)^N$. Similarly we may assume that the image of any other flat reducible under this map is disjoint from $\tilde D_k$. Now choose an $SO(3)$-invariant nonnegative function $f$ on $SU(2)^N$ which vanishes on the reducible locus, vanishes near every other flat reducible, and on the image of $D_k$ restricts to a large quadratic function $Q$ --- so large that \[Q(\text{Hol}(\exp v)) \ge 2\lambda \|v\|^2\] for all sufficiently small $v \in V_k$. Modifying the auxiliary data in this construction as appropriate (in particular, taking the radius of the thickened loops to be small) guarantees that for the corresponding perturbation $\pi$, we have \[\pi(\text{Hol}(\exp v) \ge 1.5\lambda \|v\|^2\] for all sufficiently small $v \in V_k$. It follows that \[\text{Hess}^\nu_\rho(CS+\pi)(v) \cdot v = \text{Hess}^\nu_\rho(v) \cdot v + \text{Hess}^\nu_\rho(\pi)(v) \ge -\lambda \|v\|^2 + 1.5 \lambda \|v\|^2 \ge \frac{\lambda}{2} \|v\|^2,\] and thus we have made the Hessian positive-definite on the subspace $V$, taking these $k$ non-negative eigenvalues to positive eigenvalues as desired. 

The same argument can be run with $f$ a negative function which restricts on a positive eigenspace to a negative-definite quadratic form to construct a perturbation so that the spectral flow is at least $k$.
\end{proof}

\begin{proof}[Proof of Proposition \ref{arbitrary-specshift}]
Fix an abelian flat connection $\rho$, and define 
\[
  K_\rho(\pi) := \dim_{\Bbb R}\left( \ker (\text{Hess}^\nu_{\rho}(CS+\pi))\right) \in 2\Bbb N.
\]
Notice that the function $\pi \mapsto S(\pi):=S_\pi(\rho)$ is lower semicontinuous, and the function $K$ measures its failure to be continuous: if $K_\rho(\pi) = 2m$, then for all sufficiently nearby $\pi'$, we have $S(\pi) \le S(\pi') \le S(\pi) + 2m$. In fact, any jump between $S(\pi)$ and $S(\pi) + 2m$ is realizable. Though stated under the assumption that $\pi = 0$, the argument of \cite[Proposition 3.9]{M} shows that for any $0 \le i \le m$, there exists a nearby perturbation $\pi'$ which vanishes near every reducible but $\rho$ so that $S(\pi') = S(\pi)+2i$, while $K_\rho(\pi') = 0$. That is, given a perturbation $\pi$, by a small change we may make $2i$ of the $0$-eigenvalues of $\text{Hess}^\nu_{\rho, \pi}$ negative, and $2(m-i)$ positive, leaving no $0$-eigenvalues. (The additional perturbation is chosen to be very close to a small but nondegenerate $SO(2)$-invariant quadratic form on $\Bbb R^{2m}$ with $2i$ positive eigenvalues and $2(m-i)$ negative eigenvalues.)

Now suppose we want to find a regular perturbation $\pi$ with $S_\pi(\rho) = f(\rho) = 2n$ where $n \ge 0$; the case $n < 0$ is similar. By Lemma \ref{lemma:large-specshift}, there exists some perturbation $\pi_0$ and $\pi_1$ that vanish near every reducible but $\rho$ and has 
\[
  S(\pi_0)\leq 2n \leq S(\pi_1).
\]
Now consider the family of perturbations $\pi_t:=t\pi_1+(1-t)\pi_0$. Although it might not be the case that there is some $t$ with $S(\pi_t) = 2n$, we may find $t\in [0,1]$ such that 
\[
  S(\pi_t) \le 2n \le S(\pi_t) + K_\rho(\pi_t).
\]
As we mentioned above, every jump in this range is realizable; that is, there exists a perturbation $\pi_\rho$ which vanishes near all reducibles except $\rho$ and is near to $\pi_t$, so that $S(\pi_\rho) = 2n$ and $K_\rho(\pi_\rho) = 0$. 

For each $\rho$, choose as above a perturbation $\pi_\rho$ which vanishes near all other reducible flat connections and has $S_\pi(\rho) = f(\rho)$ and $K_\rho(\pi) = 0$. Now set $\pi = \sum_\rho \pi_\rho$. Because each $\pi_\rho$ vanishes near all reducibles except $\rho$, it follows that $S_\pi(\rho) = f(\rho)$ and $K_\pi(\rho) = 0$ for all $\rho$, as desired. This perturbation is such that all reducible flat connections are non-degenerate, but it is possiible that there are irreducible critical points which are degenerate. However, nice perturbations are dense in the space of all perturbations, so choose a nice perturbation $\pi'$ near $\pi$. So long as it is chosen sufficiently close to $\pi$, this will not change the value of $S_\pi$. This is the desired nice perturbation with $S_{\pi'}(\rho) = f(\rho)$ for all $\rho$.
\end{proof}

\begin{remark}
	As discussed above, the kernel of the unperturbed extended Hessian is isomorphic to $H^1(Y;\Bbb C_{\rho})$. 
	If $0\leq S(\rho) \leq \dim_{\Bbb R} H^1(Y; \Bbb C_{\rho})$ for all $\rho$, then the function $S$ may be realized as the spectral shift $S_\pi$ 
	of an arbitrarily small perturbation. However, if $S(\rho)$ fails to lie between these bounds for any $\rho$, then $S$ is not realizable as $S_\pi$ 
	for any small perturbation. Some of the perturbations constructed above are necessarily large --- but zero on the reducible locus itself. 
	This does not make any difference for the constructions in \cite{M}. There the only reason to require that 3-manifold perturbations are small is 
	to ensure that the reducible critical set does not change. Here we enforce this with the condition on the functions defining the 
	holonomy perturbations that $f_\pi = 0$ on the reducible locus as opposed to a smallness condition.
\end{remark}

\subsection{The ASD equations on 4-manifolds}\label{cob-basics}
%!TEX root = equivariant-functoriality.tex

We are now ready to move on to the configuration spaces and moduli spaces of instantons on cobordisms. Let $(W, c): (Y, w) \to (Y', w')$ be a cobordism of weakly admissible pairs  with $\mathbf{E}_c \to W$ being the $U(2)$-bundle associated to the cycle $c$. Fix a collar neighborhood of $\partial W$ in which the embedding of $c$ is cylindrical. We will want to pass back and forth between $W$ and the manifold 
\begin{equation}\label{W*}
	W^* = (-\infty, 0] \times Y \cup_{Y} W \cup_{Y'} [0,\infty) \times Y',
\end{equation}	
obtained by adding infinite cylindrical ends to the boundary components. Notice that $W^*$ also comes equipped with an embedded surface $c^*$, which is cylindrical on the ends. In terms of the bundle $\mathbf{E}_c$, this means that its restriction to the ends is canonically isomorphic to the pullbacks of the bundles $E_{w}$, $E_{w'}$ which correspond to $w, w'$. For the definition of certain weighted Sobolev spaces below, we also fix a smooth function $\tau:W^* \to \R$ that is equal to the magnitude of the first component on the first and the third subspaces on the right hand subspace of \eqref{W*}.

For the purposes of framings, we usually suppose that $Y$ and $Y'$ are already equipped with basepoints $y$ and $y'$, which are suppressed from notation. For the purpose of parallel transport, we will also assume that $W$ is equipped with a choice of path $\gamma: [-1, 1] \to W$ so that $\gamma(-1) = y$ and $\gamma(1) = y'$, so that $\gamma$ is cylindrical in the given collar neighborhood of the boundary. There is a natural smooth extension to a map $\Bbb R \to W^*$, which we also denote by $\gamma$ in a slight abuse of notation. We consider $\gamma$ to be a `base line' on the cobordism, and $b = \gamma(0)$ to be a specific basepoint.

Now suppose that $\alpha \in \widetilde{\mathcal B}_k(Y, w)$ and $\alpha' \in \widetilde{\mathcal B}_k(Y', w')$ are fixed orbits of framed connections. Write $\mathcal A_k(\alpha)$ and $\mathcal A_k(\alpha')$ for the set of connections which lie in that gauge equivalence class. A \emph{path of connections along $(W,c)$ from $\alpha$ to $\alpha'$} (or \emph{homotopy class of path}), written $z: \alpha \to \alpha'$ where $(W,c)$ is clear from context, is an element of the set 
\[z \in \frac{\mathcal A_k(\alpha) \times \mathcal A_k(\alpha')}{\mathcal G_{k+3/2,\delta}(W,c)};\] 
here the quotient is by the group of gauge transformations of $(W^*, \mathbf{E}_{c^*})$ of Sobolev class $L^2_{k+3/2}$ which exponentially decay to constant (but not necessarily parallel) gauge transformations on the ends, and the action is $\sigma \cdot (A, B) = (\sigma(-\infty) \cdot A, \sigma(+\infty) \cdot B)$. This set is in noncanonical bijection with $\Bbb Z$.

Fix an $L^2_k$ reference connection $\mathbf{A}_0$ on $W$ which restricts to a connection in $\alpha$ and $\alpha'$ along the corresponding ends; this gives rise to a path of connections $z: \alpha \to \alpha'$ along $(W,c)$ by considering the gauge equivalence class of the pair $(\mathbf A_0(-\infty), \mathbf A_0(+\infty)) = (A, A')$. 
We then define
\[\mathcal A_{k,\delta}(W, c; A, A') = \mathbf{A}_0 + L^2_{k,\delta}\bigg(\Omega^1\big(W^*, \text{ad}(\mathbf E_c)\big)\bigg).\] 
Here $L^2_{k,\delta}$ is the space of functions of weighted Sobolev space where the weight is given by $e^{-\tau \delta}$ with $\tau$ being the function defined above. Later, we will take $\delta > 0$ to be sufficiently small with respect to choices of auxiliary data on $\partial W$. This only depends on the gauge equivalence classes of $A$, $A'$ and not the choice of $\mathbf A_0$.

There is a gauge group $\mathcal G_{k+1,\delta}(W, c; A, A')$ given by gauge transformations $\sigma$ for which $d_{\mathbf A_0} \sigma$ is of Sobolev class $L^2_{k,\delta}$. In particular, elements of this gauge group are asymptotic to stabilizers of $A$ and $A'$ on the ends. The quotient 
\[\mathcal B_{k,\delta,z}(W, c; \alpha, \alpha') = \mathcal A_k(W, c; A, A')/\mathcal G_{k+1}(W, c; A, A')\] 
depends only on the path of connections $z: \alpha \to \alpha'$, and not the choice of $\mathbf A_0$ or of $(A, A')$.
We write $\mathcal B(W; \alpha, \alpha') = \bigsqcup_{z: \alpha \to \alpha'} \mathcal B_{k,z}(W, c; \alpha, \alpha')$, suppressing all of $c$, the Sobolev index $k$, and the disjoint union over paths from notation unless they are worth drawing attention to. 

We may now follow the usual recipe to free up the action. We set 
\[\widetilde{\mathcal A}(W, c; A,A') = \mathcal A(W, c; A, A') \times \ad(\mathbf E)_{\gamma(0)},\] 
where the latter term is the space of framings of the adjoint bundle at the basepoint $\gamma(0)$. This gives rise to the quotient $\widetilde{\mathcal B}(W, c; \alpha, \alpha')$ which carries an action by $SO(3) = PU(2)$. Orbits again come in three types, corresponding to central, abelian, and irreducible connections. There is also an $SO(3)$-equivariant vector bundle $\mathcal V^+_{k,\delta} \to \widetilde{\mathcal B}_k(W, c; \alpha, \alpha')$ defined using $L^2_{k,\delta}$ sections of $\Lambda^+\otimes \ad(\mathbf E)$. There are natural completions of this space to vector bundles $\mathcal V^+_{j, \delta}$ for all $j \le k$.

\begin{remark}
The choice of framing $\Phi \in  \ad(\mathbf E)_{\gamma(0)}$ gives rise to a corresponding parallel framing $\Phi_t$ above $\ad(\mathbf{E})_{\gamma(t)}$.  If one prefers, the entire theory could be set up replacing the choice of framing at a basepoint $\Phi \in  \ad(\mathbf E)_{\gamma(0)}$ with a choice of parallel framing over a base-path, and this perspective may be better equipped to give understanding of the bubbling phenomenon.
\end{remark}

Now we introduce the ASD equations. Suppose that $(Y, w), (Y', w')$ are each equipped with metrics and regular perturbations, that $(W, c)$ is equipped with a metric which is cylindrical on the fixed collar neighborhood of the boundary, and that the orbits $\alpha$ and $\alpha'$ are $\pi, \pi'$-flat connections. At this point, we should ensure that $\delta$ is chosen smaller than any $|\lambda|$, where $\lambda$ is a non-zero eigenvalue of the perturbed extended Hessian operator associated to $\alpha$ and $\alpha'$.\footnote{If one is only interested in setting up a well-behaved Fredholm problem, one merely needs to ensure that $\delta$ is disjoint from this spectrum.}
Choose a smooth function $\beta: \R \to [0,1]$ satisfying $\beta(t)=1$ for $|t|\geq 1$ and $\beta(t)=0$ for $|t|\leq 1/2$. Then on the two cylindrical ends one may assign to each connection a self-dual 2-form 
\[(\widehat \nabla_\pi)(\mathbf A) = \beta(t) (dt \wedge \nabla_{\mathbf A(t)} \pi)^+,\hspace{1cm}(\widehat \nabla_{\pi'})(\mathbf A) = \beta(t) (dt \wedge \nabla_{\mathbf A(t)} \pi')^+.\] 
Then the map 
\begin{equation}\label{perturbed-ASD}
	\mathbf{A} \mapsto (F^+_{\mathbf{A}})^0 - \widehat \nabla_\pi(\mathbf A)- \widehat \nabla_{\pi'}(\mathbf A)
\end{equation}	
descends to a smooth $SO(3)$-equivariant section 
\[F_{\pi,\pi'}: \wt {\mathcal B}_{k,\delta}(W, c; \alpha, \alpha') \to \mathcal V^+_{k-1,\delta}.\]

We will have need to make additional perturbations on $W^*$. These are given in \cite[Definition 4.2]{M}, following the presentation in \cite{K-HigherRank}. The explicit form of these perturbations is irrelevant for our purposes, and we largely leave them as a black-box. For technical reasons in the proof of the transversality theorem\footnote{Precisely, see \cite[Proposition 4.37]{M}; the argument achieves transversality of a parameterized moduli space along the locus where the part of the perturbation contained in $W$ is zero, and uses this to show that generic small choices of interior perturbation are regular. In the gluing process, one must contend with this requirement that the perturbations are small, but this is comparatively less subtle.}, and later the gluing process, we will later need to demand that these interior perturbations are chosen sufficiently small. These define smooth $SO(3)$-equivariant sections of $\mathcal V_+$ with compact derivative. With a slight abuse of notation, we write $F_\pi$ for the sum of the section above and these new perturbations.

We say that zeroes of $F_\pi$ are $\pi$-instantons (or simply instantons, for short). Once again, these come in three types: central, abelian, irreducible. 

We write $M(W, c; \alpha, \alpha')$ for the moduli space of $\pi$-instantons. This moduli space is equipped with an $SO(3)$-equivariant map 
\[e_- \times e_+: M(W, c; \alpha, \alpha') \to \alpha \times \alpha';\] 
the map $e_-$ is defined by sending $[\mathbf A, \Phi]$ to $[\mathbf A(-\infty), \text{Hol}^{\mathbf A}_{\gamma(0) \to \gamma(-\infty)} \Phi]$, using $\mathbf A$-parallel transport along $\gamma$ to transport the framing at $\gamma(0) \in W$ to $(-\infty, y)$. The map $e_+$ is defined similarly. By equivariance and the fact that $\alpha$ and $\alpha'$ are orbits, either the left or right endpoint map separately is a fiber bundle projection. The subspace of $M(W, c; \alpha, \alpha')$ given by instantons that represent a path $z$ is denoted by $M_z(W, c; \alpha, \alpha')$.

We conclude by discussing the notion of the index, and thus expected dimension, of these moduli spaces. As with our discussion of the extended Hessian in the previous section, this discussion is easier before taking the quotient by gauge transformations.

\begin{definition}\label{def:index}
If $\mathbf A \in \mathcal A_{k,\delta}(W, c; A, A')$ is a connection in the homotopy class $z$ of paths between the $\pi$-flat connections $A$ and $A'$ at $\pm \infty$, we write $i(z)$ for the index of the elliptic complex 
\[\Omega^0_{k+1,\delta,\text{ext}}\big(W, \text{ad}(\mathbf E_c)\big) \xrightarrow{d_{\mathbf A}} \Omega^1_{k,\delta}\big(W, \text{ad}(\mathbf E_c)\big) \xrightarrow{D_{\mathbf A,\pi}} \Omega^+_{k-1, \delta}\big(W, \text{ad}(\mathbf E_c)\big),\]
where the first term is the space of sections of $\text{ad}(\mathbf E_c)$ for which $d_{\mathbf A} \sigma \in L^2_{k, \delta}$, and $D_{\mathbf A,\pi}$ is the linearization of \eqref{perturbed-ASD} after adding the secondary perturbation on the interior of $W^*$. Write $\Gamma_\alpha$ for the stabilizer of $\alpha$. For a path $z: \alpha \to \alpha'$, we also define the \emph{framed index}
\[\widetilde i(z) = i(z) + \dim \Gamma_\alpha.\] 
\end{definition}

\begin{remark}
	Because any two $\mathbf A_i$ lying in the same component of connections give rise to homotopic elliptic complexes, and any two 
	gauge-equivalent $\mathbf A_i$ give rise to isomorphic elliptic complexes, the quantities $i(z)$ and $\widetilde i(z)$ defined above only 
	depend on the path $z$. 
\end{remark}

We say $\mathbf{A}$ is a regular $\pi$-instanton if $D_{\mathbf A, \pi}$ is surjective. When $\mathbf{A}$ is abelian the ASD operator $D_{\mathbf A, \pi}$ splits as $D_{\mathbf A,\pi}^{\text{red}} \oplus D_{\mathbf A,\pi}^\nu$, a reducible part and a part normal to the reducible locus; the latter is called the \emph{normal ASD operator}. In this case, nodegeneracy of $\mathbf{A}$ splits into two hypotheses: that both $D_{\mathbf A, \pi}^{\text{red}} = D_{\mathbf A}^{\text{red}}$ and $D_{\mathbf A, \pi}^\nu$ are surjective. The first condition is independent of the perturbation because our perturbatioins vanish along the reducible locus, and always holds by Hodge theory. The second condition is called being {\it regular normal to the reducible locus}, and is more difficult to achieve. The quantities $i(z)$ and $\widetilde i(z)$ compute expected dimensions of the moduli spaces of instantons: if $M_z(W, c;\alpha, \alpha')$ consists of regular irreducible instantons, then 
\begin{align*}
	i(z) &= \dim M_z(W, c;\alpha, \alpha')\big/SO(3),\\ 
	\widetilde i(z) &= \dim \text{Fiber}\left(M_z(W, c;\alpha, \alpha') \to \alpha\right).
\end{align*} 
A more detailed discussion can be found in \cite[Section 4.2]{M}.

The framed index is additive in the following sense. Let $(W_1, c_1): (Y_1, w_1) \to (Y_2, w_2)$ and $(W_2, c_2): (Y_2, w_2) \to (Y_3, w_3)$ be cobordisms, and suppose $\mathbf{A}_1: \alpha_1 \to \alpha_2$ is a connection in the homotopy class $z$ on $(W_1, c_1)$, while $\mathbf{A}_2: \alpha_2 \to \alpha_3$ is a connection in the component $w$ on $(W_2, c_2)$. Furthermore, suppose that the $\mathbf{A}_i$ are pullback connections on the cylindrical ends and in the specified collar neighborhoods (achievable by a homotopy). Then write $W_{12} = W_1 \cup_{Y_2} W_2$; we may glue the connections $\mathbf{A}_i$ to obtain a connection $\mathbf{A}_{12}: \alpha_1 \to \alpha_3$ in the homotopy class of the path that we denote by $z \ast w$. Then we have
\[\widetilde i(z) + \widetilde i(w) = \widetilde i(z \ast w).\]
Reducing the calculation to $S^1 \times Y$, it follows quickly from the Atiyah-Singer index theorem that $\tilde i(\zeta) \in 8 \Bbb Z$ for any homotopy class $\zeta: \alpha \to \alpha$; by additivity of index, and the fact $\pi_0 \widetilde{\mathcal B}(W, c; \alpha, \alpha')$ is affine over $\pi_0 \widetilde{\mathcal B}(\Bbb R \times Y, \Bbb R \times w; \alpha, \alpha)$, it follows that the mod $8$ value of $\wt i(z)$ is independent of the choice of the path $z$, and thus similarly with $i(z)$.

For an arbitrary path $z: \alpha \to \alpha'$, we define 
\[\wt i(W, c; \alpha, \alpha') = \tilde i(z) \mod 8, \quad \quad i(W,c;\alpha, \alpha') = i(z) \mod 8.\]
Now suppose $(Y, w)$ is a weakly admissible pair. Given any two null-cobordisms $(W_i, c_i): (Y, w) \to (S^3, \varnothing)$, we can glue them to an arbitrary cobordism $(W', c'):(S^3, \varnothing) \to  (Y, w)$ to obtain cobordism $(\hat W_i, \hat c_i): S^3 \to S^3$; we have 
\[\wt i(\hat W_i, \hat c_i; \theta, \theta) = \wt i(W', c'; \theta,\alpha)+\wt i(W_i, c_i; \alpha, \theta) .\]
The index theorem for the ASD operator over closed manifolds shows that %$\wt i(\hat W, \hat c) = 3(b^1(\hat W) -b^+(\hat W)) + 2[\hat c]^2 \mod 8$.
\[
  \wt i(\hat W_i, \hat c_i; \theta, \theta) \equiv 2[\hat c_i]^2- \frac{3}{2}(\chi(\hat W_i)+\sigma(\hat W_i)) \mod 8
\]
Because Euler characteristic $\chi$ and signature $\sigma$ are additive under composition of cobordisms, we have
\[
  \wt i(W_1, c_1; \alpha, \theta)-\wt i(W_2, c_2; \alpha, \theta)\equiv -\frac{3}{2}(\chi( W_1)+\sigma(W_1))+\frac{3}{2}(\chi(W_2)+\sigma(W_2)) \mod 2
\]
This leads to the following definition of an absolute $\Bbb Z/2$ grading on $\mathfrak C(Y,w)$ \cite{Fr1}. 
\begin{definition}\label{def:mod-2-grading}
When $(Y, w)$ is a weakly admissible pair, there is an absolute $\Bbb Z/2$ grading on $\mathfrak C(Y,w)$: if $(W, c): (Y, w) \to (S^3, \varnothing)$ is a null-cobordism of $Y$ and $z: \alpha \to \theta$ a component of connections on $W$, we set 
\[\wt i(\alpha) \equiv \wt i(z)+\frac{3}{2}(\chi( W)+\sigma(W)) \mod 2.\] 
If $w = \varnothing$, taking $c = \varnothing$ the same formula defines an absolute $\Bbb Z/8$ grading.
\end{definition}

In practice, we will only need instantons to be regular up to a certain index. 

\begin{definition}\label{def:4mfd-regular}
Let $(W, c)$ be as above. We say that $\pi$ is a regular perturbation if every $\pi$-instanton in a component $z$ with $i(z) \le 7$ is regular.
\end{definition}

\subsection{The case of cylinders}
For the cylinder $(W, c) = I \times (Y,w)$, the above discussion becomes significantly simpler. The noncompact manifold $W^*$ carries a canonical diffeomorphism to $\Bbb R \times Y$ with $c^* = \Bbb R \times w$, so that the associated bundle is the pullback bundle of the bundle $E \to Y$ associated to $w$. The path $\gamma: \Bbb R \to W^*$ which we will use for parallel transport is simply given by $\gamma(t) = (t, y)$, where $y \in Y$ is the chosen basepoint.

In this case, we use simpler perturbations than the previous section; we set 
\[(\widehat \nabla_\pi)(\mathbf A) = (dt \wedge \nabla_{\mathbf A(t)} \pi)^+\] 
with no bump function term, and we use no further perturbations. In particular, the manifold $\Bbb R \times Y$ carries a natural $\Bbb R$-action via $\phi_t(s, z) = (t+s, z)$. This induces a natural $\Bbb R$-action on $\wt{\mathcal B}(\Bbb R \times Y, \Bbb R \times c; \alpha, \beta)$ by pullback and parallel transport, with 
\[\tau_t[\mathbf A, \Phi] = [\phi_t^* \mathbf A, \text{Hol}_{(0, y) \to (-t, y)} \Phi].\] 
The endpoint maps to $\alpha$ and $\beta$ are invariant under this action, as is the section $F_\pi$ itself, as both metric and perturbation are translation-invariant.

It follows that the moduli space $M(\Bbb R \times Y; \alpha, \beta)$ carries a natural $\Bbb R$-action for which the endpoint maps are invariant. We write the quotient as 
\[\breve M(Y; \alpha, \beta) = M(\Bbb R \times Y; \alpha, \beta)/\Bbb R;\] 
then $\breve M(Y; \alpha, \beta)$ is an $SO(3)$-space with equivariant endpoint maps to $\alpha$ and $\beta$. 

The set of components 
\[\pi_0 \mathcal B(\Bbb R \times Y, \Bbb R \times w; \alpha, \beta)\] 
is in canonical bijection with $\pi_0 \mathcal G(Y, w) \cong \Bbb Z.$  When $M_\zeta(\Bbb R \times Y ; \alpha, \beta)$ consists of nonconstant instantons --- that is, either $\alpha \ne \beta$ or at least $\zeta$ is not the component of the constant trajectory --- the $\Bbb R$-action is free.

When $\zeta: \alpha \to \beta$ is a given trajectory, the quantity $i(\zeta)$ depends modulo $8$ only on $\alpha$ and $\beta$. We define 
\[i(Y, w; \alpha, \beta) := i(\zeta) \mod 8, \quad \text{and similarly} \quad \wt i(Y, w; \alpha, \beta) := \wt i(\zeta) \mod 8.\]
We say that a perturbation $\pi$ on $(Y, w)$ is regular if all instantons with $i(\zeta) \le 8$ are regular. Regular perturbations exist and in fact are open and dense\footnote{If one wanted to achieve transversality for all instantons, not merely those of index bounded above by a fixed quantity, then this would rather be a countable intersection of open and dense sets} in the space of all perturbations on $(Y, w)$ by \cite[Theorem 4.49]{M}. It follows by Proposition \ref{arbitrary-specshift} that there exist regular perturbations on $(Y, w)$ with arbitrary signature data (first choose a nice perturbation with desired signature data, then choose a nearby regular perturbation).

Suppose $Y$ is a rational homology sphere and $w = \varnothing$. Then for any cobordism $(W, \varnothing): (Y, \varnothing) \to (S^3, \varnothing)$, we have $\widetilde i(W, \varnothing; \theta, \theta) = - 3 (\chi(W) + \sigma(W))/2$, and $$\widetilde i(W, \varnothing; \alpha,\theta) =  \widetilde i(Y; \alpha,\theta)+\widetilde i(W, \varnothing; \theta, \theta).$$ It follows that when $w = \varnothing$, the absolute $\mathbb Z/8$-grading on $\mathfrak C(Y, \varnothing)$ coincides with $\widetilde i(Y;\alpha, \theta)$.

%\begin{definition}\label{def:mod-8-grading}
%When $Y$ is a rational homology sphere, the set $\mathfrak C(Y, \varnothing)$ carries a natural $\Bbb Z/8$-grading $\wt i(\alpha) = \wt i(Y;\alpha, \theta)$, where $\theta$ is the trivial connection.
%\end{definition}

%This lifts the grading $\wt i$ of Definition \ref{def:mod-2-grading} and has the property that 
%\[\wt i(Y; \alpha, \beta) = \wt i(\alpha) - \wt i(\beta).\]

\subsection{Compactifications of the moduli spaces}\label{comp-moduli}
%!TEX root = equivariant-functoriality.tex

We move on to discussing the (framed) Uhlenbeck compactification of our moduli spaces. These are only well-behaved in a small dimension range (sufficient for our purposes); even then, the result is not quite a smooth manifold. We will define the compactifications and describe the extra structure they carry.

For the remainder of this section, we will abuse notation and write $\mathbf A$ to denote an equivalence class of framed connection (so, properly, an equivalence class $[\mathbf A, \Phi]$, where $\Phi$ is a framing above a chosen basepoint). We suppose that the weakly admissible pair $(Y, w)$ has been equipped with a regular perturbation $\pi$ and a basepoint $y$, so that the framing is chosen at $y$.

\begin{definition}
We say that a \textup{broken framed $\pi$-instanton} on $(Y,w)$ between $\alpha$ and $\beta$ is a sequence $\mathbf{A}_i \in \breve M(Y, w; \alpha_i, \alpha_{i+1})$ of nonconstant framed $\pi$-instantons on $\mathbb R \times Y$ modulo translation for $0 \leq i \leq n$, with $\alpha_0 = \alpha$ and $\alpha_{n+1} = \beta$, so that 
\[e_+(\mathbf{A}_i) = e_-(\mathbf{A}_{i+1})\] 
for $0 \leq i \leq n$; that is, the framings agree when parallel transported to each intermediate critical orbit. If $\mathbf A_i$ lies in the homotopy class of path $\zeta_i: \alpha_i \to \alpha_{i+1}$, we say that the homotopy class of the broken framed instanton $\mathbf A$ is $\zeta = \zeta_1 * \cdots * \zeta.$

We write $\breve M^+_\zeta(Y, w;\alpha, \beta)$ for the set of broken framed $\pi$-instantons running from $\alpha$ to $\beta$ in the homotopy class $\zeta$.
\end{definition}

Since $(Y, w)$ is fixed, we suppress them from notation for the rest of this section.

The moduli space $\breve M^+_\zeta(\alpha, \beta)$ carries the the structure of a stratified set, with strata corresponding bijectively to sequences \[\alpha \xrightarrow{\zeta_0} \alpha_1 \xrightarrow{\zeta_1} \cdots \xrightarrow{\zeta_{n-1}} \alpha_n \xrightarrow{\zeta_n} \beta;\] the corresponding stratum is given by the iterated fiber product 
\[\breve M_{\zeta_0}(\alpha, \alpha_1) \times_{\alpha_1} \cdots \times_{\alpha_n} \breve M_{\zeta_n}(\alpha_n, \beta).\]

The moduli space $\breve M^+_\zeta(\alpha, \beta)$ has a natural topology so that each stratum has the fiber product topology, given by the notion of \emph{chain-convergence}. We say that a sequence $(\mathbf{A}_k)$ of framed instantons for $k = 1, 2, \cdots$ is \emph{chain-convergent} to a broken trajectory $(\mathbf{B}_i)_{1 \leq i \leq n}$ if there is a sequence $t_k = (t^1_k, \cdots, t^n_k) \in \Bbb R^n$ with each $\lim_{k \to \infty} t^i_k - t^{i-1}_k = \infty$ and so that the translates $\tau_{t^i_k} (\mathbf{A}_k)$ converge in the appropriate Sobolev topology on compact sets to $\mathbf{B}_i$. 

There is a straightforward extension of this notion to chain-convergence of broken framed instantons (possibly with different numbers of components). If $\mathbf A_n$ is a sequence of broken instantons in the homotopy class $\zeta$ which chain-converge to $\mathbf A$, then $\mathbf A$ also lies in the homotopy class $\zeta$.

We have the following theorem on the structure of this topological space; this is stated in slightly different form as \cite[Proposition 4.43]{M}. The result can be straightforwardly obtained by either the implicit function theorem methods of \cite[Chapter 19]{KM} (the method of proof preferred in \cite{M}), or by the Banach fixed-point methods of \cite[Chapter 7.2]{DK}, which we find more convenient when studying obstructed gluing theory in Section \ref{obscob} below. 

\begin{prop}\label{prop:moduli-sss}If $\zeta: \alpha \to \beta$ is a homotopy class with $i(\zeta) \le 7$, then $\breve M^+_\zeta(\alpha, \beta)$ is compact. Further, whenever $\zeta_0, \cdots, \zeta_n$ is a sequence of homotopy classes between the critical points whose concatenation is $\zeta$, there are continuous $SO(3)$-equivariant gluing maps
\[(L, \infty]^n \times \breve M_{\zeta_1}(\alpha, \alpha_1) \times_{\alpha_1} \cdots \times_{\alpha_n} \breve M_{\zeta_n}(\alpha_n, \beta) \to \breve M^+_\zeta(\alpha, \beta)\] 
which are open embeddings onto a neighborhood of the corresponding stratum; these maps are local diffeomorphisms on each stratum separately.
\end{prop}

In particular, when $i(\zeta) \leq 7$, the moduli space of possibly broken framed instantons $\breve M^+_\zeta$ is a compact topological manifold equipped with a stratification and a smooth structure on each stratum, as well as gluing maps giving continuous charts near the corners points which are smooth on each stratum. This is nearly the structure of a smooth manifold, but not quite, and we will need to take special care when taking transverse intersections or regular values of maps defined on such objects. Before discussing that, we record the structure on $\breve M^+_\zeta$ in the following definition.

\begin{definition}\label{def-str-sm-man}
	A stratified-smooth manifold of dimension $n$ consists of the following data: 
\begin{enumerate} \vspace{-0.3cm}
\item[(i)] A partially ordered set $\Delta$ with an order-preserving {\it dimension} function $d: \Delta \to \Bbb Z_{\geq 0}$ with the maximum value $n$, 
		\item[(ii)] A topological space $P$ stratified by subsets $P_s$ for $s \in \Delta$, with $\overline{P_s} = \bigcup_{t < s} P_t$,
		\item[(iii)] A smooth structure on each stratum $P_s$ that turns it into a smooth manifold of dimension $d(s)$.
	\end{enumerate}
	We demand that for each $x \in P_s$, there exists an open, stratum-preserving embedding \[\phi: [0, \epsilon)^{n - d(s)} \times \Bbb R^{d(s)} \to P,\] with $\phi(0) = x$, so that $\phi$ induces an isomorphism of posets $\{0,1\}^{n - d(s)} \to \Delta_{\geq s}$, and so that $\phi$ is a local diffeomorphism on each stratum.

\noindent  If $M$ is a smooth manifold and $P$ is a stratified-smooth manifold, we say $f: P \to M$ is stratified-smooth if it is continuous and smooth on each stratum. We say that $f$ is stratified-transverse to a smooth submanifold $S \subset M$ if $f|_{P_s}$ is transverse to $S$ for all $s \in \Delta$. 
\end{definition}

In particular, a stratified-smooth manifold is a topological manifold with corners. Proposition \ref{prop:moduli-sss} above asserts that $\breve M^+_\zeta(\alpha, \beta)$ is a stratified-smooth $SO(3)$-manifold, where $\Delta$ is the set of sequences $s = (\zeta_0, \cdots, \zeta_n)$ of homotopy classes of paths between various critical orbits $\zeta_j: \alpha_{j} \to \alpha_{j+1}$, with $\alpha_0 = \alpha$ while $\alpha_{n+1} = \beta$ and $\zeta_0 * \cdots * \zeta_k = \zeta$. We say $t < s$ if $s$ is obtained from $t$ by concatenating some of its constituent homotopy classes.

The stratum corresponding to $s = (\zeta_0, \cdots, \zeta_k)$ is \[\breve M_{\zeta_0}(\alpha, \alpha_1) \times_{\alpha_1} \cdots \times_{\alpha_n} \breve M_{\zeta_k}(\alpha_k, \beta),\] with $d(s) = \sum_{j=0}^n \left(\wt i(\zeta_j) - 1\right) + \dim \alpha.$

\par

In our work later, we will need to look at transverse intersections of a stratified-transverse map with a smooth submanifold $S \subset M$. Unfortunately, such transverse intersections are \emph{no longer} stratified-smooth manifolds; they need not even be topological manifolds in general. Nevertheless, the resulting structure (called a \emph{stratified-smooth space}) is sufficient for all of our purposes. In the definition of the instanton chain complex, we will use a chain complex which calculates singular homology, whose generators are stratified-smooth spaces up to a certain notion of equivalence. We rapidly develop the foundations of stratified-smooth spaces in Appendix \ref{A}. These are similar to, but simpler than, the notion of $\delta$-chain; this notion has appeared in a handful of references, including \cite{KM, Lin, M}.

\begin{remark}
The notation chosen here differs from the notation chosen for the same spaces in \cite{M}, where $\breve M^+_\zeta(\alpha, \beta)$ would have been written as $\overline{\mathcal M}_{E,z,k,\delta}(\alpha, \beta)$. We prefer the current notation both for brevity and for consistency with previous literature.
\end{remark}

We now move on to a similar discussion for general cobordisms. Suppose $(W, c): (Y, w) \to (Y', w')$ is equipped with a regular perturbation $\pi_W$, extending choices of regular perturbations $\pi$ and $\pi'$ on the ends, as well as a path $\gamma$ from $y \in Y$ to $y' \in Y'$; we choose the framing in $W$ above $\gamma(0)$, and parallel transport it along the path $\gamma$ as necessary.

\begin{definition}
We say that a \textup{broken framed instanton} on $(W,c)$ between $\alpha$ and $\alpha'$ is a triple $(\mathbf A_i, \mathbf B, \mathbf A_j')$:\vspace{-0.3cm}
\begin{itemize}
\item A (possibly empty) sequence $\mathbf{A}_i \in \breve M_{\zeta_i}(\alpha_i, \alpha_{i+1})$ of nonconstant framed $\pi$-instantons (modulo translation) for $0 \leq i \leq n$, with $\alpha_0 = \alpha$ so that 
\[e_+(\mathbf{A}_i) = e_-(\mathbf{A}_{i+1})\] 
for $0 \leq i < n$;
\item A framed $\pi_W$-instanton $\mathbf B\in M_z(W, c; \alpha_{n+1}, \beta_0)$, with $e_- (\mathbf B) = e_+ (\mathbf A_n) \in \alpha_{n+1}$; 
\item A (possibly empty) sequence $\mathbf{A}'_i \in \breve M_{\zeta'_i}(\alpha'_i, \alpha'_{i+1})$ of nonconstant framed $\pi'$-instantons (modulo translation) for $0 \leq i \leq m$, with $\alpha'_{m+1} = \alpha'$, so that \[e_+(\mathbf{A}'_i) = e_-(\mathbf{A}'_{i+1})\] for $0 \le i < m$, and $e_-(\mathbf A'_0) = e_+(\mathbf B) \in \alpha'_0$.
\end{itemize}

\noindent We say this broken framed instanton represents the path of connections along $(W,c)$ given by $\zeta_1 * \cdots * \zeta_n * z * \zeta_1' * \cdots * \zeta_m'.$

\noindent We write $M^+_z(W, c;\alpha, \alpha')$ for the set of broken framed instantons on $(W, c)$ running from $\alpha$ to $\alpha'$ in the homotopy class $z$.
\end{definition}

Once again, this is a stratified set, stratified by iterated fiber products of the relevant moduli spaces; thus all strata inherit a natural smooth structure. The topology is similar to the case of a cylinder: we say that a sequence $\mathbf A_k$ of framed instantons on $(W, c)$ is chain-convergent to a broken framed instanton $(\mathbf B_i, \mathbf A, \mathbf C_j)$ if there are sequences $t_k = (t_k^1, \cdots, t_k^n) \in \Bbb R_+^n$ and $s_k = (s_k^1, \cdots, s_k^m) \in \Bbb R_-^m$ whose sequential differences go to $+\infty$ and $-\infty$, respectively, so that the framed instantons $\mathbf A_k$ converge to $\mathbf A$ on compact sets, and so that upon restricting $\mathbf A$ to the respective ends, $\tau_{t_k^i}(\mathbf A)$ converges on compact subsets\footnote{Notice that the assumptions on $t$ and $s$ imply that for sufficiently large $k$, the domain of definition of these translates contains any particular compact subset of $\Bbb R \times Y$ or $\Bbb R \times Y'$.} of $\Bbb R \times Y_1$ to $\mathbf B_i$ while $\tau_{s_k^j}(\mathbf A_k)$ converges on compact subsets of $\Bbb R \times Y_2$ to $\mathbf C_j$ (both modulo gauge, and everything in sight equipped with a framing; in restricting to the ends we parallel transport the framing along a base-path $\gamma$ to the basepoint on the end). 

Then for $i(z) \le 7$ there once again exist $SO(3)$-equivariant gluing maps, stated explicitly in \cite[Proposition 4.42]{M}. We record this as a proposition below.

\begin{prop}\label{Uhl-cob}
Suppose $(W, c)$ is a cobordism equipped with a regular perturbation. If $i(z) \le 7$, then $M^+_z(W, c; \alpha, \alpha')$ is compact with the topology induced by chain-convergence. Furthermore, with the stratification described above, $M^+_{z}(W; \alpha, \beta)$ is equipped with the structure of a stratified-smooth manifold. 
\end{prop}

\subsection{Reducibles and regularity}\label{red-reg}
%!TEX root = equivariant-functoriality.tex

We would like to construct regular perturbations $\pi_W$ on cobordisms $(W, c): (Y, w, \pi) \to (Y', w', \pi')$, where $\pi, \pi'$ are fixed regular perturbations on the boundary components. Unfortunately, this is not always possible, and in this section we will discuss the obstruction. 

Achieving regularity at irreducible connections is standard. A generic perturbation has all irreducible instantons cut out regularly; one may see this by the arguments of \cite[Section 5.5.1]{Don} except for those irreducible flat connections which are central on the ends, where one needs more care (and to introduce interior perturbations). It is more difficult to ensure regularity at the reducibles. If one of the boundary components is admissible, there are no reducible instantons whatsoever, so we may as well assume that both boundary components are rational homology spheres. 

\begin{convention}
Unless otherwise specified, for the duration of this section, $(W, c, \pi_W): (Y, w, \pi) \to (Y', w', \pi')$ is a cobordism between rational homology spheres equipped with regular perturbations $\pi, \pi'$, while $\pi_W$ is a perturbation on the cobordism extending the chosen perturbations on the ends. The signature data functions of these perturbations are written $\sigma_\pi$ and $\sigma_{\pi'}$, respectively. 
\end{convention}

First we recall the enumeration of reducible instantons on a cobordism $(W,c): (Y, w, \pi) \to (Y', w', \pi')$. Consider the set of components of the space of reducible connections on $(W,c)$; we write $\mathfrak A(W,c)$ for the set of components which include no central connections, and $\mathfrak Z(W,c)$ for the set of components which include some central connection. The following enumeration of these components is given in \cite[Proposition 2.7]{M}.

\begin{prop}\label{prop:enum-red}
Let $(W, c, \pi_W): (Y, w, \pi) \to (Y', w', \pi')$ be a cobordism between rational homology spheres. Sending a reducible connection to the pair of first Chern classes associated with the splitting into line bundles gives bijections
\begin{align*}
	\mathfrak Z(W, c) &\cong \{x \in H^2 (W;\Bbb Z) \mid 2x = {\rm PD}(c)\}, \\ 
	\mathfrak A(W, c) &\cong \big\{ \{x, y\}\subset H^2 (W;\Bbb Z) \mid x + y = {\rm PD}(c), \;\; x \neq y \big\}.
\end{align*}
\end{prop}

The central connections are regular if and only if $b^+(W) = 0$. As discussed in Section \ref{cob-basics}, for an abelian instanton $\Lambda$ to be regular, its \emph{normal ASD operator} $D^\nu_{\Lambda, \pi_W}$ must be surjective. This is a complex linear operator, so its index is an even integer; by homotopy invariance of index, this quantity depends only on $\Lambda$ and the choice of $\pi, \pi'$ on the ends.

\begin{definition}
Suppose $\Lambda$ is an abelian connection on the cobordism $(W, c, \pi_W)$. We write $N(\Lambda; \pi, \pi')$ for the index of the normal ASD operator $D^\nu_{\Lambda, \pi_W}$.
\end{definition}

If $\Lambda$ is to be cut out transversely, we need at the very least that $N(\Lambda; \pi, \pi') \ge 0$. (We may need more than this, as we will discuss later.) To make this discussion more concrete, we can compute these indices with the Atiyah-Patodi-Singer index theorem. 

To state the calculation, we set some notation. For a cohomology class $x \in H^2(W)$, write $x_i \in H^2(Y_i)$ for its restriction to each boundary component. Notice that because the boundary components are rational homology spheres, the map $j: H^2_c(W; \Bbb Q) \to H^2(W; \Bbb Q)$ is an isomorphism. Whenever $x \in H^2(W; \Bbb Z)$ we write $x^2 \in \Bbb Q$ to mean the rational number obtained as 
\[x^2 = \langle (j^{-1} x)^2, [W]\rangle,\] 
where $[W] \in H_4^{\text{lf}}(W;\Bbb Q)$ is the fundamental class in locally finite homology.

The following calculation is the main result of \cite[Section 4.5]{M}.

\begin{prop}\label{prop:normal-ind}
If $\Lambda: \alpha \to \alpha'$ is an abelian instanton in the component labeled by $\{x,y\}$, then $N(\Lambda; \pi, \pi')$ is the even integer
\begin{align}N(\Lambda; \pi, \pi') = &-2(x-y)^2 + 2(b_1(W)-b^+(W)) \nonumber\\
&+\frac{\sigma_{\pi'}(\alpha') + \rho(\textup{ad}_{\alpha'})}{2} - \frac{\sigma_\pi(\alpha) + \rho(\textup{ad}_\alpha)}{2} \nonumber\\
&+ 1 - \frac{r(\alpha) + r(\alpha')}{2},\label{N-Lambda-pi-pi'}
\end{align} 
where $\rho(\textup{ad}_\alpha)$ is the Atiyah--Patodi--Singer rho invariant of the extended Hessian operator $\widehat{\textup{Hess}}_{\alpha}(CS)$, 
$r(\alpha) = 3$ if $\alpha$ is a central connection and $r(\alpha) = 1$ if $\alpha$ is abelian. Similar comments apply to $\rho(\textup{ad}_{\alpha'})$ and $r(\alpha')$.
\end{prop}

The quantity $-2(x-y)^2$ may be understood as eight times the \emph{topological energy} 
\[\left(c_2(\mathbf E_c) - \frac 14 c_1(\mathbf E_c)^2\right) = \frac{1}{8\pi^2} \int_W \text{tr} F_{\Lambda}^2,\] 
where $\mathbf E_c$ is the bundle associated to the geometric representative $c$. 

The expression in the second line of \eqref{N-Lambda-pi-pi'} can be expressed in terms of the rho invariant of the perturbed extended Hessian $\widehat{\text{Hess}}_{\alpha}(CS + \pi)$, and the corresponding $\rho$-invariant for $\alpha'$; while this is simpler than the above expression, the statement of the proposition makes clear that $\pi$ enters in the form of its signature data. Still, we record this perturbed $\rho$ invariant as a definition.

\begin{definition}\label{def:perturbed-rho}
Let $(Y, w, \pi)$ be a rational homology sphere with an oriented $1$-cycle and a nice perturbation. Then the perturbed rho invariant $\rho_\pi: \mathfrak A(Y,w) \to \Bbb Q$ is defined as \[\rho_\pi(\alpha) = \rho(\text{ad}_\alpha) + \sigma_\pi(\alpha).\qedhere\]
\end{definition}

\begin{remark}
In Section \ref{sec:CW}, we will use $\rho_\pi$ to give an invariant of rational homology spheres which conjecturally computes twice the Casson-Walker invariant, in the spirit of Taubes' interpretation of the Casson invariant \cite{Taubes} as a signed count of irreducible flat connections. For rational homology spheres, one must use the perturbed rho invariant as a correction term, analogous to Walker's definition \cite{Walker} of his invariant.
\end{remark}

Next, we seek to find perturbations $\pi_W$ so that all instantons are regular. To guarantee that this is possible, we will need to make assumptions on the normal index. Later in this article we will seek to weaken these assumptions in two different ways, which are both recorded in the definition below.

\begin{definition}\label{def:unobs}
Suppose $\Lambda$ is a component of abelian connections on $(W,c): (Y,w,\pi) \to (Y',w',\pi')$ corresponding to the pair of cohomology classes $\{x,y\}$.
\begin{itemize}
\item We say $\Lambda$ is unobstructed if either $-2(x-y)^2 < 0$ or $b_1(W) - b^+(W) \le N(\Lambda; \pi, \pi')$; otherwise we say $\Lambda$ is obstructed.
\item We say $\Lambda$ is pseudocentral if $\Lambda$ is central on both ends and has $(x-y)^2 = 0$.
\end{itemize}

\noindent Quantifying over all components, there are three types of cobordism of interest to us.
\begin{itemize}
\item If every component of abelian connections over $(W,c)$ is unobstructed, and either $b^+(W) = 0$ or $(W,c)$ supports no central connections, then we say $(W, c)$ is an unobstructed cobordism.
\item If every component of abelian connections over $(W,c)$ is either unobstructed or pseudocentral, and either $b^+(W) = 0$ or $(W,c)$ supports no central connections, we say that $(W,c)$ is a pseudo-unobstructed cobordism.
\item If $b_1(W) = b^+(W) = 0$ and all but one component of abelian connections over $(W,c)$ is unobstructed, while the one obstructed reducible is abelian on both ends and has $N(\Lambda; \pi, \pi') = -2$, we say that $(W,c)$ is a nearly unobstructed cobordism.\qedhere
\end{itemize}
\end{definition}

\begin{remark}
The condition on $b^+(W)$ and central connections is to ensure that any central connection (which necessarily remains an instanton after perturbation) are automatically cut out transversely, so that we can focus our analysis on the more intricate case of abelian reducibles. If the perturbation $\pi_W$ has small interior part, a $\pi_W$-instanton necessarily has $\int \text{tr}(F_\Lambda^2) > -\epsilon$, so that $-2(x-y)^2 \ge 0$, which is why we ignore those components with $-2(x-y)^2 < 0$. 

We will see in Section \ref{inst-cx-unob-cob-map} that unobstructed cobordisms induce chain maps between the tilde complexes. 

A pseudocentral connection should be understood to be as close to central as an abelian connection can possibly be; a central connection on the cobordism restricts to central connections on the ends, and is also automatically flat, so $x-y$ is torsion and the quantity $(x-y)^2$ is automatically zero. We will see in Section \ref{badred} that pseudo-unobstructed cobordisms also induce maps on the tilde complexes, for which each pseudocentral connection contributes as if it were two central connections.

The last type of cobordisms --- nearly unobstructed cobordisms --- are analyzed in detail in Section \ref{obscob}, where we show that they also induce maps on instanton homology, but not between the tilde complexes (one must pass to a `suspension'). Using this, we will show that the only type of obstructed reducible we actually have to deal with are the pseudocentral connections.
\end{remark}

It follows immediately from Proposition \ref{prop:normal-ind} and the fact that $\rho(\theta) = \sigma_\pi(\theta) = 0$ for central connections $\theta$ that the pseudocentral connections have particularly simple normal index.

\begin{cor}\label{cor:pseudocentral-calc}
If $\Lambda$ is a pseudocentral connection on the cobordism $(W,c): (Y,w,\pi) \to (Y',w',\pi')$, then $N(\Lambda; \pi,\pi') = 2(b_1(W) - b^+(W)-1)$. 
\end{cor}

For each type of cobordism named in Definition \ref{def:unobs}, we can find special types of perturbations: either we can find regular perturbations, or at least perturbations which are sufficiently regular for the analysis we need to carry out. 

\begin{theorem}\label{thm:unobs-transv}Suppose we are given a cobordism $(W, c): (Y, w, \pi) \to (Y', w', \pi')$.
\begin{enumerate}[label=(\roman*)]
\item If $(W, c)$ is unobstructed, there exists a regular perturbation $\pi_W$ on $(W, c)$ extending $\pi$ and $\pi'$; the part of $\pi_W$ supported in the interior of $W$ may be taken to be arbitrarily small. Given any two regular perturbations $\pi_W$ with the same restrictions to the ends, there is a regular family $\pi_W(t)$ of perturbations interpolating between them.

\item If $(W, c)$ is pseudo-unobstructed, there exists a regular perturbation $\pi_W$ on $(W,c)$ extending $\pi, \pi'$ with arbitrarily small interior part so that all instantons are regular, except possibly for the pseudocentral connections; the same is true for 1-parameter families of perturbations.

\item If $(W, c)$ is nearly unobstructed, there exist perturbations $\pi_W$ with arbitrarily small interior part, so that all irreducible instantons are regular, all reducible instantons $\Lambda$ with $N(\Lambda; \pi, \pi') \geq 0$ are regular, and so that for the reducible instanton with $N(\Lambda; \pi, \pi' ) = -2$, we have $\textup{ker } D^\nu_{\Lambda, \pi_W} = 0$; the same is true for 1-parameter families of perturbations.
\end{enumerate}
\end{theorem}

We call perturbations as in (iii) \emph{nearly regular}.

\begin{proof}[Sketch of proof]
First, fix a perturbation $\pi_W$ with small interior part for which abelian instantons are cut out transversely within the reducible locus. Arguing along the lines of \cite[Section 3.3.1]{AB1}, we may assume all central connections are isolated and cut out transversely, and thus focus on reducible components which contain no central connections. Because a $\pi_W$-instanton has $-2(x-y)^2 = \int \text{tr}(F_\Lambda^2) \approx \|F_{\mathbf A}\|_{L^2}^2 \ge 0$ and the topological energy of an instanton lies in $\frac{2}{|H_1(\partial W)|}$, we may assume $-2(x-y)^2 \ge 0$. 

The space of reducible $\pi_W$-instantons in each abelian component has expected dimension $b_1(W)-b^+(W)$. If $\pi_W(t)$ is a one-parameter family of perturbations, the space of abelian connections which are $\pi_W(t)$-instantons for some $t$ has expected dimension $b_1(W) - b^+(W) + 1$. Now \[\Lambda \mapsto D^\nu_{\Lambda, \pi_W}\] defines a map from this space of abelian instantons to the space of complex linear Fredholm operators of index $N = N(\Lambda; \pi, \pi')$, and we want this map to miss the subspace of operators with either nontrivial kernel or nontrivial cokernel, depending on whether $N$ is negative or non-negative. By \cite[Lemma 4.35]{M}, the space of operators with non-trivial kernel and cokernel has codimension $|N|+2$, and the argument of \cite[Theorem 4.37]{M} shows that for $\pi_W$ sufficiently small, the map from the space of pairs $(\mathbf A, \pi_W)$ of perturbations and abelian $\pi_W$-instantons is transverse to this locus. 

Thus for generic small $\pi_W$ and generic small 1-parameter families $\pi_W(t)$, we can avoid the locus of operators with nontrivial kernel and cokernel so long as $b_1(W) - b^+(W) + 1 < |N| + 2$, or equivalently, $b_1(W) - b^+(W) \le |N|$. This is guaranteed for unobstructed and nearly unobstructed cobordisms by definition; further, in the unobstructed case, either $b_1(W) - b^+(W) \ge -1$ so that $N \ge 0$ and all normal ASD operators have trivial cokernel, or $b_1(W) - b^+(W) \le -2$ and reducibles do not generically appear even in one-parameter families. The same argument applies in the case of pseudo-unobstructed cobordisms.
\end{proof}

We conclude this section with a more detailed discussion of pseudocentral connections. First, we can obtain more precise control over the perturbations on the pseudocentral connections than described above. 

\begin{prop}\label{prop:pentachotomy}
If $(W,c)$ is a pseudo-unobstructed cobordism, we can obtain a perturbation $\pi_W$ (and paths of perturbations between them) which is regular at all instantons except possibly the pseudocentral instantons, with the following properties at the pseudocentral instantons. The behavior depends on the quantity $c = b_1(W) - b^+(W)$.
\begin{enumerate}
\item[(i)] If $c < -1$, then $b^+(W) > 1$, so that for a generic perturbation (and a generic path of perturbations), there are no abelian instantons whatsoever.
\item[(ii)] If $c = -1$, then $b^+(W) > 0$, so that for a generic perturbation there are no abelian instantons whatsoever; for a generic path of perturbations, there will generically be an isolated set of pseudocentral instantons, which are such that the normal Fredhom operator has trivial kernel.
\item[(iii)] If $c = 0$, then for a generic perturbation there are finitely many pseudocentral instantons, which may be assumed to have normal Fredholm operator with trivial kernel; the same is true in families.
\item[(iv)] If $c = 1$, then for a generic perturbation there is a 1-dimensional family of pseudocentral instantons, which are all regular. For a generic 1-parameter family of perturbations, there is a 2-dimensional family of pseudocentral instantons, of which finitely many are not cut out regularly.
\item[(v)] If $c > 1$, for a generic perturbation and path of perturbations, all pseudocentral instantons are regular.
\end{enumerate}
\end{prop}

In any of these cases, we call such a perturbation \emph{pseudo-regular}.

\begin{proof}
The condition needed to ensure that all abelian instantons are cut out with trivial kernel or cokernel is $c \le |N|+1$, whereas the condition needed to ensure that all abelian instantons are cut out transversely in 1-parameter families is $c \le |N|$. Corollary \ref{cor:pseudocentral-calc} asserts that $N = 2c-2$. The inequality $c \le |2c-2|+1$ holds for all $c$, whereas the inequality $c \le |2c-2|$ holds for all integers $c \ne 1$. 

So for $c \ne 1$ generically one of the kernel and cokernel will be trivial for all pseudocentral instantons appearing in 1-parameter families. As for which it is, the cokernel is trivial if $N \ge 0$, while the kernel is trivial if $N \le 0$; in terms of $c$, the former happens when $c \ge 1$ and the latter when $c \le 1$. This gives the pentachotomy above.
\end{proof}

As mentioned earlier, we will show in Section \ref{badred} that all pseudo-unobstructed cobordisms define well-defined induced maps on tilde complexes. Cases (i) and (v) require no further care; case (iv) will be shown irrelevant by index considerations; cases (ii)-(iii) will require the most work, but fall to the techniques of Section \ref{obscob}.

We conclude this section by pointing out that if one is free to change the perturbations on the ends, every cobordism is pseudounobstructed. We will use this in Section \ref{badred} to define induced maps on equivariant isntanton homology for \emph{any} cobordism $(W,c)$ for which $b^+(W) = 0$ or $(W,c)$ supports no central connections. We phrase the sattement below slightly more generally (including the case when $(Y,w)$ is non-trivial admissible, where it is tautologically true) for ease of reference.

\begin{lemma}\label{everything-is-unobs}
Let $(W, c): (Y, w) \to (Y', w')$ be a cobordism between weakly admissible pairs. Then there exist regular perturbations $\pi, \pi'$ on $(Y, w)$ and $(Y', w')$ respectively so that $(W, c): (Y, w, \pi) \to (Y', w', \pi')$ is a pseudo-unobstructed cobordism.
\end{lemma}

\begin{proof}
We saw in Proposition \ref{arbitrary-specshift} that there exist regular perturbations with arbitrary signature data. Because the normal index only depends on the signature data on the ends, this is really a lemma about signature data functions. Recall that the set of signature data functions on $Y$ is affine over the abelian group $F_{Y, w} = \text{Map}(\mathfrak A(Y,w), 4\Bbb Z)$. Choose signature data functions $\sigma_0, \sigma_0'$ on $Y$ and $Y'$. Write $r: \mathfrak A(W, c) \to \mathfrak A(Y, w)$ for the restriction and similarly for $r'$.

The normal index function computed in Proposition \ref{prop:normal-ind} may be considered as a function $N: \mathfrak A(W, c) \times F_{Y,w} \times F_{Y',w'} \to 2\Bbb Z$, with 
\[N(\Lambda, f, f') \stackrel{def}{=} N(\Lambda; \sigma_0+f, \sigma_0'+ f') = N(\Lambda; \sigma_0, \sigma_0') + \frac{f'(r'(\Lambda)) - f(r(\Lambda))}{2}.\] 
Here, if $r(\Lambda)$ or $r'(\Lambda)$ are central, we interpret $f(r(\Lambda)) = 0$ and similarly $f'\big(r'(\Lambda)\big) = 0$. 

For fixed $\pi, \pi'$, the formula in Proposition \ref{prop:normal-ind} has only one unbounded term, which is non-negative for. It follows that 
\[\{N(\Lambda, 0,0) \in \mathbb Z \mid \Lambda \in \mathfrak A(W,c) \text{ supports a reducible instanton}\}\] is bounded below. As a result, there is some nonnegative integer $n$ so that $N(\Lambda, 0, 0) - \big(b_1(W) - b^+(W)\big) \ge -2n$ for all reducible components $\Lambda$ which might contain an instanton.  

Set $f(\alpha) = -4n$ for all $\alpha \in \mathfrak A(Y,w)$ and $f'(\alpha') = 4n$ for all $\alpha' \in \mathfrak A(Y',w')$. If $\Lambda$ is not a pseudocentral reducible, then one of $r(\Lambda)$ or $r'(\Lambda)$ is non-central, $f'(r'(\Lambda)) - f(r(\Lambda))$ is either $4n$ or $8n$. As a result, $N(\Lambda, f, f') \ge \big(b_1(W) - b^+(W)\big)$ for all reducibles other than pseudocentral connections. 

Therefore $\sigma' = \sigma_0' + f'$ and $\sigma = \sigma_0 + f$ satisfy the desired hypotheses. 
\end{proof}

\subsection{Families of metrics}\label{sec:moduli-fam}
The previous discussion will be enough to show that each $(Y, w, \pi)$ gives rise to an instanton chain complex, and each unobstructed cobordism with regular perturbation $(W, c, \pi_W)$ gives rise to a chain map between the instanton complexes. Even ignoring the `unobstructed' caveat, this is not adequate. First, we want to ensure that the induced map is independent of the perturbation $\pi$ on the interior up to chain-homotopy. Second, we want to ensure that the induced maps compose correctly. 

We will need additional moduli spaces and structure theorems for them to prove these claims. Instead of dealing with the two issues separately, we introduce families of possibly \emph{broken} metrics to handle both cases at once. Because we will be focused on independence and we are not worried about our homotopies being well-defined up to homotopy, we only discuss the case of metrics and perturbations parameterized by the interval, which we choose to write as $I = [0, \infty]$. If one chooses any regular metric and perturbation $g(0), \pi(0)$ and $g(\infty), \pi(\infty)$ with the interior part of the perturbation sufficiently small, there is always a regular family of metrics and perturbations $g(t), \pi(t)$ extending these (regulular in the sense that a certain parameterized derivative operator is surjective for each $\pi(t)$-instanton $\mathbf A$; it need not be the case that each $\pi(t)$ is a regular perturbation).

Suppose 
\[(W,c): (Y_0, w_0) \to (Y_1, w_1), \quad (W',c'): (Y_1,w_1) \to (Y_2,w_2)\] 
are cobordisms. We may define a family of composites 
\[W \cup_L W' = W \cup_{Y_1} [-L, L] \times Y_1 \cup W'\] 
where $L \in [0, \infty)$ is a length parameter. The resulting manifold has a natural smooth structure, metric, and perturbation, which agrees on the left end with that of $W$ and on the right end with that of $W'$. So long as $L$ is sufficiently large, if $\pi$ and $\pi'$ are regular perturbations, so is the composite perturbation $\pi \cup_L \pi'$. 

We model the case $L = \infty$ slightly differently, as $W \cup_{Y_1} [0, \infty) \times Y_1 \sqcup (-\infty, 0] \times Y_1 \cup_{Y_1}  W'$, and we call this a broken metric on the composite $W \cup_{Y_1} W'$. In this way the family of composites above defines a smooth family of (possibly broken) cylindrical end metrics on the composite, indexed by the compact interval $[0, \infty]$. We say that a $\pi(\infty)$-instanton from $\alpha_0$ to $\alpha_2$ is an element of 
\[\bigsqcup_{\alpha_0 \xrightarrow{z} \alpha_1 \xrightarrow{z'} \alpha_2} M_{z}(W,c; \alpha_0, \alpha_1) \times_{\alpha_1} M_{z'}(W',c'; \alpha_1, \alpha_2),\] 
and its homotopy class is $z * z'$. More generally, a broken $\pi(\infty)$-instanton consists of a sequence $(\mathbf A^0_i, \mathbf A_W, \mathbf A^1_j, \mathbf A_{W'}, \mathbf A^2_k)$, where $\mathbf A^j_i$ are broken framed instantons on $(Y_j, w_j)$ for $1 \le i \le n_j$ (possibly $n_j = 0$), and $\mathbf A_W$ and $\mathbf A_{W'}$ are instantons on $W$ and $W'$, respectively. As usual, we demand that for each framed instanton in this sequence, the parallel transport of each framing to $+\infty$ coincides with the parallel transport of the next framing to $-\infty$. 

Now if $W$ is equipped with a family $g(t)$ of metrics (with $g(\infty)$ possibly broken) and $\pi(t)$ of perturbations (with $\pi(\infty)$ possibly broken), one may define a moduli space $M_z(W, c, I; \alpha, \beta)$, whose elements are pairs $(t, \mathbf A)$, where $t \in [0,\infty]$ and $\mathbf A$ is a (gauge equivalence class of) framed $\pi(t)$-instanton in the path $z$ of connections along $(W,c)$. Precisely, the moduli space for the family is \[M^+_z(W \cup W',c \cup c',I; \alpha, \beta) = \bigcup_{L\in I} M^+_z(W \cup_L W', c \cup_L c', \pi \cup_L \pi';\alpha,\beta),\] topologized as usual using chain-convergence.

\begin{remark}
Given a stratified smooth space $X$, the subset $\partial^{\text{naive}} X \subset X$ of elements not in the top-dimensional strata is \textbf{not} typically a stratified-smooth space; the correct definition of $\partial X$ takes a disjoint union over the closures of each boundary stratum. For instance, $\partial \Delta^2$ is a disjoint union of three closed intervals, and the natural map $\partial \Delta^2 \to \Delta$ is two-to-one on the corners of $\Delta^2$. 

This observation, applied to the situation above, says that \[M^+_z(W \cup W', c \cup c', \pi(\infty);\alpha,\beta) \subset M^+_z(W \cup W', c \cup c', I; \alpha, \beta)\] is usually not a stratified space, and is not the same as the fiber product $\sqcup_\gamma M^+_{z_+}(W_+; \alpha, \gamma) \times_\gamma M^+_{z_-}(W_-; \gamma, \beta)$ --- this latter moduli space overcounts sequences $(\mathbf A^0_i, \mathbf A_W, \mathbf A^1_j, \mathbf A_{W'}, \mathbf A^2_k)$ which include $m$ total $\mathbf A^1$ terms by a factor of $m+1$. 
\end{remark}

The main result of \cite[Section 4.8]{M} is that the same structure results carry through for this parameterized moduli space. 

\begin{prop}\label{Uhl-family}
If $I$ is a regular family of (possibly broken) metrics and perturbations on $(W, c)$, and $z: \alpha \to \alpha'$ is a homotopy class with $i(z) \leq 6$, then $M^+_{z}(W, c, I; \alpha, \beta)$ carries the natural structure of a compact stratified-smooth $SO(3)$-manifold of dimension $\wt i(z) + \dim \alpha + 1$. 
\end{prop}

The boundary strata of this stratified-smooth manifold are as expected: two of them correspond to restricting the family to its endpoints, while the others correspond to factorizations on the two ends of $(W,c)$. In Section \ref{inst-cx-unob-cob-map} below, these will be used to construct homotopies between the induced chain maps of cobordisms.

\subsection{Orientations}\label{sec:moduli-or}
%!TEX root = equivariant-functoriality.tex

In this section we discuss orientations on the moduli spaces of instantons. %Any discussion of orientations is somewhat tedious, so we content ourselves to quoting the results of \cite[Section 5]{M} in a way which is hopefully accessible and easy to use, while making our orientation conventions clear.

If $(W,c)$ is a cobordism $(Y,w) \to (Y',w')$, a \emph{homology orientation} of $W$ is an orientation of the real vector space $H^1(W) \oplus H^{2,+}(W) \oplus H^1(Y)$.

Each critical orbit $\alpha$ on $(Y,w)$ and $\alpha'$ on $(Y',w')$ give rise to two-element sets $\Lambda(\alpha)$ and $\Lambda(\alpha')$. In Morse theory, these would correspond to an orientation of the unstable manifold at $\alpha$; the infinite-dimensional nature of our problem means this does not make sense. These are defined using orientations of determinant line bundles over the space of connections on a 4-manifold bounding $(Y, w, \alpha)$.

Given a choice of elements of both $\Lambda(\alpha), \Lambda(\beta)$ on $(Y, w)$, this induces what is called a \emph{fiber-orientation} on the moduli space $M_\zeta(Y; \alpha, \beta)$: an orientation of the fibers of the map $M_\zeta(Y; \alpha, \beta) \to \alpha$, which negates if we change either the element of $\Lambda(\alpha)$ or the element of $\Lambda(\beta)$. This choice is such that whenever we are given a choice of elements of $\Lambda(\alpha), \Lambda(\beta)$, and $\Lambda(\gamma)$ the gluing map 
\[M_\zeta(Y; \alpha, \beta) \times_\beta M_{\zeta'}(Y; \beta, \gamma) \supset U \to M_{\zeta * \zeta'}(Y; \alpha, \gamma)\]
is a fiber-orientation preserving diffeomorphism from an open subset of the domain. Notice that here we do not reduce by the $\Bbb R$-action. Our convention is that $\breve M_\zeta$ is given a fiber-orientation so that
\[\Bbb R \times \breve M_\zeta^{\text{Fiber}} \cong M^{\text{Fiber}}_\zeta\] 
is a fiber-orientation-preserving diffeomorphisms onto their image.

If $(W, c): (Y, w) \to (Y', w')$ is a cobordism equipped with a homology orientation and a chosen element of $\Lambda(\alpha)$ and $\Lambda(\alpha')$, then the moduli spaces $M_z(W; \alpha, \alpha')$ are all equipped with fiber-orientations; these fiber-orientations are negated if you swap any one of the given choices of element. These fiber-orientations have the property that given further choices of element $\Lambda(\beta)$ and $\Lambda(\beta')$, the maps \begin{align*}M_{\zeta}(Y, w; \beta, \alpha) \times_\alpha &M_z(W, c; \alpha, \alpha') \supset U \to M_{\zeta * z}(W,c; \beta, \alpha'),\\
M_z(W,c; \alpha, \alpha') \times_{\alpha'} &M_{\zeta'}(Y',w'; \alpha', \beta') \supset U \to M_{z * \zeta'}(W, c; \alpha, \beta')
\end{align*}
are fiber-orientation preserving diffeomorphisms.

In the next section, we will find it somewhat more convenient to recast the notion of fiber-orientation; a fiber-orientation is the same as an orientation of $M_\zeta(Y; \alpha, \alpha')$ which depends on a choice of orientation of $\alpha$, and which negates if we negate the orientation on $\alpha$. In this language, our convention is that all moduli spaces are oriented as $\alpha \times M_\zeta^{\text{Fiber}}$. In particular, our convention about reduction by the $\Bbb R$-action has us orient the reduced moduli spaces $\breve M_\zeta$ so that there is an orientation-preserving open embedding from an open subset 
\[\alpha \times \Bbb R \times \breve M_{\zeta}^{\text{Fiber}} \supset U \to \alpha \times M_{\zeta}^{\text{Fiber}}.\]

When determining the orientation of boundary components of $\breve M^+_\zeta(Y; \alpha, \beta)$ below, this convention introduces an additional sign of $(-1)^{\dim \alpha}$, corresponding to commuting the factor of $\Bbb R$ past the factor of $\alpha$.

Equipped with orientations, we may explicitly give boundary relations for all of the moduli spaces relevant to us, including their natural signs. Here `stratified diffeomorphism' means an order-preserving bijection between the posets of strata and a homeomorphism between the spaces which covers said bijection and which is a diffeomorphism on each stratum; `orientation-preserving' means that it is an oriented diffeomorphism on the top strata. When $M$ is a stratified-smooth manifold, the closure $\overline{M_s}$ of each face $M_s$ is again a stratified-smooth manifold; the expression $\partial M$ means the disjoint union over $\overline{M_s}$, where $M_s$ is a codimension-1 stratum equipped with the boundary orientation.

The following proposition determines the induced orientations on the boundary components of the moduli spaces that we introduced earlier. It should be understood that all moduli spaces are taken so that $i(z)$ is in the range where these are stratified-smooth manifolds (less than or equal to $7, 7, 6$ respectively in the three cases), and the sums are taken over those moduli spaces in homotopy classes which concatenate to $z$.

\begin{prop}\label{prop:bdry-reln} 
	Suppose $(Y, w)$ and $(Y', w')$ are weakly admissible pairs with regular perturbations $\pi$ and $\pi'$. For each $\alpha\in \fC_\pi(Y,w)$, suppose an element of $\Lambda(\alpha)$ and an orientation for the orbit $\alpha$ is given. 
	For each $\alpha'\in \fC_\pi(Y',w')$, pick an element of $\Lambda(\alpha')$. Suppose $(W, c):(Y,w)\to (Y',w')$ is an unobstructed cobordism.
\begin{itemize}
	\item[(i)] Given $\alpha,\beta, \gamma\in \fC_\pi(Y,w)$ and paths $\zeta_0$ from $\alpha$ to $\beta$ and $\zeta_1$ from $\beta$ to $\delta$ with $\zeta = \zeta_0 * \zeta_1$, 
	the boundary orientation and the fiber product orientation on the boundary component 
	\[
	  \breve M^+_{\zeta_0}(Y,w;\alpha, \beta) \times_\beta \breve M^+_{\zeta_1}(Y,w;\beta, \delta)
	\] 
	of $\breve M^+_{\zeta}(Y,w;\alpha, \delta)$ are related by $(-1)^{\wt i(\zeta_0)+\dim(\alpha)}$.
	\item[(ii)] Suppose a regular perturbation $\pi_W$ for the ASD equation on $(W,c)$ is fixed. 
	Given $\alpha,\beta \in \fC_\pi(Y,w)$, $\alpha'\in \fC_\pi(Y',w')$ and paths $\zeta_0$ from $\alpha$ to $\beta$ and $z_1$ from $\beta$ to $\alpha'$ (along $(W,c)$) with $z = \zeta_0 * z_1$, 
	the boundary orientation and the fiber product orientation on the boundary component
	\[
	  \breve M^+_{\zeta_0}(Y,w;\alpha, \beta) \times_\beta M^+_{z_1}(W,c;\beta, \alpha')
	\] 	
	of $M^+_{z}(W,c;\alpha, \alpha')$ are related by $(-1)^{\dim(\alpha)}$. Similarly, given $\alpha \in \fC_\pi(Y,w)$, $\alpha',\beta'\in \fC_\pi(Y',w')$ and paths $z_0$ from $\alpha$ to $\beta'$ (along $(W,c)$) and $\zeta_1'$ from $\beta'$ to 
	$\alpha'$ with $z = z_0 * \zeta_1'$, then the boundary orientation and the fiber product orientation on the boundary component
	\begin{equation}\label{fiber-prod-bdry-2}
	  M^+_{z_0}(W,c;\alpha, \beta') \times_{\beta'}\breve  M^+_{\zeta_1'}(Y',w';\beta', \alpha')
	\end{equation} 	
	of $M^+_{z}(W,c;\alpha, \alpha')$ are related by $(-1)^{\wt i(z_0)+\dim(\alpha)+1}$. Note that an orientation of $\beta'$ is necessary to define the orientation of $M^+_{\zeta_1'}(Y',w';\beta', \alpha')$ and to define 
	the fiber product orientation. However, the resulting orientation on \eqref{fiber-prod-bdry-2} is independent of the orientation of $\beta'$.
	\item[(iii)] Suppose a regular 1-parameter family of perturbations $\{\pi_W(t)\}_{t\in I}$ with $I=[0,1]$ are chosen and $M^+_{z}(W,c,I;\alpha, \alpha')$ is the corresponding family moduli space. 
	Then $M^+_{z}(W,c,\pi(0);\alpha, \alpha')$ and $M^+_{z}(W,c,\pi(1);\alpha, \alpha')$ give two boundary components of the family moduli space whose orientations respectively differs from the boundary orientation by a factor of $(-1)^{\dim \alpha - 1}$ and $(-1)^{\dim \alpha}$.
	Given $\alpha,\beta \in \fC_\pi(Y,w)$, $\alpha'\in \fC_\pi(Y',w')$ and paths $\zeta_0$ from $\alpha$ to $\beta$ and $z_1$ from $\beta$ to $\alpha'$ with $z = \zeta_0 * z_1$, 
	the boundary orientation and the fiber product orientation on the boundary component
	\[
	  \breve M^+_{\zeta_0}(Y,w;\alpha, \beta) \times_\beta M^+_{z_1}(W,c,I;\beta, \alpha')
	\] 	
	of $M^+_{z}(W,c,I;\alpha, \alpha')$ are related by $(-1)^{\wt i(\zeta_0)+\dim(\alpha)}$. Similarly, given $\alpha \in \fC_\pi(Y,w)$, $\alpha',\beta'\in \fC_\pi(Y',w')$ and paths $z_0$ from $\alpha$ to $\beta'$ and $\zeta'_1$ from 
	$\beta'$ to $\alpha'$ with $z = z_0 * \zeta_1'$, then the boundary orientation and the fiber product orientation on the boundary component
	\begin{equation}\label{fam-fiber-prod-bdry-2}
	  M^+_{z_0}(W,c,I;\alpha, \beta') \times_{\beta'} \breve M^+_{\zeta_1'}(Y',w';\beta', \alpha')
	\end{equation} 	
	of $M^+_{z}(W,c;\alpha, \alpha')$ are related by $(-1)^{\wt i(z_0)+\dim(\alpha)}$.
\end{itemize}
\end{prop}

\begin{remark}
Pick a complex line bundle $\eta$ and a base connection $A$ on $\eta$. Replacing the bundle $E_c$ with $E_c \otimes \eta$, the map $A^{\text{det}} \mapsto A^{\text{det}} + A$ identifies the corresponding moduli spaces and all of their extra structure, \emph{except for orientations}. If $\eta$ is trivialized on the ends of $W$ and $A$ is the corresponding trivial connection, the orientations of our moduli spaces change by a uniform sign of $(-1)^{c_1^2(\eta)}$; this may be seen by reducing to the K\"ahler case using excision. See \cite[Section 4.2]{K-HigherRank} for more details. It follows that up to sign, the induced map only depends on the associated $SO(3)$-bundle $E_c$ and not the geometric representative $c$.
\end{remark}

\newpage

\section{Geometric chains and instanton homology}\label{geo-chain-inst-hom}
\subsection{Stratified-smooth spaces and the geometric chain complex}\label{sec:gchain-intro}
%!TEX root = equivariant-functoriality.tex

The language we use to describe instanton homology below is the language of \emph{geometric chains}. A theory of geometric chains is a flexible machine which, for each smooth manifold $M$, outputs a chain complex $C_*^{gm}(M)$ whose homology groups are naturally isomorphic to singular homology. This chain complex is freely generated by \emph{probing maps}: smooth maps $\phi: P \to M$ so that $P$ is a space of a certain pre-determined type, modulo certain pre-determined relations. Given a smooth map $f: M \to N$, we have an induced map $C_*^{gm}(M) \to C_*^{gm}(N)$, given at the level of generators by composition; at the level of homology, this agrees with the induced map in singular homology. Finally, given a correspondence $M \leftarrow W \to N$ between closed smooth manifolds with $r: W \to M$ a submersion, we can define chain-level \emph{fiber-product maps} 
\[\times_M W: C_*^{gm}(M) \to C_{*+\dim W - \dim M}^{gm}(N).\] 
The differential of our instanton chain complex will include both the boundary operator on geometric chains as well as fiber-products with certain moduli spaces. Thus having some form of chain-level fiber product maps is absolutely crucial.

The flexibility to choose the class of probing spaces $P$ and the relations satisfied by probing maps mean that it is easy to define a great many seemingly different theories of geometric chains. Already in \cite{FukayaSum, LipGH, Lin, M}, there are four different theories of geometric chains, whose homology functors are all canonically isomorphic to singular homology. We will need another, for a technical reason outlined below.

In Appendix \ref{A}, we construct the particular class of probing spaces and equivalence relations used in this paper, called \emph{compact oriented stratified-smooth spaces modulo collapse-equivalence}. In this section, we outline the salient features of our relevant geometric homology theory. The stated facts below are proven in the appendices. First, we'll define stratified-smooth spaces, our generalization of the notion of stratified-smooth manifolds from Definition \ref{def-str-sm-man}.
\begin{definition}\label{def-n-str-sm-sp}
	A naive stratified-smooth space of dimension $n$ consists of the following data:
	\begin{enumerate}
		\item[(i)] A partially ordered set $\Delta$ with an order-preserving grading function $d: \Delta \to \Bbb Z_{\geq 0}$  with the maximum value $n$, 
		\item[(ii)] A topological space $P$ stratified by subsets $P_s$ for $s \in \Delta$, with $\overline{P_s} = \bigcup_{t < s} P_t$,
		\item[(iii)] The structure of a smooth manifold of dimension $d(s)$ on each stratum $P_s$.
	\end{enumerate}
	
An \emph{oriented} naive stratified-smooth space is the additional data of an orientation on all top-dimensional strata.
\end{definition}

When $P$ is compact, as all examples of interest to us will be, the set $\Delta$ is necessarily finite. Given a continuous map $\psi: [0,\epsilon)^{n-d(s)} \times \Bbb R^{d(s)+\ell} \to \Bbb R^\ell$ which, on each stratum, is smooth and has zero as a regular value, the zero set $Z(\psi) = \psi^{-1}(0)$ is a naive stratified-smooth space. 

\begin{definition}\label{def-str-sm-sp}
A stratified smooth space is a naive stratified-smooth space so that each point $p \in P$ has a neighborhood which is \emph{locally isomorphic} to a naive stratified-smooth space of the form $Z(\psi)$.

We say this stratified-smooth space is oriented if it is also equipped with an orientation on the top strata so that these local isomorphisms can be chosen to be orientation-preserving.
\end{definition}

The notion of `locally isomorphic' is intricate, and includes the data of the stratification. The interested reader can find a significantly more detailed discussion, as well as some examples and non-examples, in Appendix \ref{app:str-sm}. Orientations are discussed in Appendix \ref{SS-or}. Because these are mostly a \emph{technical necessity}, we do not go into further detail here. The moduli spaces of instantons $\breve M^+_\zeta(\alpha, \beta)$ are always stratified-smooth \emph{manifolds}, and the first time stratified-smooth spaces make an appearance geometrically is in Section \ref{subsec:modified-mod-space}. For now, they are only used to define a geometric chain complex whose generators include stratified-smooth manifolds, and there is little harm in imagining all `stratified-smooth spaces' below are smooth manifolds with corners. 

We shall define a \emph{geometric chain complex} using stratified-smooth maps $\phi: P \to M$ from compact, connected, oriented stratified-smooth spaces as our probes into $M$, much the same as simplicial homology is constructed by probing $M$ with continuous maps from simplices. We call such a pair $(P, \phi)$ a probe (in $M$). We will use the term {\it stratified-smooth chain (in $M$)} for a formal linear combination of probes.

When $P$ is a stratified-smooth space, we will verify in Lemma \ref{lemma:bd-ok} that the closure $\overline{P_s}$ of each stratum is also stratified-smooth. Then the boundary of a stratified-smooth probe is given as 
\[\partial P = \sum_{d(s) = \dim P - 1} \overline{P_s},\] 
where each boundary stratum is given an appropriate boundary orientation. This is a finite sum by compactness of $P$.

While the detailed foundations of stratified-smooth spaces and the geometric chain complex are carried out in the appendix, the following key points are all that is needed for the main text:
\begin{itemize}
	\item If two probes $\phi: P \to M$ and $\eta: Q \to M$ are stratum-preserving diffeomorphic by an orientation-preserving diffeomorphism $f: P \to Q$ with $\eta f = \phi$, 
	then we set $\phi = \eta$ as elements of the geometric chain complex $C_*^{gm}(M)$.

	\item If a probe $\phi: P \to M$ supports an orientation-\emph{reversing} diffeomorphism $f: P \to P$ with $\phi f = \phi$, then we set $\phi = 0$ 
	as an element of the geometric chain complex $C_*^{gm}(M)$.

	\item If a probe $\phi: P \to M$ has $\dim P > 0$ and 
\[\text{rank}(d_x \phi_s: T_x P_s \to T_{\phi(x)} M) < d(s)\] 
for all strata with $d(s) > 0$ and all $x \in P_s$, we say $\phi$ is degenerate, and we set $\phi = 0$ as an element of $C_*^{gm}(M)$. 

	\item The most interesting relation, novel to this article and crucial in our proof of the invariance of instanton homology, is \emph{collapse-equivalence}. Let $(P, \phi)$ and $(Q,\eta)$ be probes. Suppose that there are loci $Z_P \subset P$ and $Z_Q \subset Q$ which are unions of embedded stratified-smooth subspaces of positive codimension, and \emph{degenerate} boundary strata of $P$ and $Q$, and suppose that there exists an orientation-preserving diffeomorphism $f: P \setminus Z_P \to Q \setminus Z_Q$ with $\eta f = \phi$. 
	
	Then we say $(Q, \eta)$ is {\it collapse-equivalent} to $(P, \phi)$ and set them equal in $C_*^{gm}(M)$. 
	
	\item We use a truncation procedure to ignore high-dimensional chains, so that $C_*^{gm}(M)$ is concentrated in degrees $[0, \dim M]$; in degrees $d < \dim M$, this is simply $C_*^{gm}(M)$, while in degree $d = \dim M$, this is the quotient by the subgroup $\partial C_{\dim M + 1}^{gm}$.

	\item Given smooth manifolds $M$ and $N$, a probe $(f,g): W \to M \times N$ is called an {\it oriented correspondence from $M$ to $N$}. If $f$ is a submersion and $M$ is oriented, then the assignment 
	\[
	  (P, \phi) \mapsto (P \times_M W, g \pi_2)
	\]  
	gives a map $\times_M W: C^{gm}_*(M) \to C^{gm}_{* + \dim W - \dim M}(N)$ which satisfies the relation 
	\[
	  \partial(\phi \times_M W) = (\partial \phi) \times_M W + (-1)^{\dim P + \dim M} \phi \times_M (\partial W).
	\] 
	If either one of the orientation on $M$ or $W$ is negated, then the induced map $\times_M W$ is negated. Given a pair of correspondences $V \to L \times M$ and $W \to M \times N$ (with $V \to L$ and $W \to M$ submersions), the composite $(\phi \times_L V) \times_M W$ is equal to $\phi \times_L (V \times_M W)$. That is, the fiber product operation is associative.

	\item If $(f,g): W \to M \times N$ is zero as an element of $C_*^{gm}(M \times N)$, then $\phi \times_M W = 0$ identically on the chain level.  In particular, if $\dim W > \dim M + \dim N$, then $\phi \times_M W= 0$, and if $W$ is collapse-equivalent to $W'$, then $\phi \times_M W = \phi \times_M W'$. %It also follows that if $V$ and $W$ are composable correspondences (given by geometric chains in $L \times M$ and $M \times N$, respectively) so that 
%	\[\dim V + \dim W - \dim M > \dim L  + \dim N,\] %	then $\phi \times_L (V \times_M W)$ is the zero operator at the chain level.
\end{itemize}

\begin{example}
	The most important example of collapse-equivalence for our applications is obtained by real blowup, discussed formally in Appendix \ref{BlowupChain}. If $\psi: P \to E$ is a stratified-smooth section of a vector bundle which 
	is transverse to the zero section, one may take the real blowup $B(\psi)$ of this section, replacing the zero locus with a sphere bundle over said zero locus. There is a natural map $B(\psi) \to P$ which is a collapse-equivalence. 
	This will play an important role in Section \ref{obscob}, and we introduce collapse-equivalence so that we can consider a chain equal to its blowup.
\end{example}

In Appendix \ref{CheckAxioms}, we verify that stratified-smooth spaces satisfy all the necessary requirements to give a suitable theory of geometric chains, canonically isomorphic to singular homology. We use the corresponding geometric chain complex freely below. In trying to understand the arguments of this article, it is almost entirely harmless to imagine that `stratified-smooth space' may be read `smooth manifold with corners'. There is a perfectly good notion of collapse-equivalence for smooth manifolds with corners, for instance. In Appendix \ref{gm-trunc} we add the truncation mentioned above.

Henceforth when $M$ is a smooth manifold, we suppress $gm$ from notation and write $C_*(M)$ for $tC_*^{gm}(M)$, the truncated geometric chains of Appendix \ref{gm-trunc}.

\begin{remark}
It seems plausible that with significant analytic work, we could prove that the instanton moduli spaces (and the modifications we discuss later) can be given the structure of smooth manifolds with corners. Assuming this, the theory of stratified-smooth chains would be wholly unnecessary, and we would be able to carry out our constructions with a more standard notion of geometric chains (like smooth manifolds with corners). However, all that we've actually proved is that these moduli spaces carry the structure of stratified-smooth manifolds. If we do not want to attempt to construct smooth atlases for the instanton moduli spaces, the development of the stratified-smooth theory is then forced upon us by our need for transverse intersections. This is a tradeoff, but we believe that in this case accepting stratified-smooth chains is the path of least resistance.
\end{remark}
\subsection{Flow categories and bimodules}\label{FlowCatS}
%!TEX root = equivariant-functoriality.tex

With a suitable notion of geometric chains in hand, we will now define \emph{flow categories and bimodules of geometric chains}. These are almost precisely the data rising from the moduli spaces of flowlines of a Morse--Bott function, though we make the drastic assumption that all endpoint maps are submersions. This is not necessary, but simplifies discussions of orientation and of the domains of various maps; it is also necessarily true in all of our cases of interest (which arise from $G$-equivariant Morse theory, where the critical submanifolds are orbits and endpoint maps are equivariant).

A flow category is precisely what we need to construct a Morse--Bott type complex in geometric chains.

\begin{definition}\label{def:flowcat}
A finitely-graded \textup{flow category} $\mathcal C$ in geometric chains is the following data, subject to the following restrictions: 
\begin{itemize}
\item A finite set of objects $\mathsf{Ob}(\mathcal C)$, each equipped with the structure of a closed, connected, orientable smooth manifold. (Orientability is not strictly necessary, and is completely irrelevant over $\Bbb F_2$, but makes our discussion of orientations on the moduli spaces simpler.)
\item For each object $X \in \mathsf{Ob}(\mathcal C)$ a 2-element set $\Lambda(X)$ called the \emph{orientation set}.
\item A grading function $\wt i: \mathcal C \times \mathcal C \to \Bbb Z/2N$ so that $\wt i(X,Y) + \wt i(Y,Z) = \wt i(X, Z)$; when multiple flow categories are under consideration we will write $\wt i_{\mathcal C}$ for this function. This is called the `fiber-dimension', or `fiber-index'.
\item For each pair $X, Y$ of objects, for each orientation of $X$, and for each element of $\Lambda(X)$ and $\Lambda(Y)$, we have a geometric chain, written as 
\[\mathcal C(X,Y) = \mathcal C_{XY} \in C_{\dim X + \wt i(X,Y) - 1}(X \times Y),\] 
so that the map $\mathcal C(X,Y) \to X$ is a submersion. We require that if we reverse any one of the orientation on $X$, the element of $\Lambda(X)$, or the element of $\Lambda(Y)$, the resulting geometric chain is the previous geometric chain with its orientation reversed.
\item Given a sequence $X, Y, Z$ of objects, we demand that we have equality 
\[(-1)^{\dim X}\partial \mathcal C_{XY} = \sum_{Z \in \mathsf{Ob}(\mathcal C)} (-1)^{\wt i(X,Z)} \mathcal C_{XZ} \times_Z \mathcal C_{ZY}\] 
as elements of $C_*(X \times Y)$.\qedhere
\end{itemize}
\end{definition}

Let us try to clarify a few points about this definition.

In Morse--Bott theory, to orient the moduli spaces of flowlines from $X$ to $Y$, one needs to orient a number of things: $X$, the unstable bundle of $X$, and the stable bundle of $Y$. In the second bullet point, the sets $\Lambda(X)$ should be understood as analogous to the set of orientations on the stable bundle of $X$. While these are not canonical, once one makes a choice for each $X$ (and orients each $X$) we obtain induced orientations on all moduli spaces of flowlines, with predictable change as we change any of our choices. 

By the orientation assumption, and the way fiber product maps $\times_X \mathcal C_{XY}$ behave under orientation-reversal of $X$, one obtains a well-defined map 
\[\times_X \mathcal C_{XY}: C_*(X) \otimes_{\Bbb Z/2}\Bbb Z[\Lambda(X)] \to C_*(Y)[\wt i(X,Y) - 1] \otimes_{\Bbb Z/2} \Bbb Z[\Lambda(Y)],\] 
where $\Bbb Z/2$ acts on the complex $C_*(X)$ and $C_*(Y)$ by negation/orientation-reversal and on $\Lambda(X)$ and $\Lambda(Y)$ by swapping the two elements; this map is independent of the orientation on $X$ and does not require any choice of an element of $\Lambda(X)$ or $\Lambda(Y)$.\footnote{Equivalently but less canonically, one may suppress the sets $\Lambda(X)$ by choosing an element of each from the getgo.}

In the fourth bullet point, the index should be read as an integer mod $2N$, so that $\mathcal C(X, Y)$ is an element of 
\[\bigoplus_{k \equiv_{2N} \dim X + \wt i(X,Y) - 1} C_k(X \times Y),\] 
noting that this group vanishes when $k > \dim X + \dim Y$ so there are only finitely many nonzero terms in this sum. Here and in what follows $\equiv_{2N}$ is shorthand for `equal modulo $2N$'. 

\begin{remark}
This definition is a weakening of \cite[Definition 2.8]{Zhou} of Morse--Bott flow category in three ways. First, we only pay attention to these moduli spaces $\mathcal C_{XY}$ up to a certain dimension, $\dim X + \dim Y$. As a second appearance of the same phenomenon, when looking at the fiber product $\mathcal C_{XZ} \times_Z \mathcal C_{ZY}$, all fiber product chains larger than $\dim X + \dim Y + 1$ are identically zero. Third, usually one demands that the term $\bigsqcup_Z \mathcal C_{XZ} \times_Z \mathcal C_{ZY}$ is oriented diffeomorphic to $\partial \mathcal C_{XZ}$. Here, we only ask that they are equal as geometric chains; this means that if one is diffeomorphic to a \emph{collapse} of the other, or differs from another by a degenerate chain, they still give the same element of the geometric chain complex. On the other hand, our submersive condition is more strict than Zhou's definition. If desired, it is not hard to remove this assumption in what follows, but all flow categories arising from $G$-invariant Morse theory (for compact Lie group $G$) satisfy this condition.
\end{remark}

The point of Definition \ref{def:flowcat} is that, using the fiber product constructions above with an appropriate theory of geometric chains, a flow category gives rise to a chain complex, the `flow complex'. 

\begin{definition}
Suppose one has a finitely-graded flow category $\mathcal C$ in geometric chains. There is an associated \emph{flow complex} $CM_*(\mathcal C)$ which comes equipped with a relative $\Bbb Z/2N$-grading
\[CM_*(\mathcal C) = \bigoplus_{X \in \mathsf{Ob}(\mathcal C)} C_*(X)[\wt i(X)] \otimes_{\Bbb Z/2} \Bbb Z[\Lambda(X)],\] 
where the shift denotes that if $(P, \phi) \in C_*(X)$ and $(Q, \eta) \in C_*(Y)$ then their relative grading is 
\[\text{gr}[(P, \phi), (Q, \eta)] = \dim P - \dim Q + \wt i(X,Y).\] 
One imagines that $\wt i(X,Y) = \wt i(X) - \wt i(Y)$ for some illusory function $\wt i(X)$. 

\noindent The boundary operator on $CM_*(\mathcal C)$ is $d:=\partial + d_M$, where $\partial$ is the boundary operator on geometric chains and
\[d_M \phi = \sum_Y (-1)^{\dim \phi} \phi \times_X \mathcal C_{XY}.\qedhere\] 
\end{definition}

When $\mathcal C$ arises as the flow category of a Morse--Bott-Smale function on a Riemannian closed smooth manifold, this complex should be called the Morse--Bott complex, and its homology recovers the homology of said manifold. The idea for this complex originates from a differential forms version used to develop equivariant cohomology in \cite{AB2}. Not long after, Fukaya \cite{FukayaSum} used a variant defined over the integers to discuss Morse--Bott homology and develop a connected-sum theorem for instanton homology; our version here resembles Fukaya's, who also discussed geometric chains. The notation $CM_*(\mathcal C)$ is chosen to reflect the Morse-theoretic origin of this complex.

\begin{prop}\label{prop:flow-cpx}
The above construction defines a chain complex.
\end{prop}
\begin{proof}
We will suppress the sets $\Lambda(X)$ from the proof below by, for instance, picking an element of each from the start.

Observe that the differential has relative grading $-1$: the geometric boundary operator $\partial$ has degree $-1$, and because $\times_X \mathcal C_{XY}: C_*(X) \to C_*(Y)$ has degree $\wt i(X,Y) - 1$, it has degree $-1$ when considered as a map 
\[C_*(X)[\wt i(X)] \to C_*(Y)[\wt i(Y)],\] 
so that $d_M$ has degree $-1$.

Next we should verify that the differential squares to zero. We have 
\begin{align*}
d^2 \phi = \partial^2 \phi + &\sum_Y (-1)^{\dim \phi} \partial\(\phi \times_X \mathcal C_{XY}\) + (-1)^{\dim \phi - 1} \partial \phi \times_X\mathcal C_{XY}\\
+ &\sum_{Y,Z} (-1)^{\dim \phi + \wt i(X,Z) - 1} (-1)^{\dim \phi} \phi \times_X \mathcal C_{XZ} \times_Z \mathcal C_{ZY}.
\end{align*}
The first term is zero, and the sign in the last term can be simplified to $\sum (-1)^{\wt i(X,Z) - 1}$. Using the formula for the boundary of a fiber product  boundary relation in the definition of flow category, the second term may be expanded out as 
\begin{align*}(-1)^{\dim \phi} &\left(\partial \phi \times_X \mathcal C_{XY} + (-1)^{\dim \phi + \dim X} \phi \times_X \partial \mathcal C_{XY}\right) \\
&= (-1)^{\dim \phi} \partial \phi  \times_X \mathcal C_{XY} + \sum_Z (-1)^{\wt i(X, Z)} \phi \times_X \mathcal C_{XZ} \times_Z \mathcal C_{ZY},
\end{align*}
using our formulas for the boundary of a fiber product and for the boundary of $\mathcal C_{XY}$. Thus the second term of the sum is precisely the negative of the third plus fourth terms, and therefore $d^2 \phi = 0$.
\end{proof}

\begin{remark}\label{rmk:periodic-filtration}
When we work with $\Bbb Z$-graded flow categories, the filtration by index $\wt i$ gives rise to a spectral sequence computing $HM_*(\mathcal C)$. In the $\mathbb Z/2N$-graded case, the discussion is slightly more subtle.

Again we suppress $\Lambda(X)$. Choose (arbitrarily) a lift of the relative grading to an absolute grading $\wt i(X)$. This gives us a $\Bbb Z$-graded lift 
\[\wt{CM}_*(\mathcal C) = \bigoplus_{\substack{X \in \mathsf{Ob}(\mathcal C) \\ n \equiv_{2N} \wt i(X)}} C_*(X)[n],\] 
with the appropriate differential.

We may equip $\wt{CM}_*(\mathcal C)$ with a \emph{periodic filtration}: a filtration so that 
\[F_{k+2N} \wt{CM}_{*+2N} = F_k \wt{CM}_*\] 
This filtration is obtained by taking the sum only over those components with $n \leq k$. Because $\mathcal C_{XY} \to X$ is a submersion, $\dim \mathcal C_{XY} \geq \dim X$ for all $X, Y$, so that every nonzero term in the above sum has $\wt i_{\mathcal C}(X,Y) \geq 1$. It follows that the differential preserves the filtration, and the associated graded differential is the geometric boundary operator (with no fiber product terms).\footnote{If we drop the submersive assumption, this spectral sequence fails to exist unless $\dim \mathcal C_{XY} \geq \dim X - 1$ for all $X, Y$, and the $E^1$ page has a twisted differential unless $\dim \mathcal C_{XY} \ge \dim X$ for all $X, Y$.}

As discussed in \cite[Section A.8]{M}, this gives rise to a $\Bbb Z/2N \times \Bbb Z$-graded (relatively graded in the first coordinate, absolutely graded in the second) spectral sequence $(E^r, d^r) \Rightarrow HM^*(\mathcal C)$ with 
\[E^1 \cong \bigoplus_{X \in \mathsf{Ob}(\mathcal C)} H_*(X).\] The differential on the $E^1$ page is given by fiber product with those moduli space components $\mathcal C_{XY}$ of dimension $\dim X$ (that is, where $\wt i_{\mathcal C}(X,Y) = 1$, so that $\text{Fiber}(\mathcal C_{XY} \to X$) is a finite set). This filtration is degreewise finite, so the spectral sequence collapses on a finite page --- more precisely, on $E^r$ where $r \ge \dim X + 2$ for all $X \in \mathsf{Ob}(\mathcal C)$ --- and furthermore  is strongly convergent, meaning that its $E^\infty$ page is isomorphic to the associated graded group of the periodic filtration on $HM_*(\mathcal C)$. 
\end{remark}

In practice, one constructs finitely-graded flow categories by providing actual stratified-smooth spaces (not their representatives in the chain complex, which involves a quotient by additional relations), and only doing so through a given dimension range. This use of a dimension range is precisely why we truncate our complex of geometric chains, so that the following lemma holds.

\begin{lemma}\label{lemma:flowcat-from-objects}
Suppose one has a set of objects $X$, a grading function $\wt i$, and a collection of compact oriented stratified-smooth spaces $\mathcal C_{XY}$, satisfying the first four bullet points of Definition \ref{def:flowcat} and so that the boundary relation holds when $\dim \mathcal C_{XY} \le \dim X + \dim Y + 1$. Then the construction of Proposition \ref{prop:flow-cpx} still defines a chain complex.
\end{lemma}

\begin{proof}
We suppress signs from this discussion, and write $\mathcal C_{XZ,k}$ to refer to the components of dimension $k$. The proof that $d^2 = 0$ amounts to showing that for all $X, Z$, the following two maps $C_*(X) \to C_{*+ k - \dim X}(Z)$ coincide: the first given by fiber product with $\partial\mathcal C_{XZ,k}$, and the second given by fiber product with $\sum_{Y,i} \mathcal C_{XY,i} \times_Y \mathcal C_{YZ,k+\dim Y - i}$. When $k- \dim X \le \dim Z+1$, these coincide by hypothesis. When $k - \dim X > \dim Z + 1$, these maps are identically zero, as our model for $C_*(Z)$ is supported in degrees $[0, \dim Z]$. 
\end{proof}

In the instanton setting, $\dim X + \dim Y$ is at most $6$ and $2N = 8$, so this means we need those moduli spaces of dimension between $0$ and $7$, and $\mathcal C(X,Y)$ is a manifold of some fixed dimension $0 \le d \le 7$. 

Just as flow categories give rise to chain complexes, it is natural to imagine that a similar structure gives rise to chain maps between flow complexes. This leads us to the following structure, which form the morphisms in a category $\mathsf{sFlow}$ (the `s' for `submersive').

\begin{definition}\label{def:bimod}
Suppose $\mathcal C$ and $\mathcal C'$ are both finitely graded flow categories in geometric chains. A \textup{bimodule} $\mathcal M$ between these, written $\mathcal M: \mathcal C \to \mathcal C'$, is the following data, subject to the following restrictions.

\begin{itemize} 
\item We are given a function $\wt i_{\mathcal M}: \mathsf{Ob}(\mathcal C) \times \mathsf{Ob}(\mathcal C') \to \Bbb Z/2N$, which is affine in the sense that given $X,Z \in \mathcal C$ and $Y, W \in \mathcal C'$ we have 
\[\wt i_{\mathcal C}(X, Z) + \wt i_{\mathcal M}(Z, W) + \wt i_{\mathcal C'}(W, Y) = \wt i_{\mathcal M}(X, Y).\]
\item We also have, for each $(X,Y) \in \mathsf{Ob}(\mathcal C) \times \mathsf{Ob}(\mathcal C')$ and each orientation on $X$ and choice of element of $\Lambda(X)$ and $\Lambda(Y)$, a geometric chain 
\[\mathcal M(X,Y) = \mathcal M_{XY} \in C_{\dim X + \wt i_{\mathcal M}(X, Y)}(X \times Y).\] 
We demand that $\mathcal M_{XY} \to X$ is a submersion for all $X,Y$. Once again we demand that the orientation of $\mathcal M_{XY}$ should reverse if we reverse any one of the orientation of $X$, the element of $\Lambda(X)$, or the element of $\Lambda(Y)$.
\item For all $X \in \mathsf{Ob}(\mathcal C)$ and $Y \in \mathsf{Ob}(\mathcal C')$ we have an equality at the level of geometric chains in $X \times Y$
\[(-1)^{\dim X} \partial \mathcal M_{XY} = \sum_{Z \in \mathsf{Ob}(\mathcal C)} \mathcal C_{XZ} \times_Z \mathcal M_{ZY} + (-1)^{\wt i_{\mathcal M}(X,W)+1}\sum_{W \in \mathsf{Ob}(\mathcal C')} \mathcal M_{XW} \times_W \mathcal C'_{WY}. \qedhere\]
\end{itemize}
\end{definition}

The fiber-dimension function $\wt i_{\mathcal M}$ should be understood here as the analogue of degree of a chain map, in the absence of an absolute grading.

\begin{remark}Given a generic correspondence $M \xleftarrow{f} W \xrightarrow{g} N$ between closed oriented smooth manifolds with Morse--Bott functions, there are moduli spaces of points in $W$ which, projected to $M$ (resp $N$), flow upward (downward) to a given critical submanifold. The collection of such moduli spaces $\mathcal M_W(X, Y)$ form a bimodule in the above sense. In a moment, we will see that a bimodule induces a chain map between the flow complexes. In the Morse--Bott case described here, the induced map on homology is $g f^!$, where $f^!$ is the umkehr map defined using Poincare duality on $M$ and $W$. See, for instance, \cite[Section 2.8]{KM} for the Morse case.
\end{remark}

As expected, a bimodule induces a chain map between the flow complexes.

\begin{prop}
Given a bimodule $\mathcal M: \mathcal C \to \mathcal C'$, there is a chain map 
\[CM_*(\mathcal M): CM_*(\mathcal C) \to CM_*(\mathcal C')\]
induced by the map that sends $\phi \in C_*(X)$ for $X \in \mathsf{Ob}(\mathcal C)$ to
\[CM_*(\mathcal M)\phi = \sum_{Y \in \mathsf{Ob}(\mathcal C')} \phi \times_X \mathcal M_{XY}.\]
\end{prop}

The proof that this is a chain map is a straightforward definition-push akin to the proof of Proposition \ref{prop:flow-cpx}, and we leave it to the reader. 

\begin{remark}
In terms of the periodic filtration before, the chain map induced by a bimodule is a periodically filtered map; its behavior on relative gradings is described by $i_{\mathcal M}$. The induced map on the $E^1$ page of the spectral sequence is given by the induced map of $\times_X \mathcal M(X,Y)_{\dim X}$ on homology; here the subscript indicates that we only refer to those $\mathcal M(X,Y)$ of dimension $\dim X$, so that $\mathcal M(X,Y)$ has fiber-dimension zero.
\end{remark}

We have the morphisms in our category of flow categories; for this to actually form a category, we should define identity and composite morphisms, and check the associativity and identity axioms.

\begin{definition}
Given a flow category $\mathcal C$, we say its \emph{identity bimodule} $\mathbf 1_{\mathcal C}$ has $\wt i_{\mathbf 1}(X, Y) = \wt i(X,Y)$, and 
\[\mathbf 1_{\mathcal C}(X, Y) = \begin{cases}0 & X \neq Y \\
X & X = Y 
\end{cases}.\]

\noindent Given two bimodules $\mathcal M: \mathcal C \to \mathcal C'$ and $\mathcal N: \mathcal C' \to \mathcal C''$, their composition $\mathcal N \circ \mathcal M: \mathcal C \to \mathcal C''$ is the bimodule with 
\[\wt i_{\mathcal N\circ \mathcal M}(X, Z) = \wt i_{\mathcal M}(X,Y) + \wt i_{\mathcal N}(Y, Z)\] 
for any $Y \in \mathsf{Ob}(\mathcal C'')$ and 
\[(\mathcal N \circ \mathcal M)(X, Z) = \sum_{Y \in \mathsf{Ob}(\mathcal C'')} \mathcal M(X, Y) \times_Y \mathcal N(Y, Z).\qedhere\]
\end{definition}

\begin{remark}\label{identity-morse}
One gets a slightly different prediction for the identity bimodule from the perspective of Morse theory. There, the identity bimodule should arises as the moduli spaces corresponding to the identity correspondence, or equivalently, the compactified trajectory space of flowlines on $M$ (the space of parameterized flowlines, including constant flowlines and broken flowlines). This would assign $\mathbf 1_{\text{Morse}}(X,X) = X$, while $\mathbf 1_{\text{Morse}}(X,Y) = [0,1] \times \mathcal C(X,Y)$ when $X \neq Y$. However, note that $\mathbf 1_{\text{Morse}}(X,Y)$ is degenerate as a chain in $X \times Y$: it supports an orientation-reversing isomorphism $(t,\gamma) \mapsto (1-t, \gamma)$, and hence is identically zero on the chain level for $X \ne Y$. So the identity bimodule above agrees with the one predicted by Morse theory, even at the chain level. 

\noindent The phrasing above avoids discussion of this degenerate chain. 
\end{remark}

It is straightforward to verify the bimodule relations for the bimodules introduced above; it is easy to see that composition is associative and that the identity bimodule acts as the identity, thus proving the following proposition.

\begin{prop}
Submersive flow categories, with bimodules-up-to-isomorphim as morphisms, form a category $\mathsf{sFlow}$ with identity and composition as described above. The identity and composite bimodules act as the identity and composite morphisms on the flow complex, so that $CM_*$ defines a functor.
\end{prop}

The sorts of bimodules (and thus chain maps) that arise from Morse--Bott theory depend on auxiliary data, such as metric structures. However, the induced maps on homology are well-defined. The usual argument chooses a 1-parameter family of auxiliary data, and uses this to construct a chain homotopy between the two chain maps. The definition below attempts to abstract this situation to the setting of flow categories.

\begin{definition}\label{def:htpy}
Suppose we are given two bimodules $\mathcal M_0, \mathcal M_1: \mathcal C \to \mathcal C'$ so that $\wt i_{\mathcal M_0} = \wt i_{\mathcal M_1}$ (that is, these two bimodules have the same `relative degree'). A \textup{homotopy} $\mathcal H: \mathcal M_0 \Rightarrow \mathcal M_1$ is the data of, for each $(X,Y) \in \mathsf{Ob}(\mathcal C) \times \mathsf{Ob}(\mathcal C')$, and each triple of an orientation on $X$ and element of $\Lambda(X)$ and $\Lambda(Y)$, a geometric chain 
\[\mathcal H(X, Y) = \mathcal H_{XY} \in C_{\dim X + \wt i_{\mathcal M}(X,Y) + 1}(X \times Y)\] 
so that $\mathcal H_{XY} \to X$ is a submersion for all $X, Y$. As usual we demand that the orientation of $\mathcal H_{XY}$ should reverse if we reverse any one of the orientation of $X$, the element of $\Lambda(X)$, or the element of $\Lambda(Y)$.

Finally, we demand an equality at the level of geometric chains in $C_*(X \times Y)$
\begin{align*}(-1)^{\dim X} \partial \mathcal H(X,Y) =& \;\mathcal M_1(X,Y) - \mathcal M_0(X,Y) \\
&+ \sum_{Z \in \mathsf{Ob}(\mathcal C)} (-1)^{\wt i_{\mathcal C}(X, Z)} \mathcal C_{XZ}\times_Z \mathcal H_{ZY} \\
&+ \sum_{W \in \mathsf{Ob}(\mathcal C')} (-1)^{\wt i_{\mathcal M}(X,W)} \mathcal H_{XW} \times_W \mathcal C'_{WY}.\qedhere
\end{align*}
\end{definition}

Indeed, a straightforward calculation shows that this is precisely what is needed to define a chain homotopy between induced maps of bimodules.

\begin{prop}
A homotopy $\mathcal H: \mathcal M_0 \Rightarrow \mathcal M_1$ of bimodules induces a chain homotopy 
\[CM_*(\mathcal H): CM_*(\mathcal M_0) \Rightarrow CM_*(\mathcal M_1).\]
It follows that the induced map of a bimodule on homology is well-defined up to bimodule homotopy.
\end{prop}

The formula for the chain homotopy follows the same line as for the induced map, together with an extra sign factor: for $\phi \in C_*(X)$ we write 
\[CM_*(\mathcal H)\phi = \sum_{Y \in \mathsf{Ob}(\mathcal C')} (-1)^{\dim \phi} \phi \times_X \mathcal H_{XY}.\]
If one likes, it follows from formal properties of this notion of homotopy that there is a well-defined homotopy category $\mathsf{hsFlow}$ and that taking the flow complex gives a functor 
\[CM_*: \mathsf{hsFlow} \to \mathsf{hKom}_{2N}(\Bbb Z)\] 
to the homotopy category of $\Bbb Z/2N$-graded complexes. 

Suppose all of the flow categories and bimodules and so on above carry the action of a compact Lie group $G$, acting smoothly on all objects $X$ and all moduli spaces $\mathcal C_{XY}$ so that the correspondence map $\mathcal C_{XY} \to X \times Y$ is $G$-equivariant. Then $C_*(G)$ acts on all of the $C_*(X)$'s, compatibly with the fiber product maps, so that $CM_*(\mathcal C)$ naturally carries the structure of a dg $C_*(G)$-module. Further, a bimodule $\mathcal W$ with a $G$-action compatible with those on the flow categories induces a $C_*(G)$-equivariant chain map, and a $G$-equivariant homotopy between bimodules induces a $C_*(G)$-equivariant chain homotopy between these maps.

The comment above goes through without change: there is a functor $\mathsf{hsFlow}^G \to \mathsf{hKom}^G_{2N}(\Bbb Z)$ from the homotopy category of $G$-equivariant submersive flow categories and to the homotopy category of $\Bbb Z/2N$-graded dg $C_*(G)$-modules.
\subsection{The instanton chain complex and unobstructed cobordism maps}\label{inst-cx-unob-cob-map}
%!TEX root = equivariant-functoriality.tex

Let $(Y, w, \pi)$ be a weakly admissible pair with regular perturbation $\pi$ and basepoint $y \in Y$. Then Proposition \ref{prop:bdry-reln}(i) shows that there is an associated flow category $\mathcal{I}(Y,w,\pi)$, with 
\begin{itemize}
\item object set the set of $SO(3)$-orbits of $\pi$-flat connections;
\item orientation sets the $\Lambda(\alpha)$ discussed in Section 3.7;
\item relative grading $\wt i(Y; \alpha, \beta) = \wt i(\zeta) \mod 8$, where $\zeta: \alpha \to \beta$ is any homotopy class of path of framed connections;
\item morphisms the moduli spaces $\breve M^+_\zeta(Y; \alpha, \beta)$, where $\zeta: \alpha \to \beta$ is the unique homotopy class with $0 \le \dim \alpha + \wt i(\zeta) -1 \le 7$.
\end{itemize}

As in Lemma \ref{lemma:flowcat-from-objects}, when verifying that this defines a flow category, we only need to define those moduli spaces with $\dim \breve M(Y; \alpha, \beta) \le \dim \alpha + \dim \beta + 1$. Because $\dim \breve M_\zeta(Y; \alpha, \beta) = \dim \alpha + \wt i(\zeta) - 1$, so we need those moduli spaces with $\wt i(\zeta) \le \dim \beta + 2$.

Because the dimension of an $SO(3)$-orbit is at most $3$, we need those moduli spaces with $\wt i(\zeta) \le 5$. Similarly, for the moduli spaces $M_z(W)$ on a cobordism $(W,c)$, we need those moduli spaces $M_z(W,c; \alpha, \alpha')$ with $\wt i(z) \le 4$; and for the moduli spaces $M_z(W,c, I; \alpha, \alpha')$ for a 1-parameter family of perturbations on a cobordism, we need those moduli spaces with $\wt i(z) \le 3$. 

These are promised to us by Propositions \ref{prop:moduli-sss}, \ref{Uhl-cob}, and \ref{Uhl-family} respectively, so long as $i(z) \le 7, 7, 6$ respectively. We have the relation
\[\wt i(z) = i(z) + (3 - \dim \alpha) \ge i(z),\] 
so that in particular we have well-defined stratified-smooth moduli spaces whenever $\wt i(z) \le 7, 7, 6$, regardless of what $\alpha$ is. In all cases, this means our moduli spaces are defined through fiber dimension $7 + (3 - \dim \alpha) = 10 - \dim \alpha \ge 7$, which is more than enough for our purposes.

Notice that $\mathcal I(Y,w,\pi)$ is an $SO(3)$-equivariant flow category, so its flow complex is a dg-module over $C_* SO(3)$.

\begin{definition}\label{def:tilde-cpx}
We define the {\it tilde complex} to be the associated flow complex 
\[\widetilde C_*(Y, w, \pi) = CM_*\big(\mathcal{I}(Y, w, \pi)\big),\] 
and its structure as a $C_* SO(3)$-module. 
\end{definition}

The rest of Proposition \ref{prop:bdry-reln} shows that unobstructed cobordisms $(W, c): (Y, w, \pi) \to (Y', w', \pi')$ with homology orientation give rise to bimodules between the instanton flow categories $\mathcal I(Y, w, \pi)$ and $\mathcal I(Y', w' ,\pi')$, while 1-parameter families of (possibly broken) metrics and perturbations give rise to homotopies between these bimodules (and, in the case of broken metrics, their composites).

Thus unobstructed cobordisms give equivariant chain maps between the tilde complexes, which are well-defined up to $SO(3)$-equivariant homotopy and compose as expected. We record this as a theorem.

\begin{theorem}\label{unobs-cobmap}
If $(W, c, \pi_W): (Y, w, \pi) \to (Y', w', \pi')$ is an unobstructed cobordism equipped with a homology orientation and a path $\gamma$ between the basepoints, then there is an induced $C_* SO(3)$-equivariant chain map 
\[
  \widetilde C(W, c, \pi_W, \gamma):\widetilde  C_*(Y, w, \pi) \to \widetilde C_*(Y', w', \pi').
\]
This map is independent of $\pi_W$ up to equivariant chain homotopy; the induced map of the cobordism $(\Bbb R \times Y, p^* E, p^*\pi)$ is the identity, and the composite of two unobstructed cobordism maps is equivariantly chain homotopic to the unobstructed cobordism map of the composite.
\end{theorem}

\begin{remark}The statement about the cylinder deserves comment, and is related to the situation of Remark \ref{identity-morse}. The claim follows from the fact that the moduli spaces on the cylinder are either $\breve M^+(Y,w; \alpha, \beta) \times [0,1]$ for $\alpha \neq \beta$ --- which supports an orientation-reversing isomorphism, hence is zero on the chain level --- or the set of constant trajectories $\alpha$ when $\alpha = \beta$. This means that $\mathcal{I}(\Bbb R \times Y)$ is the identity bimodule, hence induces the identity map on the flow complexes.\end{remark}

Thus the tilde complex is a sort of functor under cobordism from a certain category of manifolds $(Y, w, \pi,y)$ with morphisms the unobstructed cobordisms (with a homology orientation and path, but without a particular choice of perturbation) to a homotopy category of relatively $\Bbb Z/8$-graded $C_* SO(3)$-complexes. 

Because we need to assume the cobordisms are unobstructed, the theorem above does not even imply that the homology $\widetilde I(Y, w, \pi)$ is independent of the perturbation $\pi$. If the signature data are equal $\sigma_\pi = \sigma_{\pi'}$, then the complexes are equivariantly homotopy equivalent, but if $\sigma_\pi > \sigma_{\pi'}$ the cylinder is obstructed and we have not constructed a chain map $\widetilde C_*(Y, w, \pi) \to \widetilde C_*(Y, w, \pi')$; in fact, we will see in Corollary \ref{cor:not-h-eq} that these complexes cannot possibly be equivariantly homotopy equivalent.

We now move on to remedying this unfortunate situation; our strategy is detailed in Section 2.

\newpage

\section{Obstructed cobordism maps}\label{obscob}
%!TEX root = equivariant-functoriality.tex

The primary goal of this section is to prove Theorem \ref{intro:induced-map}. 

\begin{convention}
Throughout Section 5, we fix weakly admissible pairs $(Y,w)$ and $(Y',w')$ with regular perturbations $\pi$ and $\pi'$, and a nearly unobstructed cobordism $(W,c):(Y,w,\pi) \to (Y',w',\pi')$ as in the Definition \ref{def:unobs}. This means that $b_1(W) = b^+(W) = 0$ and $(W,c)$ supports exactly one obstructed abelian instanton, which has normal index $-2$ is asymptotic to abelian connections on $(Y,w)$ and $(Y', w')$. We call this instanton $\Lambda: \rho \to \rho'$. We fix an perturbation $\pi_W$ of the ASD equation for which all instantons but $\Lambda$ are regular and for which the deformation operator of $\Lambda$ has trivial kernel, as in Theorem \ref{thm:unobs-transv}. We write $\lambda: \rho \to \rho'$ for the path of connections over $(W, c)$ given by $\Lambda$.
\end{convention}

Due to the presence of the obstructed abelian instanton $\Lambda$, we cannot proceed as in Subsection \ref{inst-cx-unob-cob-map} to define a morphism $\widetilde C(W,c,\pi_W)$: the standard gluing theory does not apply, and the compactified moduli spaces $M^+_{z}(W,c; \alpha, \alpha')$ are not stratified-smooth spaces anymore. In Subsection \ref{obs-glu}, we use \emph{obstructed} gluing theory to analyze the moduli spaces $M^+_{z}(W,c; \alpha, \alpha')$ in the neighborhood of obstructed (broken) solutions. 

In Section \ref{subsec:modified-mod-space}, we show that \emph{modified moduli spaces} $N^+_z(W,c; \alpha, \alpha')$, obtained by deleting a neighborhood of any broken instanton factoring through $\Lambda$ and adding a `correction term', are indeed stratified-smooth spaces which satisfy an explicit boundary relation. However, they do not define bimodules between the instanton flow categories $\mathcal I$: the boundary relation is not the expected one.

Motivated by this analysis, we define the suspended flow category $\mathcal{S}_{\rho'}\mathcal{I}(Y',w',\pi')$ associated to $\mathcal{I}(Y',w',\pi')$ and $\rho'$, and then we show that the modified moduli spaces $N^+_{z}(W,c; \alpha, \alpha')$ give rise to a bimodule from $\mathcal{I}(Y,w,\pi)$ to $\mathcal{S}_{\rho'}\mathcal{I}(Y',w',\pi')$. The definition of the suspended flow category $\mathcal{S}_{\rho'}\mathcal{I}(Y',w',\pi')$ and the bimodule associated to $(W,c)$ are respectively given in Subsections \ref{cob-flow} and \ref{cob-bimods}.

With these in hand, in Section \ref{sec:invt} we apply these to $W = \Bbb R \times Y$, with $\pi_0$ and $\pi_1$ regular perturbations on $(Y, w)$ with adjacent signature data $\sigma_{\pi_0} < \sigma_{\pi_1}$, in the sense of Example \ref{adjacent signature data}. %(meaning they differ by the smallest possible amount, $4$, at exactly one reducible $\rho$). 
To show that the natural continuation map $\widetilde C(Y, w, \pi_0) \to \widetilde C(Y, w, \pi_1)$ is a quasi-isomorphism, we use the constructions above to lift this to maps to and from from the \emph{suspended complex} $S_{\rho} \widetilde C(Y, w, \pi_0) \leftrightarrow \widetilde C(Y, w, \pi_1)$, where the first complex is the flow complex of the suspended flow category $\mathcal{S}_{\rho} \mathcal I(Y, w, \pi_0)$. Using similar ideas, we conclude by showing that these two maps are chain homotopy inverse to one another. Because the continuation map described above factors as 
\[\widetilde C(Y, w, \pi_0) \to S_{\rho}(Y, w, \pi_0) \to \widetilde C(Y, w, \pi_1),\] 
and it is easy to verify that the first map is a quasi-isomorphism, this gives the desired result (and with only slightly more work, invariance of instanton homology).

\subsection{Obstructed cobordisms and obstructed gluing theory}\label{obs-glu}
Before we define the modified moduli spaces, we should understand the local structure of $M^+_z(W, c; \alpha, \alpha')$ in a neighborhood of the non-smooth locus of obstructed reducibles.

The moduli space $M^+_{z}(W,c; \alpha, \alpha')$ has three different types of obstructed elements. An {\it obstructed solution of type I} is a broken solution which belongs to the subspace 
\begin{equation}\label{obtructed-type-1}
	\breve M^+_{\zeta}(Y,w;\alpha,\rho)\times_\rho \mathcal O_\Lambda \cong \breve M^+_{\zeta}(Y,w;\alpha,\rho)\times \{\Lambda\}
\end{equation}
of a moduli space of the form $M^+_{z}(W,c; \alpha, \rho')$ with $z = \zeta * \lambda$. Similarly, an {\it obstructed solution of type II} is a broken solution in the subspace
\begin{equation}\label{obtructed-type-2}
	\mathcal O_\Lambda\times_{\rho'} \breve M_{\zeta'}^+(Y',w';\rho',\alpha')\cong  \{\Lambda\}\times \breve M_{\zeta'}^+(Y',w';\rho',\alpha')
\end{equation}
of a moduli space of the form $M^+_{z}(W,c; \rho, \alpha')$ with $z = \lambda *\zeta'$. Finally an {\it obstructed solution of type III} is a broken solution in 
\begin{equation}\label{obtructed-type-3}
	\breve M^+_{\zeta}(Y,w;\alpha,\rho)\times_\rho \mathcal O_\Lambda \times_{\rho'} \breve M_{\zeta'}^+(Y',w';\rho',\alpha')
\end{equation}
which is regarded as a subspace of $M^+_{z}(W,c; \alpha, \alpha')$ with $z = \zeta * \lambda *\zeta'$. 

In the complement of these obstructed solutions, the stratified spaces $M^+_{z}(W,c; \alpha, \alpha')$ have the structure of stratified-smooth manifolds (as in Definition \ref{def-str-sm-man}).

To study the behavior of the moduli space $M^+_{z}(W,c; \alpha, \alpha')$ in a neighborhood of an obstructed solution, we need to introduce a family of complex line bundles which are called {\it obstruction bundles}. The cokernel of the ASD operator for the elements of $\mathcal O_\Lambda$ define a complex line bundle $\mathcal H_\Lambda^+$ over $\mathcal O_\Lambda$ with an $SO(3)$ action that lifts the action on $\mathcal O_\Lambda$. %Since the restriction maps $\mathcal O_\Lambda$ to $\mathcal O_\rho$ and $\mathcal O_{\rho'}$ are isomorphisms of $SO(3)$-spaces, we may regard $\mathcal O_\Lambda$ as complex line bundles over $\mathcal O_\rho$ and $\mathcal O_{\rho'}$. 
Let $\mathcal H_{\alpha,\rho}$ be the line bundle over $\breve M^+_{\zeta}(Y,w;\alpha,\rho)\times_\rho \mathcal O_\Lambda$ pulled back from $\mathcal H_\Lambda^+$. We call $\mathcal H_{\alpha,\rho}$ the {\it obstruction bundle} associated to $\Lambda$. Similarly, we have the complex line bundle $\mathcal H_{\rho',\alpha'}$ over the space $\breve M_{\zeta'}^+(Y',w';\rho',\alpha')$.

The following proposition, which is an instance of obstructed gluing theory, gives a description for a neighborhood of an obstructed solution of type I. 
\begin{prop}\label{obs-glu-local}
	Suppose $(\bA_0,\Lambda)$ is an obstructed solution of type I in $M^+_{z}(W,c; \alpha, \rho')$ 
	where $\bA_0\in \breve M_{\zeta}(Y,w;\alpha,\rho)$ with $z = \zeta * \lambda$ given as above and 
	$i(z)\geq -1$. Then there is an $SO(3)$-invariant neighborhood $U$ of $\bA_0$ in 
	$\breve M_{\zeta}(Y,w;\alpha,\rho)$ and an $SO(3)$-equivariant section $\psi^{\textup{an}}_{\zeta\Lambda}$ of $\pi_1^*(\mathcal H_{\alpha,\rho})$ over $U\times (0, \infty]$ which vanishes on $U\times \{\infty\}$ and is transverse to the zero section over $\mathcal H_{\alpha,\rho}\times (0, \infty)$. There is also an $SO(3)$-equivariant open embedding $\Phi_{\zeta\Lambda}: \left(\psi^{\textup{an}}_{\zeta\Lambda}\right)^{-1}(0)\to M^+_{z}(W,c; \alpha, \rho')$ such that $\Phi_{\zeta\Lambda}$ maps $U\times \{\infty\}\subset \left(\psi^{\textup{an}}_{\zeta\Lambda}\right)^{-1}(0)$ to $U\times_\rho \mathcal O_\Lambda\subset M^+_{z}(W,c; \alpha, \rho')$ via the identity map on $U$.
\end{prop}

The superscript `$\text{an}$' is short for \emph{analytic}, as these sections are defined by the analytic machinery of obstructed gluing theory. %In contrast, in the next section we will \emph{choose} priveleged sections $\psi^{\text{priv}}_{\zeta \Lambda}$.

For $i(z)=-1$, the section $\psi^{\text{an}}_{\zeta\Lambda}$ is nowhere vanishing and $M^+_{z}(W,c; \alpha, \rho')$ is empty, so there is no content in this case. However, later we need compatibility of the sections $\psi^{\text{an}}_{\zeta\Lambda}$ for different $\zeta$ and $z$, and it is necessary for us to include the case that $i(z)=-1$.

\begin{proof}[Sketch of the proof]
	This proposition can be verified using standard obstructed gluing theory results for which there is 
	a large body of literature starting with the seminal work of Taubes in \cite{taubes-glue}. 
	Here we sketch the general steps in the proof and refer the reader to \cite[Chapter 7]{DK} (more precisely, the discussion following Theorem 7.2.24 and the discussion in Chapter 7.2.8) as well as 
	\cite[Chapter 4.5]{Don} 
	for more details. 
	The reader should note that the discussion of \cite{DK} is limited to the case that 
	$Y=S^3$ and the arguments are worked out in a slightly different conformal model 
	than cylindrical ends. However, working in the setup of 
	4-manifolds with cylindrical ends in fact simplifies the arguments, and the arguments of 
	\cite[Chapter 7]{DK} can be adapted to the case of a weakly admissible pair $(Y,w)$ where 
	$\fC_\pi(Y,w)$ is non-degenerate. For the ease of exposition, 
	we shall also assume that the perturbation terms in the definition of the critical sets  $\fC_\pi(Y,w)$, 
	$\fC_\pi(Y',w')$ 
	and the moduli spaces of instantons are trivial; it is straightforward to modify the arguments below to handle the general case.

	As the first step in the proof, one needs to construct approximate solutions of the ASD equation 
	by pre-gluing elements of $\breve M_{\zeta}(Y,w;\alpha,\rho)$ and $\mathcal O_\Lambda$. 
	Fix an $SO(3)$-invariant pre-compact and open neighborhood $U$ of $\bA_0$ in 
	$\breve M_{\zeta}(Y,w;\alpha,\rho)$. Fix once and for all a connection $\widehat \rho$ representing the gauge equivalence class of the reducible $\rho$, and similarly a connection $\widehat \Lambda$ representing the gauge equivalence class of obstructed reducible connection on $W$. Fix positive real numbers $T_0$ and $\epsilon_0$ which are respectively
	required to be large and small. For any $SO(3)$-orbit $\bA$ in $U$, we may find a representative connection
	$\widehat{\bA}$ such that $\widehat{\bA}$ is in temporal gauge, $\widehat{\bA}$ is asymptotic to $\hat \rho$
	and 
	\[
	  \int_{[-T_0,\infty)\times Y}\vert F_{\widehat{\bA}}\vert^2=\epsilon_0.
	\]	 
	The choice of the connection $\widehat{\bA}$ is unique up to the action of an element of the stabilizer of $\widehat \rho$.
	
	For any $T\in (0,\infty)$, we may glue the connection $\widehat{\bA}$ to $\widehat{\Lambda}$ to 
	define a connection $\widehat{\bA}\#_T \widehat\Lambda$ on $W$. To perform this pre-gluing construction, let 
	$\chi:\Bbb R\to [0,1]$ be a smooth function such that 
	$\chi(t)=0$ for $t\leq -1$ and $\chi(t)=1$ for $t\geq 0$.
	Then the connection which we call $\widehat\bA\#_T \widehat\Lambda$ is given over
	the cylindrical end $(-\infty,0)\times Y\subset W$ by
	\[
	  \chi(s+T) \widehat\Lambda+(1-\chi(s+T))\tau_{2T}^*(\widehat\bA)
	\]
	where $\tau_{T}:\Bbb R \times Y \to \Bbb R \times Y$ is the map given by $\tau_T(s,y)=(s+T,y)$. In the remaining part of $W$,  
	$\widehat\bA\#_T \widehat\Lambda$ agrees with $\widehat\Lambda$.

	Let us briefly discuss how to carry out this constructions for the whole neighborhood $U \subset \breve M_\zeta$ in the space of 
	framed connections. First, there is an open subset $\widehat U \subset \widetilde{\mathcal A}_\zeta(Y,w)$, consisting 
	of framed connections $(\widehat{\bA}, \Phi)$ so that the equivalence class $[\widehat\bA, \Phi]$ lies in $U$, and $\widehat{\bA}$ 
	satisfies the integral equality above, has asymptotic $\widehat \rho$, and is in temporal gauge. Because the framing is chosen arbitrarily, 
	this is an $SO(3)$-invariant subset of the space of framed connections; and because of the $\Gamma_{\widehat \rho}$ ambiguity in 
	the lift $\widehat{\bA}$, the projection $\pi: \widehat U \to U$ is a $\Gamma_{\widehat \rho}$-bundle.

	We write $\mathcal O_{\widehat \rho}$ for the $SO(3)$-orbit of pairs $(\widehat \rho, \Phi)$ where $\widehat \rho$ is the fixed lift above and 
	$\Phi$ is a chose of framing on $(Y,w)$; we similarly write $\mathcal O_{\widehat \Lambda}$ for pairs $(\widehat \Lambda, \Phi_W)$ where 
	$\widehat \Lambda$ is the chosen lift above and $\Phi_W$ is chosen arbitrarily. Finally, the gluing map described above defines for each $T$ 
	a map $\widehat U \times_{\mathcal O_{\widehat \rho}} \mathcal O_{\widehat \Lambda} \to \widetilde{\mathcal A}_z(W,c)$, 
	defined by 
	\[\big((\widehat\bA, \Phi), (\widehat\Lambda, \Phi_W)\big) \mapsto (\widehat\bA \#_T \widehat\Lambda, \Phi_W).\] 
	It is transparent from the construction above that this map is equivariant under the action of 
	$\Gamma_{\widehat \rho} \times_{\Gamma_{\widehat \rho}} \Gamma_{\widehat \Lambda} \cong S^1$, 
	and hence descends to an $SO(3)$-equivariant map from the quotient $U \times_{\mathcal O_{\rho}} \mathcal O_\Lambda \cong U$ to $\widetilde{\mathcal B}_z(W,c)$. 
	
	Having described how to do this construction for families, we henceforth focus on one instanton at a time; 
	the rest of the discussion can be carried out in families as well. 
	To simplify notation, we write $\bA \#_T \Lambda$ for an element of the $S^1$-orbit of framed connections produced by the construction above,
	 suppressing the framing from notation. 

	The framed connection $\bA\#_T \Lambda$ is an approximate instanton, and we want to perturb it to obtain 
	an actual solution of the ASD equation. 
	To achieve this, we need to construct right inverses for the ASD operators.
	Since the connection $\bA$ is regular, the linearization 
	of the ASD equation 
	\begin{equation}\label{dA+}
	  d_{\bA}^+:\Omega^1(\Bbb R\times Y, \text{ad}(\mathbf E_w)) \to  
	  \Omega^+(\Bbb R\times Y, \text{ad}(\mathbf E_w))
	\end{equation}
	is surjective. After taking appropriate Sobolev completions of the above function spaces, 
	we fix a right inverse $R_\bA$ for \eqref{dA+}. On the other hand, the operator $d_\Lambda^+$ 
	has non-trivial cokernel. 
	We fix a linear map 
	\[
	  \sigma_\Lambda:H^+(\Lambda) \to \Omega^+(W, \text{ad}(\mathbf E_c))
	\]
	which is supported away from the cylindrical ends of $W$ so that the operator 
	\[
	  d_\Lambda^++\sigma_\Lambda:\Omega^1(W, \text{ad}(\mathbf E_c))\oplus 
	  H^+(\Lambda) \to \Omega^+(W, \text{ad}(\mathbf E_c)),
	\]
	which maps $(a,h)$ to $d_\Lambda^+(a)+\sigma_\Lambda(h)$ is surjective. Again after taking appropriate Sobolev completions of the 
	space of 1-forms and self-dual forms on $W$, we may find a right inverse $R_\Lambda$. Linear gluing theory 
	(see, for example, \cite[Chapter 3.3]{Don}) allows us to obtain a right inverse $Q_\bA^T$ for the 
	operator 
	\[d^+_{\bA\#_T \Lambda}+\sigma_{\Lambda}:\Omega^1(W, \text{ad}(\mathbf E_c))\oplus 
	  H^+(\Lambda) \to \Omega^+(W, \text{ad}(\mathbf E_c)),
	\]
	which depends continuously on $T$ and converges to a right inverse of $d_{\mathbf A}^+ \oplus d_{\Lambda}^+ \oplus \sigma_\Lambda$ as $T \to \infty$. Note that here we use the assumption that the target of $\sigma_\Lambda$ is away from the gluing region and hence the second component of the above map is independent of $T$. We also fix the notation 
	\[
	  Q_{\bA,1}^T: \Omega^+(W, \text{ad}(\mathbf E_c)) \to \Omega^1(W, \text{ad}(\mathbf E_c)),
	  \hspace{1cm}
	  Q_{\bA,H}^T: \Omega^+(W, \text{ad}(\mathbf E_c)) \to H^+(\Lambda),
	\]
	for the components of $Q_\bA^T$. (As before, we need to consider appropriate Sobolev completions of the 
	function spaces in the domain and the co-domian of the above operators.)
	
	To perturb the connections $\bA\#_T \Lambda$ to ASD connections, first we consider 
	the weaker equation
	\[
	  F^+(\bA\#_T \Lambda+Q_{\bA,1}^T(\phi))+\sigma_\Lambda(Q_{\bA,H}^T(\phi))=0.
	\]
	A solution to this equation is the same as a fixed point of the operator
	\[
	  Z_\bA^T(\phi)=-F^+(\bA\#_T \Lambda)-(Q_{\bA,1}^T(\phi)\wedge Q_{\bA,1}^T(\phi))^+.
	\]
	One can apply the contraction mapping principle to conclude that this equation has a unique 
	solution for any $\bA\in U$ and any $T\in (0,\infty)$. (Note that we are still using the assumption 
	that $T_0$ is a large constant, so that $-F^+(\bA\#_T \Lambda)$ may be assumed to be arbitrarily small.) If $\phi_0$ is the solution of 
	this equation, then $Q_{\bA,H}^T(\phi_0)$ determines an $S^1$-equivariant map from $\widehat U\times (0,\infty)$ to $H^+(\Lambda)$. 
	This induces a section $\psi^{\textup{an}}_{\zeta\Lambda}$
	of the bundle $\mathcal H_{\alpha,\rho}$ over $U\times (0,\infty)$, and the honest solutions of 
	the ASD equation are given by the zeros of this section. 
	This section extends continuously as zero over the space $U\times \{\infty\}$.
	The desired embedding $\Phi_{\zeta\Lambda}$ is given by mapping $(\bA,T)$ in the zero set of 
	$\psi^{\textup{an}}_{\zeta\Lambda}$ to $\bA\#_T\Lambda+Q^T_{\mathbf A,1}(\phi_0)$.
\end{proof}

We need a global version of Proposition \ref{obs-glu-local}. In particular, we wish to consider an extension of the above proposition to the case that $\bA\in \breve M_{\zeta}(Y,w;\alpha,\rho)$ is replaced with a \emph{broken} trajectory in $\breve M_{\zeta}(Y,w;\alpha,\rho)$. This proposition is the first place we are obligated to use stratified-smooth spaces, as here we take the zero set of a continuous section of a vector bundle which is smooth on each stratum; a priori, there is no reason to believe the result is a stratified-smooth manifold. However, the definition of stratified-smooth spaces is chosen so that such zero sets are always stratified-smooth.

\begin{prop}\label{obs-glu-global}
	Suppose $z$ and $\zeta$ are given as above and $-1\leq i(z)\leq 6$. Then the following holds. \vspace{-0.3cm}
	\begin{itemize}
		\item[(i)] There is an $SO(3)$-equivariant continuous section $\psi^{\textup{an}}_{\zeta\Lambda}$ of $\pi_1^*(\mathcal H_{\alpha,\rho})$ over 
		$\breve M^+_{\zeta}(Y,w;\alpha,\rho)\times (0, \infty]$ which vanishes on the subspace $\breve M^+_{\zeta}(Y,w;\alpha,\rho)\times \{\infty\}$, and whose restriction to each stratum of
		$\breve M^+_{\zeta}(Y,w;\alpha,\rho)\times (0, \infty)$ is smooth and transverse to the zero section. 
		In particular, the intersection of $\left(\psi^{\textup{an}}_{\zeta\Lambda}\right)^{-1}(0)$ and $\breve M^+_{\zeta}(Y,w;\alpha,\rho)\times (0, \infty)$ is again a stratified-smooth space.
		\item[(ii)] There is an injective open continuous map $\Phi_{\zeta\Lambda}: \left(\psi^{\textup{an}}_{\zeta\Lambda}\right)^{-1}(0)\to M^+_{z}(W,c; \alpha, \rho')$ such that 
		$\Phi_{\zeta\Lambda}$ maps the subspace 
		$\breve M^+_{\zeta}(Y,w;\alpha,\rho)\times \{\infty\}\subset \left(\psi^{\textup{an}}_{\zeta\Lambda}\right)^{-1}(0)$ to $\breve M^+_{\zeta}(Y,w;\alpha,\rho) \times \{\Lambda\}\subset M^+_{z}(W,c; \alpha, \rho')$ via the identity map. 
		Moreover, the restriction of $\Phi_{\zeta\Lambda}$ to the intersection of $\left(\psi^{\textup{an}}_{\zeta\Lambda}\right)^{-1}(0)$ and $\breve M^+_{\zeta}(Y,w;\alpha,\rho) \times (0, \infty)$ 
		gives an open embedding of stratified-smooth spaces.
		\item[(iii)] For $\beta \in \fC_\pi(Y,w)$, let $\zeta_0: \alpha \to \beta$ and $\zeta_1: \beta \to \rho$ be paths of connections on $(Y, w)$ with composite $\zeta_0 * \zeta_1 = \zeta$. For any $[x,y]$ in the subspace
		\begin{equation}\label{co-dim-face}
		  \breve M^+_{\zeta_0}(Y,w;\alpha,\beta)\times_\beta \breve M^+_{\zeta_1}(Y,w;\beta,\rho)
		\end{equation}
		of $\breve M^+_{\zeta}(Y,w;\alpha,\rho)$, we have
		\begin{equation}\label{obs-sec-face}
		\psi^{\textup{an}}_{\zeta\Lambda}([x,y,\Lambda],t)=\psi^{\textup{an}}_{\zeta_1\Lambda}([y,\Lambda],t),
		\end{equation}
		and if $\psi^{\textup{an}}_{\zeta\Lambda}([y,\Lambda],t)=0$, 
		then 
		\[
		  \Phi_{\zeta\Lambda}([x, y,\Lambda],t)=(x,\Phi_{\zeta_1\Lambda}([y,\Lambda],t)).
		\]
	\end{itemize}
\end{prop}

To interpret \eqref{obs-sec-face}, note that the restriction of $\mathcal H_{\alpha,\rho}$ to the subspace of $ \breve M^+_{\zeta}(Y,w;\alpha,\rho)$ given in \eqref{co-dim-face} is induced by the pull-back of $\mathcal H_{\alpha,\rho}$. In particular, the equivariant section $\psi^{\text{an}}_{\zeta\Lambda}$ of  $\mathcal H_{\alpha,\rho}$ determines a section of $\mathcal H_{\alpha,\rho}$ restricted to \eqref{co-dim-face}.

\begin{proof}[Sketch of the proof]
	This proposition is a global version of Proposition \ref{obs-glu-local} and can be verified
	by constructing the section $\psi^{\textup{an}}_{\zeta\Lambda}$ and the map $\Phi_{\zeta\Lambda}$
	in a neighborhood of points in dimension $i$ strata of $\breve M^+_{\zeta}(Y,w;\alpha,\rho)$ 
	by induction on $i$. The induction step can be addressed by an argument similar to 
	Proposition \ref{obs-glu-local}; the bound $i(z) \le 6$ implies compactness of the relevant moduli spaces on cylinders, which is used to ensure the estimates in the proof of Proposition \ref{obs-glu-local} hold uniformly. The transversality of $\psi^{\textup{an}}_{\zeta\Lambda}$ follows from our 
	assumption on regularity of the elements of $M^+_{z}(W,c; \alpha, \rho')$ away from obstructed solutions.
	The compatibility assumption in (iii) can be verified as in \cite[Proposition 7.2.64]{DK}.
\end{proof}

There is an analogue of Proposition \ref{obs-glu-global} for obstructed solutions of type II in the moduli spaces of the form $M^+_{z}(W,c; \rho,\alpha')$. In particular, a neighborhood of such obstructed solutions in $M^+_{z}(W,c; \rho,\alpha')$ are described using an $SO(3)$-equivariant section $\psi^{\text{an}}_{\Lambda\zeta'}$ of $\mathcal H_{\rho',\alpha'}\times [-\infty, 0)$ over $\breve M_{\zeta'}^+(Y',w';\rho',\alpha')\times [-\infty, 0)$ and a map $\Phi_{\Lambda\zeta'}: \left(\psi^{\text{an}}_{\Lambda\zeta'}\right)^{-1}(0)\to M^+_{z}(W,c; \rho,\alpha')$. These satisfy the analogues of properties (i)-(iii) in Proposition \ref{obs-glu-global}.

The following proposition concerns the obstructed solutions of type III and describes the behavior of the moduli spaces $M^+_{z}(W,c; \alpha,\alpha')$ in a neighborhood of such obstructed solutions. The pull-back of the bundle $\mathcal H_\Lambda$ determines a complex line bundle over
\[X_{\zeta,\zeta'}:=\breve M^+_{\zeta}(Y,w;\alpha,\rho)\times_\rho \mathcal O_\Lambda \times_{\rho'} \breve M_{\zeta'}^+(Y',w';\rho',\alpha'),\]
which is denoted by $\mathcal H_{\alpha,\alpha'}$. It is clear from the definition of $\mathcal H_{\alpha,\alpha'}$ that it is canonically isomorphic to the pullback of the bundles $\mathcal H_{\alpha,\rho}$ and $\mathcal H_{\rho',\alpha'}$ with respect to the projection maps. Thus sections of $\mathcal H_{\alpha,\rho}$ and $\mathcal H_{\rho',\alpha'}$ can be also regarded as sections of $\mathcal H_{\alpha,\alpha'}$, and this convention is used in the statement of the following proposition about the behavior of the moduli spaces around obstructed solutions of type III. 

We omit the proof, as it follows essentially the same lines as Proposition \ref{obs-glu-global}.

\begin{prop}\label{obs-glu-global-type-3}
		Suppose $z$, $\zeta$ and $\zeta'$ are given as above and $0\leq i(z)\leq 7$. Then the following holds.
	\begin{itemize}
		\item[(i)] There is an $SO(3)$-equivariant section $\psi^{\textup{an}}_{\zeta\Lambda\zeta'}$ of $\pi_1^*(\mathcal H_{\alpha,\alpha'})$ over 
		$X_{\zeta,\zeta'}\times (0, \infty] \times [-\infty, 0)$
		whose restriction to $X_{\zeta,\zeta'}\times \big((0, \infty] \times [-\infty, 0) \setminus \{(\infty, -\infty)\}\big)$ is transverse to the zero section. For any $[x,\Lambda,x']\in X_{\zeta,\zeta'}$, we have
		\[
		  \psi^{\textup{an}}_{\zeta\Lambda\zeta'}([x,\Lambda,x'],t,-\infty)=\psi^{\textup{an}}_{\zeta\Lambda}([x,\Lambda],t),\hspace{1cm}\psi^{\textup{an}}_{\zeta\Lambda\zeta'}([x,\Lambda,x'],\infty,t)=\psi^{\textup{an}}_{\Lambda\zeta'}([\Lambda,x'],t).
		\]
		In particular, $\psi^{\textup{an}}_{\zeta\Lambda\zeta'}$ vanishes on $X_{\zeta,\zeta'}\times \{(\infty,-\infty)\}$ and the intersection of the complement of this space and $\left(\psi^{\textup{an}}_{\zeta\Lambda\zeta'}\right)^{-1}(0)$ is a stratified-smooth space.
		\item[(ii)] There is an injective open continuous map $\Phi_{\zeta\Lambda\zeta'}: \left(\psi^{\textup{an}}_{\zeta\Lambda\zeta'}\right)^{-1}(0)\to M^+_{z}(W,c; \alpha,\alpha')$ such that $\Phi_{\zeta\Lambda\zeta'}$ maps 
		$X_{\zeta,\zeta'}\times\{(\infty,-\infty)\}$ to $X_{\zeta,\zeta'}\subset M^+(W;\alpha,\alpha')$ via the identity map. Moreover, the restriction of $\Phi_{\zeta\Lambda\zeta'}$ to the intersection of 
		$X_{\zeta,\zeta'}\times \big((0, \infty] \times [-\infty, 0)\setminus \{(\infty,-\infty)\}\big)$ and $\left(\psi^{\textup{an}}_{\zeta\Lambda\zeta'}\right)^{-1}(0)$ gives an open embedding of stratified-smooth spaces.
		\item[(iii)] For $\beta \in \fC_\pi(Y,w)$, let $\zeta_0: \alpha \to \beta$ and $\zeta_1: \beta \to \rho$ be paths of connections on $(Y, w)$ with composite $\zeta_0 * \zeta_1 = \zeta$. For any $[x,y,\Lambda,x']$
		in the subspace 
		\[
		  \breve M^+_{\zeta_0}(Y,w;\alpha,\beta)\times_\beta \breve M^+_{\zeta_1}(Y,w;\beta,\rho)\times_\rho \mathcal O_\Lambda \times_{\rho'} \breve M_{\zeta'}^+(Y',w';\rho',\alpha')
		\]
		of $X_{\zeta,\zeta'}$, we have
		\[
		  \psi^{\textup{an}}_{\zeta\Lambda\zeta'}([x,y,\Lambda,x'],t,t')=\psi^{\textup{an}}_{\zeta_1\Lambda\zeta'}([y,\Lambda,x'],t,t'),
		\]		
		and if $\psi^{\textup{an}}_{\zeta\Lambda\zeta'}([y,\Lambda,x'],t,t')=0$, then
		\[
		  \Phi_{\zeta\Lambda\zeta'}([x,y,\Lambda,x'],t,t')=(x,\Phi_{\zeta_1\Lambda\zeta'}([y,\Lambda,x'],t,t')).
		\]		
		A similar claim holds for any $\beta' \in \fC_\pi(Y',w')$, paths of connections $\zeta_0': \rho' \to \beta'$ and $\zeta_1':\beta' \to\alpha'$ on $(Y', w')$ with composite $\zeta_0' * \zeta_1' = \zeta'$, and any $[x,\Lambda,x',y']$ in the following 
		subspace of $X_{\zeta,\zeta'}$:
		\[
		  \breve M^+_{\zeta}(Y,w;\alpha,\rho)\times_\rho \mathcal O_\Lambda \times_{\rho'} \breve M^+_{\zeta_0'}(Y',w';\rho',\beta')\times_{\beta'} \breve M^+_{\zeta_1'}(Y',w';\beta',\alpha')	  
		\]
		\item[(iv)] Let $[x,\Lambda,x']\in X_{\zeta,\zeta'}$. If $\psi^{\textup{an}}_{\zeta\Lambda}([x,\Lambda],t)=0$, we have
		\[
		  \Phi_{\zeta\Lambda\zeta'}([x,\Lambda,x'],t,-\infty)=[\Phi_{\Lambda\zeta'}([x,\Lambda],t),x'].
		\]
		Similarly, if $\psi^{\textup{an}}_{\Lambda\zeta'}([\Lambda,y],t)=0$, we have
		\[
		  \Phi_{\zeta\Lambda\zeta'}([x,\Lambda,x'],\infty,t)=[x,\Phi_{\Lambda\zeta'}([\Lambda,x'],t)].
		\]	
	\end{itemize}		
\end{prop}

\subsection{Modified moduli spaces}\label{subsec:modified-mod-space}
In the previous section, we investigated the local structure of the obstructed reducibles in $M^+_z(W, c;\alpha, \alpha')$. In this subsection, we modify the moduli spaces $M^+_{z}(W,c; \alpha, \alpha')$ so that the resulting moduli spaces are stratified-smooth, and so that we have explicit control over their boundary components. This modification involves two steps: removing a neighborhood of obstructed solutions and adding several new components to the moduli spaces. The first step of this process is performed using the description of neighborhoods of obstructed solutions given in the last subsection. The second step involves choosing some additional sections of the obstruction bundles, which we will fix next. Throughout this section $z$ is a path of connections along $(W,c)$ with $i(z) \le 6$.

We write $\psi^T_{\zeta\Lambda}$ and $\psi^T_{\Lambda\zeta'}$ for the restrictions of the sections $\psi^{\text{an}}_{\zeta\Lambda}$ and $\psi^{\text{an}}_{\Lambda\zeta'}$ over $\breve M^+_{\zeta}(Y,w;\alpha,\rho)\times \{T\}$ and $\breve M_{\zeta'}^+(Y',w';\rho',\alpha')\times \{-T\}$ for some large value $T \in (0, \infty)$ that will be specified later. We fix an $SO(3)$-equivariant section $\psi^{\text{priv}}_{\zeta\Lambda}$ of $\mathcal H_{\alpha,\rho}$ over $\breve M^+_{\zeta}(Y,w;\alpha,\rho)$. 

Unlike $\psi^T_{\zeta\Lambda}$, the section $\psi^{\text{priv}}_{\zeta\Lambda}$ is not determined by analysis, but is rather pre-chosen (and independent of the choice of $T$). In particular, we regard $\psi^{\text{priv}}_{\zeta\Lambda}$ as additional information attached to $(Y,w)$, rather than being determined by the cobordism $(W,c)$. 

We also need to fix a section that interpolates between $\psi^T_{\zeta\Lambda}$ and $\psi^{\text{priv}}_{\zeta\Lambda}$. Let $\psi^{\text{htpy}}_{\zeta\Lambda}$ be a section of $[0,1] \times \mathcal H_{\alpha,\rho}$ over $[0,1] \times \breve M^+_{\zeta}(Y,w;\alpha,\rho)$ whose restrictions to $\{0\} \times \breve M^+_{\zeta}(Y,w;\alpha,\rho)$ and $\{1\} \times \breve M^+_{\zeta}(Y,w;\alpha,\rho)$ are respectively equal to $\psi^T_{\zeta\Lambda}$ and $\psi^{\text{priv}}_{\zeta\Lambda}$.

Similarly, we fix $\psi^{\text{priv}}_{\Lambda\zeta'}$ and $\psi^{\text{htpy}}_{\Lambda\zeta'}$ that are respectively sections of $\mathcal H_{\rho',\alpha'}$ and $[0,1] \times \mathcal H_{\rho',\alpha'}$. Moreover, the restriction of $\psi^{\text{htpy}}_{\Lambda,\alpha'}$ to $\{0\} \times \breve M_{\zeta'}^+(Y',w';\rho',\alpha')$ and $\{1\} \times \breve M_{\zeta'}^+(Y',w';\rho',\alpha')$ are respectively equal to $\psi^T_{\Lambda\zeta'}$ and $\psi^{\text{priv}}_{\Lambda\zeta'}$. The sections fixed in this paragraph are required to satisfy the following properties:

\begin{enumerate}
	\item[(s.1)] The number $T$ is chosen such that for any $\zeta$ and $\zeta'$ the sections $\psi^T_{\zeta \Lambda}$ and $\psi^T_{\Lambda \zeta'}$, are transverse to the zero section on each stratum. We also demand that the various sections $\psi^{\text{priv}}$ and $\psi^{\text{htpy}}$ are transverse to the zero section on each stratum.
	\item[(s.2)] For any pair $\zeta$, $\zeta'$, we may regard $t\cdot \psi^{\text{priv}}_{\zeta\Lambda}+(1-t)\cdot \psi^{\text{priv}}_{\Lambda\zeta'}$ with $t\in [0,1]$ as a section of $\pi_1^*(\mathcal H_{\alpha,\alpha'})$ defined over the stratified-smooth manifold $[0,1] \times X_{\zeta,\zeta'}$. This section is required to be transverse to the zero section on each stratum.
	\item[(s.3)] For $\beta \in \fC_\pi(Y,w)$, let $\zeta_0: \alpha \to \beta$ and $\zeta_1: \beta \to \rho$ be paths of connections on $(Y, w)$ with composite $\zeta_0 * \zeta_1 = \zeta$. The restriction of $\psi^{\text{priv}}_{\zeta \Lambda}$ to 
		\begin{equation}\label{co-dim-face-2}
		  \breve M^+_{\zeta_0}(Y,w;\alpha,\beta)\times_\beta \breve M^+_{\zeta_1}(Y,w;\beta,\rho)\subset \breve M^+_{\zeta}(Y,w;\alpha,\rho)
		\end{equation}
		is given by $\psi^{\text{priv}}_{\zeta_1\Lambda}$. Similarly, the restriction of $\psi^{\text{htpy}}_{\zeta\Lambda}$ to
		\[
		  [0,1] \times \breve M^+_{\zeta_0}(Y,w;\alpha,\beta)\times_\beta \breve M^+_{\zeta_1}(Y,w;\beta,\rho)\subset [0,1] \times \breve M^+_{\zeta}(Y,w;\alpha,\rho)
		\]
	is given by $\psi^{\text{htpy}}_{\zeta_1\Lambda}$.
	\item[(s.4)] For $\beta' \in \fC_\pi(Y',w')$, let $\zeta_0': \rho' \to \beta'$ and $\zeta_1': \beta' \to \alpha'$ be paths of connections on $(Y', w')$ with composite $\zeta_0' * \zeta_1' = \zeta'$. The restriction of $\psi^{\text{priv}}_{\Lambda\zeta'}$ to 
		\begin{equation}\label{co-dim-face-3}
		  \breve M^+_{\zeta_0'}(Y',w';\rho',\beta')\times_{\beta'} \breve M^+_{\zeta_1'}(Y',w';\beta',\alpha')\subset \breve M^+_{\zeta'}(Y',w';\rho',\alpha')
		\end{equation}
		is given by $\psi^{\text{priv}}_{\Lambda\zeta_0'}$. Similarly, the restriction of $\psi^{\text{htpy}}_{\Lambda\zeta'}$ to
		\[
		   [0,1] \times \breve M^+_{\zeta_0'}(Y',w';\rho',\beta')\times_{\beta'} \breve M^+_{\zeta_1'}(Y',w';\beta',\alpha')\subset [0,1] \times \breve M^+_{\zeta'}(Y',w';\rho',\alpha')
		\]
	is given by $\psi^{\text{htpy}}_{\Lambda\zeta_0'}$.
\end{enumerate}

We shall also need some variations of the sections $\psi^{\text{an}}_{\zeta\Lambda\zeta'}$. First consider the arc 
\[
  A_T =\Big\{(t,t')\in (0, \infty] \times [-\infty, 0) \;\Big|\; \frac{1}{t'} - \frac{1}{t} = -\frac{1}{T}\Big\},
\]
which is oriented as an inverse image of $1/t' - 1/t$, so that its boundary has the point $(\infty,-T)$ with negative sign and the point $(T,-\infty)$ with positive sign. This is oriented diffeomorphic to the unit interval $[0,1]$ via the map $(t,t') \mapsto T/t = 1+T/t'$. 

Define $\psi^T_{\zeta\Lambda\zeta'}$ be the restriction of $\psi^{\text{an}}_{\zeta\Lambda\zeta'}$ to $A_T \times X_{\zeta,\zeta'}$. Let $\psi^{\text{priv}}_{\zeta\Lambda\zeta'}$ be the section over $A_T \times X_{\zeta,\zeta'}$ defined by 
\[\psi^{\text{priv}}_{\zeta\Lambda\zeta'}([x,y],t,t') = \frac{T}{t}\cdot \psi^{\text{priv}}_{\zeta \Lambda}(x)-\frac{T}{t'}\cdot \psi^{\text{priv}}_{\Lambda\zeta'}(y).\] 
Via the diffeomorphism $A_T \cong [0,1]$ above, this map is identified with the section $s\psi^{\text{priv}}_{\zeta \Lambda} + (1-s)\psi^{\text{priv}}_{ \Lambda \zeta'}.$ Our assumption (s.2) on the sections $\psi^{\text{priv}}_{\zeta\Lambda}$ and $\psi^{\text{priv}}_{\Lambda\zeta'}$ imply that $\psi^{\text{priv}}_{\zeta\Lambda\zeta'}$ is transverse to the zero section on each stratum. 

Finally, let $\psi^{\text{htpy}}_{\zeta\Lambda\zeta'}$ be a section of $\pi_1^*(\mathcal H_{\alpha,\alpha'})$ over $[0,1] \times A_T \times X_{\zeta,\zeta'}$ that satisfies the following properties:
\begin{enumerate}
	\item[(c.1)] As another requirement on $T$, we demand that $\psi^T_{\zeta\Lambda\zeta'}$ is transverse to the zero section on each stratum. We also demand that the section $\psi^{\text{htpy}}_{\zeta\Lambda\zeta'}$ is also transverse to the zero section on each stratum.
	\item[(c.2)] The restriction of $\psi^{\text{htpy}}_{\zeta,\Lambda,\zeta'}$ to $\{0\} \times A_T \times X_{\zeta,\zeta'}$ and $\{1\} \times A_T \times X_{\zeta,\zeta'}$ are respectively equal to $\psi^T_{\zeta\Lambda\zeta'}$
	and $\psi^{\text{priv}}_{\zeta \Lambda\zeta'}$.
	\item[(c.3)] For any $[x,y]\in X_{\zeta,\zeta'}$ and $t\in[0,1]$, we have 
	\[\psi^{\text{htpy}}_{\zeta\Lambda\zeta'}(t,T,-\infty,[x,y])=\psi^{\text{htpy}}_{\zeta\Lambda}(t,[x]),\hspace{1cm} \psi^{\text{htpy}}_{\zeta\Lambda\zeta'}(t,\infty,T,[x,y])=\psi^{\text{htpy}}_{\Lambda\zeta'}(t,[y]).\]
	\item[(c.4)] For $\beta \in \fC_\pi(Y,w)$, let $\zeta_0: \alpha \to \beta$ and $\zeta_1: \beta \to \rho$ be paths of connections on $(Y, w)$ with composite $\zeta_0 * \zeta_1 = \zeta$. Let $[x,y,\Lambda,x']$ be an element of the subspace 
	\[
	  \breve M^+_{\zeta_0}(Y,w;\alpha,\beta)\times_\beta \breve M^+_{\zeta_1}(Y,w;\beta,\rho)\times_\rho \mathcal O_\Lambda \times_{\rho'} \breve M_{\zeta'}^+(Y',w';\rho',\alpha') \subset X_{\zeta,\zeta'}.
	\]
	Then for any $(s,s')\in A_T$ and $t\in [0,1]$, we have \[\psi^{\text{htpy}}_{\zeta\Lambda\zeta'}(t,s,s',[x,y,\Lambda,x'])=\psi^{\text{htpy}}_{\zeta_1\Lambda\zeta'}(t,s,s',[y,\Lambda,x']).\] 
	Similarly, for any $\beta' \in \fC_\pi(Y',w')$, paths of connections $\zeta_0': \rho' \to \beta'$ and $\zeta_1': \beta' \to \alpha'$ on $(Y', w')$ with composite $\zeta_0' = \zeta_1' = \zeta'$ and any element $[x,\Lambda,x',y']$ of
	\[
	  \breve M^+_{\zeta}(Y,w;\alpha,\rho)\times_\rho \mathcal O_\Lambda \times_{\rho'} \breve M^+_{\zeta_0'}(Y',w';\rho',\beta')\times_{\beta'} \breve M^+_{\zeta_1'}(Y',w';\beta',\alpha')\subset X_{\zeta,\zeta'},
	\]
	we have 
	\[
	  \psi^{\text{htpy}}_{\zeta\Lambda\zeta'}(t,s,s',[x,\Lambda,x',y'])=\psi^{\text{htpy}}_{\zeta\Lambda\zeta_0'}(t,s,s',[x,\Lambda,x']),
	\] 
	for any $(s,s')\in A_T$ and $t\in [0,1]$. 
\end{enumerate}

Now, we are ready to define the modified moduli spaces. Fix $\alpha\in \fC_\pi(Y,w)$, $\alpha'\in \fC_\pi(Y',w')$ and a path $z: \alpha \to \alpha'$ of connections along $(W,c)$ with $i(z)\leq 6$. Recall that we have an obstructed abelian instanton $\Lambda:\rho\to \rho'$, which determines a subspace of obstructed solutions in the moduli space $M^+_{z}(W,c; \alpha, \alpha')$. If $\alpha'=\rho'$, then this subspace consists of obstructed solutions of type I given as $\breve M^+_{\zeta}(Y,w;\alpha,\rho)\times_\rho \mathcal O_\Lambda$ with $z$ being the composition of $\zeta$ and the path determined by $\Lambda$. A neighborhood $U_{\zeta\Lambda}$ of such obstructed solutions is given as
\[\Phi_{\zeta\Lambda}\(\(\psi^{\text{an}}_{\zeta\Lambda}\)^{-1}(0)\cap \(\breve M^+_{\zeta}(Y,w;\alpha,\rho) \times(T,\infty]\)\)\subset M^+_{z}(W,c; \alpha, \rho'). \]
We will verify below that the complement of $U_{\zeta\Lambda}$ in $M^+_{z}(W,c; \alpha, \rho')$ is an oriented compact stratified-smooth space; its top stratum is an open subset of $M^+_z(W,c;\alpha, \rho')$, and we give it the same orientation. 

The zeros of $\psi^{\text{htpy}}_{\zeta\Lambda}$ over $[0,1]\times \breve M^+_{\zeta}(Y,w;\alpha,\rho)$ gives another stratified-smooth space $M^{\text{htpy}}_{\zeta\Lambda}$ associated to the path $z = \zeta * \lambda$. To define an orientation on $M^{\text{htpy}}_{\zeta\Lambda}$, first fix the product orientation on $[0,1]\times \breve M^+_{\zeta}(Y,w;\alpha,\rho)$. Since $\psi^{\text{htpy}}_{\zeta\Lambda}$ is a transverse section of a complex bundle, its zero set inherits an orientation from $[0,1]\times \breve M^+_{\zeta}(Y,w;\alpha,\rho)$. We multiply this orientation by $\pm 1$ as appropriate so that $M^{\text{htpy}}_{\zeta\Lambda}$ has $\left(\psi^T_{\zeta\Lambda}\right)^{-1}(0)$ as a boundary component with the opposite orientation as the corresponding boundary component of $M^+_z(W,c;\alpha, \rho') \setminus U_{\zeta \Lambda}.$

In the case $\alpha=\rho$, the moduli space $M^+_{z}(W,c; \alpha, \alpha')$ contains a subspace of obstructed solutions of type II, given as $\mathcal O_\Lambda\times_{\rho'} \breve M_{\zeta'}^+(Y',w';\rho',\alpha')$. We define the neighborhood $U_{\Lambda\zeta'}$ of these obstructed solutions as the image of 
\[
  \Phi_{\Lambda\zeta'}\(\(\psi^{\text{an}}_{\Lambda\zeta'}\)^{-1}(0)\cap\(\breve M_{\zeta'}^+(Y',w';\rho',\alpha') \times [-\infty, T)\)\)\subset M^+_{z}(W,c;\rho,\alpha').
\] 
The compact stratified-smooth space given by the complement of $U_{\Lambda\zeta'}$ is oriented following the same convention as in the unobstructed case. We also define $M^{\text{htpy}}_{\Lambda\zeta'}$ to be the stratified-smooth space given as the zeros of $\psi^{\text{htpy}}_{\Lambda\zeta'}$ over $[0,1] \times \breve M_{\zeta'}^+(Y',w';\rho',\alpha')$. This space is again a compact stratified-smooth space, and its top stratum is oriented using the same convention as in the previous paragraph: one of its boundary components is canonically diffeomorphic to a boundary component in the complement of $U_{\Lambda \zeta'}$, and the orientations on these two boundary components should be opposite to one another.

For $\alpha\neq \rho$ and $\alpha'\neq \rho'$, the space $M^+_{z}(W,c; \alpha, \alpha')$ contains obstructed solutions of type III given as $X_{\zeta,\zeta'}$ where the composition $\zeta * \lambda * \zeta' = z$. The index assumption $i(z) \le 6$ implies that there are unique choices of $\zeta$ and $\zeta'$ such that $X_{\zeta,\zeta'}$ is non-empty. A neighborhood $U_{\zeta\Lambda\zeta'}$ of these obstructed solutions can be obtained as the image 
\[
  \Phi_{\zeta\Lambda\zeta'}\(\(\psi^{\text{an}}\)^{-1}(0)\cap \left(\left\{(s,s')\in (0, \infty] \times [-\infty, 0) \;\Big| \;\frac{1}{s'} - \frac{1}{s'}\leq -\frac{1}{T}\right\} \times X_{\zeta,\zeta'}\right)\)\subset M^+_{z}(W,c;\rho,\alpha').
\]
The complement of this open neighborhood is a compact stratified-smooth space and can be oriented following the same conventions as in the unobstructed case. Let also $M^{\text{htpy}}_{\zeta\Lambda\zeta'}$ be the stratified-smooth space given as the zero set of $\psi^{\text{htpy}}_{\zeta\Lambda\zeta'}$ over $[0,1] \times A_T \times X_{\zeta,\zeta'}$. If the domain of this section is given the product orientation, then we fix the orientation on $M^{\text{htpy}}_{\zeta\Lambda\zeta'}$ which differs from the standard orientation on the zero set of a complex line bundle by the appropriate factor $\pm 1$ so that the appropriate boundary component has the opposite sign as the corresponding boundary component of the complement of $U_{\zeta \Lambda \zeta'}.$ 

The {\it modified moduli space} $N^+_{z}(W,c; \alpha, \alpha')$ is the $SO(3)$-space given as 
\begin{equation}\label{modified-mod-space}
	N^+_{z}(W,c; \alpha, \alpha'):=\(M^+_{z}(W,c; \alpha, \alpha') \setminus U_\Lambda\) \sqcup M^{\text{htpy}}_\Lambda.
\end{equation}
For the sake of uniformity, we have dropped the notation $\zeta$ and $\zeta'$ from our notations for the spaces $U_{\zeta\Lambda}$, $U_{\Lambda\zeta'}$ and $U_{\zeta\Lambda\zeta'}$, as well as from the spaces $M^{\text{htpy}}_{\zeta\Lambda}$, $M^{\text{htpy}}_{\Lambda\zeta'}$ and $M^{\text{htpy}}_{\zeta\Lambda\zeta'}$.

\begin{prop}\label{obs-bdry-reln}
	The modified moduli space $N^+_{z}(W,c; \alpha, \alpha')$ is a compact oriented stratified-smooth space. The boundary of $N^+_{z}(W,c; \alpha, \alpha')$ is the union of the following compact orientable smooth-stratified spaces.
	\begin{itemize}
		\item[(i)] If $\beta\in \fC_\pi(Y,w)$, and $\zeta_0: \alpha \to \beta$ and $z_1: \beta \to \alpha'$ are paths of connections along $(Y,w)$ and $(W,c)$ respectively with composite $\zeta_0 * z_1 = z$, then 
		\[
		  \breve M^+_{\zeta_0}(Y,w;\alpha,\beta) \times_{\beta} N^+_{z_1}(W,c; \beta, \alpha')
		\] 
		is a boundary component of $N^+_{z}(W,c; \alpha, \alpha')$. The boundary orientation and the fiber product orientation on this boundary component are related to each other by $(-1)^{\dim(\alpha)}$.
		\item[(ii)] If $\beta'\in \fC_{\pi'}(Y',w')$, and $z_0: \alpha \to \beta'$ and $\zeta_1: \beta' \to \alpha'$ are paths of connections along $(W,c)$ and $(Y',w')$ respectively with composite $z_0 * \zeta'_1 = z$, then  
		\[
		  N^+_{z_0}(W,c; \alpha, \beta') \times_{\beta'} \breve M^+_{\zeta_1'}(Y',w';\beta',\alpha')
		\] 
		is a boundary component of $N^+_{z}(W,c; \alpha, \alpha')$. The boundary orientation and the fiber product orientation on this boundary component are related to each other by $(-1)^{\wt i(z_0)+\dim(\alpha)+1}$.
	\end{itemize}		
	Depending on the type of $\alpha$ and $\alpha'$, the remaining boundary components are given as in one of the following three cases:
	\begin{itemize}		
		\item[(l)] If $\alpha'=\rho'$ and $\zeta: \alpha \to \rho$ is a path of connections on $(Y,w)$ with composite $\zeta * \lambda = z$, then the boundary components of $N^+_{z}(W,c; \alpha, \alpha')$ contain
		\begin{itemize}
			\item[(l-i)] two copies of $\(\psi^T_{\zeta\Lambda}\)^{-1}(0)$ with opposite orientations;
			\item[(l-ii)] a copy of $\(\psi^{\textup{priv}}_{\zeta\Lambda}\)^{-1}(0)$ whose orientation differs from the standard orientation of the zero section of a complex line bundle over an oriented space by a factor of $(-1)^{\dim(\alpha)}$.
		\end{itemize}
		\item[(r)] If $\alpha=\rho$ and $\zeta': \rho' \to \alpha'$ is a path of connections on $(Y', w')$ with composite $\lambda * \zeta' = z$, then the boundary components of $N^+_{z}(W,c; \alpha, \alpha')$ contain
		\begin{itemize}
			\item[(r-i)] two copies of $\(\psi^T_{\Lambda\zeta'}\)^{-1}(0)$ with opposite orientations;
			\item[(r-ii)] a copy of $\(\psi^{\textup{priv}}_{\Lambda\zeta'}\)^{-1}(0)$ whose orientation differs from the standard orientation of the zero section of a complex line bundle over an oriented space by a factor of $-1 = (-1)^{\dim(\rho) + 1}$.
		\end{itemize}
		\item[(c)] If $\alpha\neq \rho$ and $\alpha'\neq \rho'$, and $\zeta: \alpha \to \rho, \;\zeta': \rho' \to \alpha'$ are paths of connections on $(Y,w)$ and $(Y',w')$ with composite $\zeta * \lambda * \zeta' = z$, then the boundary components of $N^+_{z}(W,c; \alpha, \alpha')$ contain
		\begin{itemize}
			\item[(c-i)] two copies of $\(\psi^T_{\zeta\Lambda\zeta'}\)^{-1}(0)$ with opposite orientations;
			\item[(c-ii)] a copy of $\(\psi^{\textup{priv}}_{\zeta\Lambda\zeta'}\)^{-1}(0)$, whose orientation differs from the standard orientation of the zero set of the section $s\psi^{\textup{priv}}_{\zeta \Lambda} + (1-s)\psi^{\textup{priv}}_{ \Lambda \zeta'}$ over $\breve M^+_\zeta(Y,w;\alpha, \rho) \times_\rho \([0,1] \times \breve M^+_{\zeta'}(Y',w';\rho',\alpha')\)$ by a factor of $(-1)^{\dim \alpha-1}$.
		\end{itemize}		
	\end{itemize}
\end{prop}

\begin{proof}[Sketch of proof.] What we will verify is that $N^+_z$ is a stratified-smooth space, and we will confirm the relevant boundary orientations; the actual enumeration of boundary strata is straightforward. The fact that $M^{\text{htpy}}_\Lambda$ is an oriented stratified-smooth space is the content of Proposition \ref{bd-fibprod-or}. 

As for the complement of $U_\Lambda$, one should define a function $L: M^+_z(W,c;\alpha, \alpha') \to (-\infty, 0)$ so that on a slight thickening of $U_\Lambda$ it takes the following form in the three different cases. In the first case, $L(x,\Lambda,t) = -1/t$. In the second case, $L(\Lambda,y,t') = 1/t'$. In the third case, $L(x,t,\Lambda,t',y) = 1/t' - 1/t$. In all cases, the set $U_\Lambda$ is defined by the property $L > -1/T$. We extend $L$ arbitrarily to the whole of $M^+_z$, so long as $L^{-1}[-1/T, 0) = \overline{U_\Lambda}$ remains true. Then Proposition \ref{bd-trunc-or} guarantees that the complement of $U_\Lambda$ is an oriented stratified-smooth space with orientation induced as a codimension-0 subspace of $M^+_z$.

What remains is to confirm the stated boundary orientations; because $M^{\text{htpy}}$ is defined so that its positive boundary component is oriented by the same conventions as $\(\psi^T\)^{-1}(0)$, it suffices to compute the boundary orientation on the stratum $L^{-1}(-1/T)$ using Proposition \ref{bd-trunc-or}. That formula contains a sign of $(-1)^{\dim X - 1}$, which one can eliminate by instead computing the orientation using the convention that the \emph{first} factor maps via $dL$ in an orientation-preserving way to $\mathbb R$. We may ignore the discussion of taking the zero set of $\psi^T$, as taking the zero set of a complex bundle induces an orientation in a canonical way which commutes with all other orientation-related operations.

Now the conventions of Section \ref{sec:moduli-or} have $M^+_\zeta(Y, w; \alpha, \rho)$ oriented as 
\[\alpha \times \Bbb R \times \breve M^{\text{Fiber}}_\zeta(Y, w; \alpha, \rho).\] 
Commuting the $\Bbb R$ factor to the front adds a sign of $(-1)^{\dim \alpha}$. The map $L$ is then given by $L(t, a, x) = -1/t$; the differential of this is $dL = t^{-2} dt$, so $dL$ gives an orientation-preserving isomorphism from the first factor to the last. Therefore $L^{-1}(-1/T)$ is oriented as 
\[(-1)^{\dim \alpha} \alpha \times \breve M^{\text{Fiber}}_\zeta \cong (-1)^{\dim \alpha} \breve M_\zeta(Y, w;\alpha, \rho).\]
This justifies the first sign. The second is similar: this space is oriented as $\rho' \times \Bbb R \times \breve M^{\text{Fiber}}_{\zeta'}(Y', w'; \rho', \alpha')$. Moving the factor of $\Bbb R$ to the front introduces a sign of $(-1)^{\dim \rho'} = 1$. The map is given by $L(t', a', A') = 1/t'$, so that $dL = -(t')^{-2} dt'$; this is orientation-reversing, which introduces an additional sign of $-1$, giving the statement in the proposition.

Finally, the set $L^{-1}(-1/T)$ over $X_{\zeta, \zeta'}$ is oriented as the zero set of a map from the following space to $\Bbb R$:
\[\alpha \times \Bbb R \times \breve M^{\text{Fiber}}_\zeta(Y, w;\alpha, \rho) \times \Bbb R \times \breve M^{\text{Fiber}}_{\zeta'}(Y', w'; \rho', \alpha').\] 
Commuting all factors of $\Bbb R$ to the front induces a sign of $(-1)^{2\dim \alpha + \wt i(\zeta) - 1}$; passing to the boundary face $L^{-1}(-1/T)$ gives us 
\[(-1)^{\wt i(\zeta)} A_T \times \alpha \times \breve M^{\text{Fiber}}_\zeta(Y, w; \alpha, \rho) \times \breve M^{\text{Fiber}}_{\zeta'}(Y', w'; \rho', \alpha'),\] 
where the extra factor of $-1$ appears because $A_T$ is oriented as the inverse image of $L^{-1}(-1/T)$, but placing the $dL$ factor \emph{first} is the opposite orientation convention on $A_T$. Finally, commuting $A_T$ back to the middle (and recalling that the oriented diffeomorphism from $A_T$ to $[0,1]$ sends our section to the one stated in the proposition) gives us the stated sign of $(-1)^{\dim \alpha -1}$.
\end{proof}

\begin{remark}
Because the map $\psi^{\text{priv}}_{\zeta \Lambda \zeta'}$ may be identified with the section $t\cdot \psi^{\text{priv}}_{\zeta\Lambda}+(1-t)\cdot \psi^{\text{priv}}_{\Lambda\zeta'}$ over $\breve M_\zeta(Y, w; \alpha, \rho) \times [0,1] \times \breve M_{\zeta'}(Y',w';\rho',\alpha')$, it can be described entirely in terms of these sections. We give an alternate description of this zero set in Lemma \ref{weird-zero-set} in the appendix, which we will use in Section 5.4 below when verifying certain boundary relations in the definition of bimodule. 
\end{remark}
\subsection{The suspended flow category}\label{cob-flow}
%!TEX root = equivariant-functoriality.tex

In the previous section, we defined modified moduli spaces $N^+_{\zeta}(W,c; \alpha, \alpha')$ attached to $(W,c)$. We would like to use the machinery of Section \ref{FlowCatS} to say that these define a bimodule (as in Definition \ref{def:bimod}) between flow categories, and hence a chain map between their flow complexes. This amounts to a boundary relation on these moduli spaces.

Inspecting the boundary relation of Proposition \ref{obs-bdry-reln}, we see that these moduli spaces do not define a bimodule $\mathcal{I}(Y,w,\pi) \to \mathcal{I}(Y',w',\pi')$ between the framed instanton flow categories. However, we still expect to be able to extract a chain map from these moduli spaces, using a \emph{slightly different flow category} as the codomain. We start by describing this slightly different flow category, via a construction called the \emph{suspended flow category}.

Before doing so, we should give a brief introduction to a construction carried out in Appendix \ref{BlowupChain}; see Construction \ref{constr:blowup} for more details. Suppose we are given a stratified-smooth space $P$ with a map to a smooth manifold $\phi: P \to M$, a complex line bundle $L \to M$, and a lift $\psi: P \to L$ of $\phi$ to a section of $L$ transverse to the zero section. Then we may define its {\it real blowup}, a stratified-smooth space $B(\psi)$ with a natural map to the sphere bundle $S(L)$; there is a natural projection map $\pi: B(\psi) \to P$, and if considered as a geometric chain in $M$, the chain $P$ is \emph{collapse-equivalent} to its real blowup $B(\psi)$. The real blowup is described in more detail in Appendix \ref{BlowupChain}; what matters here is that collapse-equivalence is preserved by fiber products, and that two collapse-equivalent stratified-smooth spaces define the same element in the geometric chain complex. This chain in $S(L)$ satisfies the boundary relation \[\partial B(\psi) = B(\psi|_{\partial P}) - (-1)^{\dim P} Z(\psi) \times_M S(L),\] where $Z(\psi)$ is the zero set of this section and the map from the second chain to $S(L)$ is projection to the second factor.

This construction allows us to lift chains in $M$ to chains in the sphere bundle $S(L)$, which we use in the construction of the suspended flow category below.

\begin{construction}\label{constr:blowup-flowcat} Let $\mathcal C$ be an $SO(3)$-equivariant, finitely graded flow category as in Definition \ref{def:flowcat}, whose objects are all $SO(3)$-orbits isomorphic to one of $SO(3), S^2$, or the singleton $\ast$. We write $\text{Ob}(\mathcal C) = \mathfrak C^* \sqcup \mathfrak A \sqcup \mathfrak Z$ for the decomposition of the set of objects into these three types of orbit; $\mathfrak C^*$ is the set of objects isomorphic to $SO(3)$, while $\mathfrak A$ is the set of objects isomorphic to $S^2$ and $\mathfrak Z$ is the set of singleton objects. 

Further, assume that each $\mathcal C(\alpha, \beta)$ is a free $SO(3)$-manifold. 

Choose an element $\rho \in \mathfrak A$. Fix an $SO(3)$-equivariant complex line bundle $\mathcal H_\rho \to \rho$, where $SO(3)$ acts freely on any nonzero vector; one may take $\mathcal H_\rho = T\rho \cong SO(3) \times_{SO(2)} \Bbb C$. 

Choose, for each $\beta \in \mathsf{Ob}(\mathcal C)$, a lift $\psi_{\rho \beta}: \mathcal C(\rho, \beta) \to \mathcal H_\rho$ and which is transverse to $0$. We demand these are compatible in the sense that, on the fiber product 
\[\mathcal C(\rho, \alpha) \times_\alpha \mathcal C(\alpha, \beta) \subset \partial \mathcal C(\rho, \beta),\] 
the restriction of $\psi_{\rho \beta}$ coincides with the pullback of $\psi_{\rho \alpha}$ to the fiber product.\footnote{Verifying that such $\psi_{\rho \beta}$ exist is a straightforward induction; the key point is that every $\mathcal C(\rho, \beta)$ is $SO(3)$-free, so there is no trouble arising from the demand for equivariant transversality.}

These choices in hand, we define the suspended flow category $\mathcal S_{\rho}(\mathcal C)$ (or simply $\mathcal S$ when the flow category $\mathcal C$ and object $\rho$ are clear from context) as follows; the construction depends on the family of sections $\psi_{\rho \beta}$, but we suppress it from notation. 

First we will describe the object space. Each object $\alpha \in \mathsf{Ob}(\mathcal C) \setminus \{\rho\}$ is left unchanged and gives an object of $\mathcal S$, with the same relative grading. The chosen orbit $\rho$ gives rise to two objects of $\mathcal S$. The first is $S_\rho = S(\mathcal H_\rho)$, the unit sphere bundle of $\mathcal H_\rho \to \rho$. This space is isomorphic to $SO(3)$, and will be given relative grading one lower than $\rho$. The second is $\rho$ itself, but it will have relative grading 2 less than $\rho$. To distinguish, we will write this as $\overline \rho$. That is, 
\[\mathsf{Ob}(\mathcal S) = \big(\mathsf{Ob}(\mathcal C) \setminus \{\rho\}\big) \sqcup S_\rho \sqcup \overline \rho,\] whe relative grading is given by the original relative grading plus a shift factor: 
\[\wt i_{\mathcal S}(\alpha, \beta) = \wt i_{\mathcal C}([\alpha], [\beta]) + s(\alpha) - s(\beta).\] 
If $\alpha \in \mathsf{Ob}(\mathcal S)$, one should interpret \[[\alpha] = \begin{cases} \rho & \alpha = S_\rho, \; \overline \rho \\ \alpha & \text{otherwise} \end{cases}\]

The shift factor is $s(\alpha) = 0$ for $\alpha \neq \overline \rho, S_\rho$ but $s(\overline \rho) = -2$ and $s(S_\rho) = -1$.

We move on to defining the morphism spaces $\mathcal S(\alpha, \beta)$ in our suspended flow category. Let us introduce some notation: 

\begin{itemize}
\item Write $Z(\rho, \beta) = Z(\psi_{\rho \beta}) = \psi_{\rho\beta}^{-1}(0)$ for the stratified-smooth space given by the zero locus of $\psi_{\rho \beta}$; this space has dimension two less than $\mathcal C(\rho, \beta)$. 
\item Write $B(\rho, \beta)$ for the real blowup $B(\psi_{\rho \beta})$ discussed above; this has the same dimension as $\mathcal C(\rho, \beta)$. This is equipped with canonical maps $\psi/\|\psi\|: B(\rho,\beta) \to S_\rho$ and $e_+ \pi: B(\rho, \beta) \to \beta$.
\end{itemize}

Then we define the morphism spaces $\mathcal S_\rho(\alpha, \beta)$ as in the following matrix:

\[\mathcal S(\alpha, \beta) = \begin{pmatrix}[c|ccc]
& \alpha \ne \rho & S_\rho & \overline \rho \\
\hline
\beta \ne \rho & \mathcal C(\alpha, \beta) & -B(\rho, \beta) & Z(\rho, \beta) \\
S_{\rho} & \mathcal C(\alpha, \rho) \times_\rho S_{\rho} & 0 & 0 \\
\overline \rho & 0 & S_\rho & 0
\end{pmatrix}\]

Here the columns correspond to the source critical orbit, and the rows correspond to the target critical orbit.
\end{construction}

\begin{remark}
The fiber product map $\times_{S_\rho} S_\rho: C_*(S_\rho) \to C_*(\overline \rho)$ arising from the term in entry $(3,2)$ of the above matrix is the induced map of $\pi: S_\rho \to \rho$; that is, if $\phi: P \to S_\rho$ is a geometric chain, then $\phi \times_{S_\rho} S_\rho = \pi \phi: P \to \overline \rho$ is the composite of $\phi$ with projection to the base. Thus, one can find a copy of the algebraic mapping cone of $\pi_*: C_*(S_\rho) \to C_*(\overline \rho)$ inside the Morse chain complex of this flow category.
\end{remark}

In a moment we will verify that this does indeed give a flow category. First, we restate the properties of real blowups for the particular spaces $B(\rho, \beta) = B(\psi_{\rho \beta})$. First, $B(\rho, \beta)$ sits in a canonical commutative diagram \[\begin{tikzcd}
	{B(\rho, \beta)} && {S_\rho \times \beta} \\
	\\
	{\mathcal C(\rho, \beta)} && {\rho \times \beta}
	\arrow["{\pi \times 1}", from=1-3, to=3-3]
	\arrow[from=1-1, to=3-1]
	\arrow["{e_- \times e_+}"', from=3-1, to=3-3]
	\arrow["{\psi/\|\psi\| \times e_+}", from=1-1, to=1-3]
\end{tikzcd}\]

Lemma \ref{blowup-equiv} asserts that when considered as a geometric chain downstairs in $\rho \times \beta$ by following the projection $\pi \times 1$,  this space is collapse-equivalent to $\mathcal C(\rho, \beta)$, so define the same elements in the geometric chain complex of $\rho \times \beta$, and induce the same chain-level fiber product maps. Lastly, by the same lemma we have 
\[\partial B(\rho, \beta) = B(\psi_{\rho \beta}|_{\partial \mathcal C(\rho, \beta)}) - (-1)^{\dim \mathcal C(\rho, \beta)} Z(\rho, \beta) \times_\rho S(\mathcal H_\rho).\]
By commuting the odd-dimensional factor $S(\mathcal H_\rho)$ to the left, the last term simplifies to $-S_\rho \times_\rho Z(\rho, \beta)$. As for the first term, expanding out using the boundary formula for $\partial \mathcal C(\rho, \beta)$ from Definition \ref{def:flowcat}, and using the fact that $\psi_{\rho \beta} = \psi_{\rho \alpha} \times 0$ over $\mathcal C(\rho, \alpha) \times_\alpha \mathcal C(\alpha, \beta)$, the whole boundary relation simplifies to 
\[\partial B(\rho, \beta) = - S_\rho \times_\rho Z(\rho, \beta) + \sum_{\alpha \in \mathfrak C \setminus \rho} (-1)^{\wt i_{\mathcal C}(\rho, \alpha)} B(\rho, \alpha) \times_\alpha \mathcal C(\alpha, \beta).\] (Recall here that $\dim \mathcal C(\rho, \alpha) = \wt i_{\mathcal C(\rho, \alpha)} - 1$.) This in hand, we are prepared to establish that $\mathcal S$ is a flow category. 

The reader should have on-hand Definition \ref{def:flowcat} while reading the argument. Below is a diagram depicting the various moduli spaces in the definition of this flow category, as well as where they run from and to; when reading the argument below, keep in mind that the goal is to show that the boundary components of $\mathcal S(\alpha, \beta)$ are a sum over $(-1)^{\wt i(\alpha, \gamma)} \mathcal S(\alpha, \gamma) \times_\gamma \mathcal S(\gamma, \beta)$, and this sum can be understood as running over all two-step paths from $\alpha$ to $\beta$ in the diagram. The moduli space $\mathcal S(\alpha, \gamma)$ decorates the arrow from $\alpha$ to $\gamma$, while the moduli space $\mathcal S(\gamma, \beta)$ decorates the arrow from $\gamma$ to $\beta$. \vspace{-0.5cm}\begin{center}\[\begin{tikzcd}
	& \alpha / \beta \\
	\\
	{S_\rho} && {\overline \rho}
	\arrow["{S_\rho}"', from=3-1, to=3-3]
	\arrow["{\mathcal C(\alpha, \rho) \times_\rho S_\rho}"'{pos=0.4}, bend right=10, from=1-2, to=3-1]
	\arrow["{-B(\rho, \beta)}"'{pos=0.4}, bend right=10, from=3-1, to=1-2]
	\arrow["{Z(\rho, \beta)}"', from=3-3, to=1-2]
	\arrow["{\mathcal C(\alpha, \beta)}", loop above, from=1-2, to=1-2]
\end{tikzcd}\]\end{center}
The term $\alpha / \beta$ means that an arrow leaving this orbit correspond to moduli spaces with source $\alpha$, while arrows running to this orbit correspond to moduli spaces with target $\beta$.

\begin{theorem}\label{thm:B-is-flowcat}
The construction of $\mathcal S_\rho(\mathcal C)$ above defines a flow category as in Definition \ref{def:flowcat}.
\end{theorem}

\begin{proof}
Because all of these spaces are $SO(3)$-spaces and all of the objects of $\mathcal S$ are $SO(3)$-orbits, it follows that all endpoint maps are still submersions. We leave it to the reader to confirm that 
\[\dim \mathcal S(\alpha, \beta) = \dim \alpha + \wt i_{\mathcal S}(\alpha, \beta) - 1\] 
for all $\alpha, \beta$. Now let us check the boundary relations in the definition of flow category, one-by-one. When proving the flow category relations below, we immediately replace every occurrence of $\mathcal S(\alpha, \beta)$ with its definition in terms of $\mathcal C$ from the matrix in Construction \ref{constr:blowup-flowcat}.

%As an aid to the reader, the following is a useful mnemonic. If we write $\mathcal S$ to mean the above matrix and matrix multiplication to mean fiber-product, the flow category relation (up to sign) is $\partial \mathcal S = \mathcal S^2$. 

\begin{itemize}
\item If $\alpha, \beta \neq \overline \rho, S_\rho$, then the flow category relations stipulate that 
\begin{align*}(-1)^{\dim \alpha}\partial \mathcal C(\alpha, \beta) &= \sum_{\gamma \ne \rho} (-1)^{\wt i_{\mathcal S}(\alpha, \gamma)} \mathcal C(\alpha,\gamma) \times_\gamma \mathcal C(\gamma, \beta)\\
&+ (-1)^{\wt i_{\mathcal S}(\alpha, S_\rho)} \mathcal C(\alpha, \rho) \times_\rho S_\rho \times_{S_\rho} (-B(\rho, \beta)).
\end{align*} 
When $\gamma \ne \rho$, the expression $\wt i_{\mathcal S}(\alpha, \gamma)$ is simply $\wt i_{\mathcal C}(\alpha, \gamma)$. Because $\wt i_{\mathcal S}(\alpha, S_\rho) = \wt i_{\mathcal C}(\alpha, \rho) - 1,$ the last term term in the sum simplifies to
\[(-1)^{\wt i_{\mathcal C}(\alpha, \rho)} \mathcal C(\alpha, \rho) \times_\rho B(\rho, \beta).\] 
Because $B(\rho,\beta)$ is collapse-equivalent to $\mathcal C(\rho, \beta)$ as a chain in $\rho \times \beta$, this fiber product is equal to $\mathcal C(\alpha, \rho) \times_\rho \mathcal C(\rho, \beta)$ in the geometric chain complex. Thus the flow category relation for $\mathcal S(\alpha, \beta)$ asks that we have an equality of elements of $C_*(\alpha \times \beta)$ 
\[(-1)^{\dim \alpha} \partial \mathcal C(\alpha, \beta) = \sum_{\gamma \in \mathsf{Ob}(\mathcal C)} (-1)^{\wt i_{\mathcal C}(\alpha, \gamma)} \mathcal C(\alpha, \gamma) \times_\gamma \mathcal C(\gamma, \beta).\] 
This is precisely what it means when we assume that $\mathcal C$ is a flow category to begin with.

\item If $\alpha = S_\rho$ and $\beta \ne \overline \rho, S_\rho$, the desired boundary relation is \begin{align*}(-1)^{\dim S_\rho} (-\partial B(\rho, \beta)) &= (-1)^{\wt i_{\mathcal S}(S_\rho, \overline \rho)} S_\rho \times_\rho Z(\rho, \beta) \\
&+ \sum_{\gamma \ne \rho} (-1)^{\wt i_{\mathcal S}(S_\rho, \gamma)} -B(\rho, \gamma) \times_\gamma \mathcal C(\gamma,\beta).
\end{align*} 

After simplifying the signs and simplifying the last terms, we are trying to show that 
\[\partial B(\rho, \beta) = -S_\rho \times_\rho Z(\rho, \beta) +  \sum_{\gamma} (-1)^{\wt i_{\mathcal C}(\rho, \gamma)} B(\rho, \gamma) \times_\gamma \mathcal C(\gamma, \beta).\] 
This is precisely the boundary relation for the real blowup recalled above.

\item The boundary relation for $\alpha = \overline \rho$ and $\beta \ne \overline \rho, S_\rho$ follows immediately from the flow category relation for $\mathcal C$ (take the zero locus of $\psi$ on both sides of the flow category relation). The case $\alpha \ne S_\rho, \overline \rho$ and $\beta = S_\rho$ is similarly tautological.

\item The boundary of $\mathcal S(\alpha, \overline \rho)$ is tautologically empty. On the right-hand side of the flow category relation, the only factorization $\mathcal S(\alpha, \beta) \times_\beta \mathcal S(\beta, \overline \rho)$ which is not tautologically empty is $\beta = S_\rho$. This means that the flow category relation in this case is 
\[0 = \mathcal C(\alpha, \rho) \times_\rho S_\rho \times_{S_\rho} S_\rho = \pm\mathcal C(\alpha, \rho) \times_\rho S_\rho,\] 
where the right-hand side is interpreted as a chain in $\alpha \times \overline \rho$. Because the map $\mathcal C(\alpha, \rho) \times_\rho S_\rho \to \overline \rho$ factors through the circle bundle projection 
\[\mathcal C(\alpha, \rho) \times_\rho S_\rho \to \mathcal C(\alpha, \rho),\] 
it follows that every stratum of this chain has small rank. So it is a degenerate chain, hence zero on the chain level, and our relation holds.
\end{itemize}

The final handful of cases follow by nearly identical arguments, occasionally recalling that $B(\rho, \beta)$ and $\mathcal C(\rho, \beta)$ define the same geometric chains in $\rho \times \beta$, as they are collapse-equivalent. 
\end{proof}

As mentioned above, there is a copy of the algebraic mapping cone 
\[\mathsf{Cone} = \mathsf{Cone}\big(C_*(SO_3)[-1] \to C_*(S^2)[-2]\big)\] 
as a summand of the flow complex $CM_*(\mathcal S_\rho)$ (though not a subcomplex!). The flow category $\mathcal S_\rho$ is a flow-categorical enrichment of this mapping cone construction. We will build a comparison map between these two flow categories by enriching the comparison map between $C_*(S^2)$ and $\mathsf{Cone}$. 

There is a degree-zero chain map $\Sigma: C_*(S^2) \to \mathsf{Cone}$, given by taking a chain $\phi: P \to S^2$ to the pullback $\phi \times_{S^2} SO(3)$ in $SO(3)$. If $P$ is a point, it is sent to a circle; if $\phi$ is the fundamental class of $S^2$ it is sent to the fundamental class of $SO(3)$. This is a chain map because $\pi(\phi \times_{S^2} SO(3)) \in C_*(S^2)$ is a degenerate chain (hence zero), so the mapping-cone portion of the differential is zero on $\phi \times_{S^2} SO(3)$. 

It is easy to see by hand that the map above induces an isomorphism on homology $H_*(S^2) \to H_*(\mathsf{Cone})$. 

Next, we use this observation to define a bimodule $\Sigma_\rho: \mathcal C \to \mathcal S_\rho$ which induces an isomorphism on the homology of their flow complexes.

\begin{construction}\label{constr:susp-map}
Define the $SO(3)$-equivariant bimodule $\Sigma_\rho$ as follows. Let $\Sigma_\rho(\alpha, \alpha) = \alpha$ for all $\alpha \ne \rho$, while $\Sigma_\rho(\rho, S_\rho) = S_\rho$ and $\Sigma_\rho(\alpha, \beta) = 0$ in all other cases. That is, 
\[\Sigma_\rho = \begin{pmatrix}[c|cc] & \alpha \ne \rho & \rho \\
\hline 
\beta \ne \rho & \delta_{\alpha = \beta} \alpha &  0 \\ 
S_\rho & 0 & S_\rho \\
\overline \rho & 0 & 0\end{pmatrix}\]
with 
\[\wt i_{\Sigma}(\alpha, \beta) = \wt i_{\mathcal C}(\alpha, \beta) - s(\beta),\] 
where $s(\beta) = 0$ if $\beta \not\in \{S_\rho, \overline \rho\}$ while $s(S_\rho) = -1$ and $s(\overline \rho) = -2$.

Here $\delta_{\alpha = \beta} \alpha$ means this term is zero unless $\alpha = \beta$, in which case it is given by $\alpha$. It is a straightforward definition-push to verify that $\mathcal S$ is a bimodule, and the nontrivial parts of the argument are similar in spirit to parts of Theorem \ref{thm:B-is-flowcat}: the key point is that $B(\rho, \beta)$ and $\mathcal C(\rho, \beta)$ define the same geometric chain in $\rho \times \beta$.
\end{construction}

This flow category and bimodule is key in our proof of invariance of equivariant instanton homology. We briefly review the definition. 

Fix a PID $R$. Associated to any $C_*(SO(3);R)$-module $M$ are the (standard, uncompleted) equivariant homology complexes $C^+(M), C^-(M), C^\infty(M)$ whose homology fits into an exact triangle; each of these are modules over $H^{-*}(BSO(3);R)$, and under certain conditions on $R$ the connecting maps in the exact triangle are known to be module homomorphisms. 

Given a $\mathbb Z/8\mathbb Z$-graded $C_*(SO(3);R)$-module $\widetilde C$ with a periodic filtration on its $\mathbb Z$-graded cover $\bar C$, following \cite[Section 5.4]{Helle} we define the completed equivariant homology groups to be \[\hat H^\bullet_{SO(3)}(\widetilde C) = H\left(\lim_{q \to -\infty} \mathop{\operatorname{colim}_{p \to \infty}} C^\bullet(F_p \bar C/F_q \bar C)\right).\] We also use the notation $\hat H^\bullet$ to include the case $H(\widetilde C)$, which needs no completion procedure. The relevant properties are enumerated as \cite[Theorem A.25]{M}, which we now review. 

Because the natural filtration on these equivariant complexes is now complete, exhaustive, and Hausdorff, the equivariant homology groups have conditionally convergent spectral sequences $H^\bullet_{SO(3)}(\text{gr } \widetilde C) \Rightarrow \hat H^\bullet_{SO(3)}(\widetilde C)$ which may be used to detect isomorphisms on equivariant homology. Because the filtration is only used to define the completion, one may in fact replace $F_p$ by any commensurate periodic filtration.

If $f: \widetilde C \to \widetilde C'$ respects filtrations in the sense that for some constant $L$ we have $f(F_p \widetilde C) \subset F_{p+L} \widetilde C'$ for all $p$, then $f$ induces a map on these completed complexes; similarly, filtered chain homotopies define chain homotopies on the completed complexes.

\begin{prop}\label{prop:suspension}
For the $SO(3)$-equivariant bimodule $\Sigma_\rho: \mathcal C \to \mathcal S_\rho(\mathcal C)$ defined above, the induced chain map 
\[CM_*(\Sigma_\rho): CM_*(\mathcal C) \to CM_*(\mathcal S_\rho)\] 
is an isomorphism on homology and on all completed equivariant homology groups $\hat H_{SO(3)}^\bullet$.
\end{prop}

\begin{proof}
The induced map on chains sends each $C_*(\alpha)$ identically onto $C_*(\alpha)$ whenever $\alpha \ne \rho$, and on the remaining factor is the map $\Sigma: C_*(\rho) \to \mathsf{Cone}$ discussed before the construction above. 

To see that $CM_*(\Sigma_\rho)$ induces an isomorphism on homology, we will use a spectral sequence argument using `periodic filtrations'; see Remark \ref{rmk:periodic-filtration} for more details on these and their associated spectral sequences. We give $CM_*(\mathcal C)$ and $CM_*(\mathcal S)$ periodic filtrations, with filtration on the unrolled complex given by
\[F_k \widetilde{CM}_*(\mathcal C) = \bigoplus_{\wt i_{\mathcal C}(\alpha) \leq k} C_*(\alpha)[\wt i_{\mathcal C}(\alpha)] \subset \widetilde{CM}_*(\mathcal C),\]
\[F_k \widetilde{CM}_*(\mathcal S) = \bigoplus_{\substack{\alpha \in \mathsf{Ob}(\mathcal S) \\ \wt i_{\mathcal C}(\alpha) \le k}} C_*(\alpha)[\wt i_{\mathcal S}(\alpha)] \subset \widetilde{CM}_*(\mathcal S),\]
where in the second expression we filter by the $\mathcal C$-index of an object of $\mathcal S$; that is, we filter by the usual $\wt i_{\mathcal C}(\alpha)$ for $\alpha \not\in \{S_\rho, \overline \rho\}$ and otherwise by
\[\widetilde i_{\mathcal C}(S_\rho) = \widetilde i_{\mathcal C}(\overline \rho) = \wt i_{\mathcal C}(\rho).\]

The associated graded complex in the first case is $\bigoplus_{\alpha \in \mathsf{Ob}(\mathcal C)} C_*(\alpha)[\wt i_{\mathcal C}(\alpha)]$ with differential the direct sum of the geometric chain complex differentials. The associated graded complex in the second case is 
\[\bigoplus_{\alpha \ne \rho} C_*(\alpha)[\wt i_{\mathcal C}(\alpha)] \oplus \mathsf{Cone}[\wt i_{\mathcal C}(\rho)].\]
Because the induced map on associated graded complexes is the identity on $C_*(\alpha)$ for $\alpha \ne \rho$, and on $C_*(\rho)$ is given by the map $\Sigma$ described above, we see that $\Sigma_\rho$ induces an isomorphism on the homology of the associated graded complexes. Hence (because the spectral sequence associated to this periodic filtration is absolutely convergent) the map $CM_*(\Sigma_\rho)$ induces an isomorphism on $HM_*$.  

Because the filtration above is commensurate to the index filtration, and our map induces an isomorphism on the associated graded of this filtration, this implies it induces an isomorphism on $H^\bullet_{SO(3)}(\text{gr} \widetilde{CM})$ and hence on $H^\bullet_{SO(3)}(\widetilde{CM})$.
\end{proof}

Therefore, for the purposes of equivariant homology groups we may freely pass between $\mathcal C$ and $\mathcal S_\rho$. The corresponding flow complexes are not, however, equivariantly homotopy equivalent. This will be discussed further in Section 7 below, where it will be used to provide a wall-crossing formula of sorts for Floer's non-equivariant instanton homology groups.

Any bimodule $(\mathcal S_\rho\mathcal C) \to \mathcal C'$ induces a bimodule $\mathcal C \to \mathcal C'$ by precomposition with $\Sigma_\rho$. This leads us to ask: when does an $SO(3)$-equivariant bimodule $\mathcal W: \mathcal C \to \mathcal C'$ enjoy a lift to an $SO(3)$-equivariant bimodule $\widehat{\mathcal W}: \mathcal S_\rho \mathcal C \to \mathcal C'$ for which $\mathcal W = \widehat{\mathcal W} \circ \Sigma_\rho$? 

This is not always the case (for instance, this is not possible for the identity bimodule). It turns out that this is true precisely when a certain index condition holds.

\begin{construction}\label{constr:Wplus-bimod}
Suppose $\mathcal C$ and $\mathcal C'$ are flow categories satisfying the assumptions at the beginning of this section, and suppose we have chosen an orbit $\rho$ of type $S^2$ in $\mathsf{Ob}(\mathcal C)$, an $SO(3)$-equivariant complex line bundle $\mathcal H_\rho \to \rho$, as well a compatible family of equivariant transverse sections $\psi_{\rho \alpha}: \mathcal C(\rho, \alpha) \to \mathcal H_\rho$ lifting the endpoint map to $\rho$. 

Suppose $\mathcal W: \mathcal C \to \mathcal C'$ is an $SO(3)$-equivariant bimodule, satisfying the following condition. For each $\alpha' \in \mathsf{Ob}(\mathcal C')$ there is a stratified-smooth subspace $\mathcal W^{\text{red}}(\rho, \alpha')$ of those points with stabilizer conjugate to $SO(2)$. 

We demand that at each $\Lambda \in \mathcal W^{\text{red}}(\rho, \alpha')$, the normal space to this locus is positive-dimensional. It follows that the normal space has dimension at least two, as this normal space is an $SO(2)$-space with weight one $SO(2)$-action.

For each $\alpha' \in \mathsf{Ob}(\mathcal C')$, choose an $SO(3)$-equivariant transverse section 
\[\psi^W_{\rho \alpha'}:  \mathcal W(\rho, \alpha') \to \mathcal H_{\rho}.\] 
On the boundary face $\mathcal C(\rho, \alpha) \times_{\alpha} \mathcal W(\alpha, \alpha')$ we demand this is given by the section $\psi_{\rho \alpha'} \times 0$ chosen above, while on the boundary face $\mathcal W(\rho, \beta') \times_{\beta'} \mathcal C'(\beta', \alpha')$ we demand this is given by $\psi^W_{\rho' \beta} \times 0$.

If there were no reducibles in these moduli spaces, these would be straightforward to construct by induction; at a reducible of type $S^2$, the construction of such sections is possible by our assumption that each reducible has non-trivial normal space. This corresponds to the fact that there exists an $SO(2)$-equivariant map $f: \Bbb C^n \to \Bbb C$ which is transverse to zero if and only if $n \ge 1$. 

We write $B(\mathcal W; \rho, \alpha')$ and $Z(\mathcal W; \rho, \alpha')$ for the blowup and zero locus, respectively, of $\psi^W_{\rho \alpha'}$ over $\mathcal W(\rho, \alpha')$.

Given this, one may define a bimodule $\widehat{\mathcal W}: \mathcal S_{\rho} \mathcal C \to \mathcal C'$ by the formula 
\[\widehat{\mathcal W}(\alpha, \alpha') =  \begin{pmatrix}[c|ccc] & \alpha \ne S_{\rho}, \overline \rho & \alpha = S_{\rho} & \alpha = \overline{\rho} \\ 
\hline 
%& & & \\
\alpha' & \mathcal W(\alpha, \alpha') & B(\mathcal W; \rho, \alpha') & Z(\mathcal W;\rho, \alpha') \end{pmatrix}\]

Here 
\[\wt i_{\widehat{\mathcal W}}(\alpha, \alpha') = \wt i_{\mathcal W}(\alpha, \alpha') - s(\alpha),\] 
where the shift $s$ is defined by the same formula as above.
\end{construction}

This defines a bimodule lifting $\mathcal W$. 

Again, here is a diagram which may be helpful when verifying the bimodule relations. The left column depicts the flow category $\mathcal S_\rho \mathcal C$, while the right column depicts the flow category $\mathcal C'$, and the arrows running left to right depict the bimodule $\mathcal W$. As before, the bimodule relation for $\partial \mathcal W(\alpha, \alpha')$ amounts to saying that the boundary of any arrow running left to right is the signed sum over all two-step paths from $\alpha$ to $\alpha'$.

\[\begin{tikzcd}
	{\alpha / \beta} \\
	\\
	{S_\rho} &&& {} & {\alpha' / \beta'} \\
	\\
	{\overline{\rho}} \\
	&&& 
	\arrow["{\mathcal W(\alpha, \alpha')}", from=1-1, to=3-5]
	\arrow["{\mathcal C(\alpha, \rho) \times_\rho S_\rho}"{pos=0.7}, from=1-1, to=3-1]
	\arrow["{S_\rho}"', from=3-1, to=5-1]
	\arrow["{Z(\rho, \beta)}", bend left=30, from=5-1, to=1-1]
	\arrow["{B(\mathcal W; \rho, \alpha')}", from=3-1, to=3-5]
	\arrow["{Z(\mathcal W; \rho, \alpha')}"', from=5-1, to=3-5]
	\arrow["{-B(\rho, \beta)}", {pos=0.6}, bend left=50, from=3-1, to=1-1]
	\arrow["{\mathcal C(\alpha, \beta)}", loop above, from=1-1, to=1-1]
	\arrow["{\mathcal C'(\alpha', \beta')}", loop above, from=3-5, to=3-5]
\end{tikzcd}\]
\begin{theorem}\label{thm:Wplus-bimod}
The formulas above define a bimodule $\widehat{\mathcal W}: \mathcal S_{\rho} \mathcal C \to \mathcal I$. This bimodule satisfies $\widehat{\mathcal W} \circ \Sigma_{\rho} = \mathcal W$.
\end{theorem}

\begin{proof}
That $\widehat{\mathcal W}$ is a bimodule follows by an analysis much like Theorem \ref{thm:B-is-flowcat}, repeatedly using three facts: that $S_{\rho} \times_{S_{\rho}} B(\mathcal W; \rho, \alpha')$ defines the same geometric chain as $\mathcal W(\rho, \alpha')$, the boundary formula for $B(\mathcal W; \rho, \alpha')$ given by Lemma \ref{blowup-equiv}, and the boundary formula for $\mathcal W(\alpha, \alpha')$ given by the definition of bimodule.

The statement about the composite bimodule follows nearly immediately from the relevant definitions. The only term in this composite where this isn't tautological is 
\[(\widehat{\mathcal W} \circ \Sigma_{\rho})(\rho, \alpha') = S_{\rho} \times_{S_{\rho}} B(\mathcal W; \rho, \alpha') = B(\mathcal W; \rho, \alpha') \textup{ as a chain in } \rho \times \alpha',\] 
which we know coincides with $\mathcal W(\rho, \alpha')$ at the level of geometric chains..
\end{proof}
\subsection{The bimodule of an obstructed cobordism}\label{cob-bimods}
%!TEX root = equivariant-functoriality.tex

We will now show that the modified moduli spaces of Section 5.2 give rise to a bimodule \emph{to} the suspended flow category, opposite to the situation of Theorem \ref{thm:Wplus-bimod}.

Let $(W, c): (Y, w, \pi) \to (Y', w', \pi')$ be a nearly unobstructed cobordism whose obstructed reducible we denote $\Lambda: \rho \to \rho'$. In Section 5.2, we constructed `modified moduli spaces' $N^+_z(W, c; \alpha, \alpha')$. Our first goal in this section will be to use these to construct a bimodule 
\[\mathcal W_-: \mathcal I(Y, w, \pi) \to \mathcal S_\rho \mathcal I(Y', w', \pi').\] 
The minus sign indicates that $\wt i(\Lambda) = -2$, so that index \emph{decreases} across this obstructed reducible.

To make the notation more parseable, let us make some abbreviations. We write $M(\alpha, \beta) = \breve M^+_\zeta(Y, w; \alpha, \beta)$ for the unique path $\zeta$ with $0 \le i(\zeta) \le 7$, while $M'(\alpha', \beta')$ is the same for $(Y', w')$. 

There are canonical equivariant diffeomorphisms $\rho \xleftarrow{e_-^\Lambda} \mathcal O_\Lambda \xrightarrow{e_+^\Lambda} \rho'$. We use these diffeomorphisms to pull back $\mathcal H_\Lambda$ and $S_\Lambda$ to bundles over $\rho'$, and we refer to the map $e_+^\rho: \rho \to \rho'$ for the map $e_+^\rho = e_+^\Lambda \circ (e_-^\Lambda)^{-1}$ obtained by inverting the first arrow in the above zig-zag. Every appearance of $M \times_\rho \mathcal O_\Lambda$ we replace with $M$ in a canonical way, and similarly we replace $\mathcal O_\Lambda \times_{\rho'} M'$ with $M'$.

As in section 5.2 and the previous section, we choose a compatible family of transverse equivariant sections 
\[\psi_{\alpha \rho}: M(\alpha, \rho) \to \mathcal H_{\rho'}, \quad \psi'_{\rho'\alpha'}: M'(\rho', \alpha') \to \mathcal H_{\rho'}.\] 
We write $B(\alpha, \rho)$ and $B'(\rho', \beta')$ for the blowups along these sections, and $Z(\alpha, \rho)$ and $Z'(\rho', \alpha')$ for the zero loci of these sections; the family of sections described here plays the role of $\psi'_{\rho'\alpha'}$ in Section 5.3. Both $B(\alpha, \rho)$ and $B'(\rho', \beta')$ come with canonical maps to $S_{\rho'}$. 

In Section 5.2, these were the sections $\psi^{\text{priv}}_{\zeta\Lambda}$, where $\zeta$ is the unique path with $0 \le i(\zeta) \le 7$.

When defining the suspended flow category, we will use $\psi'_{\rho'\alpha'} = -\psi^{\text{priv}}_{\zeta\Lambda}$. This minus sign makes no difference at the level of zero sets, but does change the map $B(\psi) \to S_{\rho'}$; it will be clear momentarily why negation is important.

Finally, we abbreviate the modified moduli spaces $N^+_{\zeta}(W, c; \alpha, \alpha')$ to $N(\alpha, \alpha').$ In indices below, we suppress the symbols $w$ in $(Y,w)$ and $c$ in $(W,c)$.

\begin{construction}\label{constr:W-minus}
We define the bimodule $\mathcal W_-: \mathcal I \to \mathcal S_{\rho'} \mathcal I'$ by the following formula.

\[\mathcal W_-(\alpha, \alpha') = \begin{pmatrix}[c|cc] & \alpha \ne \rho & \alpha = \rho \\ 
\hline 
& & \\
\alpha' \ne \rho' & N(\alpha, \alpha') & N(\rho, \alpha') \\
& &  \\
\alpha' = S_{\rho'} & N(\alpha, \rho') \times_{\rho'} S_{\rho'} - (-1)^{\wt i(Y; \alpha, \rho)} B(\alpha, \rho) & 0 \\ 
& & \\
\alpha' =\overline{\rho}' & 0 & \rho
\end{pmatrix}\]

Here we have 
\[\wt i_{\mathcal W_-}(\alpha, \alpha') = \wt i(W;\alpha, \alpha') + s(\alpha'),\] 
where $s(\overline{\rho}') = 2$ and $s(S_{\rho'}) = 1$, while $s(\alpha') = 0$ otherwise. The bottom-right entry is the correspondence $\rho \xleftarrow{1_\rho} \rho \xrightarrow{e_+^\rho} \rho'$ induced by the endpoint maps of $\mathcal O_\Lambda$ which was described above. 
\end{construction}

\begin{remark}
Entry $(2, 1)$ in this bimodule may be unexpected. What we will see in the calculation below is that this particular form is forced upon us by demanding that $\mathcal W_-(\alpha, \alpha') = N(\alpha, \alpha')$ whenever possible. 

One could place this expression in its proper context by proving a more general statement: if $\mathcal N(\alpha, \alpha')$ is a collection of moduli spaces satisfying the boundary relations obtained in Proposition \ref{obs-bdry-reln}, then it assembles to give a bimodule $\mathcal C \to \mathcal S_{\rho'} \mathcal C'$ by formulas along the lines of the above, and conversely every bimodule $\mathcal C \to \mathcal S_{\rho'} \mathcal C'$ arises in this way. We will not bother formulating or proving such a statement.
\end{remark}

Before the proof, we recall some of the important boundary relations. Proposition \ref{obs-bdry-reln} can be interpreted as giving the following boundary formula, at the level of geometric chains; the passage to geometric chains allows us to cancel out faces with opposite orientation, and equate things which are collapse-equivalent. 

When $\alpha \ne \rho$ and $\alpha' \ne \rho'$, we have the boundary relation
\begin{align*}(-1)^{\dim \alpha} \partial N(\alpha, \alpha') =& \sum_\beta M(\alpha, \beta) \times_\beta N(\beta, \alpha') \\
&+ \sum_{\beta'} (-1)^{\wt i(W; \beta, \alpha')+1} N(\alpha, \beta') \times_{\beta'} M'(\beta', \alpha') \\ &+B(\alpha, \rho)  \times_{S_{\rho'}} B'(\rho', \alpha').
\end{align*}
The final term deserves some comment. It appears in Proposition \ref{obs-bdry-reln} as a negatively-oriented copy of the zero set of the family $t \psi_{\alpha \rho} - (1-t)\psi'_{\rho' \alpha'}$ on the pullback of $\mathcal H_\Lambda$ over 
\[M(\alpha, \rho) \times [0,1] \times_{\rho'} M'(\rho', \alpha').\]
Applying the orientation-reversing isomorphism $t \mapsto 1-t$ (which flips the negative orientation) on the unit interval and negating the section (which does not change the orientation on the zero set), this section is sent to $(1-t)\psi_{\alpha \rho} - t\psi'_{\rho' \alpha'}$. We show in Lemma \ref{weird-zero-set} in the appendix that the zero set of this section is identical as a geometric chain to the fiber product $B(\alpha, \rho)  \times_{S_{\rho'}} B'(\rho', \alpha'),$ 
giving the statement written above. This is where it is important to use $\psi'_{\rho' \alpha'} = -\psi^{\text{priv}}_{\Lambda \zeta'}$; with the opposite sign, Lemma \ref{weird-zero-set} would not be applicable.

If $\alpha \ne \rho$ and $\alpha' = \rho$, the last term is instead $Z(\rho, \alpha)$. If $\alpha = \rho$ and $\alpha' \ne \rho'$, the final term is instead $-Z'(\rho', \alpha')$.

In the diagram below, the left column corresponds to the flow category $\mathcal C$, while the right column corresponds to the flow category $\mathcal S_{\rho'} \mathcal C'$. The bimodule relations for $\mathcal W_-(\alpha, \alpha')$ assert that the boundary of the arrow running from (say) $\alpha$ to $\alpha'$ is the signed sum of the fiber product over all ways to go from $\alpha$ to $\alpha'$ in two steps, one along $\mathcal W_-$ and one inside either the initial or ending flow category.  
\vspace{-0.6cm}
\[\begin{tikzcd}
	&&&& {\alpha' / \beta'} \\
	{\alpha / \beta} \\
	&&&& {S_{\rho'}} \\
	\rho \\
	&&&& {\overline \rho'}
	\arrow["{S_{\rho'}}"', from=3-5, to=5-5]
	\arrow["{\breve M'(\alpha', \rho') \times_{\rho'} S_{\rho'}}"'{pos=0.4}, from=1-5, to=3-5]
	\arrow["{Z'(\rho', \beta')}"', bend right=20, from=5-5, to=1-5]
	\arrow["{-B'(\rho', \beta')}"', bend right=40, from=3-5, to=1-5]
	\arrow["{N(\alpha, \alpha')}", from=2-1, to=1-5]
	\arrow["{\mathcal O_-}"', from=4-1, to=5-5]
	\arrow["{\breve M(\alpha, \rho)}"'{pos=0.7}, bend right=10, from=2-1, to=4-1]
	\arrow["{\breve M(\rho, \beta)}"'{pos=0.8}, bend right=10, from=4-1, to=2-1]
	\arrow["{N(\alpha, \rho') \times_{\rho'} S_{\rho'}}"{pos=0.45}, from=2-1, to=3-5]
	\arrow["{-(-1)^{\widetilde \zeta} B(\beta, \rho)}"'{pos=0.9}, from=2-1, to=3-5]
	\arrow["{N(\rho, \alpha')}"'{pos=0.1}, bend left=25, from=4-1, to=1-5]
	\arrow["{\breve M(\alpha, \beta)}", loop above, from=2-1, to=2-1]
	\arrow["{\breve M'(\alpha', \beta')}", loop above, from=1-5, to=1-5]
\end{tikzcd}\]

\begin{theorem}
For $(W, c)$ a cobordism as in the beginning of this section, the moduli spaces $\mathcal W_-$ defined above give a bimodule $\mathcal W_-: \mathcal I \to \mathcal S_{\rho'} \mathcal I'.$
\end{theorem}

\begin{proof}
The statement that $\mathcal W_-$ gives a bimodule amounts to six boundary formulas. We will only verify two of these; the rest follow by nearly identical arguments, and are furthermore similar in spirit to the verifications in Theorem \ref{thm:B-is-flowcat}. 

\begin{itemize}
\item First (filling in the values of $\mathcal W_-$), the bimodule relations for $\mathcal W_-(\alpha, \alpha')$ for $\alpha \ne \rho$ and $\alpha' \ne \overline{\rho}', S_{\rho'}$ stipulate that 
\begin{align*}
(-1)^{\dim \alpha}\partial N(\alpha, \alpha') &= \sum_{\beta \in \mathfrak C(Y,w)} M(\alpha, \beta) \times_\beta N(\beta, \alpha') \\   
&+ \sum_{\substack{\beta' \in \mathfrak C(Y',w') \\ \beta' \ne \rho'}} (-1)^{\wt i_{\mathcal W_-}(\alpha, \beta') + 1} N(\alpha, \beta') \times_{\beta'} M'(\beta', \alpha') \\ 
&+ (-1)^{\wt i_{\mathcal W_-}(\alpha, S_{\rho'})+1} N(\alpha, \rho') \times_{\rho'} S_{\rho'} \times_{S_{\rho'}} \big(-B'(\rho', \alpha')\big)\\
&+ (-1)^{\wt i_{\mathcal W_-}(\alpha, S_{\rho'})+1} (-1)^{\wt i(Y; \alpha, \rho)+1} B(\alpha, \rho) \times_{S_{\rho'}} \big(- B'(\rho', \alpha')\big).
\end{align*}

The third and last line arise from $\mathcal W_-(\alpha, S_{\rho'}) \times_{S_{\rho'}} \mathcal S_{\rho'} \mathcal I'(S_{\rho'}, \alpha')$. Our first order of business is to clean up these last few lines. First, $\wt i_{\mathcal W_-}(\rho, \overline{\rho}') = 0$. Next, $\wt i_{\mathcal W_-}(\alpha, S_{\rho'}) = \wt i(W;\alpha, \rho') + 1$ while 
\begin{align*}N(W; \alpha, \rho') \times_{\rho'} S_{\rho'} \times_{S_{\rho'}} \big(-B'(\rho', \alpha')\big) &= -N(W; \alpha, \rho') \times_{\rho'} B'(\rho', \alpha') \\
&= -N(W; \alpha, \rho') \times_{\rho'} M'(\rho', \alpha'),
\end{align*} 
because by Lemma \ref{blowup-equiv} the chains $B'(\rho', \alpha')$ and $M'(\rho', \alpha')$ are identical when considered as geometric chains in $\rho' \times \alpha'$. Lastly, the final line may be simplified with the observation that 
\[\wt i_{\mathcal W_-}(\alpha, S_{\rho'}) + 1 - \wt i(Y;\alpha, \rho) -1 = \wt i(W;\alpha, \rho') + 1 - \wt i(Y;\alpha, \rho) = \wt i(W;\rho, \rho') = -2\]
so that the sign out front is $-1$; pulling the factor of $-1$ from $-B'(\rho', \alpha')$ out front as well, this reads
\begin{align*}(-1)^{\dim \alpha}\partial N(\alpha, \alpha') &= \sum_{\beta \in \mathfrak C(Y,w)} M(\alpha, \beta) \times_\beta N(\beta, \alpha') \\   
&+ \sum_{\beta' \in \mathfrak C(Y',w')} (-1)^{\wt i(W;\alpha, \beta') + 1} N(\alpha, \beta') \times_{\beta'} M'(\beta', \alpha') \\ 
&+ B(\alpha, \rho) \times_{S_{\rho'}} B'(\rho',\alpha').
\end{align*}

This is precisely the relation written above, coming from Proposition \ref{obs-bdry-reln}(iii).

\item Let's verify the bimodule relation for the other interesting case, $\mathcal W_-(\alpha, S_{\rho'})$ where $\alpha \ne \rho$. To briefly describe what happens: $N(\alpha, \rho') \times_{\rho'} S_{\rho'}$ carries an extra boundary term $Z(\alpha, \rho') \times_{\rho'} S_{\rho'}$; the boundary relation from the moduli space $B(\alpha, \rho)$ which appears in $\mathcal W_-(\alpha, S_{\rho'})$ will `cancel it out'. To be more precise, the bimodule relations stipulate that \begin{align*}(-1)^{\dim \alpha}\partial \mathcal W_-(\alpha, S_{\rho'}) &= \sum_{\substack{\beta \in \mathfrak C(Y,w) \\ \alpha \ne \beta \ne \rho}} \mathcal I(\alpha, \beta) \times_\beta \mathcal W_-(\beta, S_{\rho'}) \\
& \sum_{\substack{\alpha' \in \mathfrak C(Y',w') \\ \alpha' \ne \rho'}} (-1)^{\wt i_{\mathcal W_-}(\alpha, \alpha')+1} \mathcal W_-(\alpha, \alpha') \times_{\alpha'} \mathcal S_{\rho'} \mathcal I'(\alpha', S_{\rho'}).\end{align*}

(The last sum does not include a factorization through $\mathcal S_{\rho'} \mathcal I(\overline \rho', S_{\rho'})$ because this moduli space is empty.) The term $\mathcal I(\alpha, \beta)$ is $M(\alpha, \beta)$; the term $\mathcal S_{\rho'}\mathcal I'(\alpha', S_{\rho'})$ is $M'(\alpha', \rho') \times_{\rho'} S_{\rho'}.$ Filling in the values of $\mathcal W_-$, the right-hand-side simplifies to 
\begin{align*}&\sum_{\beta \ne \rho} M(\alpha, \beta) \times_\beta N(\beta, \rho') \times_{\rho'} S_{\rho'} - (-1)^{\wt i(Y; \beta, \rho)} M(\alpha, \beta) \times_\beta B(\beta, \rho)\\ 
&+ \sum_{\alpha'\ne \rho'}(-1)^{\wt i(W; \alpha, \alpha')+1} N(\alpha, \alpha') \times_{\alpha'} M'(\alpha', \rho') \times_{\rho'} S_{\rho'}.
\end{align*}

The left-hand side of the bimodule relation reads 
\[(-1)^{\dim \alpha} \partial\left(N(\alpha, \rho') \times_{\rho'} S_{\rho'} - (-1)^{\wt i(Y; \alpha, \rho)} B(\alpha, \rho)\right).\]

Now use the boundary relation for $(-1)^{\dim \alpha}\partial N(\alpha, \rho')$ stated before this Proposition; most every term in that boundary relation is represented in the sum above, except for the final term corresponding to the boundary face $Z(\alpha, \rho)$. 

Subtracting $(-1)^{\dim \alpha} \partial N(\alpha, \rho') \times_{\rho'} S_{\rho'}$ from both sides, what we're left with trying to show is that 
\[(-1)^{\dim \alpha + \wt i(Y; \alpha, \rho)-1} \partial B(\alpha, \rho) = -Z(\alpha, \rho) \times_\rho S_\rho + \sum_\beta (-1)^{\wt i(Y; \beta, \rho)} M(\alpha, \beta) \times_\beta B(\beta, \rho) \times_\rho S_{\rho'}.\] 
Multiplying through by $(-1)^{\dim \alpha + \wt i(Y; \alpha, \rho)}-1$, what we assert is that 
\[\partial B(\alpha, \rho) = (-1)^{\dim \alpha + \wt i(Y; \alpha, \rho)} Z(\alpha, \rho) \times_\rho S_\rho + \sum_\beta (-1)^{\dim \alpha + \wt i(Y; \alpha, \beta)} M(\alpha, \beta) \times_\beta B(\beta, \rho).\]

Finally, this relation is given in Lemma \ref{blowup-equiv}: the final term is $B(\psi|_{\partial M(\alpha, \rho)})$, while the first term is $ (-1)^{\dim Z(\alpha, \rho)+1} Z(\alpha, \rho) \times_\rho S_\rho$, as $\dim Z(\alpha, \rho) = \wt i(Y; \alpha, \rho) + \dim \alpha - 1.$
\end{itemize}

Carrying out a similar discussion in all the other cases is similar but strictly simpler.
\end{proof}
\subsection{Invariance of instanton homology}\label{sec:invt}
%!TEX root = equivariant-functoriality.tex

Recall from Section \ref{hessian} that each weakly admissible pair $(Y, w)$ and each regular perturbation $\pi$ gives rise to a function $\sigma_\pi: \mathfrak A(Y, w) \to 2\Bbb Z$, called its \emph{signature data}. 

Recall also from Section \ref{inst-cx-unob-cob-map} that each $(Y, w, \pi)$ with basepoint $y \in Y$ gives rise to an $SO(3)$-equivariant flow category $\mathcal I(Y, w, \pi)$, whose flow complex we call $\widetilde C(Y, w, \pi)$, and that for each unobstructed cobordism $(W, c): (Y, w, \pi) \to (Y', w', \pi')$ with regular perturbation $\pi_W$, homology orientation of $W$, and path between the basepoints, we associate a bimodule $\mathcal I(W, c, \pi_W, \gamma)$ and thus an $SO(3)$-equivariant chain map 
\[\widetilde C(W, c, \pi_W, \gamma): \widetilde C(Y, w, \pi) \to \widetilde C(Y', w', \pi'),\] 
independent of $\pi_W$ up to equivariant chain-homotopy. Furthermore, these compose functorially (up to equivariant chain homotopy).

Taking $(W, c) = \Bbb R \times (Y, w)$ with trivial homology orientation and path $\gamma = \Bbb R \times \{y\}$ and applying the index computation Proposition \ref{prop:normal-ind}, this cobordism is unobstructed if and only if we have the inequality on signature data functions $\sigma_\pi \le \sigma_{\pi'}$. In that case, we obtain an $SO(3)$-equivariant continuation map 
\[\varphi_{\pi \to \pi'}: \widetilde C(Y, w, \pi) \to \widetilde C(Y, w, \pi').\]

Our goal for the remainder of this section is to prove the following.

\begin{theorem}\label{thm:main-theorem}
Whenever $\pi_0, \pi_1$ are two regular perturbations on $(Y, w)$ with $\sigma_{\pi_0} \le \sigma_{\pi_1}$, the continuation map 
\[\varphi_{\pi_0 \to \pi_1}: I^\bullet(Y, w, \pi_0) \to I^\bullet(Y, w, \pi_1)\] 
is an isomorphism on all flavors of equivariant homology.
\end{theorem}

\begin{proof}[Proof of Theorem \ref{thm:main-theorem}]
When $\sigma_{\pi_0} = \sigma_{\pi_1}$ this follows from the machinery of Section \ref{inst-cx-unob-cob-map}: there are continuation maps $\varphi_{\pi_0 \to \pi_1}$ and $\varphi_{\pi_0 \to \pi_1}$ which are equivariantly homotopy inverse, and hence each is an equivariant homotopy equivalence. The theorem above is primary interesting when $\sigma_{\pi_0} < \sigma_{\pi_1}$.

We will only need to prove this in a special case. We say that $\pi_0$ and $\pi_1$ have \emph{adjacent signature data} if there exists a reducible $\rho$ with $\sigma_{\pi_1} = \sigma_{\pi_0} + 4\delta_\rho$; that is, $\sigma_{\pi_0}$ is equal to $\sigma_{\pi_1}$ for all reducibles except $\rho$, where $\sigma_{\pi_1}$ is larger by the smallest possible amount.

To see why this is sufficient, recall that Proposition \ref{arbitrary-specshift} shows that any signature data function is realized as the signature data of some regular perturbation on $(Y, w)$. It follows that whenever $\sigma_{\pi_0} < \sigma_{\pi_1}$, we may find a sequence of intermediate perturbations $\pi_0, \pi_{1/n}, \cdots, \pi_{(n-1)/n}, \pi_1$ so that 
\[\sigma_{\pi_0} < \cdots < \sigma_{\pi_1},\] 
where each $(\pi_{i/n}, \pi_{(i+1)/n})$ have adjacent signature data. Because we can write the continuation map $\varphi_{\pi_0 \to \pi_1}$ as a composite of the individual continuation maps $\varphi_{\pi_{i/n} \to \pi_{(i+1)/n}}$, if the continuation map is an isomorphism for perturbations with adjacent signature data, it is an isomorphism whenever $\sigma_{\pi_0} \le \sigma_{\pi_1}$.

So suppose that $\pi_0$ and $\pi_1$ have adjacent signature data. Write $W_+$ for the cylinder $\Bbb R \times (Y, w)$, considered as a cobordism $(Y, w, \pi_0) \to (Y, w, \pi_1)$ equipped with a regular perturbation $\pi_+$ interpolating between these, and $W_-$ for the same cylinder considered as a cobordism in the opposite direction $(Y, w, \pi_1) \to (Y, w, \pi_0)$ equipped with an obstructed-regular perturbation $\pi_-$ interpolating between these. 

The cobordism $W_-$ is a nearly unobstructed cobordism: its reducibles correspond to the constant trajectories on the cylinder, so their normal index is given by $\frac 12(\sigma_{\pi_0}(\rho) - \sigma_{\pi_1}(\rho)$), and the one obstruced reducible $\Lambda: \rho \to \rho$ has normal index $-2$. Thus, by the machinery of the previous section we obtain a bimodule 
\[\mathcal W_-: \mathcal I(Y, w, \pi_1) \to \mathcal S_\rho \mathcal I(Y, w, \pi_0).\] We write $\varphi_-$ for the induced map on flow complexes \[\varphi_-: \widetilde C(Y, w, \pi_1) \to S_\rho \widetilde C(Y, w, \pi_0).\]

On the other hand, the cobordism $W_+$ satisfies the hypotheses of Construction \ref{constr:Wplus-bimod}: the only reducible running $\rho \to \rho$ has normal index $2$. We thus obtain a bimodule $\mathcal W_+$ running in the opposite direction, and we write $\varphi_+$ for the induced map on flow complexes.

Abbreviating $\mathcal I_0 = \mathcal I(Y, w, \pi_0)$ while $\mathcal I_1 = \mathcal I(Y, w, \pi_1)$, we obtain the following diagram of flow categories and bimodules:

\[\begin{tikzcd}
	{\mathcal I_0} && {\mathcal I_1} \\
	\\
	{\mathcal S_\rho \mathcal I_0}
	\arrow["{\mathcal W_{\pi_0 \to \pi_1}}", from=1-1, to=1-3]
	\arrow["{\Sigma_\rho}"', from=1-1, to=3-1]
	\arrow["{\mathcal W_+}"', shift right=1, from=3-1, to=1-3]
	\arrow["{\mathcal W_-}"', shift right=1, from=1-3, to=3-1]
\end{tikzcd}\]

Here $\mathcal W_{\pi_0 \to \pi_1}$ is the bimodule obtained from the usual (unmodified) moduli spaces of instantons on the unobstructed cobordism $\Bbb R \times Y: (Y, w, \pi_0) \to (Y, w, \pi_1)$; its induced map on the complex $\widetilde C$ is the continuation map we named $\varphi_{\pi_0 \to \pi_1}$ above. We established in Theorem \ref{thm:Wplus-bimod} that we have an equality of composites 
\[\mathcal W_+ \circ \Sigma_\rho = \mathcal W_{\pi_0 \to \pi_1},\] 
so that on flow complexes we have $\varphi_+ \circ \Sigma_\rho = \varphi_{\pi_0 \to \pi_1}$. We established in Proposition \ref{prop:suspension} that $\Sigma_\rho$ induces an isomorphism on equivariant homology, so we are reduced to showing that $\varphi_+$ is an isomorphism on equivariant homology. Indeed, we will prove that $\varphi_+$ is an equivariant chain homotopy equivalence.

To see this, we assert that there exist $SO(3)$-equivariant bimodule homotopies 
\[\mathcal H_1: \mathbf 1_{\mathcal I_1} \implies \mathcal W_+ \circ \mathcal W_-, \quad \mathcal H_0: \mathbf 1_{\mathcal S_\rho \mathcal I_0} \implies \mathcal W_- \circ \mathcal W_+.\]
It follows that the induced maps $\varphi_+$ and $\varphi_-$ are equivariantly homotopy inverse, as desired.

The construction of these bimodule homotopies is carried out in Section 5.5.2 below, based on the analytic setup of Section 5.5.1. The fact that they are indeed bimodule homotopies are verified in Propositions \ref{htpy1} and \ref{htpy2} in Section 5.5.2.
\end{proof}

In particular, given any two regular perturbations $\pi, \pi'$ with $\sigma_\pi \le \sigma_{\pi'}$, the argument above implies that the natural map \[\varphi_{\pi \to \pi'}: \widetilde C(Y,w,\pi) \to \widetilde C(Y,w,\pi')\] induces an isomorphism on homology and all equivariant homology groups. We now use this to derive the invariance of equivariant instanton homology.

\begin{cor}\label{welldef}
The equivariant instanton homology groups, $I^\bullet(Y, w, \pi;R)$ together with their module structures and the exact triangle relating them, are independent of regular perturbation $\pi$, 
\end{cor}

\begin{proof}
Because any signature-data function is realized as the signature data of some regular perturbation, and the maximum of two signature data functions is again a signature data function, given any two regular perturbations $\pi_0, \pi_1$, there is a regular perturbation $\pi_2$ with $\sigma_{\pi_2} = \text{max}(\sigma_{\pi_0}, \sigma_{\pi_1})$. Because the continuation maps $\varphi_{02} = \varphi_{\pi_0 \to \pi_2}$ and $\varphi_{12} = \varphi_{\pi_1 \to \pi_2}$ induce isomorphisms on all equivariant homology groups, this provides an isomorphism between 
\[\varphi_{12}^{-1} \varphi_{02}: I^\bullet(Y, w, \pi_0) \to I^\bullet(Y, w, \pi_1).\]

This isomorphism is canonical (independent of the choice of $\pi_2$). If one picks a $\pi_2'$ with $\sigma_{\pi_2} \le \sigma_{\pi_2'}$, then by functoriality of induced maps of cobordisms, we have 
\[\varphi_{22'} \varphi_{02} = \varphi_{02'}, \quad \quad \varphi_{22'} \varphi_{12} = \varphi_{12'},\]
so that 
\[\varphi_{12'}^{-1} \varphi_{02'} = \varphi_{12}^{-1} \varphi_{22'}^{-1} \varphi_{22'} \varphi_{02} = \varphi_{12}^{-1} \varphi_{02}.\]

For any two $\pi_2$ and $\pi_2'$ exceeding $\pi_0$ and $\pi_1$, there is another $\pi_2''$ exceeding all of these. Because $\varphi_{12}^{-1} \varphi_{02} = \varphi_{12''}^{-1} \varphi_{02''} = \varphi_{12'}^{-1} \varphi_{02'}$, we see that the isomorphism between $I^\bullet(Y, w, \pi_0)$ and $I^\bullet(Y, w, \pi_1)$ is independent of any choices in its construction.
\end{proof}

We use the previous corollary to give a definition of the instanton homology groups which is transparently independent of $\pi$. 

\begin{definition}\label{def:I-no-pert}
Let $(Y, w)$ be a weakly admissible pair. We set 
\[I^\bullet(Y, w; R) = \bigoplus_{\pi \textup{ regular}} I^\bullet(Y,w,\pi; R)\bigg/\varphi_{\pi \to \pi'}.\] 
Here we quotient by the additive relation generated by setting 
\[(\varphi_{\pi \to \pi'}(x), \pi') \sim (x, \pi)\] 
for any $x \in I^\bullet(Y, w, \pi; R)$ and any $\pi'$ with $\sigma_\pi \le \sigma_{\pi'}$.
\end{definition}

Notice the continuation map $\varphi_{\pi \to \pi}$ is the identity, all continuation maps are isomorphisms, we have $\varphi_{\pi' \to \pi''} \varphi_{\pi \to \pi'} = \varphi_{\pi \to \pi''}$, and for any $\pi, \pi'$ there is a perturbation $\pi''$ so that $\sigma_{\pi''} \le \sigma_\pi, \sigma_{\pi'}$. It is a straightforward exercise in group theory to see that the natural map $I^\bullet(Y, w, \pi) \to I^\bullet(Y, w)$ is an isomorphism for all regular perturbations $\pi$. This is the essential content of the fact proven in the Corollary above: not only are all $I^\bullet(Y, w, \pi)$ isomorphic, they are \emph{canonically} so.

\subsubsection{Modified moduli spaces for homotopies}
%!TEX root = invariance.tex

For the weakly admissible pair $(Y,w)$, fix two regular perturbations $\pi_0$ and $\pi_1$ with adjacent signature data. A typical element of $\fC_{\pi_0}(Y,w)$ is denoted by $\alpha$ and a typical element of $\fC_{\pi_1}(Y,w)$ is denoted by $\alpha'$. Any reducible element $\sigma$ of $\fC_{\pi_0}(Y,w)$ corresponds to a unique reducible element $\sigma'$ of $\fC_{\pi_1}(Y,w)$. There is a unique such reducible $\rho\in \fC_{\pi_0}(Y,w)$ such that 
\[
  \sigma_{\pi_1}(\rho')-\sigma_{\pi_0}(\rho)=4,
\]
and otherwise $\sigma_{\pi_0}$ and $\sigma_{\pi_1}$ are equal.

Let $(W_+,c_+,\pi_+)$ be the product cobordism from $(Y,w)$ to $(Y,w)$ together with a perturbation $\pi_+$ from $\pi_0$ to $\pi_1$; the subscript $+$ is meant to indicate that the signature data increases across the cylinder. Using our assumption on $\sigma_{\pi_0}$ and $\sigma_{\pi_1}$, we may assume that $\pi_+$ is regular and for any reducible $\sigma\in \fC_{\pi_0}(Y,w)$, there is a unique (unobstructed) reducible instanton on $(W_+,c_+,\pi_+)$ which is asymptotic to $\sigma$ and $\sigma'$ on the two ends. In particular, there is a unique reducible instanton $\Lambda_+$ on $W_+$ from $\rho$ to $\rho'$ with $\wt i(\Lambda_+)= +2$. The remaining reducible instantons on $(W_+,c_+)$ have $\wt i(\Lambda) = 0$ and are isolated.

Similarly, let $(W_-,c_-,\pi_-)$ be the product cobordism from $(Y,w)$ to $(Y,w)$ together with a perturbation $\pi_-$ from $\pi_1$ to $\pi_0$, the subscript indicating that signature data decreases along the cylinder. In this case, we may only assume that $\pi_-$ is obstructed-regular as in Theorem \ref{thm:unobs-transv}(ii). To be more precise, we may arrange for $\pi_-$ such that for any reducible $\sigma\in \fC_{\pi_0}(Y,w)$, there is a unique reducible instanton on $(W_-,c_-,\pi_-)$ which is asymptotic to $\sigma$ and $\sigma'$ on the two ends. This reducible instanton is unobstructed, with $\wt i(\Lambda) = 0$, unless $\sigma=\rho$. We will write $\Lambda_-$ for the obstructed reducible instanton that also satisfies $\wt i(\Lambda_-)= -2$, and $\lambda_-: \rho' \to \rho$ for the corresponding path of connections.

Composition of $(W_+,c_+,\pi_+)$ and $(W_-,c_-,\pi_-)$ produces the product cobordism $(W_{-}^+,c^+_{-})$ together with a family of perturbations $\pi_-^+(t)$, where $t$ denotes the length of the neck, taking values in $[T_0,\infty)$. In particular, we may regard $(W_{-}^+,c^+_{-},\pi^+_{-})$ as a cobordism from $(Y,w,\pi_0)$ to $(Y,w,\pi_0) $. We extend the parametrizing set to $I=[T_0,\infty]$, by including the broken perturbation given by $\pi_+$ and $\pi_-$. For any $\alpha,\beta \in \fC_{\pi_0}(Y,w)$ and a path $z$ from $\alpha$ to $\beta$ with $i(z)\leq 5$, consider the moduli space
\begin{equation}\label{mod-space-family}
M_z^+(W_{-}^+,I;\alpha,\beta):=\bigsqcup_{t\in I} M^+_z(W_{-}^+,\pi_-^+(t);\alpha,\beta).
%  M_z^+(W_{-}^+,c^+_{-},I;\alpha,\beta):=\bigsqcup_{t\in I} M^+_z(W_{-}^+,c^+_{-},\pi_-^+(t);\alpha,\beta).
\end{equation}
Here the disjoint union is meant in the set-theoretic sense, not the topological sense; topologically, if $\mathbf A_n$ is a sequence of $\pi_-^+(t_n)$-instantons, and $t_n \to t$ while $\mathbf A_n$ chain-converges to $\mathbf A$, then this space is topologized so that $(t_n, \mathbf A_n) \to (t, \mathbf A)$. This is completely analogous to the discussion in Section \ref{sec:moduli-fam}, and the definition of $\pi_-^+(\infty)$-instantons is identical to the definition given there.

The following lemma can be established by the same arguments which give Theorem \ref{thm:unobs-transv}.

\begin{lemma}
	The family $\pi_-^+(t)$ can be chosen such that any element in the moduli space 
	\[M_z^+(W_{-}^+,I;\alpha,\beta)\; \setminus \;M_z^+(W_-^+, \pi_-^+(\infty); \alpha, \beta)\] is regular.
\end{lemma}

Similar in spirit to Section \ref{obs-glu}, we will divide non-regular solutions of $M_z^+(W_-^+, I; \alpha, \beta)$ into two types, all of which belong to the stratum $M^+_z(W_-^+, \pi_-^+(\infty); \alpha, \beta)$. 

For any $\alpha$ and $\beta = \rho$, we say that an element of the subspace 

\begin{equation}\label{htpy-type-I}
	M_{z_+}^+(W_+,\pi_+;\alpha,\rho')\times_{\rho'} \mathcal O_{-}
\end{equation}
of $M_z^+(W_{-}^+,I;\alpha,\beta)$ is an obstructed solution of type I, where $\mathcal O_{-}$ denotes the orbit of $\Lambda_-$ in the framed moduli space. Here $z_+$ is the unique path with $z = z_+ * \lambda_-$.

For any $\alpha$ and $\beta \ne \rho$, we say that an obstructed solution of type II is an element of the subspace

\begin{equation}\label{htpy-type-II}
	X_{z_+\Lambda_- \zeta}:=M_{z_+}^+(W_+,\pi_+;\rho,\rho')\times_{\rho'} \mathcal O_{-}\times_{\rho} \breve M^+_{\zeta}(\pi_0;\rho,\beta).
\end{equation}
of $M_z^+(W_{-}^+,I;\alpha,\beta)$. Here we have $z = z_+ * \lambda_- * \zeta$.

We have analogues of Proposition \ref{obs-glu-global} for the obstructed solutions of any of the above four types that describe the behavior of the moduli space $M_z^+(W_{-}^+,I;\alpha,\beta)$ in a neighborhood of obstructed solutions. We start by analyzing the case of reducible solutions of type I. 

First note that the complex line bundle $\mathcal H^+_{{\Lambda}_-}$ over $\mathcal O_-$ pulls back to a line bundle $\mathcal H_{z_+}$ over $X_{z_+ \Lambda_-}$. Obstructed gluing theory gives an $SO(3)$-equivariant section $\psi^{\textup{an}}_{z_+\Lambda_-}$ of the pullback of $\mathcal H_{\alpha,\rho'}$  to \[(0,\infty] \times M_{z_+}^+(W_+, \pi_+; \alpha, \rho') \times_{\rho'} \mathcal O_-.\] 
This section vanishes over \begin{equation}\label{zero-psi-type-I}
	\{\infty\} \times M_{z_+}^+(W_+,\pi_+;\alpha,\rho')\times_{\rho'} \times \mathcal O_{-}
\end{equation}
and $\psi^{\textup{an}}_{z_+\Lambda_-}$ is transverse to the zero section in the complement of this subspace. 

As in part (ii) of Proposition \ref{obs-glu-global}, there is also an injective open continuous map
\[
  \Phi_{z_+\Lambda_-}: \left(\psi^{\textup{an}}_{z_+\Lambda_-}\right)^{-1}(0) \to M_{z}(W_{-}^+,I;\alpha,\rho)
\]
whose image is an open neighborhood of the subspace of obstructed solutions of type I, where $z = z_+ * \lambda_-$. 

For any $[x,y]\in M_{z_+'}(W_+,\pi_+;\alpha,\alpha')\times_{\alpha'} \breve M^+_{\zeta'}(\pi_1;\alpha',\rho')$ with $z_+=z_+'*\zeta'$, we have
\[
  \psi^{\textup{an}}_{z_+\Lambda_-}(t,[x,y])=\psi^{\textup{an}}_{\zeta'\Lambda_-}(t,y)
\]
where $\psi^{\textup{an}}_{\zeta'\Lambda_-}$ is the map to $\mathcal H_\Lambda^+$ defined in Section \ref{obs-glu}. Similarly for any $[x,y]\in \breve M^+_{\zeta}(\pi_0;\alpha,\gamma) \times _{\gamma}M_{z_+'}(W_+,\pi_+;\gamma,\rho')$ with $z_+=\zeta*z_+'$, we have 
\[
  \psi^{\textup{an}}_{z_+\Lambda_-}(t,[x,y])=\psi^{\textup{an}}_{z_+'\Lambda_-}(t,y).
\]
Just as in Section \ref{obs-glu}, there are further compatibility maps on the open embeddings $\Phi$.

For the obstructed solutions of type II, first we define the complex line bundle $\mathcal H_{z_+}$ over $X_{z_+\Lambda_- \zeta}$ as the pullback of $\mathcal H_{\Lambda}^+$ as usual. Then we have the obstruction section $\psi^{\textup{an}}_{z_+\Lambda_-\zeta}$ defined over
\[(0, \infty] \times M_{z_+}^+(W_+,\pi_+;\alpha,\rho')\times_{\rho'} \mathcal O_{-}\times_{\rho} [-\infty, 0) \times \breve M^+_{\zeta}(\pi_0;\rho,\beta).\]

This section vanishes on the strata lying over over $\{(\infty, -\infty)\} \in (0, \infty] \times [-\infty, 0)$, and is transverse to the zero section in the complement of these strata. Again, there is a map $\Phi_{z_+ \Lambda_- \zeta}$ from the zero set of $\psi^{\textup{an}}$ into an open neighborhood of the obstructed solutions of type II.

The restriction of this obstruction section to the remaining strata is determined by the obstructions sections defined on smaller moduli spaces. For any $[x,\Lambda_-,y]\in X_{z_+\Lambda_- \zeta}$, we have
\[
  \psi^{\textup{an}}_{z_+\Lambda_-\zeta}(t,[x,\Lambda_-],-\infty,y)=\psi^{\textup{an}}_{z_+\Lambda_-}(t,[x,\Lambda_-]),
  \hspace{1cm}
  \psi^{\textup{an}}_{z_+\Lambda_-\zeta}(\infty,[x,\Lambda_-],t',y)=\psi^{\textup{an}}_{\Lambda_-\zeta}(\Lambda_-,t',y).
\]
For any $[x,y]\in M_{z_+'}(W_+,\pi_+;\alpha,\alpha')\times_{\alpha'} \breve M^+_{\zeta'}(\pi_1;\alpha',\rho')$ with $z_+=z_+'*\zeta'$ and $u\in \breve M_{\zeta}(\pi_0;\rho,\beta)$, we have
\[
  \psi^{\textup{an}}_{z_+\Lambda_-\zeta}(t,[x,y,\Lambda_-],t',u)=\psi^{\textup{an}}_{\zeta'\Lambda_-\zeta}(t,[y,\Lambda_-,u],t'),
\]
and for any $[x,y]\in \breve M^+_{\zeta}(\pi_0;\alpha,\gamma) \times _{\gamma}M_{z_+'}(W_+,\pi_+;\gamma,\rho')$ with $z_+=\zeta*z_+'$ and $u\in \breve M_{\zeta}(\pi_0;\rho,\beta)$, we have
\[
  \psi^{\textup{an}}_{z_+\Lambda_-\zeta}(t,[x,y,\Lambda_-],t',u)=\psi^{\textup{an}}_{z_+'\Lambda_-\zeta}(t,[y,\Lambda_-],t',u).
\]
Finally if $x\in M_{z_+}^+(W_+,\pi_+;\alpha,\rho')$ and $[y,u]\in \breve M_{\zeta_0}(\pi_0;\rho,\gamma)\times_\gamma M_{\zeta_1}(\pi_0;\gamma,\beta)$ with $\zeta=\zeta_0*\zeta_1$, then 
\[
  \psi^{\textup{an}}_{z_+\Lambda_-\zeta}(t,[x,\Lambda_-],t',[y,u])=\psi^{\textup{an}}_{z_+\Lambda_-\zeta_0}(t,[x,\Lambda_-],t',y).
\]
For all of these obstructed solutions, there is a map $\Phi$ with appropriate decoration that identifies the zero locus of the obstruction section with an open neighborhood of obstructed solutions in the moduli space $M_z^+(W_{-}^+,I;\alpha,\beta)$.
 
As in the case of the modified moduli spaces, we first remove a neighborhood of obstructed solutions in the moduli spaces $M_z^+(W_{-}^+,I;\alpha,\beta)$ and then add a new part, which is defined in terms of a family of privileged sections of the obstruction bundles and homotopies from the analytical obstruction sections to these sections. 

To initiate this construction, we fix $T$, and let $\psi^{T}_{z_+\Lambda_-}$ denote the restriction of $\psi^{\textup{an}}_{z_+\Lambda_-}$ to the space $\{T\} \times M_{z_+}^+(W_+,\pi_+;\alpha,\rho')\times_{\rho'} \mathcal O_{-}$. Similarly, let $\psi^{T}_{z_+\Lambda_-}\zeta$ be the restriction of $\psi^{\textup{an}}_{z_+\Lambda_-\zeta}$ to the subspace of $X_{z_+ \Lambda_- \zeta}$ lying over $A_T \subset (0, \infty] \times [-\infty, 0)$; recall that this subspace $A_T$ was defined in Section \ref{subsec:modified-mod-space}.

We then defined privileged sections $\psi^{\text{priv}}_{z_+ \Lambda_-}$ and $\psi^{\text{priv}}_{z_+ \Lambda_- \zeta}$ over these subspaces, which satisfy the same compatibilty conditions that the analytic gluing sections do. We demand that the latter takes the form 
\[\psi^{\text{priv}}_{z_+\Lambda_-\zeta}(t,[x,y],t') = \frac{T}{t}\cdot \psi^{\text{priv}}_{z_+ \Lambda_-}(x)-\frac{T}{t'}\cdot \psi^{\text{priv}}_{\Lambda_-\zeta}(y);\] we assume that the two previously-chosen sections $\psi^{\text{priv}}$ (one chosen here, the other chosen in Section \ref{subsec:modified-mod-space}) are chosen so that this is a transverse section. 

We then choose homotopies of sections $\psi^{\text{htpy}}_{z_+ \Lambda_-}$ and $\psi^{\text{htpy}}_{z_+ \Lambda_- \zeta}$ once again satisfying the same compatibility and transversality properties, and interpolating between the analytic gluing sections $\psi^T$ and the priveleged sections $\psi^{\text{priv}}$.

Finally, the modified moduli spaces $N^+_z(W_-^+, I; \alpha, \beta)$ are defined by exactly the same procedure given in Section \ref{subsec:modified-mod-space}, with the same orientation conventions: they are oriented as zero sets of sections of complex line bundles, where the domains of these sections are oriented as products given the product structures written above. In the case of $\psi^{\text{priv}}_{z_+ \Lambda_- \zeta}$, we orient the zero set as the zero set of \[s \psi^{\text{priv}}_{z_+ \Lambda_-}(x) + (1-s) \cdot \psi^{\text{priv}}_{\Lambda_- \zeta}(y)\] over the space \[M_{z_+}(W_+, \pi_+; \alpha, \rho') \times_{\rho'} \mathcal O_- \times_\rho \([0,1] \times \breve M_\zeta(Y; \rho, \beta)\).\]

That these are stratified-smooth spaces follows by the same argument as in Section \ref{subsec:modified-mod-space}. Below, we record the boundary of these stratified-smooth spaces, together with their orientations; unlike Section \ref{subsec:modified-mod-space}, we use the language of geometric chains to avoid talking about faces which appear with opposite orientations.

Suppose $z = z_+ * \lambda_-$, where $z_+: \alpha \to \rho'$ is a path of connections on $W_+$.

\begin{prop}\label{prop:hard-way}
	The boundary of the stratified smooth space $N_z^+(W_{-}^+,I;\alpha,\beta)$, considered as a geometric chain in $\alpha \times \beta$, is given 
	as follows.
	\begin{itemize}
		\item[(i)] If $\beta=\rho$, then $\partial N_z^+(W_{-}^+,I;\alpha,\rho)$ is given by 
			\begin{align*}
			  (-1)^{\dim \alpha}\partial N_z^+(W_{-}^+,I;\alpha,\rho) =\; &\sum_{z = z_+ * z_-} M^+_{z_+}(W_+, \pi_+; \alpha, \alpha') \times_{\alpha'} N^+_{z_-}(W_-, \pi_-; \alpha', \rho) \\
&-N^+_z(W_-^+, \pi_-^+(T_0); \alpha, \rho) \\
&+ \sum_{z = \zeta * z'} (-1)^{\wt i(\zeta)} \breve M^+_\zeta(Y,\pi_0;\alpha, \gamma) \times_\gamma N^+_{z'}(W_-^+, I; \gamma, \rho) \\
&+ \sum_{z = z' * \zeta'} (-1)^{\wt i(z')} N^+_{z'}(W_-^+, I; \alpha, \gamma) \times_\gamma \breve M^+_{\zeta'}(Y,\pi_0; \gamma, \rho) \\
&+ \(\psi_{z_+\Lambda_-}^{\textup{priv}}\)^{-1}(0).
\end{align*}
In the final term, we have $z = z_+ * \lambda_-$, where $z_+: \alpha \to \rho'$ is a path of connections along $W_+$.

		\item[(ii)] If $\beta \ne \rho$, then $\partial N_z^+(W_-^+, I; \alpha, \beta)$ is given by 
			\begin{align*}
			  (-1)^{\dim \alpha}\partial N_z^+(W_{-}^+,I;\alpha,\beta) =\; &\sum_{z = z_+ * z_-} M^+_{z_+}(W_+, \pi_+; \alpha, \alpha') \times_{\alpha'} N^+_{z_-}(W_-, \pi_-; \alpha', \beta) \\
&-N^+_z(W_-^+, \pi_-^+(T_0); \alpha, \beta) \\
&+ \sum_{z = \zeta * z'} (-1)^{\wt i(\zeta)} \breve M^+_\zeta(Y,\pi_0;\alpha, \gamma) \times_\gamma N^+_{z'}(W_-^+, I; \gamma, \beta) \\
&+ \sum_{z = z' * \zeta'} (-1)^{\wt i(z')} N^+_{z'}(W_-^+, I; \alpha, \gamma) \times_\gamma \breve M^+_{\zeta'}(Y,\pi_0; \gamma, \beta) \\
&- \(\psi_{z_+ \Lambda_- \zeta}^{\textup{priv}}\)^{-1}(0).
\end{align*}
where in the final term we have $z = z_+ * \lambda_- * \zeta$.
	\end{itemize}
\end{prop}

\begin{remark}
The computation of boundary faces is analogous to the work concluding Section \ref{subsec:modified-mod-space}, including the discussion of signs. For the final term, the domain of this section is oriented as 
\[\alpha \times (T_0,\infty) \times M^{\text{Fiber}}_{z_+}(W_+; \alpha, \rho') \times (-\infty, 0) \times \breve M^{\text{Fiber}}_\zeta(Y, \pi_0; \rho, \beta);\]
here $t \in (T_0, \infty)$ corresponds to the parameter in $\pi_-^+(t)$, while $(-\infty, 0)$ corresponds to how far we translate the instanton on $Y$ to the left before gluing.

To compute the orientation on the zero set, one shuffles $(T_0,\infty) \times (-\infty, 0)$ to the front incurring a sign of $(-1)^{\wt i(z_+)}$, then passes to the locus $A_T \cong [0,1]$ where $T > T_0$ is large, which incurs a sign of $-1$ because the outward-normal-first convention is here negative the orientation of the regular zero set convention. Finally, one shuffles that back to 
\[(-1)^{\dim \alpha - 1}\alpha \times M^{\text{Fiber}}_{z_+}(W_+; \alpha, \rho') \times [0,1] \times \breve M^{\text{Fiber}}_\zeta(Y, \pi_0; \rho, \beta).\qedhere\]
\end{remark}

Moving on to the other composite, one may construct a family of perturbations $\pi_+^-(t)$ on the cylinder which are equal to $\pi_1$ on the ends, so that $\pi_+^-(\infty)$ is the broken perturbation $\pi_- \sqcup \pi_+: (Y, \pi_1) \to (Y, \pi_0) \to (Y, \pi_1)$. Running through the same construction, we obtain modified moduli spaces $N_z^+(W_+^-, I)$, and we enumerate the boundary strata (with signs) below. Here when doing gluing on both sides of $\Lambda_-$ our open set is indexed by $(T_0,\infty) \times (0,\infty)$, with the first term corresponding to the perturbation parameter, the second corresponding to how far left we translate the instanton on $(\Bbb R \times Y, \pi_1)$ before gluing. We instead use the function $L(t,t') = -1/t-1/t'$ to truncate, with locus $A'_T = \{(t,t') \mid 1/t + 1/t' = 1/T\}$ oriented as the regular preimage $L^{-1}(-1/T)$. The isomorphism $A'_T \cong [0,1]$ given by $t \mapsto 1-T/t = T/t'$ is orientation-preserving, and sends the section $\psi_{\zeta' \Lambda_- z_+}^{\text{priv}}$ to \[s \psi^{\text{priv}}_{\zeta' \Lambda_-}(x) + (1-s) \cdot \psi^{\text{priv}}_{\Lambda_- z_+}(y).\]

\begin{prop}\label{prop:easy-way}
	The boundary of the stratified smooth space $N_z^+(W_+^-,I;\alpha',\beta')$ as a geometric chain is given 
	as follows.
	\begin{itemize}
		\item[(i)] If $\alpha'=\rho'$, then $\partial N_z^+(W_+^-,I;\rho', \beta')$ is given by 
			\begin{align*}
			  (-1)^{\dim \rho'}\partial N_z^+(W_+^-,I;\rho',\beta') =\; &\sum_{z = z_- * z_+} N^+_{z_-}(W_-, \pi_-; \rho', \alpha) \times_{\alpha} M^+_{z_+}(W_+, \pi_+; \alpha, \beta') \\
&-N^+_z(W_+^-, \pi_+^-(T_0); \rho', \beta') \\
&+ \sum_{z = \zeta' * z'} (-1)^{\wt i(\zeta')} \breve M^+_{\zeta'}(Y,\pi_1;\rho', \gamma') \times_{\gamma'} N^+_{z'}(W_+^-, I; \gamma', \beta') \\
&+ \sum_{z = z' * \zeta} (-1)^{\wt i(z')} N^+_{z'}(W_+^-, I; \rho', \gamma') \times_{\gamma'} \breve M^+_\zeta(Y,\pi_1; \gamma', \beta') \\
&+ \(\psi_{\Lambda_- z_+}^{\textup{priv}}\)^{-1}(0).
\end{align*}
In the final term, we have $z = \lambda_- * z_+$, where $z_+: \rho \to \beta'$ is a path of connections along $W_-$. 
		\item[(ii)] If $\alpha' \ne \rho'$, then $\partial N_z^+(W_+^, I; \alpha, \beta)$ is given by 
			\begin{align*}
			  (-1)^{\dim \alpha'}\partial N_z^+(W_+^-,I;\alpha',\beta') =\; &\sum_{z = z_- * z_+} N^+_{z_-}(W_-, \pi_-; \alpha', \alpha) \times_{\alpha} N^+_{z_+}(W_+, \pi_+; \alpha, \beta') \\
&-N^+_z(W_+^-, \pi_+^-(T_0); \alpha', \beta') \\
&+ \sum_{z = \zeta' * z'} (-1)^{\wt i(\zeta')} \breve M^+_{\zeta'}(Y,\pi_1;\alpha', \gamma') \times_{\gamma'} N^+_{z'}(W_+^-, I; \gamma', \beta') \\
&+ \sum_{z = z' * \zeta} (-1)^{\wt i(z')} N^+_{z'}(W_+^-, I; \alpha', \gamma') \times_{\gamma'} \breve M^+_\zeta(Y,\pi_1; \gamma', \beta') \\
&-(-1)^{\wt i(\zeta')} \(\psi_{\zeta'\Lambda_- z_+}^{\textup{priv}}\)^{-1}(0),
\end{align*}
where in the final term we have $z = \zeta' * \lambda_- * z_+$.
	\end{itemize}
\end{prop}

\begin{remark}
For the final term, the domain of this section is oriented as 
\[\alpha' \times (T_0, \infty) \times (0,\infty) \times \breve M^{\text{Fiber}}_{\zeta'}(Y, \pi_1; \alpha', \rho') \times M^{\text{Fiber}}_{z_+}(W_+; \rho, \beta');\]
one shuffles $(T_0,\infty) \times (0,\infty)$ to the front incurring no additional sign, passes to the locus $A'_T \cong [0,1]$ (which incurs a sign of $-1$ because the outward-normal-first convention is here negative the orientation of the regular zero set convention), and shuffles that back to 
\[(-1)^{\dim \alpha' + \wt i(\zeta') - 1}\alpha' \times M^{\text{Fiber}}_{\zeta'}(Y, \pi_1; \alpha', \rho') \times [0,1] \times \breve M^{\text{Fiber}}_{z_+}(W_+;\rho, \beta').\qedhere\]
\end{remark}

\subsubsection{Homotopies of bimodules from modified moduli spaces}
We can now use the constructions of the preceding section to define the desired bimodule homotopies. To check these easily, have on-hand the definitions of $\mathcal W_+$ from Construction \ref{constr:Wplus-bimod} and $\mathcal W_-$ from Construction \ref{constr:W-minus}, as well as the definition of bimodule homotopy from Definition \ref{def:htpy}. We begin with the easier of the two cases.

For compactness of notation, write $\mathcal H_+^-(\alpha, \beta)$ for the modified moduli spaces $N^+(W,\pi_+^-(t); \alpha, \beta)$ from the previous section. Let $\pi_+^-(T_0)$ be the constant perturbation at $\pi_1$, so that $M^+(W, \pi_+^-(T_0))$ defines the identity bimodule on $\mathcal I_1 = \mathcal I(Y, \pi_1)$.

\begin{prop}\label{htpy2}
The modified moduli spaces above define a bimodule homotopy \[\mathcal H_+^-: \mathbf 1_{\mathcal I_1} \implies \mathcal W_+ \circ \mathcal W_-.\] 
\end{prop}
\begin{proof}
We hope that the following diagram indicating the relevant bimodules is helpful when parsing the argument; the composite bimodule is given by summing over the fiber products of all two-step paths running from $\alpha'$ to $\beta'$, and we want to find a homotopy from the identity bimodule to this composite bimodule. 

\[\begin{tikzcd}
	&&& {\alpha} \\
	&&& {} \\
	\alpha' & {} && {S_{\rho}} && {} & \beta' \\
	&&& {} \\
	&&& {\overline \rho}
	\arrow["{B(W_+; \rho, \beta')}"', from=3-4, to=3-7]
	\arrow["{M(W_-; \alpha', \rho) \times_{\rho} S_{\rho}}", from=3-1, to=3-4]
	\arrow["{M(W_-; \alpha', \alpha)}", bend left=20, from=3-1, to=1-4]
	\arrow["{M(W_+; \alpha, \beta')}", bend left=20, from=1-4, to=3-7]
	\arrow["{\delta_{\alpha' = \rho'} \mathcal O_-}"', bend right=20, from=3-1, to=5-4]
	\arrow["{Z(W_+; \rho, \beta')}"', bend right=20, from=5-4, to=3-7]
	\arrow["{-(-1)^{\widetilde i(\zeta')} B(\pi_1;\alpha', \rho') }"', from=3-1, to=3-4]
\end{tikzcd}\]

The boundary relations enumerated in Proposition \ref{prop:easy-way} almost immediately show that this gives a bimodule homotopy from the identity bimodule. The only interesting stratum is the last one, arising from obstructed gluing theory; we will focus on the case $\beta' \ne \rho'$, as the case $\beta = \rho$ is simpler. In the homotopy relation, this boundary stratum comes from the composite term, as 
\[\mathcal W_-(\alpha', \rho) \times_\rho \mathcal W_+(\rho, \beta')\] 
has the extra term 
\[-(-1)^{\wt i(\zeta)} B_\zeta(\alpha', \rho') \times_{S_\rho} B(W; \rho, \beta');\] 
because we use $-\psi^{\text{priv}}_{\zeta' \Lambda_-}$ in defining the blowups, Lemma \ref{weird-zero-set} identifies this with $-(-1)^{\wt i(\zeta)}$ times the zero set of $(1-s)\psi^{\text{priv}}_{\Lambda_- z_+} - s \psi^{\text{priv}}_{\zeta' \Lambda_-}$, which is identified with the zero set of the section $\psi_{\zeta' \Lambda_- z_+}^{\text{priv}}.$ So this additional term in the homotopy relation is precisely the same as the additional term in the boundary relation of Proposition \ref{prop:easy-way}, as desired.
\end{proof}

Though the obstructed-gluing analysis of the two cases is nearly identical (except for one sign), the definition of the bimodule homotopy running in the other direction, and the verification that it is a bimodule homotopy, are significantly more intricate. 

\begin{construction}
For compactness of notation, write $H_-^+(\alpha, \beta)$ for the modified moduli spaces $N^+(W,\pi_-^+(t); \alpha, \beta)$ from the previous section. Here we assume that $\pi_-^+(T_0)$ is the constant perturbation at $\pi_0$, so that $N_z^+(W_-^+, \pi_-^+(T_0); \alpha, \beta)$ is either diffeomorphic to $[0,1] \times \breve M(\pi_0; \alpha, \beta)$, with endpoint maps determined by the latter factor), or when $\alpha = \beta$ equal to $\alpha$ itself. In the former case, this is a trivial chain (it supports the orientation-reversing isomorphism given by $t \mapsto 1-t$ in the first factor), hence zero at the chain level, as in Remark \ref{identity-morse}.

Before moving on, observe that $H_-^+(\rho, \rho)^{\text{red}}$ is diffeomorphic to $[0,1] \times \rho$, and (because the endpoint maps determine an equivariant map to $\rho \times \rho$, which has two isolated reducible orbits) the endpoint maps factor through the projection to $\rho$, so this is also a trivial chain.

Over each moduli space $H_-^+(\rho, \beta)$, we will define certain equivariant sections of the complex line bundle $\mathcal H_\rho$, written $\psi^H_{\rho\beta}$. We prescribe the boundary values on the various faces enumerated in Proposition \ref{prop:hard-way}: $\psi^H_{\rho \beta}$ should be equal to the privileged sections $\psi^{\text{priv}}_{z_+ \Lambda_-}$ on the first listed face; should agree with the privileged sections $\psi_{\zeta \Lambda_-}^{\text{priv}}$ on the third face; should be pulled back from $\psi^H_{\rho \gamma}$ on the fourth enumerated face; and on the final face, diffeomorphic to $B(W_+; \rho, \rho') \times_{S_{\rho'}} B(\pi_0; \rho, \rho)$, it should be pulled back from a section $\psi^+_{\rho \rho}$ over the first factor. On the part of the fourth face given by $H_-^+(\rho, \rho)^{\text{red}} \times_\rho \breve M(\pi_0; \rho, \beta)$, we demand that the section $\psi^H_{\rho \beta}$ factor through the projection to $\breve M(\pi_0; \rho, \beta)$. This guarantees that even after blowing up, this boundary component is a trivial chain.

Finally, in the case $\beta = \rho$, we only attempt to define a transverse section over the irreducible locus of $H_-^+(\rho, \rho)$; there is no transverse section over the reducible locus, but we will (mostly) be able to ignore it, as this is a trivial chain.

Define $B_-^+(\rho, \beta)$ to be the blowup of $\psi^H_{\rho \beta}$, so a geometric chain in $S_\rho \times \beta$; let $Z_-^+(\rho, \beta)$ be the zero set of $\psi^H_{\rho \beta}$, so still a geometric chain in $\rho \times \beta$.

In the special case $\beta = \rho$, we set $\widehat B_-^+(\rho, \rho)$ to be the blowup of the section $\psi^H_{\rho \rho} \sqcup \psi^+_{\rho \rho}$ over 
\[H_-^+(\rho, \rho)^{\text{irr}} \times_\rho S_\rho - B(W_+; \rho, \rho');\] transversality dictates that the zero set $\widehat Z_-^+(\rho, \rho)$ is empty.

Then the bimodule homotopy 
\[\mathcal H_-^+: \mathbf 1_{\mathcal S_\rho \mathcal I_0} \implies \mathcal W_- \circ \mathcal W_+\] 
is defined by the formula 
\[\mathcal H_-^+(\alpha, \beta) = \begin{pmatrix}[c|ccc] & \alpha & S_\rho & \overline \rho \\
\hline 
\beta & H_-^+(\alpha, \beta) & -B_-^+(\rho, \beta) & Z_-^+(\rho, \beta) \\
S_\rho & H_-^+(\alpha, \rho) \times_\rho S_\rho - (-1)^{\wt i(z_+)} B_{z_+}(W_+; \alpha, \rho) & -\widehat B_H(\rho, \rho) & 0 \\ 
\overline \rho & 0 & 0 & 0 
\end{pmatrix}.\qedhere\] 
\end{construction}

\begin{prop}\label{htpy1}
The formula above defines a bimodule homotopy.
\end{prop}

To aid with the verification, we include a diagram indicating both of the bimodules we are composing.

\[\begin{tikzcd}
	\alpha & {} &&&&& \beta \\
	&&&&&& {} \\
	{S_{\rho}} &&& \alpha' & {} && {S_{\rho}}\\
	&&&&&&\\
	{\overline \rho} & {} &&&&& {\overline \rho}
	\arrow["{\delta_{\alpha' \ne \rho'}N(W_-; \alpha', \rho) \times_{\rho} S_{\rho}}", from=3-4, to=3-7]
	\arrow["{N(W_-; \alpha', \beta)}", bend left=20, near end, from=3-4, to=1-7]
	\arrow["{\delta_{\alpha' = \rho'} \mathcal O_-}", bend right=20, from=3-4, to=5-7]
	\arrow["{B(W_+; \rho, \alpha')}", from=3-1, to=3-4]
	\arrow["{Z(W_+; \rho, \alpha')}"', bend right=20, from=5-1, to=3-4]
	\arrow["{M(W_+; \alpha, \alpha')}", bend left=20, near start, from=1-1, to=3-4]
	\arrow["{-(-1)^{\widetilde i(\zeta')} B(\pi_1;\alpha', \rho') }"', from=3-4, to=3-7]
\end{tikzcd}\]

Next, we depict the relevant flow categories and the bimodule homotopies.

\[\begin{tikzcd}
	\alpha / \gamma & {} &&&&& \gamma / \beta \\
	&&&&&& {} \\
	{S_{\rho}} &&& {} & {} && {S_{\rho}}\\
	&&&&&& {} & {} \\
	{\overline \rho} & {} &&&&& {\overline \rho}
	\arrow["{H_-^+(\alpha, \beta)}", from=1-1, to=1-7]
	\arrow["{-\widehat B_H(\rho, \rho)}"', bend right=15, from=3-1, to=3-7]
	\arrow["{H_-^+(\alpha, \rho) \times_\rho S_\rho}", bend right=10, pos=0.18, from=1-1, to=3-7]
	\arrow["{-(-1)^{\wt i(z_+)} B(W_+; \alpha, \rho')}"', bend right=10, pos=0.7, from=1-1, to=3-7]
	\arrow["{-B_-^+(\rho, \beta)}", near start, from=3-1, to=1-7]
	\arrow["{Z^+_-(\rho, \beta)}", bend right=40, near start, from=5-1, to=1-7]
	\arrow["{S_{\rho}}", from=3-1, to=5-1]
	\arrow["{}", from=1-1, to=3-1]
	\arrow["{S_{\rho}}"', from=3-7, to=5-7]
	\arrow["{}"', from=1-7, to=3-7]
	\arrow["{-B(\pi_0; \rho, \alpha)}", bend left=40, from=3-1, to=1-1]
	\arrow["{Z(\pi_0; \rho, \alpha)}",bend left=20, from=5-1, to=1-1]
	\arrow["{-B(\pi_0; \rho, \beta)}"', bend right=40, from=3-7, to=1-7]
	\arrow["{Z(\pi_0; \rho, \beta)}"', bend right=20, from=5-7, to=1-7]
	\arrow["{\breve M(\pi_0; \alpha, \gamma)}", loop above, from=1-1, to=1-1]
	\arrow["{\breve M(\pi_0; \gamma, \beta)}", loop above, from=1-7, to=1-7]
\end{tikzcd}\]

Then the relation we are attempting to prove asserts that for every left-to-right arrow in the second diagram (from, say, $\alpha$ to $\beta$), its boundary is a sum over all two-step paths in the second diagram which start at $\alpha$ and end at $\beta$, plus the corresponding fiber product in the first diagram, minus $\delta_{\alpha = \beta} \alpha$.

\begin{proof}[Sketch of proof]
Because this computation is intricate and delicate, we will avoid discussing the signs (which would only serve to obfuscate the discussion); the careful reader should not have difficulty inserting the signs into the arguments below. That $\mathcal H_-^+(\alpha, \beta)$ satisfies the homotopy relation for $\alpha, \beta \ne \rho$ follows exactly as in the computation that $\mathcal W_-$ is a bimodule. 

The novel parts are those involving moduli spaces of the form $H_-^+(\rho, \rho)$. We will only discuss those parts of the bimodule relation which involve this term below, as all other parts of the discussion follow the same strategy as previous arguments; the reader should understand this as a guide to clarifying potentially confusing details when running through the verification.

\begin{itemize}
\item In entry $(1,2)$, the term $\mathcal H_-^+(S_\rho, S_\rho) \times_{S_\rho} \mathcal S_\rho \mathcal C(S_\rho, \beta)$ is given by 
\[B\(H_-^+(\rho, \rho)^{\text{irr}} \times_\rho S_\rho + B(W_+; \rho, \rho')\) \times_{S_{\rho}} B_\zeta(\rho, \beta);\] the left-hand blowup is determined entirely by the first factor, so we focus on this last. The latter term simplifies to give the $\(\psi^{\text{priv}}_{z_+ \Lambda_- \zeta}\)^{-1}(0)$ term, blown up on the left-hand side (this is where we use the assumption that $\psi_{\rho \beta}^H$ is pulled back from a section $\psi^H_{\rho \rho}$ on the first factor here). The former term simplifies to give 
\[B\(H_-^+(\rho, \rho)^{\text{irr}} \times_\rho \breve M_\zeta(\rho, \beta)\);\] by the assumption that the section on this term is pulled back from the first factor, this simplifies to $B_-^+(\rho, \rho)^{\text{irr}} \times_\rho \breve M_\zeta(\rho, \beta)$, which is one of the terms promised by the boundary relation Proposition \ref{prop:hard-way}. Finally, this leaves one term in the known boundary relation unaccounted for: the blowup of $H_-^+(\rho, \rho)^{\text{red}} \times_\rho \breve M_\zeta(\rho, \beta).$ However, we assumed that the section used to define the blowup factors through a projection to the latter factor (and every stratum of the first factor has positive dimension), so that this chain has small rank on each positive-dimensional stratum. In particular, the corresponding blown-up chain is degenerate, so that this term is zero on the chain level. This squares away the unexpected factors in both the bimodule relation and the boundary relation Proposition \ref{prop:hard-way}.

Entry $(1,3)$ is similar but strictly simpler because $\widehat Z_-^+(\rho, \rho) = \varnothing$. 

\item In entry $(2, 1)$, the novel terms arise from $\mathcal S_\rho \mathcal C(\alpha, S_\rho) \times_{S_\rho} \mathcal H_-^+(S_\rho, S_\rho),$ which may be expanded as 
\[\breve M(\alpha, \rho) \times_\rho S_\rho \times_{S_\rho} B\left(H_-^+(\rho, \rho)^{\text{irr}} \times_\rho S_\rho + B(W_+; \rho, \rho')\right);\] cancelling out $S_\rho \times_{S_\rho}$ and replacing the blowup with the original chain, this simplfies to 
\[\breve M(\alpha, \rho) \times_\rho H_-^+(\rho, \rho)^{\text{irr}} \times_\rho S_\rho + \breve M(\alpha, \rho) \times_\rho M(W_+; \rho, \rho'),\] which contribute the terms corresponding to factorization through $\rho$ to the boundary relations of both $H_-^+(\alpha, \rho) \times_\rho S_\rho$ and $B_{z_+}(W_+; \alpha, \rho)$. In the former case, this still leaves an absent term corresponding to $H_-^+(\rho, \rho)^{\text{red}}$; but this is a degenerate geometric chain, so zero at the chain level, so irrelevant in checking the boundary relation.

\item In entry $(3, 1)$, there are only two terms in the bimodule homotopy relation for $\partial \mathcal H_-^+(\alpha, \overline \rho)$ which are not tautologically empty: the first is 
\[\mathcal M_1(\alpha, \overline \rho) = \mathcal W_+(\alpha, \rho') \times_{\rho'} \mathcal W_-(\rho', \overline \rho) = M(W_+; \alpha, \rho') \times_{\rho'} \mathcal O_- \cong M(W_+; \alpha, \rho'),\] 
and the second is 
\[\mathcal H_-^+(\alpha, S_\rho) \times_{S_\rho} \mathcal S_\rho \mathcal C(S_\rho, \overline \rho) = \(H_-^+(\alpha, \rho) \times_\rho S_\rho + B_{z_+}(W_+; \alpha, \rho)\) \times_{S_\rho} S_\rho,\] 
which simplifies to 
\[H_-^+(\alpha, \rho) \times_\rho S_\rho + B_{z_+}(W_+; \alpha, \rho)\] 
as a chain in $\alpha \times \rho$.

Now the first term is degenerate --- it factors through the projection to $H_-^+(\alpha, \rho)$, so has small rank --- and considering the second as a chain in $\alpha \times \rho$ instead of $\alpha \times S_\rho$ leaves it equivalent to the chain $M(W_+; \alpha, \rho)$. Thus this term cancels out with the previous.

\item Entry $(2,2)$ is by far the most intricate; there is no avoiding talking about every component. First, let us point out that in 
\[\widehat B_H(\rho, \rho) = B\(H_-^+(\rho, \rho)^{\text{irr}} \times_\rho S_\rho + B(W_+; \rho, \rho')\),\] the first term has no boundary: this moduli space has fiber-dimension $1$ over $S_\rho$, so that $H_-^+(\rho, \rho)^{\text{irr}}$ is 3-dimensional, so its quotient by $SO(3)$ is hence a 0-manifold, which has no boundary. Because the section over this space is non-vanishing, taking the left blowup is equivalent to choosing a lift of the map to $\rho$ to a map to $S_\rho$; no new boundary terms appear.

So we are engaged in a computation of the boundary of the second term. Let us first write the boundary predicted by the bimodule homotopy. This consists of four terms: 
\begin{gather*}\sum_\gamma \mathcal S(S_\rho, \gamma) \times_\gamma \mathcal H(\gamma, S_\rho) \quad \quad \sum_\gamma \mathcal H(S_\rho, \gamma) \times_\gamma \mathcal S(\gamma, S_\rho) \\
\mathcal \sum_{\gamma'} \mathcal W_+(S_\rho, \gamma') \times_{\gamma'} \mathcal W_-(\gamma', S_\rho) \quad \quad \mathcal M_0(S_\rho, S_\rho) = S_\rho.
\end{gather*}

Spelling these out, this amounts to $S_\rho$ and six additional boundary terms: \begin{align}
&B(\pi_0; \rho, \gamma) \times_\gamma H^+_-(\gamma, \rho) \times_\rho S_\rho \\
&B(\pi_0; \rho, \gamma) \times_\gamma B(W_+; \gamma, \rho') \\
&B^+_-(\rho, \gamma) \times_\gamma \breve M(\pi_0; \gamma, \rho) \times_\rho S_\rho \\
&B(W_+; \rho, \gamma') \times_{\gamma'} M(W_-; \gamma', \rho) \times_\rho S_\rho \\
&B(W_+; \rho, \gamma') \times_{\gamma'} B(\pi_1; \gamma', \rho') 
\end{align}
The boundary has fiber-dimension zero, so we may discard all terms which necessarily have positive fiber-dimension: the first, third, and fourth do not appear. What we are left with is precisely the boundary relation for $B(W_+; \rho, \rho')$; the final term $S_\rho$ corresponds to $Z(W_+; \rho, \rho') \times_{\rho} S_\rho$, as this is identifiable with $\mathcal O_+ \times_\rho S_{\rho} = S_\rho$, cancelling out $\mathcal M_0(S_\rho, S_\rho) = S_\rho$. So this is indeed the boundary of $\widehat B_H(\rho, \rho)$, as desired.

\item In entry $(3,3)$, the only terms in the bimodule homotopy relation which aren't tautologically zero (or of fiber-dimension $\ge 8$) are $\mathcal M_1(\overline \rho, \overline \rho) - \mathcal M_0(\overline \rho, \overline \rho)$. The latter is a copy of $\overline \rho$ itself, with the identity map as its endpoint maps. The first is a sum over fiber products, but the only term which isn't tautologically zero is $Z(W_+; \overline \rho, S_\rho) \times_{\overline \rho} \mathcal O_-$. Both factors consist of a single reducible orbit, so that the fiber product is also a single reducible orbit; the reducible locus in $H_-^+(\rho, \rho)^{\text{red}}$, diffeomorphic to $[0,1] \times \rho$, gives a homotopy between their endpoint maps. Because the reducible orbits in $\rho \times \rho$ are isolated, this means the endpoint maps can be identified, so that $\mathcal M_1(\overline \rho, \overline \rho) - \mathcal M_-(\overline \rho, \overline \rho) = \overline \rho - \overline \rho = 0$. So this term in the bimodule homotopy relation is satisfied.
\end{itemize}

We leave boundary relations $(2,3)$ and $(3,2)$ to the interested reader; the arguments are no more difficult.
\end{proof}

This completes the proof of Theorem \ref{thm:main-theorem}.

\newpage

\section{Functoriality of equivariant instanton homology}
\subsection{Cobordism maps in the pseudo-unobstructed case}\label{badred}
%!TEX root = equivariant-functoriality.tex

The results we have thusfar are enough to prove that equivariant instanton homology is a well-defined invariant of rational homology sphere. We described in Definition \ref{def:I-no-pert} how to define a group which is transparently independent of perturbation, which we called $I^\bullet(Y,w)$: we identify the various $I^\bullet(Y, w, \pi)$ using the continuation maps $\varphi_{\pi \to \pi'}$ whenever $\sigma_\pi \le \sigma_{\pi'}$, which we have established are isomorphisms on equivariant homology. We would like to upgrade $I^\bullet(Y,w)$ to a functorial invariant. The goal of this section is to establish the following result, whose proof is given at the end of the section.

\begin{theorem}\label{thm:exists-induced-map}
If $(W,c): (Y,w) \to (Y',w')$ is a cobordism for which either $b^+(W) = 0$ or $(W,c)$ supports no central connections, there is a well-defined induced map $I^\bullet(W,c): I^\bullet(Y,w) \to I^\bullet(Y', w')$ on equivariant insanton homology, which is a module map in the appropriate sense and respects the exact triangles.
\end{theorem}
The statement that $(W,c)$ supports no central connections is equivalent to the claim that $PD(c) \in H^2(W;\Bbb Z)$ is not twice some other class, or equivalently, there exists a closed and oriented surface $\Sigma \subset W$ with $\#(\Sigma \cap c)$ odd.

We would like to argue as follows. Because $I^\bullet(Y,w)$ is independent of the choice of perturbation, we might hope to choose perturbations $\pi, \pi'$ so that $(W,c): (Y,w,\pi) \to (Y',w',\pi')$ is an unobstructed cobordism. Then we could define the induced map on equivariant instanton homology $I^\bullet(Y, w, \pi) \to I^\bullet(Y', w',\pi')$ by applying the construction of Theorem \ref{unobs-cobmap} and then composing with the equivariant homology functors. This induced map $(W, c)_{\pi \to \pi'}$ descends to a well-defined map on $I^\bullet(Y,w)$; to see this, suppose $(\pi_0, \pi_0')$ is another pair of perturbations for which $(W,c)$ is unobstructed and $\sigma_{\pi_0} \le \sigma_\pi$ while $\sigma_{\pi'} \le \sigma_{\pi_0'}$. Then the fact that unobstructed cobordism maps compose functorially implies that \[\varphi_{\pi' \to \pi_0'} (W,c)_{\pi \to \pi'} \varphi_{\pi_0 \to \pi'} = (W,c)_{\pi_0 \to \pi'_0}.\] We will discuss the structure of this argument in more detail towards the end of this section.

Unfortunately, we cannot always find perturbations $\pi, \pi'$ so that $(W,c)$ is unobstructed. The normal index $N(\Lambda; \pi, \pi')$ of a pseudocentral reducible is independent of $\pi$ and $\pi'$, and this quantity may well be negative. By Lemma \ref{everything-is-unobs}, the best we can do is find perturbations $(\pi, \pi')$ so that $(W,c): (Y,w,\pi) \to (Y',w',\pi')$ is `pseudo-unobstructed': all its reducibles are either unobstructed or pseudocentral. 

To prove Theorem \ref{thm:exists-induced-map}, we must construct cobordism maps $\widetilde C(Y,w,\pi) \to \widetilde C(Y',w',\pi')$ associated to pseudo-unobstructed cobordisms, which are well-defined up to equivariant homotopy and compose correctly with the continuation maps.

Now recall the pentachotomy of Proposition \ref{prop:pentachotomy}. In cases (i) and (v), the relevant moduli spaces are cut out transversely (even in families), so that we have already constructed the relevant cobordism maps and verified that they are well-defined up to equivariant homotopy, by Theorem \ref{unobs-cobmap}. In cases (ii) and (iv), we have an induced map associated to a regular perturbation $\pi_W$, but one has to argue carefully to show that the induced map is well-defined up to equivariant homotopy. In case (iii) one has to argue carefully to show that the induced map exists at all.

First we will handle case (iii) of pseudo-unobstructed cobordisms with $b_1(W) - b^+(W) = 0$, explaining both how to guarantee that an induced map exists and why it is well-defined. 

Let $(W, c): (Y, w, \pi) \to (Y', w', \pi')$ be a pseudo-unobstructed cobordism with $b_1(W) - b^+(W) = 0$; by Proposition \ref{prop:pentachotomy}(iii), we may choose a perturbation $\pi_W$ on $(W, c)$ so that all instantons are cut out transversely except these pseudocentral reducibles, each of which corresponds to a 2-dimensional orbit $\mathcal O_\Lambda$ of framed instantons on $W$ whose linearized ASD operator has trivial kernel and cokernel isomorphic to $\Bbb C$.

Given $\theta \in \mathfrak Z(Y,w)$ and $\theta' \in \mathfrak Z(Y', w')$, we will write $\mathcal{P}(W,c;\theta, \theta')$ for the space of pseudocentral instantons on $(W,c)$ running from $\theta$ to $\theta'$. The $SO(3)$-space $\mathcal{P}$ is a disjoint union of orbits $\mathcal O_\Lambda$, each of which is diffeomorphic to $S^2$. There is an obstruction bundle $\mathcal H_{P} \to \mathcal{P}$, an $SO(3)$-equivariant complex line bundle, assigning to each instanton $\Lambda \in \mathcal{P}$ the cokernel of the associated linearized ASD operator. 

We will write $\mathfrak P(W,c;\theta, \theta') = \pi_0 \mathcal{P}(W,c;\theta, \theta')$ for a corresponding discrete space of components. When $(W,c)$ or $(\theta, \theta')$ are clear from context, we will suppress them from notation.

Let us briefly discuss orientations; the naturally oriented object is the $0$-dimensional space $\mathfrak P(W,c;\theta, \theta')$. To see this, recall the conventions of Section \ref{sec:moduli-or}. So long as one chooses an element of $\Lambda(\theta)$ and $\Lambda(\theta')$ and a homology orientation on the cobordism, the moduli space $\mathfrak P(W,c;\theta, \theta')$ is canonically oriented. More precisely, these choices induce an orientation on the determinant line of the complex $T_\Lambda \mathcal{P} \xrightarrow{0} \mathcal H_\Lambda$; perturbing this so that the map is a nonzero $\Gamma_\Lambda$-equivariant map, the determinant line corresponds to the orientation of a generic zero set of a $\Gamma_\Lambda$-equivariant section of $\mathcal H_\Lambda$. 

Such a generic section has two zeroes, which are canonically oriented with the same orientation. We thus consider this an orientation on the set $\pi_0 \mathcal{P} = \mathfrak P$. 

\begin{remark}\label{rmk:signs}
Each element of $\mathfrak P(W,c;\theta, \theta')$ is oriented with the same sign, because the corresponding deformation complexes are homotopic. This sign is the same as the sign on the set of central connections $\mathfrak Z(W,c;\theta, \theta')$, for the same reason. This is mostly interesting when $b_1(W) = b^+(W) = 0$; when $b^+(W) > 0$ our standing assumption is that $(W,c)$ supports no central connections, and one can perturb away the abelian connections.
\end{remark}

Fix an orbit $\mathcal O_\Lambda \subset \mathcal P(W,c;\theta, \theta')$, and write $\lambda: \theta \to \theta'$ for the corresponding path of connections along $(W,c)$. When carrying out the obstructed gluing theory of Section 5.1, one finds that the ends of the moduli spaces $M^+_z(W, c; \alpha, \alpha')$ corresponding to gluing along $\mathcal O_\Lambda$ are modeled on the zero set of an appropriate section: 
\begin{itemize}
\item[(l)] If $\alpha' = \theta'$ and $z = \zeta * \lambda$, the corresponding end is modeled on the zero set of an $SO(3)$-equivariant section of the pullback bundle over $(0, \infty) \times \breve M^+_\zeta(Y, w; \alpha, \theta) \times \mathcal O_\Lambda$.  
\item[(r)] If $\alpha = \theta$ and $z = \lambda * \zeta'$, the corresponding end is modeled on the zero set of an $SO(3)$-equivariant section of the pullback bundle over $\mathcal O_\Lambda \times (-\infty, 0) \times \breve M^+_{\zeta'}(Y', w'; \theta', \alpha')$.
\item[(c)] If $\alpha \ne \theta$ and $\alpha' \ne \theta'$ and $z = \zeta * \lambda * \zeta'$, the corresponding end is modeled on the zero set of an $SO(3)$-equivariant section of the pullback bundle over
\[(0, \infty) \times \breve M^+_\zeta(Y, w; \alpha, \theta) \times \mathcal O_\Lambda \times (-\infty, 0) \times \breve M^+_{\zeta'}(Y', w'; \theta', \alpha').\]
\end{itemize}

Applying the modified moduli space ideas of Section 5.2, we interpolate between these sections and a pre-chosen section. In choosing these pre-chosen sections, we will want them to satisfy some special properties. 

For each pseudocentral reducible $\Lambda$, we choose privileged $SO(3)$-equivariant sections 
\begin{gather*}
\psi^{\text{priv}}_{\zeta \Lambda}: \breve M^+_\zeta(Y, w; \alpha, \theta) \times \mathcal O_\Lambda \to \mathcal H_\Lambda \\
\psi^{\text{priv}}_{\Lambda \zeta'}: \mathcal O_\Lambda \times \breve M^+_{\zeta'}(Y', w'; \theta', \alpha') \to \mathcal H_\Lambda \\
\psi^{\text{priv}}_{\zeta \Lambda \zeta'}: \breve M^+_\zeta(Y, w; \alpha, \theta) \times \mathcal O_\Lambda \times [0,1] \times \breve M^+_{\zeta'}(Y', w'; \theta', \alpha') \to \mathcal H_\Lambda 
\end{gather*}
which satisfy conditions precisely analagous to those enumerated in Section 5.2. In small degrees, we can understand the zero sets of these sections explicitly.

\begin{itemize}
\item When $\wt i(\zeta) = 1$, the first map has domain equivariantly diffeomorphic to a disjoint union of copies of $SO(3) \times S^2$; an $SO(3)$-equivariant section of this bundle may be identified with a section of the complex line bundle over $S^2$ of Euler class $2$. Equivariance implies that the projection map $\pi: \psi^{-1}(0) \to SO(3)$ is a local diffeomorphism, and $\psi^{-1}(0)$ consists of a number of copies of $SO(3)$. The Euler class computation shows that the signed count of these is two, so that $\pi: \psi^{-1}(0) \to SO(3)$ identifies $\psi^{-1}(0)$ as twice the fundamental class in $C_3^{gm} SO(3) = \Bbb Z$. 

\item When $\wt i(\zeta') = \dim \alpha + 1$, the second map has domain equivariantly diffeomorphic to a disjoint union of copies of $S^2 \times SO(3)$, and a similar discussion applies: $\psi^{-1}(0)$ can be identified with twice the fundamental class of $SO(3)$.
\end{itemize}

Now $\(\psi^{\text{priv}}\)^{-1}(0)$ takes a particularly simple form. In what follows, $\{\Lambda\}$ is oriented as the corresponding element of $\mathfrak P(W,c;\theta, \theta')$.

\begin{itemize}
\item In the first case, $\(\psi^{\text{priv}}_{\zeta \Lambda}\)^{-1}(0)$ has dimension $\wt i(\zeta) + \dim \alpha - 1$. As a geometric chain in $C_*(\alpha \times \theta')$, this is zero unless $\wt i(\zeta) = 1$, in which case this geometric chain coincides with $2\breve M_\zeta \times \{\Lambda\}$, as discussed above.
\item In the second case, $\(\psi^{\text{priv}}_{\Lambda \zeta'}\)^{-1}(0)$ has dimension $\wt i(\zeta') - 1$; this is always at least $\dim \alpha'$, and this chain is zero as an element of $C_*(\theta \times \alpha')$ unless $\wt i(\zeta') = 1 + \dim \alpha$, in which case it coincides with $-2\{\Lambda\} \times \breve M_{\zeta'}$ as above. The sign arises because $\{-\infty\}$ is a negatively-oriented boundary point in $[-\infty, 0)$.
\item In the third case, $\(\psi^{\text{priv}}_{\zeta \Lambda \zeta'}\)^{-1}(0)$ has dimension $\wt i(\zeta) + \wt i(\zeta') + \dim \alpha - 1$. The term $\wt i(\zeta) + \dim \alpha - 1$ is the dimension of $\breve M_\zeta(\alpha, \theta)$, which is always at least $\dim \alpha$; the term $\wt i(\zeta') - 1$ is the dimension of $\breve M_{\zeta'}(\theta', \alpha')$, which is always at least $\dim \alpha'$. It follows that 
\[\wt i(\zeta) + \wt i(\zeta') + \dim \alpha - 1 = (\wt i(\zeta) + \dim \alpha - 1) + (\wt i(\zeta') - 1) + 1 \ge \dim \alpha + \dim \alpha' + 1,\] 
and the chain complex $C_*(\alpha \times \alpha')$ is supported in degrees $* \le \dim \alpha + \dim \alpha'.$ So this chain is zero on the chain level. 
\end{itemize}

All this being said, one may follow the recipe of Section 5.2 to produce \emph{modified moduli spaces} associated to $(W, c, \pi_W)$: these will be moduli spaces $N^+_z(W, c; \alpha, \alpha')$ with the usual boundary terms together with, for each $\Lambda$, two oppositely-oriented boundary terms corresponding to the section of $\mathcal H_\Lambda$ arising from obstructed gluing theory, and additional boundary terms corresponding to the zero set of $\psi^{\text{priv}}_{\zeta \Lambda}$ or $\psi^{\text{priv}}_{\Lambda \zeta'}$ or $\psi^{\text{priv}}_{\zeta \Lambda \zeta'}$ as appropriate. 

\begin{prop}\label{prop:bdry-pseudo}
Let $N^+_z(W, c; \alpha, \alpha')$ be the modified moduli spaces constructed above. As geometric chains in $\alpha \times \alpha'$, we have the following equalities. 
\begin{itemize}
\item If $\alpha' = \theta'$ is central, then \begin{align*}
(-1)^{\dim \alpha} \partial N^+(W, c; \alpha, \theta') &= \sum_\beta \breve M^+(Y, w; \alpha, \beta) \times_\beta N^+(W, c; \beta, \theta') \\
&+ \sum_{\beta'} (-1)^{\wt i(W; \alpha, \beta')} N^+(W, c; \alpha, \beta') \times_{\beta'} \breve M^+(Y', w'; \beta', \theta') \\
&+ \sum_\theta 2\breve M^+(Y, w; \alpha, \theta) \times \mathfrak P(W,c;\theta, \theta').
\end{align*}
\item If $\alpha = \theta$ is central, then \begin{align*}
\partial N^+(W, c; \theta, \alpha') &= \sum_\beta \breve M^+(Y, w; \theta, \beta) \times_\beta N^+(W, c; \beta, \alpha') \\
&+ \sum_{\beta'} (-1)^{\wt i(W; \theta, \beta')} N^+(W, c; \theta, \beta') \times_{\beta'} \breve M^+(Y', w'; \beta', \alpha') \\
&+ \sum_{\theta'} -2\mathfrak P(W,c;\theta, \theta') \times \breve M^+(Y', w'; \theta', \alpha').
\end{align*}
\item If neither $\alpha$ nor $\alpha'$ is central, then \begin{align*}
(-1)^{\dim \alpha}\partial N^+(W, c; \alpha, \alpha') &= \sum_\beta \breve M^+(Y, w; \alpha, \beta) \times_\beta N^+(W, c; \beta, \alpha') \\
&+ \sum_{\beta'} (-1)^{\wt i(W; \alpha, \beta')} N^+(W, c; \alpha, \beta') \times_{\beta'} \breve M^+(Y', w'; \beta', \alpha').
\end{align*}
\end{itemize}
\end{prop}

Given this, we may define a bimodule $\mathcal I(W, c)$ associated to the cobordism $(W, c)$, together with the necessary additional data of a homology orientation and a path between the basepoints on the boundary components. We define this by the formula 
\begin{equation}\label{eq:bimod-pseudo}\mathcal I(W,c)(\alpha, \alpha') = N^+(W, c; \alpha, \alpha') + 2\mathfrak P(W,c;\alpha, \alpha').
\end{equation} 
Here we interpret the final term $\mathfrak P(W,c;\alpha, \alpha')$ as a discrete set of points which is empty unless $\alpha, \alpha'$ are both central, in which case it is the set of pseudocentral reducibles running from $\alpha$ to $\alpha'$, with the orientation discussed at the start of the section. Notice that if $\Lambda: \alpha \to \alpha'$ is pseudocentral, we have $\wt i(W; \alpha, \alpha') = 0$, so this has the correct expected dimension.

This gives an induced map between instanton chain complexes for a chosen pseudo-regular perturbation $\pi_W$ on $(W,c)$. When arguing that the induced map is well-defined, one may carry out essentially the same construction for the moduli spaces $M^+_z(W, c, I; \alpha, \alpha')$ for a 1-parameter family of metrics and perturbations: one needs to delete a neighborhood of the stratum $[0,1] \times \mathcal{P}(W,c)$ from all the compactified moduli spaces for 1-parameter families of metric and perturbation. 

Thus the boundary relations above (and the corresponding version in the case of 1-parameter families) establishes the following.

\begin{cor}
Whenever $(W, c): (Y,w,\pi) \to (Y',w',\pi')$ is a pseudo-unobstructed cobordism with $b_1(W) - b^+(W) = 0$, there is a natural map $\widetilde C(Y,w,\pi) \to \widetilde C(Y', w',\pi')$, independent of the choice of perturbation on $(W,c)$ up to equivariant chain homotopy.
\end{cor}

Before moving on to the remaining cases, we take this opportunity to point out that the central components of $\widetilde C(W,c)$ are determined by explicit cohomological information. Write $r_W: H^2(W;\Bbb Z) \to H^2(Y;\Bbb Z) \oplus H^2(Y';\Bbb Z)$. Recall here from Lemma \ref{lemma:enum-red-3d} that a central connection $\theta$ on $(Y,w)$ corresponds to a cohomology class $x \in H^2(Y;\Bbb Z)$ with $2x = PD(w)$, and similarly on $(Y', w')$. The statement of the following proposition uses the set of pairs $\tilde Z(W, c; \theta, \theta')$ defined as \[\big\{(y,y') \in H^2(W;\Bbb Z) \times H^2(W;\Bbb Z) \mid y + y' = PD(c), \;\; y-y' \text{ is torsion}, \;\; r_W(y) = r_W(y') = (x,x')\big\}.\]

\begin{prop}
Suppose $(W,c): (Y,w,\pi) \to (Y',w',\pi')$ is a pseudo-unobstructed cobordism with $b_1(W) = b^+(W) = 0$. If $\theta \in \mathfrak Z(Y,w)$ and $\theta' \in \mathfrak Z(Y',w')$ and $x,x'$ are the corresponding cohomology classes, the component of the map $\widetilde C(W,c)$ running from $C_*(\theta) = \Bbb Z$ to $C_*(\theta') = \Bbb Z$ is multiplication by $\pm \big|\tilde Z(W,c; \theta, \theta')|$.
\end{prop}
\begin{proof}
The only component of $N^+(W,c; \theta, \theta')$ which is 0-dimensional is the component of central connections. By the definition of $\mathcal I(W,c)(\theta, \theta')$ in \eqref{eq:bimod-pseudo}, the induced map $\theta \to \theta'$ is given by the signed count of $\mathfrak Z(W,c;\theta, \theta') + 2\mathfrak P(W,c;\theta, \theta')$. Finally, the sign on every point is the same by Remark \ref{rmk:signs}, so we just need to determine the number of points in this count.

Recall from Proposition \ref{prop:enum-red} that a central connection on $(W,c)$ extending $\theta$ and $\theta'$ corresponds to a cohomology class $y \in H^2(W;\Bbb Z)$ with $2y = PD(c)$ and $r_W(y) = (x, x')$, where $x$ and $x'$ are the cohomology classes corresponding to $\theta$ and $\theta'$. An abelian connection running between these is identified with a pair $\{y,y'\}$ with $y+y' = PD(c)$ and $y \ne y'$ for which $r_W(y) = r_W(y') = (x,x')$. Finally, Definition \ref{def:unobs} says that such an abelian connection is called pseudocentral precisely when $(y'-y)^2 = 0$. Because $b^+(W) = 0$, this is true if and only if $y' - y$ is torsion. 

$\widetilde Z$ supports an involution $\iota(y,y') = (y',y)$; the pseudocentral reducibles correspond to free orbits of $\iota$, while the central reducibles correspond to fixed points of $\iota$. It follows that \[\big|\tilde Z(W,c;\theta,\theta')\big| = \big|\mathfrak Z(W,c;\theta, \theta')\big| + 2\big|\mathfrak P(W,c;\theta, \theta')\big|.\qedhere\] 
\end{proof}

\begin{remark}
Consider the special case that $PD(c)$ is the sum of a torsion class and an even class which vanishes on the boundary: $PD(c) = 2a+t$ where $r_W(a) = 0$. Then the set $\tilde Z(W,c;\theta, \theta')$ is in bijection with $r_W^{-1}(x,x')$, via the map $z \mapsto (a+z,a+t-z)$. When $w,w',$ and $c$ are furthermore empty and $\theta, \theta'$ are trivial connections, this set is in bijection with \[\{x \in H^2(W;\mathbb Z) \mid x\text{ torsion, } r_W(x) = 0\},\] and is in particular guaranteed to be nonempty.
\end{remark}

We move on to verifying that the induced map is well-defined in the remaining cases. 
\begin{prop}\label{prop:c-is-minus-one}
Suppose $(W,c): (Y,w,\pi) \to (Y',w',\pi')$ is a pseudo-unobstructed cobordism with $b_1(W) - b^+(W) = -1$. For a pseudo-regular path of perturbations $\pi^t_W$, the modified moduli spaces $N_z^+(W,c,I;\alpha, \alpha')$ define an $SO(3)$-equivariant bimodule homotopy. 
\end{prop}
\begin{proof}
Recall that being pseudo-regular means here that the ASD operator associated to $\Lambda$ has trivial kernel and cokernel isomorphic to $\Bbb C^2$, which assembles into an $SO(3)$-equivariant rank 2 complex vector bundle $\mathcal H_\Lambda$ over $\mathcal O_\Lambda$. 

Given a pseudocentral $\Lambda: \theta \to \theta'$ in the family moduli space, consider a neighborhood of $\breve M_\zeta(\alpha, \theta) \times \mathcal O_\Lambda$ in $M_z^+(W,c,I;\alpha,\theta')$; here $\wt i(z) = \wt i(\zeta) + \wt i(\lambda) = \wt i(\zeta) - 3$. Truncating and modifying as usual, we see that the modified moduli spaces will have an additional boundary component corresponding to the zero set of a privileged section \[\breve M_\zeta(\alpha, \theta) \times \mathcal O_\Lambda \to \mathcal H_\Lambda.\] The zero set is a free $SO(3)$-space with dimension $\dim \alpha + \wt i(\zeta) - 3$, so can only be nonempty of this quantity exceeds $3$; it can only be non-trivial in $C_*(\alpha \times \theta')$ if this quantity is at most $\dim \alpha$. These conditions are only simultaneously possible when $\alpha$ is irreducible and $\wt i(\zeta) = 3$. 

In this case, one may restrict this $SO(3)$-equivariant section to a non-equivariant map $\breve M_\zeta^{\text{Fiber}}(\alpha, \theta) \times \mathcal O_\Lambda \to \mathcal H_\Lambda$, and conversely induce an $SO(3)$-equivariant section on the total space from a non-equivariant one on this subspace. Because this rank 2 complex vector bundle is pulled back from the second factor, there exists a nowhere-vanishing section (for instance, the composite of projection to the second factor and a nowhere-vanishing section $\mathcal O_\Lambda \to \mathcal H_\Lambda$, which exists for dimension reasons). Choose $\psi_{\zeta \Lambda}^{\text{priv}}$ to be such a section. Running the same argument on the other side, we may do the same for $\psi_{\Lambda \zeta'}$.

Then the boundary relations analogous to Proposition \ref{prop:bdry-pseudo} have no factors corresponding to factorization through $\Lambda$ in the dimension range where such factors could be nonzero, so the modified moduli spaces define a bimodule homotopy.
\end{proof}

This establishes that the induced map is well-defined up to equivariant homotopy in the case $b_1(W) - b^+(W) = -1$. The only remaining case in the pentachotomy is $b_1(W) - b^+(W) = 1$. This case is in some sense more irritating than the others, because in 1-parameter families there may exist pseudocentral reducibles which have both non-trivial kernel and cokernel. The preceding gluing analysis does not apply here. Fortunately, this also doesn't matter: the obstructed reducibles only appear in moduli spaces of large dimension.

\begin{lemma}
If $(W,c): (Y, w, \pi) \to (Y', w', \pi')$ is a pseudo-unobstructed cobordism with $b_1(W) - b^+(W) = 1$, and $M^+_z(W,c,I; \alpha, \alpha')$ is a family moduli space which contains a broken trajectory factoring through a pseudocentral reducible $\Lambda: \theta \to \theta'$, we have $\wt i(z) + 1 \ge \dim \alpha' + 4$, so that this moduli space has dimension at least $\dim \alpha + \dim \alpha' + 4$.
\end{lemma}

\begin{proof}
Supposing this moduli space contains the broken instanton $(\mathbf A, \Lambda, \mathbf A')$, write $z = \zeta * \lambda * \zeta'$, where $\zeta: \alpha \to \theta$ and $\zeta': \theta' \to \alpha'$ are the paths of connections along $Y$ (resp $Y'$) given by $\mathbf A$ (resp $\mathbf A'$); here $\zeta$ and $\zeta'$ may be constant if, for instance, this moduli space runs from $\theta$ to $\alpha'$ or $\alpha$ to $\theta'$. For any transversely cut out trajectory on a cylinder $\zeta: \alpha \to \theta$ we have $\wt i(\zeta) \ge 0$ (with $\wt i(\zeta) \ge 1$ if this is a nonconstant trajectory) and for $\zeta': \theta \to \alpha'$ we have $\wt i(\zeta') \ge \dim \alpha'$ because the fiber over $\theta'$ is the moduli space itself, which maps submersively onto $\alpha'$ (with this framed index being larger than $\dim \alpha' + 1$ if the trajectory is nonconstant). 

Now in general for a reducible $\Lambda: \theta \to \theta'$ we have \[\wt i(\lambda) = (b_1(W) - b^+(W) - \dim \Gamma_\Lambda) + N(\Lambda) + 3;\] the first term is the index of the deformation complex internal to the $SO(2)$-reducible locus, the second the index of the deformation complex normal to that locus, and the last term $\dim \Gamma_\theta$, this contribution arising from Definition \ref{def:index}. If $\Lambda$ is pseudocentral, $N(\Lambda) = 2(b_1(W) - b^+(W) - 1)$ and $\dim \Gamma_\Lambda = 1$. so that this reduces to $\wt i(\lambda) = 3(b_1(W) - b^+(W)) = 3$ in our case of interest. 
f
Thus in this case $\wt i(\lambda) = 3$. By additivity, we have $\wt i(z) = \wt i(\zeta) + \wt i(\lambda) + \wt i(\zeta') \ge 3 + \dim \alpha'$, which gives the inequality stated in the Lemma.
\end{proof}

The analogue of Lemma \ref{lemma:flowcat-from-objects} for bimodule homotopies shows that to define the relevant bimodule homotopy, we only need those moduli spaces $M^+_z(W,c,I;\alpha, \alpha')$ of dimension $\le \dim \alpha + \dim \alpha' + 1$. In this dimension range, all of our moduli spaces are stratified-smooth and satisfy the expected boundary relations. It follows that even in this case the induced map $\wt C(W,c,\pi_W)$ is independent of $\pi_W$ up to equivariant homotopy, which we have now established in all cases.

\begin{proof}[Proof of Theorem \ref{thm:exists-induced-map}]
Observe that when $(W_1, c_1)$ and $(W_2, c_2)$ are composable pseudo-unobstructed cobordisms, the composite is also pseudo-unobstructed, and the bimodules described above compose functorially up to equivariant homotopy (by a small variation on the standard argument). 

For each regular perturbation $\pi$ we have a well-defined equivariant instanton homology group $I^\bullet(Y, w; \pi)$. We write $\pi \le \pi'$ to mean that their signature data satisfy $\sigma_\pi \le \sigma_{\pi'}$. For each $\pi \le \pi'$ we have a continuation map $\varphi_{\pi \to \pi'}$ between these, which is an isomorphism by Theorem \ref{thm:main-theorem}, and these compose as you expect. The homology group $I^\bullet(Y, w)$ is defined by identifying all of these groups for varying $I^\bullet(Y, w; \sigma)$ using the continuation maps. 

To rephrase this, we have a preordered set $\Pi(Y,w)$, the set of regular perturbations $\pi$ on $(Y, w)$. When we have $\pi \le \pi'$, we can identify the associated instanton homology groups via the continuation maps. Finally, for any two perturbations, there is always a larger one, so we can identify them both by passing to that larger one. We will use this same recipe for cobordism maps.

For a cobordism, there is a set $\Pi(W,c) \subset \Pi(Y, w) \times \Pi(Y', w')$ of pairs $(\pi, \pi')$ so that $(W,c): (Y, w, \pi) \to (Y', w', \pi')$ is pseudo-unobstructed. This set is nonempty by Lemma \ref{everything-is-unobs}. For each such, we have defined a module map $I^\bullet(W,c)_{\pi \to \pi'}: I^\bullet(Y, w, \pi) \to I^\bullet(Y', w', \pi')$ respecting the exact triangles by applying the equivariant instanton homology functors to the map $\widetilde C(W,c)$ described above, which we have now verified are well-defined up to equivariant homotopy. We abbreviate this to $(W,c)_{\pi \to \pi'}$ for convenience.

To see that this gives a well-defined map $I^\bullet(Y,w) \to I^\bullet(Y',w')$, we need to show that for any two pairs $(\pi, \pi') \in \Pi(W, c)$, the maps $(W,c)_{\pi \to \pi'}$ are related by the relevant continuation maps.

By the index formula, if $(\pi_1, \pi'_1) \in \Pi(W,c)$ and we have $\pi_0 \le \pi_1$ and $\pi'_1 \le \pi'_2$, then $(\pi_0, \pi_2') \in \Pi(W, c)$ as well. By the functoriality mentioned at the beginning of this proof, we have 
\[\varphi_{\pi_1' \to \pi_2'} (W, c)_{\pi_1 \to \pi_1'} \varphi_{\pi_0 \to \pi_1} = (W,c)_{\pi_0 \to \pi_2'}.\]

Given any $(\pi_0, \pi_0')$ and $(\pi_1, \pi_1')$ in $\Pi(W,c)$, pick a perturbation $\pi_2$ with $\sigma_{\pi_2} = \text{min}(\sigma_{\pi_0}, \sigma_{\pi_1})$ and a perturbation $\pi_2'$ with $\sigma_{\pi_2'} = \text{max}(\sigma_{\pi_0}', \sigma_{\pi_1}')$. We have $(\sigma_2, \sigma') \in \Pi(W,c)$; the procedure above identifies both maps $(W,c)_{\pi_i \to \pi_i'}$ with $(W,c)_{\pi_2 \to \pi_2'}$, and hence with each other, using the various continuation maps. It is not much more work to show that the identification induced by any other pair $(\pi_3, \pi_3')$ with $\pi_3 \le \pi_2$ and $\pi_2' \le \pi_3'$ agrees with the identification induced by $(\pi_2, \pi_2')$.

Altogether, this shows that the cobordism maps given above assemble into a well-defined map between $I^\bullet(Y, w) \to I^\bullet(Y', w')$.
\end{proof}

\begin{remark}\label{rmk:inducedmap}
In practice, there are two ways to describe the induced map. If one feels free to choose the perturbations $\pi, \pi'$, then one can say the induced map $I^\bullet(W,c)$ is defined by choosing perturbations $(\pi, \pi')$ so that $(W,c)$ is pseudo-unobstructed, and defining this to be the map $I^\bullet(Y,w,\pi) \to I^\bullet(Y',w',\pi')$ given by applying the equivariant homology functor to the construction described earlier in this section.

On the other hand, if one wishes to define the induced map $I^\bullet(W,c)_{\pi_0 \to \pi_0'}: I^\bullet(Y,w,\pi_0) \to I^\bullet(Y',w',\pi_0')$ between perturbations where $(W,c)$ is not necessarily pseudo-unobstructed, instead one should choose new perturbations $\pi_1, \pi_1'$ with $(W,c): (Y,w,\pi_1) \to (Y',w',\pi_1')$ pseudo-unobstructed, and with $\sigma_{\pi_1} \le \sigma_{\pi_0}$ and $\sigma_{\pi_0'} \le \sigma_{\pi_1'}$; then the construction above defines \[I^\bullet(W,c)_{\pi_0 \to \pi_0'} = \varphi_{\pi_0' \to \pi_1'}^{-1} I^\bullet(W,c)_{\pi_1 \to \pi_1'} \varphi_{\pi_1 \to \pi_0}^{-1},\] using that the continuation maps are known to be isomorphisms on equivariant instanton homology by Theorem \ref{thm:main-theorem}.
\end{remark}

\subsection{Composites}\label{sec:composite}
Theorem \ref{thm:exists-induced-map} establishes that each cobordism $(W,c): (Y,w) \to (Y', w')$ between admissible pairs induces a well-defined cobordism map. What remains is to establish that these compose functorially. If $(W_1, c_1): (Y_0, w_0) \to (Y_1, w_1)$ and $(W_2, c_2): (Y_1, w_1) \to (Y_2, w_2)$ are cobordisms, write $(W_{12}, c_{12}): (Y_0, w_0) \to (Y_2, w_2)$ for their composite, obtained by pasting the two cobordisms along $(Y_1, w_1)$.

As mentioned at the beginning of the proof of Theorem \ref{thm:exists-induced-map}, if $\pi_i$ are regular perturbations on the respective admissible pairs $(Y_i, w_i)$ for which $(W_i, c_i)$ are pseudo-unobstructed, the composite bimodule $\mathcal I(W_2, c_2) \circ \mathcal I(W_1, c_1)$ is homotopic to $\mathcal I(W_{12}, c_{12})$. Unfortunately, it is not always possible to choose $\pi_i$ so that all three of these cobordisms are pseudo-unobstructed, even in favorable-seeming situations. For instance, suppose $b_1(W_i) = b^+(W_i) = 0$ for $i = 1, 2, 12$, and that $W_{12}$ supports a pseudocentral reducible $\Lambda: \theta_0 \to \theta_2$ whose restriction to $(Y_1, w_1)$ is an \emph{abelian} connection $\rho_1$. 

In this situation, the normal index is additive; writing $\Lambda_i$ for the restriction of $\Lambda$ to $W_i$, we see that $N(\Lambda_1) + N(\Lambda_2) = N(\Lambda) = -2$. It follows that \emph{no matter what perturbations we choose on the $Y_i$, one of $W_1$ or $W_2$ must be obstructed:} neither $\Lambda_i$ is a pseudocentral reducible, as it is abelian on one end, but one of them must have negative normal index. 

The best we can prove is the following.

\begin{lemma}\label{lemma:close-to-good}
Suppose $(W_i, c_i)$ are composable cobordisms with $b^+(W_i) - b_1(W_i) \ge 0$, and for which either $b^+(W_i) = 0$ or $(W_i,c)$ supports no central connections. Then there exist perturbations $\pi_0, \pi_1^\pm, \pi_2$ on $(Y_i, w_i)$ so that \[(W_1, c_1): (Y_0, w_0, \pi_0) \to (Y_1, w_1, \pi_1^+),\] \[(W_2, c_2): (Y_1, w_1, \pi_1^-) \to (Y_2, w_2, \pi_2),\] \[(W_{12}, c_{12}): (Y_0, w_0, \pi_0) \to (Y_2, w_2, \pi_2)\] are all pseudo-unobstructed, and $0 \le \sigma_{\pi_1^+} - \sigma_{\pi_1^-} \le 4$.
\end{lemma}

\begin{proof}
This is really a lemma about signature data functions, so we speak in terms of signature data functions $\sigma_i$ instead of perturbations $\pi_i$; choose one on each $(Y_i, w_i)$ to start with. Write $F_{Y, w} = \text{Map}(\mathfrak A(Y,w), 2\Bbb Z)$. Just as in the proof of Lemma \ref{everything-is-unobs}, we will consider the normal index of a reducible as a function 
\[N_i: \mathfrak A(W_i, c_i) \times F_{Y_{i-1}, w_{i-1}} \times F_{Y_i, w_i} \to 2\Bbb Z\] 
with 
\[N_i(\Lambda, f_{i-1}, f_i) = N_i(\Lambda, \sigma_{i-1}, \sigma_i) + \frac{f_i(r_i \Lambda) - f_{i-1}(r_{i-1} \Lambda)}{2}\]
for $i = 1,2$, and similarly for $i=12$. Here $r_i$ is restriction to $Y_i$. For convenience below, $\rho$ refers to an abelian flat connection on a 3-manifold, while $\theta$ to a central flat connection.

If $b_1(Y_1) > 0$ then $Y_1$ is admissible and the cobordisms $(W_i, c_i)$ support no reducible instantons; the relevant unobstructedness conditions hold tautologically for both $W_1$ and $W_2$. The composite $W_{12}$ need not be unobstructed, but Lemma \ref{everything-is-unobs} shows that one may choose the signature data on the ends appropriately so that it is. A similar discussion applies in the case that $Y_0$ or $Y_2$ has positive Betti number, so for the remainder of the proof we assume all $Y_i$ are rational homology spheres.

If $N_{12}$ is the corresponding normal index function for the composite cobordism, we have an additivity property (thanks, in part, to the assumption that $b_1(Y_1) = 0$): 
\[N_{12}(\Lambda, f_0, f_2) = N_1(\Lambda|_{W_1}, f_0, f_1) + N_2(\Lambda|_{W_2}, f_1, f_2) + c(\Lambda|_{Y_1}),\] 
where $c(\alpha) = 0$ if $\alpha$ is abelian and $c(\alpha) = 2$ if $\alpha$ is central. Our goal is to choose $f_0, f_1, f_2$ so that $N_1 \ge b_1(W_1) - b^+(W_1)$ and $N_{12} \ge b_1(W_{12}) - b^+(W_{12})$ while $N_2 \ge b_1(W_2) - b^+(W_2) - 2$; then the desired signature data will be given by $\sigma_0 + f_0$ and $\sigma_2 + f_2$, while $\sigma_1^+ = \sigma_1 + f_1$ and $\sigma_1^- = \sigma_1 + f_1 - 2$.

To ensure that thorese inequalities hold for every reducible running $\Lambda_1: \theta_0 \to \rho_1$ along the first cobordism and $\Lambda_2: \rho_1 \to \theta_2$ along the second, we are obligated to choose $f_1$ first. In this case, we have 
\[N_1(\Lambda_1, f_0, f_1) = N_1(\Lambda_1; \sigma_0, \sigma_1) + \frac{f_1(\rho_1)}{2} \quad \text{and} \quad N_2(\Lambda_2, f_1, f_2) = N_2(\Lambda_2; \sigma_1, \sigma_2) - \frac{f_1(\rho_1)}{2}.\] 
The desired inequalities on $N_1$ and $N_2$ are equivalent to
\[2N_2(\Lambda_2; \sigma_1, \sigma_2) - 2\big(b_1(W_2) - b^+(W_2)\big) + 4 \ge f_1(\rho_1) \ge 2\big(b_1(W_1) - b^+(W_1)\big) - 2N_1(\Lambda_1; \sigma_0, \sigma_1)\]
for all abelian flat connections $\rho_1$ and all reducible instantons $\Lambda_i$ which are asymptotic to $\rho_1$ at the appropriate end and a central connection on the other end.

Such a function exists if and only if \[2N_1(\Lambda_1; \sigma_0, \sigma_1) + 2N_2(\Lambda_2; \sigma_1, \sigma_2) - 2\big(b_1(W_1) + b_1(W_2) - b^+(W_1) - b^+(W_2)\big) + 4 \ge 0\] for all $\Lambda_1: \theta_0 \to \rho_1$ and $\Lambda_2: \rho_1 \to \theta_2$. 

By additivity of index and the fact that $Y_1$ is a rational homology sphere, this is equivalent to the inequality \[N_{12}(\Lambda_{12}; \sigma_0, \sigma_2) \ge b_1(W_{12}) - b^+(W_{12}) - 2.\] Finally, the index formula Proposition \ref{prop:normal-ind} implies $N_{12} \ge 2b_1 - 2b^+ - 2$, and the assumption that $b_1 - b^+ \ge 0$ implies that this is at least $b_1 - b^+ - 2$. Therefore one may choose $f_1$ so that the given inequality holds. (If $b_1 - b^+ \ge 2$, it follows that in fact choose $\pi_1, \pi_2, \pi_3$ so that $N_i(\Lambda_i) \ge 0$ for all $i$.) 

This establishes the desired inequality for all abelian instantons $\Lambda_1: \theta_0 \to \rho_1$ and $\Lambda_2: \rho_1 \to \theta_2$ with this choice of $f_1$; this choice is irrelevant to $W_{12}$.

This is the only difficult point; the rest of the argument follows as in Lemma \ref{everything-is-unobs}. In fact, with these choices of signature data, the only obstructed non-pseudocentral reducibles $\Lambda_2$ on $(W_2, c_2)$ run from an abelian connection to a central connection.
\end{proof}

The situation in the preceding Lemma is subtle, and causes difficulty in the proof of functoriality, which requires slightly more general versions of our previous constructions. 

\begin{theorem}
Suppose $(W_1, c_1): (Y_0, w_0) \to (Y_1, w_1)$ and $(W_2, c_2): (Y_1, w_1) \to (Y_2, c_2)$ are cobordisms between admissible pairs which have either $b^+(W_i) = 0$ or that $(W_i, c_i)$ supports no central connections. Then the same is true for the composite $(W_{12}, c_{12})$, and the cobordism maps defined in Theorem \ref{thm:exists-induced-map} satisfy \[I^\bullet(W_2, c_2) \circ I^\bullet(W_1, c_1) = I^\bullet(W_{12}, c_{12}).\]
\end{theorem}

\begin{proof}[Sketch of proof]
If $(W_{12}, c_{12})$ supports central connections, then by restriction so too do both $(W_i, c_i)$ (and neither end is admissible, so all boundary components are rational homology spheres); in this case, we must have $b^+(W_1) = b^+(W_2) = 0$ by assumption, but because $Y_1$ is a rational homology sphere this implies $b^+(W_{12}) = 0$ as well, proving the first claim about the composite. 

As for the cobordism maps, to take advantage of Lemma \ref{lemma:close-to-good} we want to reduce to the case $b_1(W_i) - b^+(W_i) \ge 0$ for both $i = 1,2$. First we handle the case where at least one of these is negative. By an argument along the lines of Lemma \ref{lemma:close-to-good}, we may choose perturbations $\pi_0, \pi_1, \pi_2$ so that one of $(W_i, c_i)$ is a pseudo-unobstructed cobordism and the other --- while not pseudo-unobstructed --- has $b_1(W_i) - b^+(W_i) < 0$. Because $b^+(W_i) < 0$, this implies that there exists a regular perturbation $\pi^W_i$ on this cobordism with no reducible connections whatsoever; the same is necessarily true for $\pi^W_1 \cup_L \pi^W_2$ when the length of the neck $L$ is sufficiently large. The stated result follows in this case by a standard neck-stretching argument which implies that the induced map on the composite is the same as the map induced by $\pi^W_1 \cup_L \pi^W_2$ for sufficiently large $L$.

So we assume $b_1(W_i) - b^+(W_i) \ge 0$ for both $i = 1,2$, and suppose that we have chosen perturbations $\pi_0, \pi_1^\pm, \pi_2$ as in the statement of Lemma \ref{lemma:close-to-good}. 

Now the situation is much more complicated. If we can establish the existence of the diagram below, as well as show that it is homotopy-commutative, we will have proved the desired theorem. In this diagram, we write for instance $\widetilde C(\pi_0)$ in place of the full notation $\widetilde C(Y_0, w_0, \pi_0; R)$ for the sake of space, and similarly we suppress the geometric representatives $c$. On the other hand, we include the perturbations on the ends as part of the notation for clarity.

\[\begin{tikzcd}
	&& {\widetilde C(\pi_1^+)} \\
	\\
	{\widetilde C(\pi_0)} && {S\widetilde C(\pi_1^-)} && {S \widetilde C(\pi_2)} \\
	&& {\widetilde C(\pi_1^-)} && {\widetilde C(\pi_2)} \\
	&& {} && {}
	\arrow["{\widetilde{\text{Cyl}}^{1^+ \to 1^-}}", from=1-3, to=3-3]
	\arrow["{W_1^{0 \to 1^+}}", from=3-1, to=1-3]
	\arrow["{\widetilde{W_1}^{0 \to 1^-}}"', from=3-1, to=3-3]
	\arrow["{(\Sigma W_2)^{1^- \to 2}}", from=3-3, to=3-5]
	\arrow["{\Sigma_1}", from=4-3, to=3-3]
	\arrow["{\Sigma_2}"', from=4-5, to=3-5]
	\arrow["{W_2^{1^- \to 2}}", from=4-3, to=4-5]
	\arrow["{W_{12}^{0 \to 2}}"', bend right=40, from=3-1, to=4-5]
\end{tikzcd}\]

A handful of notational points need to be explained here. The existence of a handful of these $SO(3)$-equivariant chain complexes and chain maps are slight generalizations of work carried out earlier in this article.

\begin{itemize}
\item Here $S\widetilde C(\pi_1^-)$ and $S \widetilde C(\pi_2)$ are defined by a variation on the blowup construction of Section \ref{cob-flow}, where we blow up simultaneously at \emph{every} abelian reducible. This variation can be constructed by hand or by iteratively applying the version of the construction where we blow up at a single abelian reducible.
\item The maps $\Sigma_1, \Sigma_2$ are defined by a variation on the map defined in Construction \ref{constr:susp-map}, with the difference being that $\Sigma_i(\rho, S_\rho) = S_\rho$ for \emph{every} abelian reducible $\rho$. Again, these induce an isomorphism on homology, and therefore equivariant homology.
\item The maps $\widetilde{\text{Cyl}}^{1^+ \to 1^-}$ and $\widetilde W_1^{0 \to 1^-}$ are defined by a variation on the construction of Section \ref{cob-bimods} of the bimodule of an obstructed cobordism, applied to the cylinder $\Bbb R \times Y$ and the cobordism $W_1$; this variation allows for more than one obstructed abelian connection $\Lambda$, which may run from a central connection to an abelian connection. A careful reading shows that the construction and proofs of facts regarding these maps go through with minimal change; the only crucial points are that all obstructed abelian connections run $\alpha \to \rho$ where $\rho$ is abelian (but $\alpha$ may be abelian or central), and each has normal index $-2$. 
\item The bimodule inducing the map $(\Sigma W_2)^{1^- \to 2}$ is defined by applying Construction \ref{constr:Wplus-bimod} to the bimodule $\Sigma_2 \circ W_2^{1^- \to 2}$; it can also be constructed by hand with little difficulty.
\end{itemize}

Now we assert that this diagram commutes up to homotopy. In particular, we assert the existence of equivariant homotopies between the following pairs of maps: 
\begin{align*}
(\Sigma W_2)^{1^- \to 2} \circ \Sigma_1 &\sim \Sigma_2 \circ W_2^{1^- \to 2} \\
\widetilde{\text{Cyl}}^{1^+ \to 1^-} \circ W_1^{0 \to 1^+} &\sim \widetilde{W_1}^{0 \to 1^-}\\
(\Sigma W_2)^{1^- \to 2} \circ \widetilde W_1^{0 \to 1^-} &\sim \Sigma_2 \circ W_{12}^{0 \to 2}.
\end{align*}

The first statement is tautological: these maps are in fact identically equal. 

The second statement follows by obstructed gluing theory in the spirit of Section \ref{obscob}, now applied to the obstructed cobordism given by the composite of $W_1$ and the obstructed cylinder $\Bbb R \times Y_1: (Y_1, \pi_1^+) \to (Y_1, \pi_1^-)$; the modified moduli space ideas apply with minimal change.

The final statement is more difficult, but is rather similar to the arguments carried out in Section \ref{sec:invt}. Again, this amounts to using obstructed gluing theory to construct the appropriate modified moduli spaces, and recognizing that their boundary relations establish a bimodule homotopy between the two composite bimodules represented by this relation.

Now recall our goal. Abusing notation, we write $W$ also for the induced map on equivariant homology $I^\bullet$; for consistency with previous sections we write $\varphi_{1^- \to 1^+}$ as the map induced by the unobstructed continuation map and $\widetilde \varphi_{1^+ \to 1^-}$ for the map induced by the obstructed continuation map. As discussed in Remark \ref{rmk:inducedmap}, the composite \[I^\bullet(W_2, c_2) \circ I^\bullet(W_1, c_1): I^\bullet(\pi_0) \to I^\bullet(\pi_2)\] is defined to be the composite \[W_2^{1^- \to 2} \circ H^\bullet(\Sigma_1)^{-1} \circ W_1^{0 \to 1^+},\] and we want to show that this is equal to the map $W_{12}^{0 \to 2}$. 

This follows purely formally using the identities above. Applying each of these in turn (occasionally rearranging using that $\Sigma_i$ and $\varphi$ are isomorphisms on equivariant homology), we find the desired equality.

%:
%\begin{align*}W_2^{1^- \to 2} \varphi_{1^- \to 1^+}^{-1} W_1^{0 \to 1^+} &= W_2^{1^- \to 2} \Sigma_1^{-1} \widetilde \varphi_{1^+ \to 1^-} W_1^{0 \to 1^+} \\
%&= W_2^{1^- \to 2} \Sigma_1^{-1} \widetilde{W_1}^{0 \to 1^-} \\
%&= \Sigma_2^{-1} (\Sigma W_2)^{1^- \to 2} \widetilde{W_1}^{0 \to 1^-} \\
%&= \Sigma_2^{-1} \Sigma_2 W_{12}^{0 \to 2} = W_{12}^{0 \to 2}.\qedhere
%\end{align*}
\end{proof}

Altogether, we have proven Theorem \ref{thm:Functor-w}(i): $I^\bullet$ is functorial on a category whose objects are admissible pairs $(Y, w)$ with a basepoint on $Y$, and whose morphisms are cobordisms $(W,c)$ equipped with a homology orientation and a path between the basepoints, where either $b^+(W) = 0$ or $(W,c)$ admits no central connections.

\subsection{The case $b_2^+ > 0$ via the blowup trick}\label{MM-blowup}
Thusfar we have established functoriality for negative-definite cobordisms and cobordisms which support no central connections. Our techniques do not apply to the case $b_2^+ > 0$ without change; now we have an obstructed \emph{central} connection on the cobordism, while we have explained only how to resolve (mildly) obstructed abelian connections.

It is likely that a modification of the previous techniques (a suspended flow category construction for central critical orbits instead) will allow us to extend functoriality to this setting, and this will be necessary for an investigation of the $\Bbb F_2$-theory. We leave this to future work, instead opting to use a trick which requires inverting $2$ in the coefficient ring. This is founded on an idea appearing in \cite{MM}: by taking the connected sum with a bundle over $\overline{\Bbb{CP}}^2$ with odd first Chern class, we can remove central connections and thus define induced maps associated to all cobordisms.

We briefly recall the relevant definitions and results. See \cite[Section 5.2]{DK} for a more careful discussion.

Let $(W, c)$ be a cobordism, and $\Sigma \subset W$ a closed oriented surface embedded in $W$; arbitrarily choose a spin structure on this surface. The index bundle of the dual Dirac operator $\slashed{\partial}^*_{A,\Sigma}$ paired to the restriction of $A \in \widetilde{\mathcal B}(W,c)$ to $\Sigma$ induces an $SO(3)$-equivariant complex line bundle $\widetilde L_\Sigma$ over the configuration space of framed connections on $W$. Given a reducible orbit $\mathcal O_\Lambda \subset \widetilde{\mathcal B}(W,c)$, the stabilizer $\Gamma_\Lambda$ acts on the fiber of $\widetilde L_\Sigma$ over $\Lambda$. If $\Lambda$ is central, this action is necessarily trivial, as every homomorphism $SO(3) \to SO(2)$ is trivial. If $\Lambda$ is abelian, respecting a parallel splitting $E_c \cong L_1 \oplus L_2$ into complex line bundles with first Chern classes $\{x,y\}$, the weight of the $SO(2) \cong \Gamma_\Lambda \subset SO(3)$ action is given by $|\langle x-y, [\Sigma]\rangle|$. The first Chern class of the associated complex line bundle over $\mathcal O_\Lambda$ is $-2|
\langle x-y, [\Sigma]\rangle|$.

In favorable situations (where there are few reducibles, or at least those reducibles that do exist appear with positive normal index) one may choose a generic $SO(3)$-equivariant section of $\widetilde L_\Sigma$ which is transverse to zero on all moduli spaces $M^+_z(W,c;\alpha, \alpha')$. In this situation, the cut-down moduli spaces again define a bimodule $\mathcal I(W,c, \Sigma)$ of relative degree two smaller than $\mathcal I(W,c)$.

We will apply this construction in a special case.

\begin{construction}
Suppose $(W,c)$ is a cobordism equipped with a base-path. Choose a point $p \in W$ disjoint from this base-path and the surface $c$. Take the blowup at $p$, producing the connected-sum $(W_\#, c_\#) = (W \# \overline{\Bbb{CP}}^2, c \cup \overline{\Bbb{CP}}^1)$. If $z$ is a component of connections on $W$, we write $z_\#$ for the corresponding component on $W_\#$ arising as connected sum with the component of connections on $(\overline{\Bbb{CP}}^2, \overline{\Bbb{CP}}^1)$ of index $-1$ (framed index $2$), so that $\wt i(z_\#) = \wt i(z) + 2$. We will be investigating the cut-down moduli spaces associated to the surface $\Sigma = \overline{\Bbb{CP}}^1$.

The reducibles on $(W_\#, c_\#)$ are labeled by pairs $\{x' = (x,i), y' = (y,j)\}$, where $x,y \in H^2(W;\Bbb Z)$ and $i,j$ are integers such that  $x + y = PD(c)$ and $i+j = 1$. Because it is never possible that $x' = y'$ (as $i+j$ is odd), this cobordism supports no central reducibles; if $(W,c)$ supports a central reducible, then the corresponding reducible of least index on $(W_\#, c_\#)$ is pseudo-central. If $\Lambda$ is an abelian reducible on $W$ with cohomology classes $\{x,y\}$, the reducible $\Lambda_i$ on $(W_\#,c_\#)$ with cohomology classes $\{(x,i), (y, 1-i)\}$ has normal index $N(\Lambda_i) = N(\Lambda) + 2(2i-1)^2 \ge N(\Lambda) + 2$. 

It follows that if $(W,c)$ is pseudo-unobstructed, then $(W_\#, c_\#)$ is unobstructed. When $b^+(W) > 0$, it follows that for a generic choice of perturbation $(W_\#, c_\#)$ supports no reducibles whatsoever. Choose such a perturbation, and then a section of $\widetilde L_\Sigma$ so that the moduli spaces $M^+_{z_\#}(W_\#, c_\#; \alpha, \alpha')$ are always transverse to zero. We write $M^+_{z_\#}(W_\#, c_\#, \overline{\Bbb{CP}}^1; \alpha, \alpha')$ for the moduli spaces cut-down by these sections. These define a bimodule $\mathcal I(W_\#,c_\#, \pi_\#)$. This bimodule is transparently well-defined up to equivariant homotopy if $b^+(W) > 1$.

%When $b^+(W) = 0$, one should first assume that $(W,c)$ is pseudo-unobstructed (in which case $(W_\#, c_\#)$ is as well). If $(W,c)$ supports central connections, it is not actually possible to choose a section of $\widetilde L_\Sigma$ which is transverse to all moduli spaces: if $\Theta$ is a central reducible on $(W,c)$ with cohomology classes $\{x,x\}$ and $\Theta_1$ is the corresponding abelian reducible with cohomology classes $\{(x,1), (x,0)\}$, we have $N(\Theta_1) = 0$. Thus the orbit $\mathcal O_{\Theta_1} \cong S^2$ is isolated in the moduli spaces of instantons, but the weight of $\widetilde L_\Sigma$ over this orbit is one. It is isomorphic to $TS^2$ with its standard orientation, and any $SO(3)$-equivariant section $S^2 \to TS^2$ is zero.

Well-definedness in the case $b^+(W) = 1$ is somewhat more subtle, but entirely analogous to Proposition \ref{prop:c-is-minus-one}; if $\Theta$ is a central connection on $(W,c)$, then in 1-parameter families of perturbations on $(W_\#, c_\#)$, the abelian connection $\Theta_0$ may appear as an isolated obstructed reducible of normal index $-2$ on $(W_\#, c_\#)$. Over such a reducible there is no transverse section of $\widetilde L_\Sigma \to \mathcal O_{\Theta_0}$; these play the role that the obstructed reducibles of normal index $-4$ played in Proposition \ref{prop:c-is-minus-one}, though the sign on the zero locus is opposite, as the associated complex line bundle has first Chern class $-2$ instead of $+2$. 

The same modified moduli space construction from Proposition \ref{prop:c-is-minus-one} implies that the resulting bimodule $\mathcal I(W_\#, c_\#, [\overline{\Bbb{CP}}^1], \pi_\#)$ is independent of $\pi_\#$ up to equivariant homotopy.

We write $I^\bullet(W_\#,c_\#,[\overline{\Bbb{CP}}^1])$ for the induced map on equivariant instanton homology.
\end{construction}

Before moving on to discussing functoriality, we should compare this cobordism map to the ones we've already defined.

\begin{prop}\label{same-map}
If $(W,c)$ is a pseudo-unobstructed cobordism with $b^+(W) > 0$, and $R$ is a field, then $-2I^\bullet(W,c;R) = I^\bullet(W_\#,c_\#,[\overline{\Bbb{CP}}^1];R)$.
\end{prop}

\begin{proof}[Sketch of proof.]
Choose a perturbation on $(W,c)$ with no reducibes; because $(W,c)$ is pseudo-unobstructed with $b^+(W) > 0$, it supports no central connections. Write $(W_L, c_L)$ for the connected sum $(W,c) \# (\overline{\Bbb{CP}}^2, \overline{\Bbb{CP}}^1)$ above with a neck of length $L$. We may perform this connected sum in such a way that when $L = \infty$, the resulting space is conformally equivalent to $W \setminus p \sqcup \overline{\Bbb{CP}}^2 \setminus \{q\}$, and therefore finite-energy on $(W_\infty, c_\infty)$ instantons can be identified with pairs of instantons on $(W,c)$ and instantons on $(\overline{\Bbb{CP}}^2, \overline{\Bbb{CP}}^1)$. The framed moduli space on $(\overline{\Bbb{CP}}^2, \overline{\Bbb{CP}}^1)$ has dimension $2$ modulo $8$, so the only relevant moduli space is the 2-dimensional one, which consists of a single abelian instanton. 

It follows that \[M^+_{z_\#}(W_\infty, c_\infty; \alpha, \alpha') = M^+_z(W, c; \alpha, \alpha') \times S^2.\] Furthermore, the complex line bundle $\widetilde L_\Sigma$ is determined entirely by the $\overline{\Bbb{CP}}^2$ factor, so is pulled back from (a negatively oriented) $TS^2$ over the latter factor. \emph{If we ignore equivariance}, we could choose our section to be a section $S^2 \to TS^2$ with exactly two negatively-oriented zeroes; this would identify $M^+_{z_\#}(W_\infty, c_\infty, \overline{\Bbb{CP}}^1; \alpha, \alpha')$ with $-2$ copies of $M^+_z(W,c; \alpha, \alpha')$. This observation, first appearing (with different orientation conventions) as \cite[Theorem 6.9]{Kotschick} and also appearing as \cite[Proposition 3.2]{MM}, says that we have not produced anything new. 

To give an identification of the bimodules $-2\mathcal I(W,c)$ and $\mathcal I(W_\#,c_\#, [\overline{\Bbb{CP}}^1])$, we will use the variation on geometric chains described in Remark \ref{rmk:disconn-equiv}, where one quotients by the relation of collapse-equivalent between arbitrary (possibly disconnected) chains, rather than quotienting connected chains by collapse-equivalence and then taking disjoint unions. We may use this variation on geometric chains so long as $R$ is a field. In this theory, it suffices to construct an equivariant section of $\widetilde L_\Sigma$ whose zero set is collapse-equivalent to $-2M^+_z(W,c;\alpha, \alpha')$ in a way preserving the endpoint maps. We construct such a section inductively, as follows, using the fact that $(W,c)$ supports no reducible instantons so that its moduli spaces are free over $SO(3)$. 

The construction is general. To start, let $\psi: \mathbb R^3 \times S^2 \to TS^2$ be the map given by $\psi(v, p) = \Pi_{T_p S^2}(v)$. If $A \in SO(3)$ we have $\psi_{Av}(Ap) = A_* \psi_v(p)$; in particular, $\psi_v$ is equivariant under the group of rotations which fix $v$, while $\psi_0 = 0$. 

Suppose $X$ is a free $SO(3)$-space. Choose an $SO(3)$-equivariant map $v: X \to \mathbb R^3$ which has zero as a regular value; that this is possible follows because $X$ is $SO(3)$-free. If such a map is already chosen on $\partial X$, we may extend it arbitrarily to the whole of $X$, so long as it remains $SO(3)$-equivariant and has zero as a regular value. Writing $L = \pi_2^* TS^2$ for the pullback to $X \times S^2$, we may define a section $\Psi: X \times S^2 \to L$ by $\Psi(x, p) = \psi_{v(x)}(p)$. Writing $X_{SL} = v^{-1}(0) \times S^2$, the projection map $\pi_1: \Psi^{-1}(0) \setminus X_{SL} \to X \setminus v^{-1}(0)$ is seen to be a double-cover: if $x \in X \setminus v^{-1}(0)$, then exactly two points in $\Psi^{-1}(0)$ map to $x$: the pairs $(x, \pm v(x)/\|v(x)\|)$. That this double-covering map is negatively-oriented follows because $TS^2$ has opposite the usual orientation. 

Therefore, $\Psi^{-1}(0)$ agrees with a negatively-oriented double cover of $X$ at the level of geometric chains, modulo the more general collapse-equivalence relation between disconnected chains. Finally, every prientation-reversing double-covering map is collapse-equivalent to two copies of the base: if $f: B \to \mathbb{RP}^\infty$ represents the double covering map, the double cover trivializes in the complement of the small locus $f^{-1}(\mathbb{RP}^{\infty-1})$, and hence every orientation-reversing double cover of $X$ is collapse-equivalent to $-2X$.

Now construct the map $v$ inductively over $M^+_z(W,c; \alpha, \alpha') \times S^2$, with the boundary values specified by inductive hypothesis.
\end{proof}

\begin{remark}
The statement above extends to the case $b^+(W) = 0$, but one needs to explain more clearly how to define the right-hand-side in the case that $b^+ = 0$ (using a setup analogous to that of Section \ref{badred}); we do not see substantial benefit to doing so. 

Further, this statement should hold with respect to all coefficient rings, not simply those in which $2$ is invertible. (In particular, over $R = \Bbb F_2$, the induced map $I^\bullet(W_\#, c_\#, [\overline{\Bbb{CP}}^1])$. One proof would define a finite-dimensional chain complex model using an appropriate cell structure on $SO(3)$. We will not carry it out here.
\end{remark}

This in hand, we can prove Theorem \ref{thm:Functor-w}(ii) from the guide: 

\begin{theorem}
When $R$ is a field with $\textup{char}(R) \ne 2$, then the functors $I^\bullet(Y, w;R)$ extend to the category $\mathsf{Cob}_{3+1}^w$ of weakly admissible pairs, allowing for indefinite cobordisms.
\end{theorem}

\begin{proof}
Let $(W,c): (Y, w, \pi) \to (Y, w, \pi')$ be a cobordism (with the necessary additional data). We define 
\[I^\bullet(W, c) = \frac{1}{2} I^\bullet_\#(W, c).\]

We have already verified in Proposition \ref{same-map} that this agrees with the definition of $I^\bullet(W,c)$ from Theorem \ref{thm:exists-induced-map}. That this map is functorial follows from the following chain of equalities. Write $W^{12} = W^1 \cup_{Y} W^2$, where $W^1: Y_1 \to Y$ and $W^2: Y \to Y_2$ are cobordisms (with suitable additional data, which is suppressed from notation for convenience). Then \begin{align*} \frac 12 I^\bullet(W^{12}_\#, \mathbf{E}^{12}_\#, \overline{\Bbb{CP}}^1; R) &= \frac 14 I^\bullet(W^{12}_\#, (\mathbf{E}^{12}_\#)_\#, \overline{\Bbb{CP}}_1^1 \sqcup \overline{\Bbb{CP}}_2^1; R) \\
&= \frac 14 I^\bullet(W^1_\#, \mathbf{E}^1_\#, \overline{\Bbb{CP}}^1_1) \circ I^\bullet(W^2_\#, \mathbf{E}^2_\#, \overline{\Bbb{CP}}^1_2) \\
&= \left(\frac 12 I^\bullet(W^1_\#, \mathbf{E}^1_\#, \overline{\Bbb{CP}}^1_1)\right) \circ \left(\frac 12 I^\bullet(W^2_\#, \mathbf{E}^2_\#, \overline{\Bbb{CP}}^1_2)\right).
\end{align*} 

Here the subscripts on $\overline{\Bbb{CP}}^1$ indicate which of the two blowups they correspond to. The second-to-last equality follows from the fact that instanton homology is functorial for unobstructed cobordisms by the usual neck-stretching argument. 
\end{proof}

\newpage

\section{Irreducible Floer homology and the Casson-Walker invariant}\label{sec:CW}
%!TEX root = equivariant-functoriality.tex

In Floer's original definition of the instanton homology of integer homology spheres, the invariant $I_*(Y)$ was obtained by ignoring the central connection and constructing a chain complex generated by the irreducible connections. 

One may define such a chain complex for rational homology spheres, but it is no longer independent of the choice of perturbation. Instead, we will see it depends on a choice of signature data function $\sigma$, and use the results of Section \ref{sec:invt} to see how it changes as we pass between adjacent signature data.

\subsection{Floer's complex from algebra}
Given an $SO(3)$-equivariant flow category $\mathcal C$ as in Section \ref{FlowCatS}, we can extract a chain complex called the \emph{irreducible flow complex}. By default, it is relatively graded; if $\mathcal C$ is given an absolute grading, then the irreducible flow complex inherits an absolute grading as well.

We say an object $\alpha$ of $\mathcal C$ is \emph{irreducible} if it is diffeomorphic to $SO(3)$. We define a complex
\[C_*^{\text{irr}}(\mathcal C) = \bigoplus_{\substack{\alpha \in \mathsf{Ob}(\mathcal C) \\ \alpha \text{ irreducible}}} \Bbb Z[\wt i(\alpha)] \times_{\Bbb Z/2} \Lambda(\alpha),\] 
where as in Section \ref{FlowCatS} the two-element set $\Lambda(\alpha)$ corresponds to orientation data, while $\wt i(\alpha)$ is the absolute grading if one exists, but otherwise is meant to indicate that the relative grading between the $\alpha$ and $\beta$ summands is $\wt i(\alpha, \beta)$.  

When $\widetilde i(\alpha, \beta) = 1$, the moduli space $\mathcal C(\alpha, \beta)$ has dimension $3 = \dim \alpha$, and the map $\mathcal C(\alpha, \beta) \to \alpha$ is an $SO(3)$-equivariant local diffeomorphism. The fiber above any point in $\alpha$ is canonically identified with any other, and we write $\mathcal C^{\text{irr}}(\alpha, \beta)$ for this fiber. When one chooses an element of $\Lambda(\alpha)$ and an element of $\Lambda(\beta)$, this fiber is given an orientation, which reverses upon changing either of those choices. To suppress discussion of these orientation elements, let us suppose we have prechosen such elements for each irreducible orbit.

We may thus define a differential on $C_*^{\text{irr}}(\mathcal C)$ via 
\[d^{\text{irr}}(\alpha) = \sum_\beta \# \mathcal C^{\text{irr}}(\alpha, \beta) \beta.\]

Instead of verifying that this squares to zero (which is straightforward), it is useful to identify the construction above with a purely algebraic avatar. Write $u = [SO(3)] \in C_3 SO(3)$ for the representative of the fundamental class. Consider the image of $u: CM_*(\mathcal C) \to CM_{*+3}(\mathcal C)$, which we write as $uCM_*(\mathcal C)$. This is a subcomplex of $CM_*(\mathcal C)$, because the Leibniz rule gives is $\tilde d(ux) = (du) x - u \tilde d(x)$ and $d(u) = 0$. Next, we will show that this complex is closely related to $C_*^{\text{irr}}(\mathcal C)$.

\begin{prop}\label{prop:Floer-from-alg}
The complex $\big(C_*^{\textup{irr}}(\mathcal C), d^{\textup{irr}}\big)$ is isomorphic to the complex $\big(u CM_*(\mathcal C)[-3], -\tilde d\;\big)$.
\end{prop}

\begin{proof}
First, let us recall some structural properties of the complex $CM_*(\mathcal C)$, and in particular of geometric chains. This complex has underlying relatively graded abelian group given by 
\[CM_*(\mathcal C) = \bigoplus_{\alpha \in \mathsf{Ob}(\mathcal C)} C_*(\alpha)[\wt i(\alpha)] \otimes_{\Bbb Z/2} \Bbb Z[\Lambda(\alpha)].\]
Because we use the truncated geometric chain complex defined in Appendix \ref{gm-trunc}, and this complex has $C_*(\alpha)$ supported in degrees $[0, \dim \alpha]$, we see that if $\alpha$ is diffeomorphic to a point or $S^2$, we must have $u C_*(\alpha) = 0$, while if $\alpha \cong SO(3)$, we have a noncanonical isomorphism $u C_*(\alpha) \cong \Bbb Z[3],$ with the fundamental class with respect to some orientation $[\alpha]$ being a cyclic generator. 

By choosing the class $u$ at the start we fixed an orientation on $SO(3)$; orient each irreducible $\alpha$ so that the map $SO(3) \to \alpha$ given by acting on a point is an orientation-preserving diffeomorphism. This makes the above isomorphism canonical, taking the cyclic generator to be the fundamental class with respect to the orientation described above; equivalently, by sending $u \cdot \text{pt}$ to $1 \in \Bbb Z$.

Thus, at the level of graded abelian groups, we have an identification 
\[C_*^{\text{irr}}(\mathcal C) = u CM_*(\mathcal C)[-3]\] 
sending $\alpha$ to $[\alpha]$ We will see this in fact their differentials coincide, so long as one interprets the differential on the right-hand-side as $-\tilde d$. To see this, observe that 
\[\tilde d [\alpha] = \sum_\beta (-1)^{\dim \alpha} [\alpha] \times_\alpha \mathcal C(\alpha, \beta) = \sum_{\beta \in \mathsf{Ob}^{\textup{irr}}(\mathcal C)} (-1)^{\dim \alpha} \mathcal C(\alpha, \beta);\] 
if $\wt i(\alpha, \beta) > 1$, the space $\mathcal C(\alpha, \beta)$ has dimension at least 4, so is zero in $C_*(\beta)$. Using also that $\dim \alpha = 3$, we may thus write this sum as 
\[\tilde d[\alpha] = \sum_{\substack{\beta \in \mathsf{Ob}^{\textup{irr}}(\mathcal C) \\ \wt i(\alpha, \beta) = 1}} -\mathcal C(\alpha, \beta).\] 
Now each $\mathcal C(\alpha, \beta)$ here is three-dimensional, so each component projects diffeomorphically onto $\beta$; these components are in one-to-one bijection with $\mathcal C^{\text{irr}}(\alpha, \beta)$, the points in a fiber above some $x \in \alpha$. 

To compare orientations, recall that $\mathcal C^{\text{irr}}(\alpha, \beta)$ is oriented as a fiber over $x \in \alpha$, so that 
\[\mathcal C(\alpha, \beta) \cong SO(3) \times \mathcal C^{\text{irr}}(\alpha, \beta)\]
is an oriented diffeomorphism, and the endpoint map $SO(3) \to \beta$ is equivariant, hence an oriented diffeomorphism by our initial choice of orientations. It follows that 
\[\tilde d[\alpha] = -\sum_\beta \# \mathcal C^{\text{irr}}(\alpha, \beta) [\beta].\qedhere\]
\end{proof}

Now observe that if $\tilde C$ is a dg $C_* SO(3)$-module, then $u \tilde C$ is again a chain complex; if $f: \tilde C \to \tilde C'$ is a $C_* SO(3)$-equivariant map, then because $f(ux) = u f(x)$ we see that $f$ restricts to a chain map $u \tilde C \to u \tilde C'$. Finally, if $H: \tilde C \to \tilde C'[1]$ is a $C_* SO(3)$-equivariant chain homotopy, this inducues a chain homotopy $u\tilde C \to u \tilde C'[1]$. 

\begin{cor}
Given a weakly admissible pair with choice of signature data $(Y, w, \sigma)$, there is a well-defined `irreducible instanton homology group' $I_*(Y, w, \sigma)$. If $(W, c): (Y, w, \sigma) \to (Y', w', \sigma')$ is pseudo-unobstructed, then there is a well-defined induced map
\[I(W,c): I_*(Y, w, \sigma) \to I_*(Y', w', \sigma'),\]
which composes functorially.
\end{cor}

\begin{proof}
As a corollary of Theorem \ref{unobs-cobmap}, whenever $\pi_0$ and $\pi_1$ are two regular perturbations on $(Y, w)$ with the same signature data, there is a $C_* SO(3)$-equivariant continuation map $\tilde C(Y, w, \pi_0) \to \tilde C(Y, w, \pi_1)$ which is an $SO(3)$-equivariant homotopy equivalence. Passing to the homology of $u \tilde C[-3]$, one may construct a canonically defined $I_*(Y, w, \sigma)$ which depends only on $\sigma$ by the mechanism of Definition \ref{def:I-no-pert}. 

We showed in Section \ref{badred} that when $(W, c): (Y, w, \pi) \to (Y', w', \pi')$ is pseudo-unobstructed, there is a corresponding $SO(3)$-equivariant chain map $\tilde C(Y, w, \pi) \to \tilde C(Y', w', \pi')$ which is well-defined up to equivariant chain homotopy and which composes functorially up to equivariant chain homotopy. Applying $u \tilde C[-3]$ and passing to homology, we get the claimed functoriality results.
\end{proof}

While the groups $I_*(Y, w, \sigma)$ do depend on the signature data $\sigma$, we can take advantage of the proof of Theorem \ref{thm:main-theorem}. If $\sigma_0 < \sigma_1$ are adjacent signature data on $(Y, w)$ which differ on $\rho \in \mathfrak A(Y, w)$, we proved that theorem by establishing the existence of an $SO(3)$-equivariant homotopy equivalence 
\[CM_*\mathcal S_\rho \mathcal I(Y, w, \sigma_0) \simeq CM_* \mathcal I(Y, w, \sigma_1).\]
Writing $S_\rho \tilde C(Y, w, \sigma_0)$ for the former complex and $\tilde C(Y, w, \sigma_1)$ for the latter, we thus see that 
\[I_*(Y, w, \sigma_1) \cong u\widetilde C(Y, w, \sigma_1)[-3] \cong u S_\rho \tilde C(Y, w, \sigma_0)[-3].\] 
To understand how $I_*$ changes in passing from $\sigma_0$ to $\sigma_1$, we merely need to understand the latter complex.

\begin{prop}\label{prop:cone}
Let $\mathcal C$ be an $SO(3)$-equivariant flow category of the type considered above (for instance, any framed instanton flow category). The Morse complex $C^{\textup{irr}}_*(\mathcal S_\rho \mathcal C)$ is obtained as an algebraic mapping cone. Precisely, there is an isomorphism of (relatively graded) abelian groups 
\[uCM_*(\mathcal S_\rho \mathcal C)[-3] = u CM_*(\mathcal C)[-3] \oplus \Bbb Z[\wt i_{\mathcal C}(\rho)-1]\] 
so that 
\[d^{\textup{irr}}_\rho = \begin{pmatrix} d^{\textup{irr}} & 0 \\ -e_2 & 0 \end{pmatrix},\] 
where $\langle e_2(\rho), \beta \rangle$ is the signed count of isolated irreducible components of $\mathcal C(\rho, \beta)$.
\end{prop}

\begin{proof}
We have 
\[u CM_*(\mathcal S_\rho \mathcal C) = \left(\bigoplus_{\alpha \in \mathsf{Ob}^{\textup{irr}}(\mathcal C)} \Bbb Z\langle \alpha\rangle[\wt i(\alpha)]\right) \oplus \Bbb Z\langle S_\rho\rangle[\wt i(\rho)-1];\]
we obtain the additional generator from $uC_*(S_\rho)[-3]$, as $\wt i_{\mathcal S}(S_\rho) = \wt i_{\mathcal C}(\rho) - 1$. This gives the stated isomorphism at the level of relatively graded abelian groups.

Now it remains to compute the differential, which we may write as 
\[d^{\text{irr}} = \begin{pmatrix} d_{CC} & d_{C\rho} \\ d_{\rho C} & d_{\rho \rho} \end{pmatrix}.\]
Each component is given by computing $\# \mathcal S^{\text{irr}}(\alpha, \beta)$ where $\alpha, \beta$ are irreducible objects of $\mathcal S_\rho \mathcal C$; equivalently, this is the signed count of the number of components. To enumerate these, recall the definition of the moduli spaces $\mathcal S(\alpha, \beta)$ from Construction \ref{constr:blowup-flowcat}. When neither $\alpha$ nor $\beta$ is $S_\rho$, we have an equality $\mathcal S(\alpha, \beta) = \mathcal C(\alpha, \beta)$, so that $d_{CC} = d^{\text{irr}}$ is the irreducible differential coming from $\mathcal C$ itself. 

Because $\mathcal S(\alpha, S_\rho) = \mathcal C(\alpha, \rho) \times_\rho S_\rho$, and $\mathcal C(\alpha, \rho)$ consists of irreducible trajectories, this space has dimension at least 4 and hence does not contribute to the differential; so $d_{C \rho} = 0$. Even more easily, $\mathcal S(S_\rho, S_\rho) = \varnothing$ by definition, so that $d_{\rho \rho} = 0$.

The only remaining interesting component is $\mathcal S(S_\rho, \beta) = -B(\rho, \beta)$, where $B(\rho, \beta)$ is the real blowup of the moduli space $\mathcal C(\rho, \beta)$ along a certain $SO(3)$-equivariant section. In the case we're interested in, $\mathcal C(\rho, \beta)$ is 3-dimensional --- a disjoint union of free orbits of $SO(3)$ --- this section may be assumed nonzero everywhere (define it to be a nonzero value at a point in each component, and extend the definition everywhere by equivariance). The components of $-B(\rho, \beta)$ are thus oriented diffeomorphic to the components of $-\mathcal C(\rho, \beta)$, and we defined $e_2(\rho) = \sum \# \mathcal C^{\text{irr}}(\rho, \beta).$ This shows that $d_{\rho C} = e_2$, as desired.
\end{proof}

Because $I_*(Y, w, \sigma_1)$ is isomorphic to the homology of the irreducible flow complex attached to $\mathcal S_\rho \mathcal I(Y, w, \sigma_1)$, the preceding mapping cone formula immediately implies the following.

\begin{cor}\label{cor:wall-crossing}
Suppose $\sigma_0 < \sigma_1$ are adjancent signature data on $(Y, w)$, differing at $\rho$. Write the continuation map as $\varphi: I_*(Y, w, \sigma_0) \to I_*(Y, w, \sigma_1)$, and $e_2: \Bbb Z[\wt i(\rho)] \to I_*(Y, w, \sigma_0)$ for the map sending $1 \mapsto [\sum_\alpha \# M(\rho, \alpha)]$. Then there is an exact triangle 
\[\cdots \to \Bbb Z[\wt i(\rho)] \xrightarrow{e_2} I_*(Y, w, \sigma_0) \xrightarrow{\varphi} I_*(Y, w, \sigma_1) \to \cdots\]
\end{cor}

Recall from Definition \ref{def:mod-2-grading} that $I_*(Y, w, \sigma)$ always carries a canonical $\Bbb Z/2$ grading. When $\rho$ is reducible, the absolute $\Bbb Z/2$-grading $\wt i(\rho)$ is always even; this can be seen by an application of the index formula to a 4-manifold with reducible connection bounding $(Y, w, \rho)$. 

It follows that $I_*(Y, w, \sigma)$ has a well-defined Euler characteristic, and as the signature data increases, the Euler characteristic decreases: 
\[\chi\big(I_*(Y, \sigma_0)\big) = \chi\big(I_*(Y, \sigma_1)\big) + 1.\]
In the final section, we will use this to define a Casson--Walker type invariant. Through the corollary below, it also shows that our technical passage to the suspended flow complex $S_\rho \widetilde C_*(Y, w, \sigma)$ was in some sense necessary; the continuation map is merely a quasi-isomorphism, and we couldn't have constructed an equivariant homotopy inverse.

\begin{cor}\label{cor:not-h-eq}
If $\sigma_0 \le \sigma_1$ are signature data on $(Y, w)$, the continuation map 
\[\varphi_{\sigma_0 \to \sigma_1}: \widetilde C(Y, w, \sigma_0) \to \widetilde C(Y, w, \sigma_1)\] 
is an $SO(3)$-equivariant homotopy equivalence if and only if $\sigma_0 = \sigma_1$, even though it is always an equivariant quasi-isomorphism.
\end{cor}
\begin{proof}
When $\sigma_0 = \sigma_1$ the cylinder is unobstructed running in either direction, and the corresponding cobordism maps give equivariant chain homotopy inverses. 

If $C$ is a dg-module over $C_* SO(3)$, and $f: C \to C'$ is an $SO(3)$-homotopy equivalence, then the corresponding map $f_*: uC \to uC'$ is a homotopy equivalence. So by Proposition \ref{prop:Floer-from-alg}, it follows that if $\varphi_{\sigma_0 \to \sigma_1}$ is an equivariant homotopy equivalence, then $\varphi_*: I_*(Y, w, \sigma_0) \to I_*(Y, w, \sigma_1)$ is an isomorphism. 

Now by Corollary \ref{cor:wall-crossing}, if $\sigma_0 \le \sigma_1$, then
\[\chi(I_*(Y, \sigma_0)) - \chi(I_*(Y, \sigma_1)) = \frac 14 \sum_{\alpha \in \mathfrak A(Y,w)} \sigma_1(\alpha) - \sigma_0(\alpha) \ge 0\] 
with equality if and only if $\sigma_0 = \sigma_1$. So if the continuation map $\varphi$ is an equivariant homotopy equivalence, we have $\sigma_0 = \sigma_1$.
\end{proof}

\subsection{An instanton Casson--Walker invariant}
If $\alpha \in \mathfrak A(Y, w)$, there is a corresponding representation $\text{ad}_\alpha: \pi_1(Y) \to SO(2)$. Recall from Definition \ref{def:perturbed-rho} that the perturbed rho invariant of $\alpha$ with respect to a regular perturbation $\pi$ is the rational number 
\[\rho_\pi(\alpha) = \sigma_\pi(\alpha) + \rho(\text{ad}_\alpha),\] 
with $\rho$ being the relative $\eta$ invariant of Atiyah--Patodi--Singer, and $\sigma_\pi$ is the signature data function of Definition \ref{def:sigdata}.

Notice two things. First, this function depends only on the signature data function $\sigma_\pi$, and does not otherwise depend on the perturbation. Second, if $\sigma_0 < \sigma_1$ are adjacent signature data differing precisely at $\alpha$, we have $\rho_{\sigma_1}(\alpha) = \rho_{\sigma_0}(\alpha) + 4$. We can thus use this as a correction term to counteract the change in Euler characteristic of $I_*$ as we move between chambers and obtain a well-defined rational number.

\begin{definition}
Let $(Y, w)$ be a weakly admissible pair. We say that the {\it instanton Casson--Walker invariant} is the number 
\[\lambda_I(Y,w) = \frac{1}{|\textup{Tors } H_1 Y|}\left(\chi \big(I_*(Y, w,\pi)\big) + \frac 14 \sum_{\alpha \in \mathfrak A(Y)} \rho_\pi(\alpha)\right),\]
for any choice of regular perturbation $\pi$. 
\end{definition}

This is independent of $\pi$ because the change in the displayed quantities cancel out as we move between adjacent signature data. We expect this to be an instanton-theoretic definition of the Casson--Walker invariant.

\begin{conjecture}\label{conj:I-eq-CW}
For any rational homology sphere $Y$, the instanton and Casson--Walker invariants coincide; we have $\lambda_I(Y, w) = \lambda_{CW}(Y)$ for all $w$.
\end{conjecture}
In particular, this would mean that $\lambda_I(Y, w)$ is independent of $w$. 

\begin{remark}
The Casson invariant and its relatives have a long history. Soon after Casson defined his invariant of integer homology spheres, Taubes \cite{Taubes} gave a definition in terms of signed counts of flat connections modulo gauge; the signs amount to saying that the Casson invariant computes the Euler characteristic of Floer's instanton homology. Around the same time, Walker published his generalization to rational homology spheres \cite{Walker}; a signed count of points no longer gives an invariant, and instead one must add canonical rational-valued correction terms (as we do above) to get a well-defined invariant. 

The Seiberg--Witten analogue of the Casson--Walker invariant is defined in \cite{Chen:SW}. There, even in the case of integer homology spheres (indeed, even in the case of $S^3$) one needs to introduce correction terms to get a well-defined invariant; these correction terms only fail to be integers for rational homology spheres. The resulting quantity is indeed equal to the Casson--Walker invariant, as proved in \cite{Lim:SW} for integer homology spheres and \cite{MWCW} for rational homology spheres.

Finally, a version of the Casson--Walker invariant in terms of a precursor to Heegaard Floer homology was provided in \cite{OS:casson}; here again one uses a topological correction term, and their Casson--Walker invariant is given as a sum over spin$^c$ structures.
\end{remark}

As yet we have neither a surgery exact triangle or a connected sum theorem, so we cannot say much presently outside of some specific examples; but we can at least verify it in those few examples. To set convention, the lens space $L(p,q)$ is $-p/q$ surgery on the unknot.

\begin{theorem}
Conjecture \ref{conj:I-eq-CW} is true when $Y$ is an integer homology sphere or a lens space $L(p,q)$. 
\end{theorem}
\begin{proof}
When $Y$ is an integer homology sphere, the sum over abelian flat connections is tautologically zero and the choice of signature data is vacuous, so we simply have $\lambda_I(Y) = \chi\left(I_*(Y)\right)$. The main result of \cite{Taubes} is that we have an equality $\chi\left(I_*(Y)\right) = 2\lambda_C(Y)$ with half of the Casson invariant of $Y$. Walker's normalization conventions have $\lambda_{CW}(Y) = 2\lambda_C(Y)$, giving us the desired result.

For lens spaces, the trivial perturbation is a regular perturbation, so we may take $\sigma = 0$. Then the invariant $\lambda_I(L(p,q), w)$ is given as \[\frac 1{4p} \sum_{\alpha \in \mathfrak A(L(p,q), w)} \rho(\alpha).\]
Write $H^2\big(L(p,q)\big) \cong \Bbb Z/p$. In the case $w = \varnothing$, the abelian reducibles can be identified with those integers $\ell$ with $0 < \ell < p/2$, and have adjoint representation given by $\text{ad}_\ell(i) = e^{2\pi i (2\ell)/p}$. In the case $p$ is even and $[w]$ is nonzero in $\Bbb Z/2$-homology, the abelian reducibles can be identified with the odd integers $0 < 2\ell+1 < p$, and $\text{ad}_\ell(i) = e^{2\pi i (2\ell+1)/p}$.

Let us handle the first case. Using the formula of \cite[Proposition 2.18]{HK}, we have 
\[\rho(\text{ad}_\ell) = \frac 4p \sum_{k=1}^{p-1} \cot\left(\frac{\pi k}{p}\right) \cot \left(\frac{\pi kq}{p}\right)  \sin^2 \left(\frac{2\pi k \ell}{p}\right).\]
Summing over $0 < \ell < p/2$ and multiplying by $1/4p$, we see that $\lambda_I(L(p,q), \varnothing)$ is 
\[\frac{1}{p^2} \sum_{0 < \ell < p/2} \sum_{k=1}^{p-1}  \cot\left(\frac{\pi k}{p}\right) \cot \left(\frac{\pi kq}{p}\right)  \sin^2 \left(\frac{2\pi k \ell}{p}\right).\]
Let us compute the identical sum 
\[\frac 1{2p^2} \sum_{\ell=0}^{p-1} \sum_{k=1}^{p-1} \cot\left(\frac{\pi k}{p}\right) \cot \left(\frac{\pi kq}{p}\right)  \sin^2 \left(\frac{2\pi k \ell}{p}\right);\] 
this coincides with the sum above because the $\sin^2$ term is symmetric with respect to swapping $\ell$ and $p-\ell$, and when $p$ is even the term corresponding to $\ell = p/2$ automatically vanishes. 

When $p$ is even, notice that every term corresponding to $k=p/2$ automatically vanishes because $\cot(\pi/2) = 0$. Further, observe that $\sum_{\ell=0}^{p-1} \sin^2\left(\frac{2\pi k \ell}{p}\right) = p/2$ whenever $0 < k < p$ and $k \ne p/2$: this follows by the double-angle formula $\sin^2(x) = \frac 12 - \frac 12 \cos(2x)$ and the fact that $\sum_{\ell=0}^{p-1} \cos(4\pi k \ell/p) = 0$, because this latter expression is the real part of the sum $0 = 1 + \zeta + \cdots + \zeta^{p-1}$. Here $\zeta = e^{4\pi i k/p} \ne 1$ because $0 < k < p$ and $k \ne p/2$. 

Combining these two facts, we may write 
\[\lambda_I(L(p,q), \varnothing) = \frac{1}{4p} \sum_{k=1}^{p-1}\cot\left(\frac{\pi k}{p}\right) \cot \left(\frac{\pi kq}{p}\right),\] 
which is a Dedekind sum traditionally denoted $s(q;p)$. 

Finally, we have $\lambda_{CW}(L(p,q)) = s(q;p)$, by \cite[Proposition 6.3]{Walker}; in comparing to the reference, note that Walker's $L_{p/q}$ is what we call $L(-p, q)$, and hence $\lambda_{CW}(L(p,q))$ is written in Walker's book as $\lambda(-L_{p/q}) = - \lambda(L_{p/q})$, as desired.

The computation for the case $[w] \ne 0 \in H_1(Y; \Bbb F_2)$ is similar, with the relevant $\rho$ invariant being 
\[\rho(\text{ad}_\ell) = \frac 4p \sum_{k=1}^{p-1} \cot\left(\frac{\pi k}{p}\right) \cot \left(\frac{\pi k q}{p}\right) \sin^2 \left(\frac{\pi k \ell}{p}\right),\]
for odd $0 < \ell < p$, so that 
\[\lambda_I(L(p,q), w) = \frac{1}{p^2} \sum_{j=0}^{p/2-1} \sum_{k=1}^{p-1} \cot\left(\frac{\pi k}{p}\right) \cot \left(\frac{\pi k q}{p}\right) \sin^2 \left(\frac{\pi k (2j+1)}{p}\right).\]
This simplifies to the same quantity $s(q;p)$. The argument is nearly identical: eliminate the terms with $k = p/2$, and sum over $j$ to get rid of the final factor; here one uses that when $\zeta$ is a $p$'th root of unity for $p$ even we have $\zeta + \zeta^3 + \cdots + \zeta^{2p-1} = 0$.
\end{proof}

\begin{remark}
Using the computations in \cite{Helle}, together with further computations of the rho invariants for abelian connections on spherical manifolds, one can also verify the conjecture for any 3-manifold of the form $SU(2)/\Gamma$ (with $w = \varnothing$), where $\Gamma$ is a finite subgroup of $SU(2)$. We will not carry out this computation in detail here. 
\end{remark}

\newpage

\begin{appendices}
\section{Stratified-smooth spaces}\label{A}
In this appendix, we construct the model for singular chains $C_*(X)$ we use in the main body of the text, where $X$ is a smooth manifold, which we will call \emph{geometric chains}; this complex should be functorial for smooth maps and have homology groups canonically isomorphic to singular homology. Let us briefly outline the desired construction. 

This geometric chain complex should be free on a certain generating set, the \emph{nondegenerate, nontrivial probes}, where these probes are something analogous to $(P, \phi)$ a compact oriented smooth manifold with corners $P$ with a smooth map $\phi$ to $X$; the boundary map should be analogous to the passage to $\partial P$ in the case of smooth manifolds with corners. 

For our purposes, these probes must be more general than smooth manifolds with corners: they must include the moduli spaces of instantons $\breve M^+_\zeta(\alpha, \beta)$ discussed in Section \ref{comp-moduli}; in Definition \ref{def-str-sm-man}, we named the kind of object we obtain a \emph{stratified-smooth manifold}. When proving the excision property, one needs to be able to take transverse intersections, which leads to the more general notion of \emph{stratified-smooth spaces} below, roughly the kind of space obtained by taking transverse intersections of stratified-smooth manifolds.

The first two sections of this appendix develop the bare minimum of foundations and orientation theory of stratified-smooth manifolds needed to define geometric chains. In Section 5 of the main text, we also need a certain blowup construction on these stratified-smooth spaces, which we set up in Appendix \ref{BlowupChain}, and carry out some necessary computations. In Appendix \ref{CheckAxioms} we construct the geometric chain complex and verify an appropriate version of the Eilenberg-Steenrod axioms, and therefore verify that the resulting homology groups are canonically isomorphic to singular homology. In the final section, Appendix \ref{gm-trunc}, we perform a truncation procedure on the geometric chain complex; this is used in Section \ref{geo-chain-inst-hom} to ensure we only need to construct moduli spaces within a certain dimension range.
\subsection{Foundations of stratified-smooth spaces}\label{app:str-sm}
Smooth manifolds with corners should be special cases of stratified-smooth spaces, so first we must decide on a preferred definition. There are a number of inequivalent definitions (see \cite{Joyce} for a survey). The more general definitions have more complicated boundary operators, which cannot be defined solely in terms of the codimension-1 strata. This is harder to work with effectively when defining and computing the boundary of a stratified-smooth space; we prefer to be more restrictive. 

With eyes towards the moduli spaces in the main text, we instead use a version of manifolds with corners and stratified-smooth spaces where one is stratified by a set of faces, which need not be connected, whose closures are all again manifolds with corners (resp. stratified-smooth spaces). To do so, we record the stratification with an appropriate kind of partially ordere set.

\begin{definition}
A \textup{graded poset} is a finite partially ordered set $(\Delta, <)$ equipped with a `dimension' function $d: \Delta \to \Bbb Z_{\geq 0}$ satisfying the following: if $s < t$ implies $d(s) < d(t)$, and if furthermore $\{r \mid s < r < t\} = \varnothing$ then $d(t) = d(s) + 1$.

We say that the largest integer in the image of $d$ is the \textup{dimension} of $\Delta$. 
\end{definition}

The simplest example of a graded poset, and one which will be very useful for us, is the cube poset $\{0,1\}^n$ for some natural number $n$. Writing elements as functions $s: \{1, \cdots, n\} \to \{0,1\}$, the ordering has $s \le t$ if $s(i) \le t(i)$ for all $i$, and the dimension function is $d(s) = \sum s(i)$.

For the objects of interest, the graded poset $\Delta$ will have the property that the sub-poset \[\Delta_{\geq s} = \{t \in \Delta \mid t \ge s\}\] is always isomorphic to a cube poset $\{0, 1\}^{\dim X - d(s)}$, with dimension function shifted up by $d(s)$. 

Partially ordered sets are relevant because the natural relation between the different strata of a stratified space is a partial ordering. Below, we use a definition of stratified space which explicitly records the relevant partially ordered set. 

\begin{definition}\label{def:stratspace}
A \textit{stratified space} $(X, \Delta)$ is a pair of a poset $\Delta$ and a topological space $X$ equipped with, for each $s \in \Delta$, a locally closed subset $X_s \subset X$; these subsets are required to satisfy the property $\overline{X_s} = \bigcup_{t \le s} X_t.$

An isomorphism $(X, \Delta_X) \cong (Y, \Delta_Y)$ of stratified spaces is a pair $(f, \bar f)$, where $\bar f: \Delta_X \to \Delta_Y$ is a poset isomorphism and $f$ is a homeomorphism with $f(X_s) = Y_{\bar f(s)}$.
\end{definition}

Next we record the naive version of a smooth version of stratified spaces, where one merely equips each stratum with a smooth structure.

\begin{definition}\label{NSS-def}
A \textit{naive stratified-smooth space} $(X, \Delta)$ of dimension $n$ is a pair of a graded poset $\Delta$ and a second-countable Hausdorff space $X$ stratified by $\Delta$, where each stratum $X_s$ is given the structure of a smooth manifold of dimension $d(s)$. 

A \textit{stratified-smooth map} $f: X \to M$ from $X$ to a smooth manifold $M$ is a continuous function for which each $f_s: X_s \to M$ is smooth; we say $m \in M$ is a regular value if it is a regular value for all maps $f_s: X_s \to M$.

An \textit{isomorphism} $(f, \bar f)$ between two naive stratified-smooth spaces $(X,\Delta_X)$ and $(Y, \Delta_Y)$ is an isomorphism of stratified spaces so that the induced maps between strata $f_s: X_s \to Y_{\bar f(s)}$ are all diffeomorphisms.
\end{definition}

Let us discuss one standard example of a naive stratified-smooth space, as well as two important constructions which will be used in the definition of the local models below.

\begin{example}
The space $[0,\infty)^n$ can be given the structure of a naive stratified-smooth space. The stratification is by the cube poset $\Delta = \{0,1\}^m$; the faces are products of $\{0\}$ and $(0, \infty)$. More precisely, write \[F(i) = \begin{cases} \{0\} & i = 0 \\ (0,\infty) & i = 1 \end{cases}\] Then the face associated to the stratum $s \in \{0,1\}^m$ is the subspace \[(0,\infty)^s := \prod_{i=1}^n F(s(i));\] that is, the $i$'th factor is $\{0\}$ if $s(i) = 0$ and is $(0,\infty)$ if $s(i) = 1$.
\end{example}

\begin{example}
Let $(X, \Delta)$ be a naive stratified-smooth space, and let $U \subset X$ be an open subset. Then $U$ also inherits the structure of a naive stratified-smooth space, taking $\Delta(U) \subset \Delta$ to be those $s$ for which $U \cap X_s \ne \varnothing$ with strata $U_s = U \cap X_s$.
\end{example}

\begin{example} 
Let $(X, \Delta)$ be a naive stratified-smooth space of dimension $n$, and let $\psi: X \to \Bbb R^\ell$ be a stratified-smooth map with zero as a regular value. Write $\Delta(\psi) \subset \Delta$ for the sub-poset \[\Delta(\psi) = \{s \in \Delta \mid \psi^{-1}(0) \cap X_s \ne \varnothing\},\] with grading function $d_\psi(s) = d(s) - \ell$. Then $(\psi^{-1}(0), \Delta(\psi))$ is a naive stratified-smooth space of dimension $n-\ell$, with $\psi^{-1}(0)_s = \psi_s^{-1}(0)$.
\end{example}

We combine these into one type of model space; a stratified-smooth space will be a naive stratified-smooth space locally isomorphic to one of the model spaces.

\begin{definition}\label{def:cutout-locus}
Let $V \subset [0,\infty)^m \times \Bbb R^n$ be an open subset containing zero, and let $\psi: V \to \Bbb R^\ell$ be a naive stratified-smooth map with zero as a regular value. Then $Z(\psi) = \psi^{-1}(0)$ carries the natural structure of a naive stratified-smooth space of dimension $m+n-\ell$. We call such a naive stratified-space $Z(\psi)$ a {\it cutout locus}.
\end{definition}
 
Notice that the assumption $0 \in V$ implies that $\Delta(V) = \{0,1\}^m$; every face has a point near $0$ and hence a point in $V$. Thus the stratification poset only changes when passing to $\Delta(\psi)$.

Before moving on to the definition of stratified-smooth space, we make a brief observation on the combinatorics of this poset $\Delta(\psi)$. The following argument is inspired by \cite[Section 2c]{SeSm}, which we follow again when discussing orientations later.

\begin{lemma}\label{lemma:upward-closed}
Suppose $\psi: V \to \mathbb R^\ell$ is a stratified-smooth map with $0$ as a regular value.  Then the poset $\Delta(\psi) \subset \{0,1\}^m$ is an upward-closed sub-poset; that is, if $s \in \Delta(\psi)$ and $t \in \{0,1\}^m$ has $t > s$, then $t \in \Delta(\psi)$ as well.
\end{lemma}

\begin{proof}
Suppose that $Z(\psi)_s \subset (0,\infty)^s \times \Bbb R^n$ is nonempty.\footnote{Recall that here $s \in \{0,1\}^m$ and the exponent means the $i$'th factor is either $\{0\}$ or $(0, \infty)$, depending on $s(i)$.} Because $Z(\psi)_s$ is the preimage of a regular value, it is a smooth manifold of codimension $\ell$ inside $(0,\infty)^s$. Choose a small sphere $c: S^{\ell - 1} \to \big((0,\infty)^s \times \Bbb R^n\big) \setminus Z(\psi)_s$ which has linking number one with $Z(\psi)_s$, so that in particular $\psi c: S^{\ell - 1} \to \Bbb R^\ell - \{0\}$ is not null-homotopic.

It follows that in fact $c: S^{\ell - 1} \to \big([0,\infty)^m \times \Bbb R^n\big) \setminus Z(\psi)$ is not null-homotopic (for otherwise $\psi c$ would be null-homotopic, which it is not).

Now for $s < t$, we may homotop $c$ into $c': S^{\ell - 1} \to (0,\infty)^t \times \Bbb R^n$ by an arbitrarily small homotopy; taking this homotopy sufficiently small, its image will not meet $Z(\psi)$, as $Z(\psi)$ is closed and the domain is compact. In particular, $\psi c'$ is homotopic to $\psi c$ as a map to $\Bbb R^n \setminus \{0\}$, and therefore $\psi c'$ is not null-homotopic; thus $c'$ is not null-homotopic in $\big((0,\infty)^t \times \Bbb R^n\big) \setminus Z(\psi)_t$. In particular, $Z(\psi)_t$ is nonempty, which is what we wanted to show.
\end{proof}

\begin{remark}\label{rmk:upward-closed}
It follows that $\Delta(\psi)$ is a \emph{cubical poset}, which means that for each $s \in \Delta(\psi)$, we have $\Delta(\psi)_{\ge s} \cong \{0,1\}^m_{\ge s} \cong \{0,1\}^{m - d(s)}$; the poset of strata lying above a given stratum is isomorphic to the cube poset of appropriate dimension. Soon after we define the boundary operator, we will see that this observation underlies the relation $\partial^2 = 0$.\end{remark}

We now define the geometric objects of interest.

\begin{definition}\label{def:strat-smooth}
A \emph{stratified-smooth space} is a naive stratified-smooth space $(X, \Delta)$ so that, for each $x \in X$, there exists an open subset $x \in U \subset X$ and a cutout locus $Z(\psi)$ so that $U$ is isomorphic to $Z(\psi)$ as a naive stratified-smooth space. We call a choice of open set $U$ and isomorphism $U \cong Z_\psi$ a {\it cutout chart} for $x$.
\end{definition}

Note that in particular we do not have an `atlas' of such charts in any sense; we merely demand that they exist.

\begin{example}
In \cite{Joyce}, a sequence of increasingly more strict notions of manifolds with corners are defined, and {\it $\langle N\rangle$-manifolds} defined there give examples of stratified-smooth spaces. Our Definition \ref{def-str-sm-man} of stratified-smooth manifold is very closely related to Joyce's definition $\langle N\rangle$-manifold, though ours is a weaker object, as the smooth structure on the top stratum of a stratified-smooth manifold does not extend up to the boundary.

In particular, the simplex $\Delta^k$ is a stratified-smooth space for all $k$, as are all products of simplices. 
\end{example}

\begin{example}
The definition of stratified-smooth space is more restrictive than it may seem at first, because of the combinatorial nature of an isomorphism between naive stratified-smooth spaces. This is best exemplified by the {\it teardrop}, defined in \cite[Remark 2.11]{Joyce} as an example of a naive topological manifold with corners which does not satisfy any of the more stringent definitions of manifold with corners. His teardrop is the set \[\{(x,y) \mid x \ge 0, \; y^2 \le x^2 - x^4\},\] which has three strata: a corner stratum, one boundary stratum, and one open stratum. The stratification poset $\Delta$ can be identified with $\{0, 1, 2\}$ with the total order $0 < 1 < 2$.

This is not locally isomorphic as a naive stratified-smooth space to $[0,\infty)^2$ near the corner point. A neighborhood of the corner point in the teardrop has \emph{one} (disconnected) codimension-1 stratum, while a neighborhood of $0$ in $[0,\infty)^2$ has two.

In fact, by Lemma \ref{lemma:upward-closed} above, the teardrop is not a stratified-smooth space: if it were, $\Delta_{\ge 0} = \{0, 1, 2\}$ would be isomorphic to a cubical poset, which it is not.

\begin{center}
\begin{figure}[h]
\centering
\includegraphics[width=6cm]{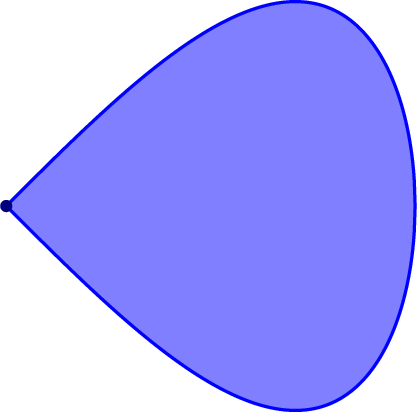}
\caption*{A picture of the teardrop with its three strata.} \qedhere
\end{figure}
\end{center}
\end{example}

Now that we have a definition of stratified-smooth space, we can define its boundary. Part of the reason we choose the combinatorial language in the definition of stratified-smooth space is that it makes the definition of boundary very simple, using the following lemma; notice that this lemma fails for the boundary stratum of the teardrop.

\begin{lemma}\label{lemma:bd-ok}
Let $X$ be a stratified-smooth space of dimension $n$, and let $X_s$ be a stratum of $X$. Then $\overline{X_s}$ naturally carries the structure of a stratified-smooth space of dimension $d(s)$.
\end{lemma}

\begin{proof}
The underlying topological space is of course the closure $\overline{X_s} \subset X$ with the subspace topology. Because $\overline{X_s} = \cup_{t \leq s} X_t$ we see that $\overline{X_s}$ carries the natural structure of a stratified space with smooth structure on each stratum, with face poset $\Delta_{\leq s}$.

We need to construct cutout charts for this subspace. If $x \in X_t \subset \overline{X_s}$, choose a cutout chart for $X$ near $x$. To spell out all this data: choose an open set $x \in U \subset X$, an open subset $V \subset [0,\infty)^m \times \Bbb R^n$, a naive stratified-smooth map $\psi: V \to \Bbb R^\ell$ with zero as a regular value, an isomorphism of posets $\overline \phi: \Delta(U) \to \Delta(\psi)$, and $\phi: U \cong Z_\psi$ an isomorphism of naive stratified-smooth spaces lifting $\overline \phi$. Then $\overline \phi(s) \in \{0,1\}^m$; we abbreviate this to $\overline s$.

Write $[0,\infty)^{\overline s}$ for the product of $\{0\}$ and $[0,\infty)$, with a $0$ in the $i$'th coordinate if $\overline s(i) = 0$ and $[0,\infty)$ otherwise. Write \[V_{\le \overline s} = V \cap \left([0,\infty)^{\overline s} \times \Bbb R^n\right),\] and $\psi_{\le \overline s}$ for the restriction of $\psi$ to this subspace. 

Then $\phi$ and $\overline \phi$ restrict to give an isomorphism of naive stratified-smooth spaces between $U \cap \overline{X_s} = U_{\le s}$ and $Z(\psi_{\le \overline s})$. This gives the desired cutout chart on $\overline{X_s}$.
\end{proof}

This is enough to define the boundary of a stratified-smooth space.

\begin{definition}
Let $X$ be a stratified-smooth space. The \textup{boundary} of $X$ is the stratified-smooth space \[\partial X = \bigsqcup_{\substack{s \in \Delta \\ d(s) = \dim X - 1}} \overline{X_s}.\qedhere\]
\end{definition}

Next, we will use Lemma \ref{lemma:upward-closed} to obtain an understanding of $\partial^2 X$.

\begin{prop}\label{del-sq}
The boundary $\partial^2 X$ is isomorphic to $\{0,1\} \times S$ for some stratified-smooth space $S$ of dimension $\dim X - 2$.
\end{prop}

\begin{proof}
First we show that for all $s \in \Delta$ with$d(s) = \dim X - 2$, we have $\Delta_{\ge s} \cong \{0,1\}^2$; in particular, for each such $s$, there are exactly two $t$ with $d(t) = \dim X + 1$ and $s < t$. This follows by picking a point $x \in X_s$ and a cutout chart $x \in U \subset X$ with $U \cong Z_\psi$. Then we have \[\Delta_{\ge s} = \Delta(U)_{\ge s} \cong \Delta(Z_\psi)_{\ge s} \cong \{0,1\}^2;\] the first isomorphism follows because if $s < t$ we have $X_t \cap U \ne \varnothing$ (because $X_s \subset \overline{X_t}$ and $U$ is open), the second by the definition of isomorphism of stratified-smooth spaces, and the last by Lemma \ref{lemma:upward-closed} and the discussion in Remark \ref{rmk:upward-closed}.

\noindent Next, observe that because \[\partial{\overline X_t} = \bigsqcup_{\substack{s \in \Delta \\ s < t \\ d(s) = \dim X - 2}} \overline{X_{s}},\] we have \[\partial^2 X = \bigsqcup_{\substack{t,s \in \Delta \\ s < t \\ d(s) + 1 = d(t) = \dim X - 1}} \overline{X_s}.\] Because there are exactly two such $t$ for each $s$ siwth $d(s) = \dim X - 2$, this disjoint union is isomorphic to \[\partial^2 X = \{0,1\} \times \bigsqcup_{\substack{s \in \Delta \\ d(s) = \dim X - 2}} \overline{X_s},\] as desired.
\end{proof}

We move on to foundations of smooth maps and fiber products. 

\begin{lemma}\label{lemma:strat-smooth}
Suppose $(X, \Delta)$ is a stratified space, and $f: X \to M$ a stratified-smooth map. Then given any $x \in X$ and any sufficiently small neighborhood $x \in U \subset X$, $f$ extends to a stratified-smooth map on the total space $V$ of any cutout locus $U \cong Z(\psi) \subset V$ around $x$. If $f$ is transverse to $S \subset M$, we may ensure this is true also for the extension to $V$.
\end{lemma}
\begin{proof}
Because this is a local question, we may assume $X = Z(\psi)$ and $M = \mathbb R^\ell$. Given a stratified-smooth map $f: Z(\psi) \to \mathbb R^\ell$, extend $f$ inductively to $V^{(k)} \cup Z(\psi)$ as a function which is overall continuous smooth on each stratum of $V^{(k)}$. To extend this to $V^{(k+1)}$ consistently with the values on $Z(\psi)^{(k)}$, first choose a continuous extension by the Tietze extension theorem, say $F_0$. On each open stratum of $V^{(k+1)}$, diffeomorphic to $\mathbb R^{k+1}$, we now have a continuous function $F_0: \mathbb R^{k+1} \to \mathbb R^\ell$ which restricts to $f$ along $Z(\psi)_k$. Choose a smooth function $F'$ which coincides with the smooth function $f$ along $Z(\psi)_k$ and for which $\|F'(x) - F(x)\|$ goes to zero as $\|x\| \to \infty$.  Redefining $F$ to be equal to $f$ on $V^{(k)} \cup Z(\psi)$ and $F'$ along the open faces of $V^{(k+1)}$, now $F'$ is continuous and smooth on each stratum, as desired.
\end{proof}

Then as desired the fiber product of two stratified-smooth spaces is again stratified-smooth. This includes, for instance, preimages of regular values of stratified-smooth functions to smooth manifolds (or more generally, intersection with submanifolds, hence in particular a generic zero set of a section of a vector bundle).

\begin{prop}\label{fibprod-ssmooth}
Suppose $f: X \to B$ and $g: Y \to B$ are stratified-smooth maps to the smooth manifold $B$ so that the restrictions to each pair of open strata $f_s: X_s \to B$ and $g_t: Y_t \to B$ are transvers. Then there is a natural stratified-smooth space structure on the fiber product $X \times_B Y$, whose underlying topological space is the usual fiber product. 
\end{prop}

\begin{proof}
The face poset of $X \times_B Y$ is a subset \[\Delta_B \subset \Delta \times \Delta',\] consisting of those pairs $(s, s')$ so that $X_s \times_B Y_{s'}$ is nonempty, equipped with the grading $d_B(s,s') = d(s) + d'(s') - \dim B$. That the ordering induced from $\Delta \times \Delta'$ coincides with the closure-inclusion ordering amounts to the relation $\overline{X_s \times_B Y_{s'}} = \overline{X_s} \times_B \overline{Y_{s'}}.$ The forward inclusion is obvious, while the reverse is not; still, it follows along lines analogous to the proof of Lemma \ref{lemma:upward-closed}.

Pick $(x,y) \in X_s \times_B Y_{s'}$. Write $k = \dim X - d(s)$ and $k' = \dim Y - d'(s')$; we have open sets \[V \subset [0,\infty)^k \times \Bbb R^{d(s) + m}, \;\;\; V' \subset [0,\infty)^{k'} \times \Bbb R^{d'(s') + n}\] and open neighborhoods \[x \in U \subset X, \;\;\; y \in U' \subset Y,\] and maps $\psi: V \to \Bbb R^m$ and $\psi': V' \to \Bbb R^n$, with isomorphisms of naive stratified-smooth spaces \[\phi: Z(\psi) \cong V, \;\;\; \phi': Z(\psi') \to V'\] sending $x$ and $y$ to zero, respectively.

Further, by Lemma \ref{lemma:strat-smooth}, we may assume $V$ and $V'$ are equipped with extensions $F$ and $G$ of $f$ and $g$, which are only defined on $Z(\psi)$ and $Z(\psi')$. We suppose these charts are chosen sufficiently small that the images of $F$ and $G$ lie in a single chart of $B$; for convenience, we thus suppose $B = \Bbb R^b$.

Then consider the open set \[\overline U = U \times U' \subset [0,\infty)^{k} \times \Bbb R^{d(s) + m} \times [0,\infty)^{k'} \times \Bbb R^{d'(s')+n},\] and the map $\xi = (\psi,F-G,\psi'): \overline U \to \Bbb R^{m + b + n}$. Because the restriction of $F-G$ to the zero set of $(\psi, \psi')$ is precisely $f-g$, the zero set of $\xi$ is the fiber product $Z(\psi) \times_{\Bbb R^b} Z(\zeta)$. Further, because $\psi$ and $\psi'$ have zero as a regular value, and $F-G$ has zero as a regular value on the zero locus of $\psi$ and $\zeta$, it follows that $\xi$ has zero as a regular value on each stratum as well.

Thus $U \times_B U' = Z(\xi)$, and $(\phi, \phi')$ restricts to give a homeomorphism $Z(\xi) \cong V \times_B V'$ which sends strata to strata and is a diffeomorphism on each stratum. We have constructed the necessary cutout chart.
\end{proof}

This fiber product operation is associative: if we have maps $f: X \to A$, $g: Y \to A \times B$ and $h: Z \to B$ satisfying the appropriate transversality requirements, then there is a canonical diffeomorphism of stratified-smooth spaces \[X \times_A (Y \times_B Z) \cong (X \times_A Y) \times_B Z.\] Further, using the description of its face poset in the previous argument, we may explicitly describe its boundary.

\begin{prop}\label{bd-fibprod}
The fiber product $X \times_B Y$ has boundary \[\partial(X \times_B Y) = (\partial X) \times_B Y \sqcup X \times_B \partial Y.\]
\end{prop}

\begin{proof}
The left-hand side is \[\bigsqcup_{\substack{s \in \Delta_X, t \in \Delta_Y \\ d(s) + d'(t) - \dim B = \dim X + \dim Y - \dim B - 1}} \overline{X_s} \times_B \overline{Y_t};\] now $d(s) + d'(t) = \dim X + \dim Y - 1$ if and only if one of $\dim X - d(s)$ and $\dim Y - d'(t)$ is one and the other is zero. The disjoint union simplifies to \[\left(\bigsqcup_{\substack{s \in \Delta_X, t \in \Delta_Y \\ d(s) = \dim X \\ d(t) = \dim Y - 1}} \overline{X_s} \times_B \overline{Y_t}\right) \sqcup \left(\bigsqcup_{\substack{s \in \Delta_X, t \in \Delta_Y \\ d(s) = \dim X - 1 \\ d(t) = \dim Y}} \overline{X_s} \times_B \overline{Y_t}\right).\] Because the disjoint union $\bigsqcup_{d(s) = \dim X} \overline{X_s}$ is $X$ itself, while $\bigsqcup_{d(s) = \dim X - 1} \overline{X_s} = \partial X$, the above sum simplifies to \[\partial X \times_B Y \sqcup X \times_B \partial Y,\] as desired. 
\end{proof}

In addition to stratified-smooth spaces defined by equations (such as the general fiber products described above), we will also find need in Section \ref{subsec:modified-mod-space} for stratified-smooth spaces defined by \emph{inequalities}; we record the relevant statement below.

\begin{prop}\label{trunc-still-smooth}
If $(X,\Delta)$ is a stratified-smooth space and $L: X \to \Bbb R$ is stratified-smooth, then for a regular value $c$ of $L$, the set $X^{\le c} = L^{-1}(-\infty, c]$ is a stratified-smooth space with boundary \[\partial X^{\le c} = L^{-1}(c) \sqcup (\partial X)^{\le c}.\]
\end{prop}

\begin{proof}
The stratification poset for $X^{\le c}$ is denoted $\Delta^{\le c}$, and is a sub-poset of $\Delta \times \{-1, 0\}$; we say $(s, -1) \in \Delta^{\le c}$ if $L^{-1}(c) \cap X_s$ is nonempty, while $(s, 0) \in \Delta^{\le c}$ if $L^{-1}(-\infty, c) \cap X_s \ne \varnothing$. The notation for $\Delta^{\le c}$ chosen so that $d(s,0) = d(s)$ and $d(s,-1) = d(s)-1$. Once we verify that $X^{\le c}$ is indeed a stratified-smooth space, this gives the stated boundary formula.

It is straightforward to see that this is a stratified topological space; the strata $(X^{\le c})_{(s,0)}$ are open subsets of $X_s$, hence carry the natural structure of smooth manifolds, whereas the strata $(X^{\le c})_{(s,-1)}$ are given as regular values of a smooth map $L: X_s \to \Bbb R$, hence also in a natural way smooth manifolds. 

The existence of cutout charts for $x \in X^{\le c}_{(s,0)}$ follows immediately by restricting cutout charts for $x \in X_s$. The interesting case is $x \in X^{\le c}_{(s,-1)}$. In this case, one starts by picking a cutout chart for $x \in X_s$; that is, there is a choice of open set $V \subset [0,\infty)^m \times \Bbb R^n$, a stratumwise-smooth map $\psi: V \to \Bbb R^\ell$, and an isomorphism $Z(\psi) \cong U$ of naive stratified-smooth spaces, where $U$ is an open neighborhood of $x$. By the assumption that $L$ is stratified-smooth, we may assume $L$ extends to a function on $V$, which we denote by the same symbol.

We may extend this to a map $\hat \psi: V \times [0,\infty) \to \Bbb R \oplus \Bbb R^\ell$ defined by \[\hat \psi(v,s) = (\psi(v),L(v) + s - c).\] For $(s,v)$ to be in the zero set of this map means precisely that $\psi(v) = 0$ and that $L(v) \le c$. The map $\pi: Z(\hat \psi) \to Z(\psi)$ is a topological embedding, and a smooth embedding on each stratum; composing with the map $\varphi: Z(\psi) \cong U$ gives us the desired stratified-smooth isomorphism from $Z(\hat \psi)$ to $U_{\le c}$.
\end{proof}
\subsection{Orientations on stratified-smooth spaces}\label{SS-or}
%!TEX root = equivariant-functoriality.tex

In this section, we will define and analyze the notion of an \emph{oriented} stratified-smooth space. Briefly, this is an orientation on the top strata which is consistent with the orientations induced by some cutout chart. We will also use cutout charts to define the boundary orientation. 

The first of these requires that we first understand the orientation on a cutout chart; for the second, we must verify that the induced boundary orientation does not depend on the choice of cutout chart.

\begin{construction}
To set conventions: suppose $f: X \to B$ is a smooth map between oriented smooth manifolds and $b \in B$ is a regular value. Given $x \in f^{-1}(b)$, we write $T_x Z$ for the tangent space to $f^{-1}(b)$; write $di: T_x Z \to T_x X$ for the differential of the inclusion map. If $s$ is a splitting of the surjective map $df_x: T_x X \to T_b B$, then we orient $T_x Z$ by the assumption that $(di, s): T_x Z \oplus T_b B \to T_x X$ is an orientation-preserving isomorphism.

Now let $\psi: [0,\infty) \times \Bbb R^n \to \Bbb R^m$ be a continuous map which, on each stratum, is smooth and has zero as a regular value. Write $\psi^\circ$ for its restriction to $(0,\infty) \times \Bbb R^n$, and write $\psi^\partial$ for its restriction to $\Bbb R^n \cong \{0\} \times \Bbb R^n$, equipped with the boundary orientation (here, negative the standard orientation). 

Then the zero sets $Z(\psi^\circ)$ and $Z(\psi^\partial)$ are canonically oriented, being regular level sets of smooth maps between oriented smooth manifolds.
\end{construction}

For convenience in what follows below, it is worth mildly generalizing the definition of cutout chart. Instead of frontloading the half-line coordinates and demanding that the cutout function $\psi$ be defined on $V \subset [0,\infty)^m \times \Bbb R^n$ an open subset, we allow $V$ to be an open subset of any product of copies of $[0,\infty)$ and $\Bbb R$ in any order. This is convenient primarily for easy discussion of products where we would otherwise have to shuffle coordinates to move the half-line coordinates to the front, introducing a sign. When discussing the orientation on the boundary of $Z(\psi)$, one first passes to the boundary (with its standard boundary orientation) and then passes to the regular zero set of $\psi^\partial$ on the corresponding oriented boundary strata.

This in hand, we can define an oriented stratified-smooth space.

\begin{definition}
Let $(X, \Delta)$ be a stratified-smooth space. An \emph{orientation} of $(X, \Delta)$ is a choice of orientation for each top stratum so that for each $x \in X$, there \emph{exists} a cutout chart $(U, V, \psi, \phi)$ around $x$ so that the induced diffeomorphism on the top stratum, $\phi: U^\circ \to Z(\psi)^\circ$, is orientation-preserving.
\end{definition}

In understanding this definition, a useful example to keep in mind is the $\theta$-graph: there is a continuous function on $[0,1]^2$ (which, on each stratum, is smooth and has zero as a regular value) whose zero set on the boundary is $\{1/2\} \times \{0, 1\}$  but on the interior is three arcs running between these points. An orientation on this stratified-smooth space orients two of the arcs in one direction and one arc in the opposite direction. For a given cutout chart, only one of these six possible orientations will be oriented compatibly with that cutout chart, but the definition only asks that \emph{some} cutout chart be oriented compatibly.

\begin{figure}[h]
\centering
{\color{blue}\[\begin{tikzcd}
	\bullet &&& \bullet
	\arrow["{>}"{marking}, no head, start anchor={[xshift=-1.5ex, yshift=-0.4ex]}, end anchor={[xshift=1.5ex, yshift=-0.4ex]}, bend left=20, from=1-1, to=1-4]
	\arrow["{>}"{marking}, no head, start anchor={[xshift=-1.5ex, yshift=0.4ex]}, end anchor={[xshift=1.5ex, yshift=0.4ex]}, bend right=20, from=1-1, to=1-4]
	\arrow["{<}"{marking}, no head, start anchor={[xshift=-1.5ex]}, end anchor={[xshift=1.5ex]}, from=1-1, to=1-4]
\end{tikzcd} \quad \quad \quad \begin{tikzcd}
	\bullet &&& \bullet
	\arrow["{>}"{marking}, no head, start anchor={[xshift=-1.5ex, yshift=-0.4ex]}, end anchor={[xshift=1.5ex, yshift=-0.4ex]}, bend left=20, from=1-1, to=1-4]
	\arrow["{<}"{marking}, no head, start anchor={[xshift=-1.5ex, yshift=0.4ex]}, end anchor={[xshift=1.5ex, yshift=0.4ex]}, bend right=20, from=1-1, to=1-4]
	\arrow["{>}"{marking}, no head, start anchor={[xshift=-1.5ex]}, end anchor={[xshift=1.5ex]}, from=1-1, to=1-4]
\end{tikzcd}\quad \quad \quad \begin{tikzcd}
	\bullet &&& \bullet
	\arrow["{<}"{marking}, no head, start anchor={[xshift=-1.5ex, yshift=-0.4ex]}, end anchor={[xshift=1.5ex, yshift=-0.4ex]}, bend left=20, from=1-1, to=1-4]
	\arrow["{>}"{marking}, no head, start anchor={[xshift=-1.5ex, yshift=0.4ex]}, end anchor={[xshift=1.5ex, yshift=0.4ex]}, bend right=20, from=1-1, to=1-4]
	\arrow["{>}"{marking}, no head, start anchor={[xshift=-1.5ex]}, end anchor={[xshift=1.5ex]}, from=1-1, to=1-4]
\end{tikzcd} \]
\[\begin{tikzcd}
	\bullet &&& \bullet
	\arrow["{<}"{marking}, no head, start anchor={[xshift=-1.5ex, yshift=-0.4ex]}, end anchor={[xshift=1.5ex, yshift=-0.4ex]}, bend left=20, from=1-1, to=1-4]
	\arrow["{<}"{marking}, no head, start anchor={[xshift=-1.5ex, yshift=0.4ex]}, end anchor={[xshift=1.5ex, yshift=0.4ex]}, bend right=20, from=1-1, to=1-4]
	\arrow["{>}"{marking}, no head, start anchor={[xshift=-1.5ex]}, end anchor={[xshift=1.5ex]}, from=1-1, to=1-4]
\end{tikzcd} \quad \quad \quad \begin{tikzcd}
	\bullet &&& \bullet
	\arrow["{<}"{marking}, no head, start anchor={[xshift=-1.5ex, yshift=-0.4ex]}, end anchor={[xshift=1.5ex, yshift=-0.4ex]}, bend left=20, from=1-1, to=1-4]
	\arrow["{>}"{marking}, no head, start anchor={[xshift=-1.5ex, yshift=0.4ex]}, end anchor={[xshift=1.5ex, yshift=0.4ex]}, bend right=20, from=1-1, to=1-4]
	\arrow["{<}"{marking}, no head, start anchor={[xshift=-1.5ex]}, end anchor={[xshift=1.5ex]}, from=1-1, to=1-4]
\end{tikzcd}\quad \quad \quad \begin{tikzcd}
	\bullet &&& \bullet
	\arrow["{>}"{marking}, no head, start anchor={[xshift=-1.5ex, yshift=-0.4ex]}, end anchor={[xshift=1.5ex, yshift=-0.4ex]}, bend left=20, from=1-1, to=1-4]
	\arrow["{<}"{marking}, no head, start anchor={[xshift=-1.5ex, yshift=0.4ex]}, end anchor={[xshift=1.5ex, yshift=0.4ex]}, bend right=20, from=1-1, to=1-4]
	\arrow["{<}"{marking}, no head, start anchor={[xshift=-1.5ex]}, end anchor={[xshift=1.5ex]}, from=1-1, to=1-4]
\end{tikzcd} \]
}\\
\caption*{The six ways to make the $\theta$-graph into an oriented stratified-smooth space.}
{\color{red}
\[\begin{tikzcd}
	\bullet &&& \bullet
	\arrow["{<}"{marking}, no head, start anchor={[xshift=-1.5ex, yshift=-0.4ex]}, end anchor={[xshift=1.5ex, yshift=-0.4ex]}, bend left=20, from=1-1, to=1-4]
	\arrow["{<}"{marking}, no head, start anchor={[xshift=-1.5ex, yshift=0.4ex]}, end anchor={[xshift=1.5ex, yshift=0.4ex]}, bend right=20, from=1-1, to=1-4]
	\arrow["{<}"{marking}, no head, start anchor={[xshift=-1.5ex]}, end anchor={[xshift=1.5ex]}, from=1-1, to=1-4]
\end{tikzcd}\quad \quad \quad \quad \quad  \begin{tikzcd}
	\bullet &&& \bullet
	\arrow["{>}"{marking}, no head, start anchor={[xshift=-1.5ex, yshift=-0.4ex]}, end anchor={[xshift=1.5ex, yshift=-0.4ex]}, bend left=20, from=1-1, to=1-4]
	\arrow["{>}"{marking}, no head, start anchor={[xshift=-1.5ex, yshift=0.4ex]}, end anchor={[xshift=1.5ex, yshift=0.4ex]}, bend right=20, from=1-1, to=1-4]
	\arrow["{>}"{marking}, no head, start anchor={[xshift=-1.5ex]}, end anchor={[xshift=1.5ex]}, from=1-1, to=1-4]
\end{tikzcd} \]
}\\
\caption*{The two orientations of the top strata which do not define an oriented stratified-smooth space.}
\end{figure}

There is no harm in the fact that this connected stratified-smooth space admits six distinct orientations. That the two `orientations' in which all three arcs go the same direction are not compatible with any cutout chart is shown by the following lemma, which is the same argument as in \cite[Section 2c]{SeSm}. This lemma will be the foundation of our proof that the boundary of a 1-dimensional stratified-smooth space is zero points (counted with sign), and that the boundary orientation is well-defined.

\begin{lemma}\label{cut-count}
Suppose $\psi: [0,\infty) \times \Bbb R^n \to \Bbb R^n$ is continuous and smooth with zero as a regular value on each stratum. Suppose furthermore that $(\psi^\partial)^{-1}(0) = \{0\}$, oriented negatively. Write $\psi_{t,r}: \{t\} \times D^n_r \to \Bbb R^n$ for the restriction to a slice at time $t$ and the ball of radius $r>0$ around the origin. For sufficiently small $t,r$, the map $\psi_{t,r}$ has zero as a regular value and $\psi_{t,r}^{-1}(0)$ is a finite collection of points whose orientations sum to one.
\end{lemma}

\begin{proof}To briefly remark on the signs, remember that $\psi^\partial$ has domain $\{0\} \times \Bbb R^n$ with the boundary orientation, which is negative the standard orientation; thus $\psi_{0,r}^{-1}(0)$ is a single point, oriented oppositely to $(\psi^\partial)^{-1}(0)$. In particular, $\psi_{0,r}^{-1}(0)$ is a single \emph{positively-oriented} point.

By the Sard lemma, for a generic set of $t,r \in (0,\infty)^2$, the given map has zero as a regular value and has no zeroes on the boundary of $D^n_r$. Because $D^n_r$ is compact, the desired count is finite. Further, for $r,t$ sufficiently small, the subset $S^{n-1}_r \times [0,t]$ does not meet $\psi^{-1}(0)$ (by continuity of $\psi$ and its behavior on the boundary). Write \[\overline \psi_{t,r} = \psi i_{t,r}: S^{n-1} \to \Bbb R^n \setminus 0\] for the composite of the inclusion $S^n_r \hookrightarrow \{t\} \times \Bbb R^n$ and the map $\psi$. Then the degree of $\overline \psi_{t,r}$ is precisely the number of intersection points of $\{t\} \times D^n_r$ with $\psi^{-1}(0)$. 

By assumption, $\deg \overline \psi_{0,r} = 1$ for all $r$. Further, because $S^{n-1}_r \times [0,t]$ does not meet $\psi^{-1}(0)$, the maps $\overline \psi_{t,r}$ are homotopic for all sufficiently small $t \geq 0$ and $r > 0$. Thus the degree of $\overline \psi_{t,r}$ for small $t > 0$ is also equal to one, giving the desired result.
\end{proof}

We may understand these counts intrinsically. Suppose $X$ is an oriented 1-dimensional stratified-smooth space. We classify the connected components of the open stratum of $X$ as three types: components which are equal to their own closure are called \emph{isolated}. Otherwise, either $\overline C = C \cup \{x\}$ or $\overline C = C \cup \{x,y\}$. In either case, $C \cong \Bbb R$; if both ends of $C$ accumulate (and accumulate to the same point) we say $C$ is \emph{heteroclinic}, whereas if exactly one end of $C$ accumulates to $x$ we say $C$ is \textit{homoclinic}. In the latter case, we say $C$ \emph{leaves $x$} or \emph{approaches $x$} depending on whether or not $x$ is an accumulation point of the negatively-oriented end of $C$ or the positively-oriented end, respectively.

\begin{lemma}
In the situation of Lemma \ref{cut-count}, for sufficiently small $t, r > 0$ the count $\# \psi_{t,r}^{-1}(0)$ is the number of homoclinic components of $Z(\psi^\circ)$ entering $0$ minus the number of homoclinic orbits of $Z(\psi^\circ)$ leaving $0$. 
\end{lemma}

\begin{proof}
This amounts to counting the contribution of each component to the count $\psi_{t,r}^{-1}(0)$. Consider a neighborhood of $[0,\epsilon] \times D^n_r$ of $0$ with $\epsilon, r>0$ sufficiently small that \[Z(\psi) \cap \partial \left([0,t] \times D^n_r\right) = Z(\psi) \cap \left(\{0, t\} \times D^n_r\right),\] and so that this intersection is transverse. Write $Z'(\psi) = Z(\psi) \cap \left([0,t] \times D^n_r\right)$; this is a compact 1-dimensional stratified-smooth space. Now consider what happens to the components of $Z(\psi)$ as we pass to $Z'(\psi)$: 

\begin{itemize}
\item An isolated component of $Z(\psi)$ either intersects $Z'(\psi)$ trivially, or in an arc whose boundary lies on $\{t\} \times D^n_r$. 
\item A heteroclinic component intersects with $Z'(\psi)$ to give two arcs running from $0$ to $\{t\} \times D^n_r$ (one oriented forward, the other oriented backward) and a collection of arcs running from $\{t\} \times D^n_r$ back to $\{t\} \times D^n_r$.
\item A homoclinic component intersects with $Z'(\psi)$ to give one arc running from $0$ to $\{t\} \times D^n_r$ (oriented forward if the component leaves $x$, and backward if it approaches $x$) and a handful of arcs running from $\{t\} \times D^n_r$ to $\{t\} \times D^n_r$.
\end{itemize}

Now $\# \psi_{t,r}^{-1}(0)$ is the same as the oriented sum of points in the intersection of $Z'(\psi)$ with the submanifold $\{t\} \times D^n_r$. The intersection number of this submanifold with any arc which runs from $\{t\} \times D^n_r$ back to itself is zero, while the intersection with an arc running from $0$ to $\{t\} \times D^n_r$ is $+1$ and the intersection with an arc running $\{\epsilon\} \times D^n_r$ to $0$ is $-1$. This means isolated components contribute nothing to the count; heteroclinic components contribute $1-1 = 0$ to the count; and homoclinc components contribute either $1$ or $-1$ to the count, depending on whetherh they leave $x$ or approach $x$. 
\end{proof}

In particular, this explains why only six of the eight apparent orientations of the $\theta$-graph are compatible with a cutout chart. What's more, this description of the signed count is \emph{intrinsic to the stratified space}: if $Z(\psi)$ is 1-dimensional, we may determine whether or not $\partial Z(\psi)$ is compatibly oriented with the orientation on the interior by investigating the orientation on each arc, without ever passing to a cutout chart.

\begin{cor}\label{cor:dim-one-orpres}
If $Z(\psi_i)$ are 1-dimensional and $\phi: Z(\psi_1) \to Z(\psi_2)$ is an isomorphism of stratified-smooth spaces which is orientation-preserving on the interior, then $\phi$ is orientation-preserving on the boundary, as well.
\end{cor}

We can now use this to confirm that the boundary orientation may be intrinsically defined in general (so is independent of the choice of a cutout chart).

\begin{lemma}\label{lemma:well-def-bdry-or}
Suppose $\psi_i: [0,\infty) \times \Bbb R^{n_i} \to \Bbb R^{m_i}$ are a pair of continuous maps which are smooth and have zero as a regular value on each stratum; further suppose $\phi: Z(\psi_1) \to Z(\psi_2)$ is an isomorphism of stratified-smooth spaces which is orientation-preserving on the top stratum. Then $\phi$ is also orientation-preserving on the boundary strata.
\end{lemma}

\begin{proof}
We reduce to the 1-dimensional case. If $(X, \Delta)$ is a naive stratified-smooth space, and $f: X \to \Bbb R^m$ is a stratitified-smooth map, we say $f$ is \emph{special} if $f^{-1}(0)$ is a stratified-smooth space (that is, it admits cutout charts). Notice that if $\phi: X \to Y$ is an isomorphism of naive stratified-smooth spaces, and $f: Y \to \Bbb R^m$ is special, then $f\phi$ is also special; $\phi$ induces an isomorphism $f^{-1}(0) \cong (f \phi)^{-1}(0)$, and the property of being a stratified-smooth space --- that is, admitting cutout charts --- is preserved under isomorphisms of naive stratified-smooth spaces.

Now consider the situation in the statement of the lemma. Pick $x \in \partial Z(\psi_1)$; we aim to show that $d\phi_x$ is orientation preserving as a map $T_x \partial Z_1 \to T_{\phi(x)} \partial Z_2$. Passing to a smaller open neighborhood of $x$ if necessary, choose a stratified-smooth function $f: Z(\psi_2) \to \Bbb R^{n_2 - m_2}$ so that $f$ induces an isomorphism $T_x \partial Z_2 \to \Bbb R^{n_2 - m_2}$. Explicitly, one may choose $f$ to be a linear map restricted to $Z \subset [0,\infty) \times \Bbb R^{n_2}$. Because this function extends to a cutout chart, it follows that $f^{-1}(0) \subset Z(\psi_2)$ again carries the structure of a stratified-smooth space, so $f$ is special. Call this zero locus $Z(\psi_2, f)$. 

We thus have a stratification-preserving homeomorphism which is a diffeomorphism on each stratum $\phi_f: Z(\psi_1, f\phi) \to Z(\psi_2, f)$. Further, because $\phi$ was orientation-preserving and the construction defining orientations of zero sets is natural, $\phi_f$ is orientation-preserving on the top (1-dimensional) stratum and is orientation-preserving on the boundary stratum if and only if $\phi$ is orientation-preserving at $x$.

But by Corollary \ref{cor:dim-one-orpres}, we see that $\phi_f$ is indeed orientation-preserving on the boundary; thus $\phi$ is orientation-preserving at an arbitrary boundary point, as desired.
\end{proof}

Now we are prepared to define the boundary of oriented stratified-smooth spaces and check the basic properties.

\begin{definition}
Let $(X, \Delta)$ be an oriented stratified-smooth space. Its \emph{boundary} is the disjoint union \[\partial X = \bigsqcup_{d(s) = \dim X - 1} \overline{X_s},\] where each $X_s$ is oriented with the \emph{boundary orientation} induced by a choice of cutout chart at an arbitrary point $x \in X_s$.
\end{definition}

Note that this boundary orientation is well-defined by the previous lemma: it is independent of the choice of cutout chart. Because an orientation on a stratified-smooth space is induced from an orientation on $[0,\infty)^m \times \Bbb R^n$ via a cutout chart, the boundary operation on stratified-smooth spaces shares many formal properties with the usual boundary operator on oriented smooth manifolds with corners.

\begin{prop}\label{del-sq-or}
Let $(X, \Delta)$ be an oriented stratified space. Then the double boundary $\partial^2 (X, \Delta)$ supports an orientation-reversing involution. More precisely, \[\partial^2(X,\Delta) \cong \bigsqcup_{d(s) = \dim X - 2} \overline{X_s} \sqcup -\overline{X_s}.\]
\end{prop}

\begin{proof}
We saw that \[\partial^2 X = \bigsqcup_{\substack{s \in \Delta \\ d(s) = \dim X - 2}} \overline{X_s} \sqcup \overline{X_s}\] as unoriented stratified-smooth spaces in the proof of Proposition \ref{del-sq}. We need to verify that these two copies of $\overline{X_s}$ appear with the opposite orientation.

Because the definition of boundary orientation is local, it suffices to verify this claim for a given codimension-2 face of $Z(\psi)$. But the (iterated) boundary orientation on a given stratum of $Z(\psi)$ is induced in a uniform way by the (iterated) boundary orientation on the corresponding stratum of $[0,\infty)^m \times \Bbb R^n$, and the stated property holds for Euclidean space with corners.
\end{proof}

We have an analogue of the usual result that the boundary of a compact oriented 1-manifold consists of a finite number of points whose signs sum to zero. The proof is essentially identical to \cite[Lemma 6.7]{M} and to the argument opening \cite[Section 2c]{SeSm}, and it is an immediate application of Lemma \ref{cut-count}. 

\begin{prop}\label{1d-ss-bdry-zero}
Let $(X, \Delta)$ be a compact connected oriented 1-dimensional stratified-smooth space. Then $\partial (X, \Delta)$ is a finite set of points which are canonically oriented with a sign $\pm 1$; the sum of these signs is zero, and hence $\partial(X,\Delta)$ supports an orientation-reversing diffeomorphism.
\end{prop}

\begin{proof}
First, $\Delta$ is a finite set, so there are finitely many 0-dimensional strata. Each bottom stratum is a closed set, hence compact because $X$ is compact. Thus the 0-dimensional strata in $\Delta$ form a compact smooth manifold of dimension 0. 

Now each zero-dimensional stratum has a cutout chart $(U_s, V_s, \psi_s, \phi_s)$. If we delete each neighborhood $U_s$, what remains is a compact 1-dimensional manifold with boundary $X' \subset X$, and it is well-known that the boundary of a 1-dimensional manifold with boundary consists of zero points, counted with sign. 

To be more precise, $\partial X' = \sqcup_s \partial U_s$. Then to conclude, one needs to verify that the signed count of points in $\partial U_s$ agrees with the signed count of points in $X_s$; but (choosing the $U_s$ appropriately) this is precisely the content of Lemma \ref{cut-count}.
\end{proof}

We conclude with a discussion of orientations on fiber products. 

\begin{prop}\label{bd-fibprod-or}
If $X$ and $Y$ are oriented stratified-smooth spaces and $f: X \to B$ and $g: Y \to B$ are stratified-smooth maps to an oriented smooth manifold $B$ without boundary, so that $f$ and $g$ are stratified-transverse, then $X \times_B Y$ has the natural structure of an oriented stratified-smooth space and \[\partial(X \times_B Y) = (\partial X) \times_B Y \sqcup (-1)^{\dim X + \dim B} X \times_B (\partial Y).\]
\end{prop}
\begin{proof}
Recall that the fiber-product of smooth manifolds is oriented with `$B$ in the middle'. Precisely, at the level of vector spaces, suppose $f_i: V_i \to W$ are transverse. Write $I_j = \text{im}(f_j)$ and $K_j = \text{ker}(f_j)$; additionally write $I_1 \cap I_2 = I_{12}$. 

Split $I_1 \cong I_1' \oplus I_{12}$ while $I_2 \cong I_{12} \oplus I_2'$. Further split 
\[V_1 \cong K_1 \oplus I_1 \cong K_1 \oplus I_1' \oplus I_{12}\] 
and 
\[V_2 \cong I_2 \oplus K_2 \cong I_{12} \oplus I_2' \oplus K_2.\] 
Lastly, we may now split 
\[W \cong I_1' \oplus I_{12} \oplus I_2'.\] 

If $V_1, V_2,$ and $W$ are oriented, choose orientations on $I_1', I_2', I_{12}, K_1, K_2$ arbitrarily so that the three above isomorphisms are orientation-preserving; there are four ways to do so, say \begin{align*}
&(o_{I_1'}, o_{I_2'}, o_{I_{12}}, o_{K_1}, o_{K_2}), \\
&(-o_{I_1'}, o_{I_2'}, -o_{I_{12}}, o_{K_1}, -o_{K_2}),\\
&(o_{I_1'}, -o_{I_2'}, -o_{I_{12}}, -o_{K_1}, o_{K_2}), \\
&(-o_{I_1'}, -o_{I_2'}, o_{I_{12}}, -o_{K_1}, -o_{K_2}).
\end{align*}

Now orient $V_1 \times_W V_2 \cong K_1 \oplus I_{12} \oplus K_2$ so that this isomorphism is orientation-preserving; this defines the same orientation on the fiber product in all four cases above, so is well-defined. Applying this construction to the appropriate tangent spaces gives rise to our fiber product orientation on $X \times_B Y$.

The cutout charts constructed in the proof of Proposition \ref{fibprod-ssmooth} are not orientation-preserving. However, we can analyze the difference from the expected orientation. 

Suppose the cutout chart used for $X$ in the above proposition is $\psi: V \to \Bbb R^{\ell}$, and write $b = \dim B, y = \dim Y$; a straightforward and tedious linear algebra argument shows that the orientation of the cutout chart in Proposition \ref{fibprod-ssmooth} disagrees with the fiber product orientation by a factor of $(-1)^{\ell b + \ell y + by}$, and in particular this is consistent over the entire cutout chart. 

It follows that the fiber product orientation on $X \times_B Y$ does indeed give the fiber product the structure of an oriented stratified-smooth space; one may simply flip the orientation on the cutout charts constructed before, if necessary, to make them orientation-preserving.

The desired boundary orientation follows from a tedious definition chase and the above sign computation; notice that it agrees with the expected boundary orientations on a fiber product when this is a transverse fiber product of smooth manifolds with boundary.
\end{proof}

\begin{remark}
Suppose the maps $X \to B$ and $Y \to B$ are stratumwise submersions. Then the coordinate-swapping diffeomorphism $\phi: X \times_B Y \to Y \times_B X$ has sign $(-1)^{(x-b)(y-b)}$, with $x = \dim X$ and so on. Because this is computed in the top stratum, this result follows from the corresponding fact for smooth manifolds, which is straightforward from the definition of fiber product orientation. 
\end{remark}

We conclude with a similar computation for the final proposition of Appendix \ref{app:str-sm}.

\begin{prop}\label{bd-trunc-or}
Suppose $(X, \Delta)$ is an oriented stratified-smooth space and $L: X \to \Bbb R$ a stratified-smooth map, with $c$ a regular value. Then $X^{\le c}$ is also an oriented stratified-smooth space, with boundary \[\partial X^{\le c} = (-1)^{\dim X-1} L^{-1}(c) \sqcup (\partial X)^{\le c},\] where $L^{-1}(c)$ is oriented as a regular zero set.
\end{prop}
\begin{proof}
The orientations on the top strata are inherited as an open subset of $X$. Recall that we orient zero sets $Z(\psi)$ for $\psi: V \to \Bbb R^\ell$ by the rule that $TZ \oplus \Bbb R^\ell \cong TV$.

It is straightforward linear algebra and definition-pushing to see that the cutout charts constructed in Proposition \ref{trunc-still-smooth} are orientation-preserving with respect to the existing orientation on the top stratum of $X$; this amounts to saying that \[(TZ \oplus \Bbb R^\ell) \oplus \Bbb R \cong TV \oplus \Bbb R\] is orientation-preserving, where the map is upper-triangular and the identity on the bottom-right entry, and we know that the top-left entry is orientation-preserving.

In computing the sign on the boundary stratum $L^{-1}(c)$, the additional factor arises because boundaries are oriented with an outward-normal-first convention, while zero sets place the $d/dt$ factor as the final basis element; comparing them requires commuting a factor of $\Bbb R$ across $TL^{-1}(c)$, introducing the stated sign.
%In computing the sign on the boundary stratum $L^{-1}(T)$, recall the construction of boundary orientations. We start with the cutoff chart $V \times [0,\infty) \to \Bbb R^\ell \oplus \Bbb R$ constructed in Proposition \ref{trunc-still-smooth}. We pass to the boundary of the domain, $(-1)^{\dim V + 1} V \to \Bbb R^\ell \oplus \Bbb R$, and compute the orientation on the zero set. 
%
%Write $Z' \cong L^{-1}(T) \cap U$ for the subset of $Z(\psi) \subset V$ on which $L(v) = T$. Then these are oriented by the rules \[TZ \oplus \Bbb R^\ell \cong TV, \quad \quad TZ' \oplus \Bbb R^\ell \oplus \Bbb R \cong (-1)^{\dim V + 1} TV \oplus \Bbb R, \quad \quad TZ' \oplus \Bbb R \cong TZ.\] 
%
%To compare the three of these, in the second expression we should commute the $\Bbb R$ factor in the domain to the middle, which induces an sign factor of $(-1)^\ell$. This establishes that the boundary face which is isomorphic to $L^{-1}(T)$ is isomorphic by an isomorphism with sign $(-1)^{\dim V  \ell + 1} = (-1)^{\dim X + 1}$.
\end{proof}

\subsection{The real blowup of a stratified-smooth space}\label{BlowupChain}
%!TeX root = equivariant-functoriality.tex

In Section \ref{obscob} of the main text, we will find ourselves in the following situation. 

Suppose we have a stratified-smooth space $P$, a vector bundle $q: E \to M$ over a smooth manifold, and a stratified-smooth map $\psi: P \to E$ transverse to the zero section. We will want to construct, in a canonical way, a stratified-smooth space $P'$ with a map $\psi': P' \to S(E)$, where $S(E)$ is the unit sphere bundle of $E$. Over the locus $P \setminus \psi^{-1}(0)$, the map $\psi/\|\psi\|$ gives a canonical such lift; however, this lift cannot be extended over the locus where $\psi = 0$.

To remedy this, we can replace our stratified-smooth space with its \emph{real blowup} (or simply `blowup') along $\psi = 0$; the definition of this space furnishes it with a canonical such lift. 

In this section, we first define this blowup construction and verify that the result is again a stratified-smooth space. After that, we define an equivalence relation on stratified-smooth spaces with maps to $M$ so that a stratified-smooth space is equivalent to any of its real blowups. We conclude by verifying some technical facts which we will make use of in Sections \ref{cob-flow}--\ref{sec:invt}.

\begin{construction}\label{constr:blowup}
Let $P$ be an oriented stratified-smooth space. Let $M$ be a closed smooth manifold equipped with an oriented vector bundle $q: E \to M$, and a stratified-smooth map $\psi: P \to E$ transverse to the zero section. Give $E$ a fiberwise metric; its only purpose is to pick out a sphere bundle of $E$.

We define $B(\psi) \subset  P \times [0,\infty) \times S(E)$ to be the subset given by 
\[B(\psi) = \{(p,t,v) \mid \psi(p) = tv\}.\] 
Both maps $\pi_1: B(\psi) \to P$ and $\pi_3: B(\psi) \to S(E)$, given by projection to first and last coordinate, will be important to us.

Let us give $B(\psi)$ the structure of a stratified-smooth space.

If $P$ has face poset $\Delta$, then call $B(\psi)$'s face poset $\Delta_B(\psi) \subset \Delta \times \{-1,0\}$. To be explicit about this, we have a subset $\Delta_Z(\psi) \subset \Delta$ given by the faces on which $\psi^{-1}(0) \cap P_s \neq \varnothing$. Then $B(\psi)$'s face poset is 
\[\Delta_B(\psi) = \Delta_Z(\psi) \times -1 \sqcup \Delta \times 0,\] 
with grading function $d_\psi(s, -1) = d(s) - 1$ and $d_\psi(s, 0) = d(s)$, and with order induced from $\Delta \times \{-1,0\}$ with the lexicographic order. 

For any $s \in \Delta$, the strata $(s, 0)$ are given by 
\[B(\psi)_{s,0} = \pi_1^{-1}\left(P_s \setminus \psi^{-1}(0)\right).\] 
Being the inverse image of locally closed subsets by a continuous map, these are locally closed; a point in this stratum takes the form $\left(p, \|\psi(p)\|, \psi(p)/\|\psi(p)\|\right)$ for $p \in P_s$, and we give $B(\psi)_{s,0}$ the unique smooth structure so that $\pi_1: B(\psi)_{s,0} \to P_s \setminus \psi^{-1}(0)$ is a diffeomorphism.

On the other hand, if $s \in \Delta_Z(\psi)$, then $B(\psi)_{s,-1} = \pi_1^{-1}(Z(\psi)_s);$ a point in this stratum takes the form $(p, 0, v)$ for $p \in \psi^{-1}(0) \cap P_s = Z(\psi)_s$ and any $v \in S(E_{q\psi(p)})$. In this case, $\pi_3$ gives $B(\psi)_{s,-1}$ the structure of a sphere bundle over $Z(\psi)_s$, canonically identified with the pullback of the sphere bundle $E \to M$ under $f|_{Z(\psi)_s}$. We give $B(\psi)_{s,-1}$ the unique smooth structure so that this identification is a diffeomorphism.

Because \[\dim Z(\psi)_s = d(s) - \text{rank}(E),\] and the sphere bundle has fiber of dimension $\text{rank}(E) - 1$, we see that $\dim B(\psi)_{s,-1} = d(s) - 1 = d_\psi(s,-1)$.

It is straightforward to confirm that \begin{gather*}\overline{B(\psi)_{s,-1}} = \bigcup_{\substack{t \in \Delta_Z \\ t \leq s}} B(\psi)_{t,-1},\\
\overline{B(\psi)_{s,0}} = \bigcup_{\substack{i \in \{-1,0\} \\ t \in \Delta_Z \\ t \leq s}} B(\psi)_{t,i}.
\end{gather*}

%The only part of this statement that is not obvious is that every element $(p,0,v) \in B(\psi)_{t,0}$ arises as a limit of elements $(p_n, r_n, v_n) \in B(\psi)_{s,1}$ for $t \leq s$. But because $\psi: P_t \to L$ is transverse to zero, we see that there is a small disc $c: D(E_p) \to P_t \setminus \psi^{-1}(0)$ with $(\psi c)(v) = \epsilon v$ for some small $\epsilon$. It follows that there is a sequence $(x_n, \epsilon/n, v) \in B(\psi)_{t,1}$ with $(x_n, \epsilon/n, v) \to (p, 0,v)$ --- given by $x_n = \frac 1n c(p,0,v).$ Now one may choose sequences $x_{nm} \in P_s \setminus \psi^{-1}(0)$ with $\lim_m x_{nm} = x_n$, and then correspondingly \[\lim_n (r_n, x_n, v_n) = \lim_n \left(x_{nn}, \|\psi(x_{nn})\|, \frac{\psi(x_{nn})}{\|\psi(x_{nn})\|}\right) = (0, x, v),\] as desired.

Because we have an isomorphism of naive stratified-smooth spaces $\pi_1: B(\psi)_{\ge s, 0} \cong P_{\ge s}$, a cutout chart around any $(p,r,v) \in B(\psi)_{s,0}$ is induced by a choice of cutout chart around $p \in P_s$ (simply lift it along the isomorphism $\pi_1$).

Cutout charts around $(p, 0, v) \in B(\psi)_{s,-1}$ are more intricate. Suppose $(p, 0, v) \in B(\psi)_{s,-1}$. Pick a cutout chart \[\big(U \subset P, V, f: V \to \Bbb R^k, \phi: U \to Z(f)\big)\] for $p \in P_s$. We may choose $U$ small enough to lie inside a trivializing chart for the vector bundle $E$; trivializing the bundle and setting $m = \dim M, e = \text{rank}(E)$, we may write $\psi$ as a map \[\psi: Z(f) \to \Bbb R^m \times \Bbb R^e.\] By definition of stratified-smooth map, we may choose this cutout chart so that the section $\psi$ extends to a map $\widehat \psi = (\sigma, g): V \to \Bbb R^m \times \Bbb R^e$. 

Now choose a small chart $j: \Bbb R^{e-1} \to S^{e-1} \cong E_{q\psi(p)}$ around $v$. Recall here that $q: E \to M$ is the bundle projection. Our cutout chart for $B(\psi)$ around $p$ has domain $V' = V \times [0,\epsilon) \times \Bbb R^{e-1}$ and smooth map \[h(x,t,v) = (f(x), g(x) - tj(v)) \in \Bbb R^k \times \Bbb R^e.\] The isomorphism $\phi': h^{-1}(0) \to B(\psi)$ is given by \[(x,t,v) \mapsto (\phi^{-1}(x), t, j(v)) \in B(\psi) \subset P \times [0,\infty) \times S(E);\] this is a homeomorphism onto its image, the set of $(t,x,v) \in B(\psi)$ with $\|\psi(x)\| < \epsilon$ and $x \in U$, which is open. Because this map $\phi'$ also induces a bijection between the natural stratifications and is a difeomorphism on each stratum, it gives a cutout chart near $(p, 0, v)$.

Now orient $B(\psi)$ is oriented so that the map $\pi_1: B(\psi) \to P$ is an orientation-preserving diffeomorphism off of $Z$. It is easy to see the cutout charts constructed above are orientation-preserving, so that this gives $B(\psi)$ the structure of an oriented stratified-smooth space. 
\end{construction}

In the situation above, $(P, \phi)$ and $(B(\psi), \phi \pi_1)$ are both stratified-smooth spaces equipped with stratified-smooth maps to $M$; if $P$ is connected, then $B(\psi)$ is as well. For the sake of our technical arguments in the main text, we wish to construct an equivalence relation on connected stratified-smooth spaces mapping to $M$ so that these are considered equivalent. To get there --- and to define the geometric chain complex in the next section --- we need a handful of definitions.

\begin{definition}
Let $M$ be a smooth manifold. The set of \emph{stratified-smooth probes} into $M$, written $\mathsf P(M)$, is the set of pairs $(P, \phi)$, where $P$ is a compact connected oriented stratified-smooth space, $\phi: P \to M$ is a stratified-smooth map, and we consider $(P, \phi) = (Q, \eta)$ if there exists an orientation-preserving isomorphism $f: P \to Q$ with $\eta = \phi f$. This set carries an involution, $\iota(P, \phi) = (-P, \phi)$, orientation-reversal of $P$. We call a probe which is fixed under this involution \emph{trivial}.

We say that a probe $(P, \phi)$ is \emph{degenerate} if $\dim P > 0$ and, for all strata $P_s$ of dimension $d(s) > 0$ and all $p \in P_s$, we have \[\text{rank}(d\phi_s: T_p P_s \to M) < d(s).\] We write $\mathsf D(M) \subset \mathsf P(M)$ for the set of degenerate probes in $M$; this subset is $\iota$-invariant.
\end{definition}

\begin{remark}\label{rmk:compare-degen}
There are two other notions of degenerate probe in the literature. In \cite{LipGH}, a probe is said to be \emph{small image} if its image can be covered by the countable union of images of probes of smaller dimension, while in \cite{FMMS2} a probe is said to be \emph{small rank} if $d\phi$ has non-maximal rank on the top stratum. Both authors then say that a degenerate probe is a small probe whose boundary is the disjoint union of trivial probes, pairs of isomorphic but oppositely-oriented probes, and small probes.

Our definition of `degenerate probe' includes fewer probes (for instance, the projection map $[0,1]^2 \to [0,1]$ is not degenerate according to our definition, whereas it is for both of these authors), and in particular our $C_*^{gm}(M)$ will usually be supported in infinitely many degrees; for these authors, it is supported in degrees $* \in [0, \dim M + 1]$. However, when analyzing the notion of collapse-equivalence below, it is difficult to work with the condition that the boundary be a sum of a small rank probe and a trivial probe. The definition chosen above is better suited for studying collapse-equivalence.
\end{remark}

We will now define \emph{collapse-equivalence} of stratified-smooth spaces (though we will only apply it under the assumption that the given stratified-smooth spaces are connected), and then verify that $(P, \phi)$ and $(B(\psi), \phi \pi_1)$ are collapse-equivalent. 

\begin{definition}
Let $(P,\Delta)$ be a stratified-smooth space. We say a \emph{stratified-smooth subspace of codimension $k$} is a stratified-smooth space $(Z, \Delta_Z)$, an order-preserving injection $\overline j: \Delta_Z \to \Delta$ with $d(\overline j(s)) = d_Z(s) + k$, and a topological embedding $j: Z \to P$ so that $j$ restricts on each stratum $Z_s$ to a smooth embedding \[j: Z_s \to P_{\overline j(s)}.\] We identify $Z$ with its image $j(Z)$ in the following discussion.

Let $(P, \phi) \in \mathsf P(M)$. We say that a closed subset $P_{SL} \subset P$ is a \emph{small locus} if $P_{SL} = P_{\codim > 0} \cup P_{D}$, where $P_{\codim > 0}$ is the union of finitely many codimension $k>0$ stratified-smooth subspaces of $P$ (not necessarily compact), and $P_D \subset \partial P$ is the union of degenerate strata of positive codimension.

Suppose there exist small loci $P_{SL} \subset P, Q_{SL} \subset Q$, and an orientation-preserving diffeomorphism of stratified smooth spaces $f: P \setminus P_{SL} \to Q \setminus Q_{SL}$ so that $f\phi = \eta$. Then we say $(P, \phi)$ and $(Q, \eta)$ are {\it collapse-equivalent}, and write this as $(P, \phi) \sim_C (Q, \eta)$.
\end{definition}

This is indeed an equivalence relation, and plays well with the notion of degeneracy and boundary operators.

\begin{lemma}\label{lemma:degen-collapse}
Collapse-equivalence gives an equivalence relation $\sim_C$ on $\mathsf P(M)$. If $(P, \phi)$ is degenerate, then so is $(Q, \eta)$. If $(P, \phi) \sim_C (Q, \eta)$ then $\iota(P, \phi) \sim_C \iota (Q, \eta)$ and $\partial (P, \phi) \sim_C \partial (Q, \eta) + (R, \sigma_{D})$, where $(R, \sigma_D)$ is degenerate.
\end{lemma}

\begin{proof}
Reflexivity and symmetry are obvious, as is the claim about orientation-reversal. 

For transitivity, observe that the union of two small loci is again a small locus; if $(P, \phi) \sim_C (Q, \eta) \sim_C (R, \psi)$ and we have small loci $P_{SL} \subset P,$ $Q^1_{SL}, Q^2_{SL} \subset Q$ and $R_{SL} \subset R$ so that \begin{align*}f: (P \setminus P_{SL}, \phi) \cong (Q \setminus Q^1_{SL}, \eta)  \\
g: (Q \setminus Q^2_{SL}, \eta) \cong (R \setminus R_{SL}, \psi),
\end{align*}
Each of these small loci is written as a union $P_{SL} = P_{\codim > 0} \cup P_D$ of positive-codimension subspaces and degenerate strata. Write 
\[P'_{\codim > 0} = P_{\codim > 0} \cup f^{-1}\(Q^2_{\codim > 0}\) \quad \quad R'_{\codim > 0} = R_{\codim > 0} \cup g\(Q^1_{\codim > 0}\);\] take the locus $P'_D$ to be the union of $P_D$ and the strata intersecting $f^{-1}(Q^2_D)$, similarly for $R'_D$. Notice that the components meeting $f^{-1}(Q^2_D)$ are degenerate: those which aren't already in $P_D$ have small rank on the open dense set $(P \setminus P_{SL})_s$ of each positive-dimension stratum, so by semicontinuity of rank they have small rank on each positive-dimension stratum. So $P'_D$ and similarly $R'_D$ are both unions of degenerate strata of $P$ and $R$.

Then $P'_{SL} = P'_{\codim > 0} \cup P'_D$ is a small locus, as is $R'_{SL} = R'_{\codim > 0} \cup R'_D$. Finally, $fg$ gives an isomorphism $(P \setminus P'_{SL}, \phi) \cong (R \setminus R'_{SL}, \psi)$. This proves transitivity of $\sim_C$.

As implicitly mentioned above, degeneracy off of a small locus implies degeneracy on the whole probe, so collapse-equivalence preserves degeneracy.

The claim about boundaries is also clear: for all strata not contained in the small locus, the isomorphism $f: P \setminus P_{SL} \to Q \setminus Q_{SL}$ descends to an isomorphism $\overline{P_s} \setminus P_{SL} \to \overline{Q_{f(s)}} \setminus Q_{SL}$. Any strata contained in the small locus are degenerate, by definition of small locus. So $f$ provides a collapse-equivalence between the non-degenerate parts of $\partial P$ and the non-degenerate parts of $\partial Q$, and the remaining pieces are collapse-equivalent and form the chain $(R, \sigma_D)$ mentioned in the statement of the proposition.
\end{proof}

Henceforth, we write $\mathsf Q(M) = \mathsf P(M)/\sim_C$. We write $\mathsf D_Q(M) = \mathsf D(M)/\sim_C \subset \mathsf Q(M)$. 

\begin{lemma}\label{blowup-equiv}
Let $\phi: P \to M$ be an oriented stratified-smooth probe equipped with a lift $\psi: P \to E$ of $\phi$ to a section of a vector bundle $q: E \to M$, which is stratified-transverse to the zero section. Then the real blowup $B(\psi)$ satisfies the following properties: \begin{enumerate}
\item[(i)] $B(\psi)$ is an oriented stratified-smooth probe in $S(E)$ via the projection map $\pi_3: B(\psi) \to S(E)$, with boundary relation \[\partial (B(\psi), \pi_3) = (B(\psi|_{\partial P}), \pi_3) - (-1)^{\dim Z_\psi} (Z_\psi \times_B S(E), \pi_2).\] 
\item[(ii)] When composed with the projection map $q: S(E) \to M$, the blowup defines a probe $(B(\psi), q\pi_3) = (B(\psi), \phi \pi_1)$ in $M$ which is collapse-equivalent to $(P, q\psi)$.
\end{enumerate}
\end{lemma}

\begin{proof}
The boundary relation may be read off from the description of the cutout charts. 

The projection map $\pi_1: B(\psi) \to P$ is an isomorphism away from $Z = \psi^{-1}(0) \subset P$ (a small locus), and $(\pi^{-1}(Z), \phi \pi_1)$ is a subset of $\partial B(\psi)$ which factors through the fiber bundle projection $\pi_1^{-1}(Z) \to Z$ with positive-dimensional fibers. It follows that every face of $(\pi^{-1}(Z), \phi \pi_1)$ is small rank, and hence this is a degenerate probe. Thus $(B(\psi), \phi \pi_1)$ is collapse-equivalent to $(P, \sigma)$.
\end{proof}

In the main text, we will need a variation on this which identifies two stratified-smooth probes related to this blowup construction. We include this construction and lemma in this section, as it's in the spirit of the work carried out above.

\begin{construction}\label{constr:weird-zero-set}
Let $(P, \sigma)$ and $(Q, \eta)$ be transverse stratified-smooth probes mapping submersively to $B$, equipped with lifts $\psi: P \to E$ and $\phi: Q \to E$ lifting $\sigma$ and $\eta$, respectively (where $E \to B$ is an oriented vector bundle of rank $\text{rk}(E) > 1$). Assume further that the section $(p,t,q) \mapsto  (1-t)\psi(p) - t\phi(q)$ on $P \times [0,1] \times Q$ is transverse to the zero section. 

We write $Z_{\psi \phi} \subset P \times_B ([0,1] \times Q)$ for the zero set of $(p,t,q) \mapsto (1-t)\psi(p) - t\phi(q)$; by the assumption that this is a transverse zero set, $Z_{\psi \phi}$ is an oriented stratified-smooth space. If $P$ and $Q$ are furthermore equipped with maps to $A$ and $C$, respectively, $Z_{\psi \phi}$ defines a probe in $A \times C$.

There is another way to use these maps to construct a stratified-smooth probe in $A \times C$. We write $B(\psi)$ and $B(\phi)$ for the blowups, considered as probes in $A \times S(E)$ and $S(E) \times C$, and consider the probe $B(\psi) \times_{S(E)} B(\phi)$ in $A \times C$. This fiber product is clearly transverse over blowup strata (which map submersively onto the fibers of $S(E)$), and as will be seen in the proof below, on the open strata the fact that this is a transverse fiber product is equivalent to the statement that $(1-t)\psi -t\phi$ is transverse to zero.
\end{construction}

In Section \ref{obscob}, we will see stratified-smooth chains of the form $Z_{\psi \phi}$ naturally appear from analytic considerations, whereas stratified-smooth chains of the form $B(\psi) \times_{S(E)} B(\phi)$ appear in the same context from algebraic considerations. To square these two perspectives, we will need to see that these two stratified-smooth chains are collapse-equivalent. 

\begin{lemma}\label{weird-zero-set}
In the situation of the above construction, the fiber product $B(\psi) \times_{S(E)} B(\phi)$ is collapse-equivalent to $Z_{\psi \phi}$ as probes in $A \times C$. If $B$ and $E$ are both even-dimensional, then this collapse-equivalence is orientation-preserving (otherwise there is an overall sign).
\end{lemma}

\begin{proof}
One should intuit this claim as follows. On the open stratum where $\psi, \phi \ne 0$, the space $B(\psi) \times_{S(E)} B(\phi)$ is diffeomorphic to the subset of the fiber product $P \times_B Q$ of those pairs $(p, q)$ so that $\psi(p)$ and $\phi(q)$ are positive scalar multiples of one another. On the other hand, supposing $\psi(p), \phi(q)$ are both nonzero, there exists a unique $t$ with $(1-t)\psi(p) = t\phi(q)$ if and only if the two vectors $\psi(p)$ and $\phi(q)$ are positive scalar multiples of one another. 

Indeed, this suggests a reasonably simple map from one to the other. Points in the fiber product $B(\psi) \times_{S(E)} B(\phi)$ are sextuples $(p, t, v, q, s, w)$ with \[\sigma(p) = \eta(q), \;\;\; \psi(p) = tv, \;\;\; \phi(q) = sw, \;\;\; v = w.\] We would like to map $F: B(\psi) \times_{S(E)} B(\phi) \to Z_{\psi \phi}$ via \[f(p, t, v, q, s, w) = (p, \frac{t}{s+t}, q).\] Notice that because the maps to $A$ and $C$ are defined by projecting to $P \times Q$, and $F$ preserves the projection to $P \times Q$, the map $F$ preserves the map to $A \times C$. 

However, $F$ is not defined everywhere, and is not surjective. 

The map $F$ is ill-defined on is the set where $s = t = 0$; these are the points in the blowup corresponding to $\psi = \phi = 0$. This is a union of degenerate strata \[(Z_\psi \times_B S(E)) \times_{S(E)} (S(E) \times_B Z_\phi) = Z_\psi \times_B S(E) \times_B Z_\phi\] in the domain. This is degenerate because each stratum takes the form $Z_s \times_B S(E) \times_B Z_t$, which factors through the projection to $Z_s \times_B Z_t$ with positive-dimensional closed fibers, and hence has small rank on every stratum. 

Next, $F$ misses the codimension $\text{rank}(E)$ subset of $Z_{\psi \phi}$ defined by $t\psi - (1-t)\phi = 0$; this condition implies that $\psi = \phi = 0$, because the points in $Z_{\psi \phi}$ have $(1-t)\psi - t\phi = 0$ to begin with; this subset can be identified with $Z_\psi \times_B [0,1] \times Z_\phi$. Thus $F$ defines an isomorphism away from small loci in the domain and codomain, which preserves the map to $A \times C$, and hence is a collapse-equivalence. 

What remains is a discussion of orientations. We will freely use global symbols (such as $M$, or $B$) to mean small open subsets of these as appropriate; because everything maps submersively to $B$ and $B$ is even-dimensional, we may effectively ignore it by passing to fibers, so we suppress it below, writing $E$ for an oriented vector space. We orient $S(E)$ with the outward-normal-first convention, so that the scaling map $\Bbb R_+ \times S(E) \cong E - 0$ is an orientation-preserving isomorphism.

An open locus in $B(\psi) \times_{S(E)} B(\phi)$ may be identified with \[P^{\langle v\rangle} \times S(E) \times Q^{\langle v\rangle},\] 
where the superscript indicates that we take a fiber above $\Bbb R^+ v \subset E$, where $v$ is fixed; in the first case, we orient this so that the natural locally-defined map $P^{\langle v\rangle} \times S(E) \cong P$ is an oriented diffeomorphism; in the second case $S(E) \times Q^{\langle v\rangle} \cong Q$, as is our convention for orientations on the fiber product. Then an open subset of the fiber product is identified with \[B(\psi) \times_{S(E)} B(\phi) \sim P^{\langle v\rangle} \times S(E) \times Q^{\langle v\rangle}.\]

On the other side, there is a natural identification of $Z_{\psi \phi}$ with the zero locus of the map $(p, t, q) \mapsto (1-t)\psi(p) - t\phi(q)$ from $P \times [0,1] \times Q$ to $E$. Using the identifications above near a point $(p,t,q)$ with $\psi(p)$ and $\phi(q)$ positive scalar multiples of $v \in S(E)$, we write this as the zero set of a map on $P^{\langle v\rangle} \times S(E) \times [0,1] \times S(E) \times Q^{\langle v\rangle}$. Further, we identify a neighborhood of $v \in S(E)$ with a neighborhood of zero in $\langle v\rangle^\perp$, oriented so that $\langle v\rangle \oplus \langle v\rangle^\perp \cong E$. If $\psi(p) = v_p$ and $\psi(q) = v_q$ with $(1-t_0) v_p = t_0 v_q$, we can identify of $Z_{\psi \phi}$ near $(p, t_0, q)$ as the zero set of the map \[P^{\langle v\rangle} \times \langle v\rangle^\perp \times (0,1) \times \langle v\rangle^\perp \times Q^{\langle v\rangle} \to E \quad \text{given by} \quad (p,w,t,u,q) \mapsto (1-t)w - tu + (t-0-t)v_p + (t_0-t)v_q.\] For now, we ignore the outside terms and focus on the zero set, considered as a function on $\langle v\rangle^\perp \times (0,1) \times \langle v\rangle^\perp$.

The first part of the given expression lies in $\langle v\rangle^\perp$, while the latter lies in $\langle v\rangle$; they both must be zero. Because $v_p$ and $v_q$ are positive scalar multiples of one another, $-(t-t_0)(v_p+v_q)$ is only zero if $t = t_0$, and then the first equation gives $u = \frac{1-t_0}{t_0} w$, so by projection to the $w$-factor, this zero set is diffeomorphic to $\langle v\rangle^\perp$. 

To justify that this projection is an oriented diffeomorphism, we check that for fixed $w_0$, the map \[(t,u) \mapsto (1-t)w_0 - tu -(t-t_0)(v_p+v_q)\] is an oriented diffeomorphism $(0,1) \times \langle v\rangle^\perp \to E$ between open subsets. To see this, take the derivative; if $(t_0, u_0)$ is fixed and $\{u_2, \cdots, u_n\}$ is an oriented basis for $\langle v\rangle^\perp$, at $(t_0, u_0)$ the maps the oriented basis $\{t, u_2, \cdots, u_n\}$ to $\{-u_0-w_0-v_p-v_q, -tu_2, \cdots, -tu_n\}$; this is equivalent to the oriented basis $\{-v, -u_2, \cdots, -u_n\}$, which is equivalent to $(-1)^{\text{rank}(E)}$ the standard basis. By the assumption that $E$ has even rank, this proves the claim that $Z_{\psi \phi}$ is locally oriented diffeomorphic to $P^{\langle v\rangle} \times \langle v\rangle^\perp \times Q^{\langle v\rangle}$. 

Finally, with respect to this description in the domain and codomain, the map $F$ is isotopic to the identity. Thus $F$ is orientation-preserving, as desired.
\end{proof}
\subsection{Geometric homology from stratified-smooth spaces}\label{CheckAxioms}
In this section, we use the set of stratified-smooth probes modulo collapse-equivalence to define a chain complex associated to smooth manifolds without boundary which is natural under smooth maps. Using a variation on the Eilenberg--Steenrod axioms, we will verify that the resulting homology theory is naturally isomorphic to singular homology on closed manifolds (or, more precisely, on pairs $(X, A)$ which are the interiors of a handlebody-pair).

We first outline the construction of the underlying graded group of our chain complex.

\begin{construction}
Write $\mathsf Q_d(X) = \mathsf P_d(X)/\sim_C$ for the set of $d$-dimensional probes in $X$ modulo collapse-equivalence.

Set \[K_d(X) = \Bbb Z[\mathsf Q_d(X)]/\Bbb Z[\mathsf Q_d(X)]^\iota.\] The purpose of collapsing the $\iota$-invariant subgroup is to set $\iota (P, \phi) = -(P, \phi)$, so that orientation-reversal corresponds to negation; however, quotienting by this more naive relation would introduce 2-torsion into the complex, corresponding to probes which support orientation-reversing isomorphisms. Collapsing the entire $\iota$-invariant subgroup avoids this issue. 

If one writes $\mathsf T_d(X) = \mathsf Q_d(X)^\iota$ for the `trivial probes' (those which are collapse-equivalent to themselves with the opposite orientation), then $K_d(X)$ can be identified with the free abelian group on $(\mathsf Q_d \setminus \mathsf T_d)/\iota$.

Next, collapse the subgroup of $K_d(X)$ spanned by the image of the degenerate probes. The degenerate probes $\mathsf D(X) \subset \mathsf P(X)$ define a $\sim_C$-invariant subset; we define the quotient by this relation to be $\mathsf D_Q(X) \subset \mathsf Q(X)$. We set \[C^{gm}_d(X) = K_d(X)/\langle \mathsf D_Q(X)\rangle.\] We call this the complex of \emph{geometric chains}, and an element a geometric chain. Again, the group $C^{gm}_d(X)$ is free abelian: one may give it a basis in bijection with $\big(\mathsf Q_d \setminus (\mathsf T_d \cup (\mathsf D_Q)_d)\big)/\iota.$

Then $C_*^{gm}(X)$ is the underlying graded abelian group of what we will call the geometric chain complex. 
\end{construction}

To justify the name `geometric chain complex', we should construct a well-defined boundary operator which squares to zero. 

\begin{lemma}\label{bdry-well-def}
The operation $\partial (P, \phi) = \sum_{d(s) = \dim P - 1} (\overline{P_s}, \phi)$ sending a stratified-smooth probe to its geometric boundary gives a well-defined map $\partial: C_d(X) \to C_{d-1}(X)$ with $\partial^2 = 0$.
\end{lemma}

\begin{proof}
To be slightly more precise, recall that the faces $P_s$ are not assumed to be connected, whereas probes are; when we write $(\overline{P_s}, \phi)$ above, we mean the (finite) sum over its connected components. The compactness of $P$ guarantees that the stratification poset $\Delta$ is finite, and the compactness of each $\overline{P_s}$ guarantees there are only finitely many such components, so this is indeed a finite sum. As for orientations, each $\overline{P_s}$ is equipped with its boundary orientation.

By Proposition \ref{lemma:degen-collapse}, the boundary operator gives a well-defined map $\partial: \mathsf{Q}(X) \to \mathbb Z[\mathsf{Q}]/\mathbb Z[\mathsf{D}_Q]$, which may be extended linearly to all of $\mathbb Z[\mathsf{Q}]$. Because $\partial \iota = \iota \partial$, this descends to a well-defined operator \[K_*(X) = \frac{\mathbb Z[\mathsf Q]}{\mathbb Z[\mathsf Q]^\iota} \to \frac{\mathbb Z[\mathsf Q]}{\mathbb Z[\mathsf Q]^\iota + \mathbb Z[\mathsf D_Q]} = C_*^{gm}(X).\] 
The result will follow from the fact that observe that $\partial \mathsf{D}_Q \subset \mathbb Z[\mathsf{D}_Q] + \mathbb Z[\mathsf Q]^\iota$. In degree $d > 1$, it is obvious that the boundary of a degenerate chain is degenerate. When $d = 1$, observe that the assumption `$\phi$ has small rank' means that $d\phi = 0$ on each connected component of a 1-dimensional stratum, and hence $\phi$ is constant on each such stratum. It follows that $\phi$ is constant on the closure of each connected component of the 1-dimensional stratum; because $P$ is connected and the $0$-dimensional strata lie in the closure of the 1-dimensional strata, it follows that $\phi$ is also constant on the 0-dimensional stratum and thus that the constant value $\phi$ takes on different components all coincide. So $\phi$ is constant. By Proposition \ref{1d-ss-bdry-zero}, the boundary $\partial(X, \phi)$ is an even number of points all mapping to the same point, whose signed sum is zero, and hence is an element of $\mathbb Z[\mathsf{Q}]^\iota$.

The fact that this operation has $\partial^2 = 0$ is Proposition \ref{del-sq-or}.
\end{proof}

So for each smooth manifold without boundary $X$, we obtain a chain complex $C_*^{gm}(X)$ and a corresponding homology group $H_*^{gm}(X)$. This defines a functor on the category of smooth manifolds: if $f: X \to Y$ is a smooth map, simply set $f_* (P, \phi) = (P, f\phi)$. It is obvious both that this is functorial and that it respects the boundary operator.

We would like to extend this to a functor on the category of manifold pairs, which enjoys a relative long exact sequence on homology. This requires a technical lemma.

\begin{lemma}
When $i: A \to X$ is an injective immersion, the map $i_*: C_*^{gm}(A) \to C_*^{gm}(X)$ is injective.\vspace{-0.2cm}
\end{lemma}
\begin{proof}
We will verify that $i_*$ is injective on a basis, by showing that $i: \mathsf{Q}(A) \to \mathsf{Q}(X)$ is injective, with $i^{-1}(\mathsf D_Q(X)) = \mathsf D_Q(A)$ and $i^{-1} \mathsf T(X) = \mathsf T(A)$. 

Suppose $(P, \phi), (Q, \eta) \in \mathsf P(A)$ and $f: P \setminus P_{SL} \to Q \setminus Q_{SL}$ is an isomorphism away from a small locus for which $i \eta = i \phi f$. Because $di$ is injective, a small locus relative to $i \eta$ is a small locus relative to $\eta$, and because $i$ is injective we have $\eta = \phi f$. Then $f$ defines a collapse-equivalence between $(P, \phi)$ and $(Q, \eta)$, so $i: \mathsf{Q}(A) \to \mathsf{Q}(X)$ is injective. By injectivity and the equation $\iota i = i \iota$, if $(P, i\phi) \sim_C \iota (P, i\phi)$, then $(P, \phi) \sim_C \iota (P, \phi)$. Conversely, the pushforward of a trivial chain is always trivial. Therefore, $i^{-1}(\mathsf T(X)) = \mathsf T(A)$.

Suppose $(P, i\phi)$ is degenerate. Because $i\phi$ has small rank on every stratum of $P$ of positive dimension, and $di$ is injective, $\phi$ also has small rank on every stratum of $P$. Thus $(P, \phi)$ is degenerate. Conversely, the pushforward of a degenerate chain is degenerate. Because collapse-equivalence preserves degeneracy, we see that $i^{-1}(\mathsf D_Q(X)) = \mathsf D_Q(A)$.
\end{proof}

Suppose $A \subset X$ is an open submanifold. We may thus define a homology group \[C_*^{gm}(X, A) = C_*^{gm}(X)/i_* C_*^{gm}(A),\] and this is obviously natural for maps of manifold pairs $(X,A) \to (Y, B)$. Further, because the map $C_*^{gm}(A) \to C_*^{gm}(X)$ is injective, the standard homological algebra shows that there exists a long exact sequence \[\cdots \xrightarrow{[-1]} H_*^{gm}(A) \to H_*^{gm}(X) \to H_*^{gm}(X, A) \xrightarrow{[-1]} \cdots\] 

\begin{lemma}\label{ES-additivity}
For two pairs $(X, A)$ and $(Y, B)$, the natural maps $i: (X, A) \to (X \sqcup Y, A \sqcup B)$ and $j: (Y, B) \to (X \sqcup Y, A \sqcup B)$ induce an isomorphism \[H^{gm}_*(X,A) \oplus H^{gm}_*(Y, B) \cong H^{gm}_*(X \sqcup Y, A \sqcup B).\] 
\end{lemma}
\begin{proof}
Because probes $(P, \phi)$ are defined with $P$ connected, a probe in $X \sqcup Y$ factors through either $X$ or $Y$. Thus there are canonical bijections $\mathsf Q(X \sqcup Y) = \mathsf Q(X) \sqcup \mathsf Q(Y)$, and the boundary operator respects the obvious decomposition $C^{gm}_*(X \sqcup Y) \cong C^{gm}_*(X) \oplus C_*^{gm}(Y)$. The inverse to this decompositon is $i_* \oplus j_*$. The result follows.
\end{proof}

\begin{lemma}\label{ES-dimension}
Let $\star$ be the one-point manifold. There is a canonical isomorphism $C_*^{gm}(\star) \cong \Bbb Z[0]$, and hence a canonical isomorphism $H_*^{gm}(\star) \cong \Bbb Z[0].$
\end{lemma}
\begin{proof}
A probe $(P, \phi)$ to the point has $\phi$ constant. By definition, a probe is degenerate if $\dim P > 0$ and for all positive-dimensional strata, $d\phi$ has rank less than the dimension of the stratum. Because $\phi$ is constant, this is tautological, and thus all positive-dimensional probes are zero. Because there is only one zero-dimensional connected stratified-smooth space, $C_*^{gm}(\star)$ is free abelian on one degree-zero generator; further, the unique zero-dimensional stratified-smooth probe has a canonical orientation: the orientation $+1 \in \Bbb R = \det(\Bbb R^0)$. This gives a canonical generator of $C_0^{gm}(\star)$, and hence a canonical isomorphism to $\Bbb Z[0]$.
\end{proof}

The homotopy axiom is only slightly more complicated. The usual proof follows by constructing a map to $C_*(X) \to C_*(I \times X)$ which satisfies an appropriate boundary formula, but we choose to only define $C_*^{gm}(X)$ for open manifolds. Instead, we construct an analogous map to $\Bbb R \times X$.

\begin{lemma}\label{ES-homotopy}
There is a natural map \[\mathsf{Cyl}: C_*^{gm}(X) \to C_{*+1}^{gm}(\Bbb R \times X);\] writing $i^0: X \to \Bbb R \times X$ for the inclusion of $\{0\} \times X$ and similarly for $i^1$, we have the boundary formula \[\partial \mathsf{Cyl} + \mathsf{Cyl} \partial = i^1_* - i^0_*.\] 
Given a smooth homotopy $H: I \times X \to Y$, one may extend this noncanonically to a smooth map $\hat H: \Bbb R \times X \to Y$. Then $\hat H_* \mathsf{Cyl}$ gives a chain homotopy between $H_0$ and $H_1$.
\end{lemma}
\begin{proof}
We set $\mathsf{Cyl}(P, \phi) = (I \times P, \text{id} \times \phi)$, where $I = [0,1]$. Because this operation respects orientation-reversal, collapse-equivalence, and sends degenerate chains to degenerate chains, it induces a map on the underlying chain complexes. The boundary formula here is a special case of the fiber-product boundary formula in Proposition \ref{bd-fibprod-or}. 
\end{proof}

We have now verified almost all of the Eilenberg-Steenrod axioms. The excision axiom is the most subtle, and an order of magnitude more difficult to prove. 

Suppose $(X, A)$ is a manifold pair, and $B$ is a closed subset of $X$ contained in $A$. We want to show that $H_*(X-B, A-B) \to H_*(X, A)$ is an isomorphism. Choose a separating codimension-1 closed submanifold $C$ which is contained in $A$ and so that $B$ is contained on one side of $C$; write $X = Y \cup_C Z$, with $B \subset Y$.

Given a stratified-smooth probe in $X$ whose boundary lies in $A$, then `intersection with $Z$' should produce a stratified-smooth chain in $X-B$ whose boundary lies in $A-B$. In fact, if $\phi$ is the chain in $X$ and $\phi_Y, \phi_Z$ are intersections with $Y$ and $Z$ respectively, then $\phi$ should be homologous to $\phi_Y + \phi_Z$. 

To see why, let us be more precise. Using the Urysohn lemma, choose a smooth function $f: X \to [-1,1]$ with $B \subset f^{-1}(-1)$ and $A^c \subset f^{-1}(1)$; further suppose that $0$ is a regular value. We write $Y = f^{-1}[-1,0]$ and $Z = f^{-1}[0, 1]$, and $C = f^{-1}(0)$ for the codimension 1 submanifold arising as their common boundary. We write \[W = \{(t,x) \in [-1, 1] \times X \mid t \leq 1 - |f(x)|\}\] for the `creased cylinder', or `silo'. It is a manifold with corners, whose boundary is given by $Y \sqcup Z \sqcup -X$, where $-X$ is $X$ equipped with the opposite orientation.

Then intersection with $W$ should cobound $\phi$ and $\phi_Y + \phi_Z$. This, together with some formal manipulation, should prove the excision property. The details are somewhat subtle, as one needs to make these chains transverse to $C$ in a way compatible with degeneracy of chains and the collapse-equivalence relation.

\begin{lemma}\label{ES-excision}
If $(X,A)$ is a manifold pair and $B \subset X$ is a closed subset which is contained in $A$, then the map $i_*: H_*(X - B, A - B) \to H_*(X, A)$ induced by inclusion is an isomorphism.
\end{lemma}

\begin{proof}
We will show that $i_*$ is both injective and surjective. 

To start, pick a homology class $[\phi] \in H_*(X,A)$; pick a representative \[\phi = \sum [(P_i, \phi_i)] \in C_*^{gm}(X)\] where $(P_i, \phi_i)$ are stratified-smooth probes, with $\partial \phi = \phi_A + \phi_{\text{cancel}}$ given by a sum of two terms: $\phi_A$, a chain whose image lies in $A$, and $\phi_{\text{cancel}}$, a chain which is zero in $C_*^{gm}(X)$. This means that $\phi_{\text{cancel}}$ is the disjoint union of a degenerate chain and a chain which is collapse-equivalent to itself with the opposite orientation, neither assumed connected. We refer to a small locus appearing in these collapse-equivalences as a relevant small locus.

Write $j_C: C \hookrightarrow X$ for the inclusion of the submanifold $C$ constructed before this proposition. It need not be the case that the $\phi_i$ are transverse to $j_C$. However, by standard transversality arguments, one may construct an arbitrarily small smooth homotopy $H: [0,1] \times X \to X$ so that $H_0 = \text{Id}$ and $H_1 \phi_i$ are all transverse to $j_C$. Furthermore, one can ensure that the restriction of $H_1 \phi_{\text{cancel}}$ to each relevant small locus is also transverse to $C$. If $H$ is chosen sufficiently small, then $H(I \times \partial \phi)$ has image in $A$, and $[H_1 \phi] = [\phi] \in H_*(X, A)$. By abuse of notation, we write $\phi$ for the representative $H_1 \phi$ with the desired transversality properties. 

Next, set $\phi_Z = Z \times_X \phi$ and similarly for $\phi_C, \phi_Y, \phi_W$. Because \[\partial \phi_Z = \phi_C + (\partial \phi)_Z,\] and $C$ is contained in $A$, and $\partial \phi = \phi_A + \phi_{\text{cancel}}$ with the first term contained in $A$, we see that $\partial \phi_Z$ lies in $C_*^{gm}(A)$; this requires the observation that $Z \times_X \phi_{\text{cancel}}$ is still the sum of a degenerate chain and a chain collapse-equivalent to itself with the opposite orientation. The first fact follows because $j_Z$ is transverse to $\phi_{\text{cancel}}$, and the second because $j_Z$ is transverse to the restriction of $\phi_{\text{cancel}}$ to each small locus, so that the fiber product of the small locus with $Z$ again defines a small locus. One may argue similarly to determine the boundary of the chains $\phi_Y$ and $\phi_W$. 

In particular, $[\phi_Z] \in H_*(X, A)$ is well-defined. Further, because $Z \subset X-B$, this in fact defines an element of $H_*(X-B, A-B)$. We wish to show that $[\phi_Z] = [\phi]$, so that the map $i_*$ is surjective. 

To see this, observe that \[\partial(\phi_W) = \phi_Y + \phi_Z - \phi + (\partial \phi)_W.\] Because $\phi_Y$ lies in $C_*^{gm}(A)$, and $(\partial \phi)$ lies in $C_*^{gm}(A)$, so too does $(\partial \phi)_W$; it follows that in $C_*^{gm}(X, A)$ we have $\partial \phi_W = \phi_Z - \phi$, and thus $[\phi] = [\phi_Z]$ in $H_*(X, A)$, as desired.

Injectivity follows essentially the same lines, now using that we may construct $H$ so small that chains which are already transverse to $j_C$ remain transverse to $j_C$.
\end{proof}

We are now ready to assemble this into the main theorem of this Appendix.

\begin{theorem}
The geometric chain complex $C_*^{gm}(X)$ is a functorial chain complex on the category of smooth manifolds and smooth maps. When $X$ is the interior of a compact manifold, the resulting homology groups are canonically isomorphic to singular homology $H_*(X)$.
\end{theorem}

\begin{proof}
We have now seen that $H_*^{gm}(X,A)$ defines a functor on the category $\mathsf{MP}$ of pairs $(X, A)$, where $X$ is a smooth manifold without boundary and $A$ is an open submanifold, and that this functor satisfies an appropriate version of the Eilenberg-Steenrod axioms for this category. Let us now explain how to modify the usual Eilenberg-Steenrod argument to this setting. 

There is a closely related category, the category $\mathsf{HP}$ of \emph{handlebody pairs} $(Y, B)$. Here, $Y$ is a compact manifold (possibly with boundary) equipped with a handlebody decomposition, where the index of newly attached handles is monotonically increasing and $B$ is a submanifold (possibly with boundary) obtained by only attaching some of these handles. We demand that $Y$ is obtained by attaching handles of index no smaller than the largest index of a handle in $B$. Then $(Y,B)$ carries a natural filtration $(Y^{(k)}, B^{(k)}$ of the submanifolds obtained by attaching only handles up to index $k$. 

We say that a smooth map $(Y, B) \to (Y', B')$ is handle-preserving if it preserves this filtration. There is a functor $\mathsf{HP} \to \mathsf{MP}$ sending $(Y,B)$ to $(Y^\circ, B^\circ)$. 

Both $\mathsf{MP}$ and $\mathsf{HP}$ have a natural notion of homotopy; in the former category, this just means smooth homotopy; in the latter category, we say $f_0, f_1: (Y, B) \to (Z, C)$ are homotopic if there exists a handle-preserving smooth map $(I \times Y, I \times B) \to (Z, C)$ interpolating between them. A cellular approximation-type lemma implies that the corresponding functor on homotopy categories $\mathsf{hHP} \to \mathsf{hMP}$ is an equivalence of categories. 

Now by passing to $\mathsf{MP}$, the geometric homology functor $H_*^{gm}$ defines a functor on the category $\mathsf{hHP}$. If we can show that $H_*^{gm}$ is naturally isomorphic to singular homology on $\mathsf{hHP}$, the result follows for $\mathsf{hMP}$ and hence $\mathsf{MP}$ as well. But what remains follows essentially the standard argument; one determines a canonical identification $H_*^{gm}(X^{(k)}, A^{(k)}) \cong H_*(X^{(k)}, A^{(k)})$ inductively on the index $k$, using the Eilenberg-Steenrod axioms.
\end{proof}

\begin{remark}
It seems probable that the assumption that $X$ is the interior of a compact manifold could be dropped, but we will only actually use this theorem for the manifolds $X = SO(3), S^2, \star$, so it seems rather pointless to extend the result more generally.
\end{remark}

\begin{remark}\label{rmk:disconn-equiv}
One could generalize the notion of collapse-equivalence to disconnected stratified-smooth spaces, and define a geometric chain complex by quotienting $\mathbb Z[\mathsf{P}(X)]$ by this more general relation, then collapsing degenerate chains and trivial chains. The arguments above go through with little change, except that it is no longer clear that the geometric chain complex is free as an $R$-module. This is undesirable for us, as we only obtain results using this complex when $R$ is a field.
\end{remark}
\subsection{Truncating the geometric chain complex}\label{gm-trunc}
In some presentations of geometric homology theory, the chain complex $C_*^{gm}(M)$ is concentrated in degrees $[0, \dim M+1]$. This is a desirable feature for two reasons. First, the periodic homological algebra of \cite[Appendix A]{M} only works for bounded dg-algebras. Second, this phenomenon guarantees we only need the moduli spaces $\breve M^+_\zeta(\alpha, \beta)$ up through a certain dimension range, which allows us to ignore the bubbling phenomenon.

This is not the case for our model $C_*^{gm}(M)$ whenever $\dim M > 0$. Pick any nondegenerate $\phi: P \to M$ of positive degree (for instance, take $(P, \phi) = (M, 1_M)$ to be the identity probe). Then $(P \times I^k, \phi \pi_1)$ is nondegenerate, as the positive-dimensional stratum $(P \times \{0\}, \phi)$ does not have small rank. So $C_*^{gm}(M)$ is nonzero in arbitrarily large degrees. This is because our rather strict notion of nondegeneracy includes fewer nondegenerate probes than other presentations; see Remark \ref{rmk:compare-degen} for more details.

We briefly explain how to correct for this failure.

Let $M$ be a smooth manifold. Set the \emph{truncated geometric chain complex} to be $$tC_k^{gm}(M) = \begin{cases}C_k^{gm}(M) & k < \dim M \\
C_k^{gm}(M)/\partial C_{k+1}^{gm}(M) & k = \dim M \\
0 & k > \dim M\end{cases}.$$

That is, $tC_*^{gm}(M)$ is the complex obtained by collapsing the acyclic subcomplex $$D_*(M) = \begin{cases} 0 & k < \dim M \\ \partial C_{k+1}^{gm}(M) & k = \dim M \\
C_k^{gm}(M) & k > \dim M.
\end{cases}$$

\begin{theorem}\label{thm:trunc-gm}
The complex $tC_*^{gm}(M;R)$ satisfies the following properties. \begin{enumerate}
\item[(i)] It is concentrated in degrees $[0, \dim M]$, and free over $R$ if $R$ is a PID and $M$ is orientable.
\item[(ii)] Its homology groups are canonically isomorphic to singular homology.
\item[(iii)] If $f: M \to N$ is a smooth map with $\dim M \ge \dim N$, then there is an induced chain map $f_*: tC_*^{gm}(M) \to tC_*^{gm}(N)$.
\item[(iv)] Suppose $(W, f \times g) \in \mathsf P(M \times N)$ is a probe, and suppose $f: W \to M$ is a submersion. Then there is an induced map $\times_M W: tC^{gm}_*(M) \to tC^{gm}_*(N)$ of degree $\dim W - \dim M$, sending $(P, \phi) \mapsto (P \times_M W, \pi_2 g)$, and satisfying the boundary relation $$\partial (P \times_M W) = (\partial P) \times_M W + (-1)^{\dim P + \dim W} P \times_M (\partial W).$$ 
\end{enumerate}
\end{theorem}
\begin{proof}
It is obvious that this complex is concentrated in degrees $[0, \dim M]$. Further, because $D_*(M)$ is an acyclic complex, the canonical map $C_*^{gm}(M) \to tC_*^{gm}(M)$ is a quasi-isomorphism. Because $H_*^{gm}(M)$ was already seen to be canonically isomorphic to singular homology, so too is $tH_*^{gm}(M)$. The quotient complex agrees with $C_*^{gm}(M)$ through degree $\dim M - 1$. There is a short exact sequence $H_{\dim M}(M) \to tC_{\dim M}^{gm}(M) \to \partial C_{\dim M}^{gm}(M) \subset C_{\dim M - 1}^{gm}(M)$. The quotient module is a submodule of a free module, hence free when $R$ is a PID, so the short exact sequence splits. Because $M$ is orientable, $H_{\dim M}(M)$ is free, and the claim follows.

The claimed maps are induced by maps on $C_*^{gm}$. We have already discussed the induced map $f_*$ on $C_*^{gm}$; because $f_*$ is a chain map of degree zero, it obviously has $f_* D_*(M) \subset D_*(N)$, and thus induces a map on the truncated complexes. 

To see that $\times_M W$ is well-defined at the level of $C_*^{gm}$, one needs to see that it respects all of: orientation-reversal of probes; degenerate probes; and collapse-equivalence of probes. It is easy to verify that it does. The stated boundary relation is given in Proposition \ref{bd-fibprod-or}. We now need to verify it descends to $tC$. 

Because $f: W \to M$ is a submersion, the map $\times_M W$ has non-negative degree. If $\deg(\times_M W) = \dim W - \dim M > 0$ it is obvious that $F(W) D_*(M) \subset D_*(N)$ for degree reasons. If $\deg F(W) = 0$, then $\dim W = \dim M$; for $f: W \to M$ is a submersion it must be on all boundary strata as well, but because boundary strata have dimension less than $\dim M$, it follows that $\partial W = \varnothing$. It follows that $F(W)$ is a chain map of degree zero and any such map descends to $tC^{gm}_*$ as desired.
\end{proof}

As a corollary to Theorem \ref{thm:trunc-gm}(iii), if $G$ is a Lie group, then $tC^{gm}_*(G)$ is a dg-algebra and if $G$ acts on $M$ then $tC^{gm}_*(M)$ is a dg-module.
\end{appendices}

\bibliography{references}
\bibliographystyle{hplain}
%\Addresses
\end{document}